\documentclass[11pt]{amsart}
\usepackage{pdfsync}
\usepackage{hyperref}
\usepackage{array} % for the table at the end
\usepackage{amsfonts}
\usepackage{amsmath}
\usepackage{amssymb}
\usepackage{amsthm}
\usepackage{graphicx}
\usepackage{color}

\usepackage[
  hmarginratio={1:1},     % equal left and right margins
  vmarginratio={1:1},     % equal top and bottom margins
  textwidth=17cm,        % new text width
  textheight=21cm,
  heightrounded,          % always useful
]{geometry}

\newtheorem{thm}{Theorem}[section]
\newtheorem{cor}[thm]{Corollary}
\newtheorem{lem}[thm]{Lemma}

\newtheorem{defn}[thm]{Definition}
\newtheorem{rem}[thm]{Remark}
\newtheorem{ex}[thm]{Example}
\newtheorem*{thm*}{Theorem}
\newtheorem*{defn*}{Definition}
\numberwithin{equation}{section}

\def\supp{\mathrm{supp}} %per il supporto
\def\dim{\mathrm{dim}}
\def\dist{\mathrm{dist}}

\newcommand{\BB}{\mathbb{B}}
\newcommand{\RR}{\mathbb{R}}
\newcommand{\NN}{\mathbb{N}}
\newcommand{\QQ}{\mathbb{Q}}
\newcommand{\ZZ}{\mathbb{Z}}
\newcommand{\G}{\mathcal{G}}

\newcommand{\si}{\sigma}
\newcommand{\kk}{\kappa}

\def\dist{\mathrm{dist}} %per la distanza
\def\qed{\,\unskip\kern 6pt \penalty 500d
\raise -2pt\hbox{\vrule \vbox to8pt{\hrule width 6pt
\vfill\hrule}\vrule}\par}
\begin{document}

\title[On the nodal set of solutions to nonlocal parabolic equations]{On the nodal set of solutions to a class of nonlocal parabolic equations}

\author[A. Audrito]{Alessandro Audrito}\thanks{}
\address{Alessandro Audrito \newline \indent
Institut f\"ur Mathematik
Universit\"at Z\"urich, \newline \indent
Winterthurerstrasse 190
CH-8057 Z\"urich, Switzerland
 }
\email{alessandro.audrito@gmail.com, alessandro.audrito@math.uzh.ch}

\author[S. Terracini]{Susanna Terracini}\thanks{}
\address{Susanna Terracini \newline \indent
 Dipartimento di Matematica ``Giuseppe Peano'', Universit\`a di Torino, \newline \indent
Via Carlo Alberto, 10,
10123 Torino, Italy}
\email{susanna.terracini@unito.it}

\date{\today} %%  this cancels date in article format


%
%
%
%
%
%
%%%%%%%%%%%%%%%%%%%%%%%%%%%%%%%%%%%%%%%%%%%%%%%%%%%%%%%%%%%%%%%%%%%%%%%%%%%%%%%%%%%%%%%%%%%%%%%%%%%%%%%%
%%%%%%%%%%%%%%%%%%%%%%%%%%%%%%%%%%%%%%%%%%%%%%%%%%%%%%%%%%%%%%%%%%%%%%%%%%%%%%%%%%%%%%%%%%%%%%%%%%%%%%%%
%
%
%
%
%
%%%%%%%%%%%%%%%%%%%%%%%%%%%%%%%%%%%%%%%%%%%%%%%%%%%%%%%%%%%%%%%%%%%%%%%%%%%%%%%%%%%%%%%%%%%%%%%%%%%%%%%%

\thanks{Work partially supported by the ERC Advanced Grant 2013 n.~339958 Complex Patterns for Strongly Interacting Dynamical Systems - COMPAT. The first author has been partially funded by Projects MTM2011-24696 and  MTM2014-52240-P (Spain), by local projects ``Equazioni differenziali lineari e non lineari'', ``Equazioni differenziali non lineari e applicazioni'' (Italy), and by the INDAM-GNAMPA projects ``Equazioni diffusive non-lineari in contesti non-Euclidei e disuguaglianze funzionali associate'' and ``Ottimizzazione Geometrica e Spettrale'' (Italy).}

\begin{abstract}
We investigate the  \emph{local properties}, including the \emph{nodal set} and the \emph{nodal properties} of solutions to the following parabolic problem of Muckenhoupt-Neumann type:
\begin{equation*}
\begin{cases}
\partial_t \overline{u} - y^{-a} \nabla \cdot(y^a \nabla \overline{u}) = 0 \quad &\text{ in } \mathbb{B}_1^+ \times (-1,0) \\
-\partial_y^a \overline{u} = q(x,t)u \quad &\text{ on } B_1 \times \{0\}  \times (-1,0),
\end{cases}
\end{equation*}
where $a\in(-1,1)$, is a fixed parameter $\mathbb{B}_1^+\subset \RR^{N+1}$ is the upper unit half ball and $B_1$ is the unit ball in $\RR^N$.   Our main motivation comes from its relation with a class of nonlocal parabolic equations involving the fractional power of the heat operator
\begin{equation*}
H^su(x,t) = \frac{1}{|\Gamma(-s)|} \int_{-\infty}^t \int_{\RR^N} \left[u(x,t) - u(z,\tau)\right] \frac{G_N(x-z,t-\tau)}{(t-\tau)^{1+s}} dzd\tau.
\end{equation*}

We characterise the possible blow-ups and we examine the structure of the nodal set of solutions \emph{vanishing with a finite order}. More precisely, we prove that the nodal set has at least parabolic Hausdorff codimension one in $\RR^N\times\RR$, and can be written as the union of a locally smooth part and a singular part, which turns out to possess remarkable stratification properties. Moreover, the asymptotic behaviour of general solutions near their nodal points is classified in terms of a class of explicit polynomials of Hermite and Laguerre type, obtained as eigenfunctions to an Ornstein-Uhlenbeck type operator. Our main results are obtained through a fine blow-up analysis which relies on the monotonicity of an Almgren-Poon type quotient and some new Liouville type results for parabolic equations, combined with more classical results including Federer's reduction principle and the parabolic Whitney's extension.
\end{abstract}

\maketitle

\tableofcontents

%
%
%
%
%
%
%
%
%%%%%%%%%%%%%%%%%%%%%%%%%%%%%%%%%%%%%%%%%%%%%%%%%%%%%%%%%%%%%%%%%%%%%%%%%%%%%%%%%%%%%%%%%%%%%%%%%%%%%%%%%%%
%
%
%%%%%%%%%%%%%%%%%%%%%%%%%%%%%%%%%%%%%%%%%%%%%%%%%%%%%%%%%%%%%%%%%%%%%%%%%%%%%%%%%%%%%%%%%%%%%%%%%%%%%%%%%%%
%
\section{Introduction}
The main purpose of this paper is to give a detailed description of the \emph{local properties}, including, besides regularity,  the \emph{nodal set} and the \emph{nodal properties} of solutions to the following nonlocal parabolic equation involving the fractional power of the heat operator
\begin{equation}\label{eq:NONLOCALHEPOTENTIAL}
(\partial_t - \Delta)^su = q(x,t)u \quad \text{ in } B_1\times(-1,0),
\end{equation}
where $s \in (0,1)$ is a fixed exponent and the potential $q = q( x,t)$ satisfies suitable assumptions \eqref{eq:ASSUMPTIONSPOTENTIAL} (for the moment, we can assume that $q$ is a bounded and smooth function on $B_1\times(-1,0)$).  Our study will be performed through the analysis of the solutions to the auxiliary local problem of Muckenhoupt-Neumann type:

\begin{equation}\label{eq:ExtensionEquationLocalForward}
\begin{cases}
\partial_t \overline{u} - y^{-a} \nabla \cdot(y^a \nabla \overline{u}) = 0 \quad &\text{ in } \mathbb{B}_1^+ \times (-1,0) \\
-\partial_y^a \overline{u} = q(x,t)u \quad &\text{ on } B_1 \times \{0\}  \times (-1,0),
\end{cases}
\end{equation}
where $a=2s-1\in(-1,1)$ and $\mathbb{B}_1^+\subset \RR^{N+1}$ is the upper unit half ball; $u$ will be the trace of $\overline{u}$ at the characteristic hyperplane $y=0$ (see Subsection \ref{subsec:extension}).  While equation \eqref{eq:ExtensionEquationLocalForward} can be easily interpreted in a classical or distributional sense, the definition of fractional powers of the heat operator appearing in \eqref{eq:NONLOCALHEPOTENTIAL} is not straightforward and requires some thorough explanation.

\subsection{Fractional powers of the heat operator}

The fractional power of the heat operator $H := \partial_t - \Delta$ was introduced by M. Riesz in \cite{Riesz1938:art} (see also \cite{Riesz1949:art}) and it is a notable example of pseudo-differential operator with high relevance in pure mathematics as well as in applied fields such as biology, elasticity and finance (see for instance  \cite{AthCaffaMilakis2016:art,BanGarDanPetr2018:art,DainelliGarofaloPetTo2017:book,DuvantLions1972:art,
Hilfer2000:book,MetzlerKlafter2000:art}). As showed in \cite[Section 28]{SamkoKilbasMarichev:1993}, it can be written through the following integral form
\begin{equation}\label{eq:IntegralRepresentationFHO}
H^su(x,t) = \frac{1}{|\Gamma(-s)|} \int_{-\infty}^t \int_{\RR^N} \left[u(x,t) - u(z,\tau)\right] \frac{G_N(x-z,t-\tau)}{(t-\tau)^{1+s}} dzd\tau,
\end{equation}
for smooth bounded functions $u = u(x,t)$ (see also \cite{BanGarofalo2017:art,NystromSande2016:art,StingaTorrea2015:art}), where
\begin{equation}\label{eq:EXPRESSIONFORGAUSSIANINRNHEATEQUATION}
G_N(x,t) =  \frac{1}{(4\pi t)^{N/2}} e^{-\frac{|x|^2}{4t}},
\end{equation}
is the standard Gaussian, the fundamental solution of the heat equation. The integral representation \eqref{eq:IntegralRepresentationFHO} has several interesting aspects. First, it clearly shows that the nonlocal character of  $H^s$ holds both in space and time, since the value $H^su(x,t)$ depends on the values $u(z,\tau)$ for all $z \in \mathbb{R}^N$ and $\tau < t$. Secondly, it allows us to recover the fractional laplacian
\[
(-\Delta)^su(x) = C_{N,s} \int_{\mathbb{R}^N} \frac{u(x)-u(z)}{|x-z|^{N+2s}} dz,
\]
when $u(x,t) = u(x)$, and the Marchaud-Weyl derivative \cite{Marchaud1927:art,Weyl1917:art}
\[
(\partial_t)^su(t) = \frac{1}{|\Gamma(-s)|} \int_{-\infty}^t \frac{u(t) - u(\tau)}{(t-\tau)^{1+s}}d\tau,
\]
when $u(x,t) = u(t)$ (see \cite{StingaTorrea2015:art}), obtained as the inverse of the  Riemann-Liouville integral operator (also called Weyl integral operator)
\[
I^s u(t) = \frac{1}{\Gamma(s)} \int_{-\infty}^t \frac{u(\tau)}{(t-\tau)^{1-s}}d\tau,
\]
for a suitable class of functions (see, for instance, \cite{Hilfer2000:book,SamkoKilbasMarichev:1993}).

Finally, we stress that such operators play a crucial role in many applications (see, for instance, \cite[Chapter 8]{SamkoKilbasMarichev:1993} and/or \cite[Chapter VI-IX]{Hilfer2000:book}) and possess interesting probabilistic interpretations in the context of  Continuous Time Random Walks (CTRW) (see  for instance \cite{BaeumerEtAl2005:art}, \cite[Chapter 3]{Hilfer2000:book}, \cite{KenkreEtAl:1973} \cite[Section 3]{MetzlerKlafter2000:art}, \cite{MontrollWeiss1965:art}). We present below some fundamental ideas.

We consider a test particle moving on a $N$-dimensional lattice, and we denote by $u = u(x,t)$ its probability density function, that is,  $u(x,t)$ represents the probability of finding the test particle at the space-time position $(x,t)$. It is a classical fact that if the particle moves according to a Brownian motion law, then the function $u$ satisfies the classical Heat equation
\[
\partial_tu - \Delta u = 0, \quad t > 0.
\]
We recall that in the standard diffusion model described by the Heat equation, the particle moves of a unit step in space along a random direction, at each unit time. Nonlocal equations are obtained when the nature of this simple diffusion law is  modified and all the lattice sites can be reached by the particle at each time step, though with different probabilities. For example, if, at each unit time, the particle is allowed to make \emph{long jumps}, the probability density function $u = u(x,t)$ is a L\'evy distribution. When in addition $u$ is symmetric, stable and, roughly speaking, the probability of jumping to a distance $r$ decays as $r^{-N-2s}$ for some $s \in (0,1)$, then it satisfies the nonlocal equation
\[
\partial_tu + (-\Delta)^s u = 0, \quad t > 0,
\]
see for instance \cite[Section 3.5]{MetzlerKlafter2000:art} and also \cite{BonSirVaz2016:art,Valdinoci2009:art}, where $(-\Delta)^s$ is the fractional Laplacian. A second option is obtained by assuming that the test particle is allowed to make \emph{long rests} at any position $x \in \mathbb{R}^N$ (see   \cite{KenkreEtAl:1973,MetzlerKlafter2000:art,MontrollWeiss1965:art}), which gives the so called CTRW model. In particular, if the probability of having a waiting time $t$ before the next jump falls of like $t^{-s}$, the probability density function $u$ satisfies the equation
\begin{equation}\label{eq:TimeNonlocalModel}
D_t^su -\Delta u = \frac{t^{-s}}{\Gamma(1-s)} u_0(x), \quad t > 0,
\end{equation}
where $D_t^s$ is the Caputo derivative (see \cite{BaeumerEtAl2005:art} and \cite[Section 3 and Section 5]{Hilfer2000:book}) and $u_0$ is the distribution at time $t = 0$ (in this presentation we are implicitly assuming $u=0$ whenever $t < 0$). This is a first interesting example of sub-diffusive equation (see   \cite{BaeumerEtAl2005:art} and \cite[Section 3.4]{MetzlerKlafter2000:art}). Looking at the right hand side, we see that the equation keeps memory of the initial data $u_0$, for any $t > 0$ and it decays as a negative power of $t$ for large times, in sharp contrast with the classical diffusion framework. Combining the last two models i.e., if the particle can jump and rest, the resulting equation for $u$ turns out to be
\begin{equation}\label{eq:DoubleNonlocalModel}
D_t^{s_0}u + (-\Delta)^{s_1} u = \frac{t^{-s_0}}{\Gamma(1-s_0)} u_0(x), \quad t > 0,
\end{equation}
for some $s_0,s_1 \in (0,1)$. We notice that in such a case, the distribution governing the waiting times and the distribution governing the jumps are \emph{independent}, in the sense that the length of a jump from a position $x$ to a position $z$ does not depend the time spent at the location $x$ (see \cite{Zaslavsky1994:art}).

We are now ready to combine these two behaviours in a model where the probability density function $u = u(x,t)$ satisfyies a parabolic equation involving $(\partial_t-\Delta)^s$. More precisely, let us denote by $T$ the waiting time spent by the particle at some spatial location and assume that $\mathbb{P}(T > t) \sim t^{-s}$ for $t$ large. Furthermore, let $J$ be the random variable governing the jumps lenght. If the two distribution are \emph{not} independent, the equation for $u$ is neither \eqref{eq:TimeNonlocalModel} nor \eqref{eq:DoubleNonlocalModel}. In particular, assuming that  $J = J_t \sim \mathcal{N}(0,2t)$ when $T = t$ (i.e. $J$ is normally distributed with zero average and variance $2t$, whenever $T = t$), then we find the equation
\begin{equation}\label{eq:DoubleCoupledNonlocalModel}
(D_t - \Delta)^s u = \frac{t^{-s}}{\Gamma(1-s) } u_0(x), \quad t > 0,
\end{equation}
where, as above, $u_0$ denotes the initial probability density function. Comparing with \eqref{eq:DoubleNonlocalModel}, the waiting times and the jumps are now \emph{coupled}: if the time spent to a location $x$ is short, then the size of the jump is expected to be short as well, while if the waiting time at $x$ is long, one equally expects a long jump. In other words, the distance covered by the particle in a jump depends on how long the particle has rested at the location before jumping. Notice that in \eqref{eq:TimeNonlocalModel}, \eqref{eq:DoubleNonlocalModel} and \eqref{eq:DoubleCoupledNonlocalModel},  we have used the notation $(D_t - \Delta)^s$, instead of $(\partial_t -\Delta)^s$: this is to emphasize  the fact that the above models reduce to a fractional ODE of Caputo type in time (see again \cite[Section 3 and Section 5]{Hilfer2000:book}), when the distribution is uniform in space. We refer the interested reader to the extensive literature on such models  (\cite{BaeumerEtAl2005:art,KenkreEtAl:1973,MetzlerKlafter2000:art,MontrollWeiss1965:art}), while we quote \cite{Hilfer2000:book,Samko2001:book,SamkoKilbasMarichev:1993} for the general theory of fractional derivates of different type (Caputo, Riemann-Liouville, Marchaud,\ldots) and the links between them.

\medskip

Coming back to our presentation, it is worthwhile stressing that the integral representation \eqref{eq:IntegralRepresentationFHO} is not the unique way to define $H^s$, and there are other options, corresponding to different functional settings. Perhaps the most direct one is to define $H^s$ in terms of its Fourier transform, namely
\begin{equation}\label{eq:FourierDefHs}
\widehat{H^su}(\eta,\vartheta) := (i\vartheta + |\eta|^2)^s\,\widehat{u}(\eta,\vartheta),
\end{equation}
for all functions $u = u(x,t)$ belonging to the natural domain
\begin{equation}\label{eq:NatDomHs}
\text{dom}(H^s):=\left\{u\in L^2(\RR^{N+1}): (i\vartheta + |\eta|^2)^s \, \widehat{u} \in L^2(\RR^{N+1}) \right\},
\end{equation}
and recover $H^su$ through the Plancherel's theorem. As highlighted in \cite{StingaTorrea2015:art} (see also \cite{BanGarofalo2017:art}), the above definition admits an equivalent formulation based on the theory of semigroups. The main idea is to consider a solution $p = p(x,t,\tau)$ to the auxiliary problem
\[
\begin{cases}
\partial_{\tau}p = -Hp \quad &\text{in } \mathbb{R}^{N+1}\times \mathbb{R}_+ \\
p(x,t,0) = u(x,t)     \quad &\text{in } \mathbb{R}^{N+1}\times \{0\},
\end{cases}
\]
where $u \in \text{dom}(H^s)$, i.e. $p(x,t,\tau) = P_{\tau}^Hu(x,t)$ where $P_{\tau}^H = \left\{e^{-\tau H}\right\}_{\tau > 0}$ is the semigroup with generator $H$. From \cite{BanGarofalo2017:art,StingaTorrea2015:art} we know that
\[
P_{\tau}^H u(x,t) = (G_N(\cdot,\tau) \star \Lambda_{-\tau}u)(x,t) = \int_{\mathbb{R}^N} G_N(x-z,\tau) u(t-\tau,z) dz,
\]
where $\Lambda_{\tau} u(x,t) = u(x,t+\tau)$, with
\[
\widehat{P_{\tau}^H u}(\xi,\vartheta) = e^{-\tau(i\vartheta + |\xi|^2)}\widehat{u}(\xi,\vartheta).
\]
Now we recall that the classical formula
\begin{equation}\label{eq:ClassicalSemigroupFormula}
\lambda^s = \frac{1}{\Gamma(-s)} \int_0^{\infty} \left( e^{-\tau\lambda} - 1 \right) \frac{d\tau}{\tau^{1+s}},
\end{equation}
holding for $r > 0$, can be extended to a much larger class of nonnegative closed operators $L$ defined on a Banach space:
\[
L^s f(x) = \frac{1}{\Gamma(-s)} \int_0^{\infty} \left[ P_{\tau}f(x) - f(x) \right] \frac{d\tau}{\tau^{1+s}},
\]
as established by Balakrishnan \cite{Balakrishnan1960:art}, where $P_{\tau} = e^{-\tau L}$ is the semigroup with infinitesimal generator $L$. The main idea in \cite{StingaTorrea2015:art} is to extend formula \eqref{eq:ClassicalSemigroupFormula} to any complex number $\lambda \in \mathbb{C}$ with positive real part, and define the fractional power of $H = \partial_t - \Delta$ through an adaptation of Balakrishnan's result. In particular, the fractional power of the Heat operator can be defined as
\begin{equation}\label{eq:SemigroupDefHs}
H^su(x,t) := \frac{1}{\Gamma(-s)} \int_0^{\infty} \left[ P_{\tau}^Hu(x,t) - u(x,t) \right] \frac{d\tau}{\tau^{1+s}}.
\end{equation}
In \cite{StingaTorrea2015:art} it has been shown that the definitions \eqref{eq:FourierDefHs} and \eqref{eq:SemigroupDefHs} are in facts equivalent and the integral formulation \eqref{eq:IntegralRepresentationFHO} has been deduced as a main consequence. Actually, the three options are equivalent for functions belonging to the Schwartz class $\mathcal{S}(\mathbb{R}^{N+1})$.

In this paper we will always assume that $H^s$ is defined through the integral representation \eqref{eq:IntegralRepresentationFHO} and that the function $q : \RR^{N+1} \to \RR$ satisfies
\begin{equation}\label{eq:ASSUMPTIONSPOTENTIAL}
\begin{cases}
\|q\|_{C^1} \leq K \quad &\text{ if } 1/2 \leq s < 1 \\
\|q\|_{C^2}, \; \|x\cdot\nabla_x q\|_{L^{\infty}} \leq K \quad &\text{ if } 0 < s < 1/2,
\end{cases}
\end{equation}
for some universal constant $K > 0$ (see  \cite{BanGarofalo2017:art} for the definition of the spaces $C^k = C^k(\RR^{N+1})$). These assumptions where already introduced in \cite{BanGarofalo2017:art}: they are needed to prove the regularity of solutions to the extended problem and the monotonicity of an Almgren-Poon quotient. Actually, the majority of our results hold true under weaker assumptions (see  \eqref{Rem:WeakAssPotential}) but, since our results are new even in the case $q = 0$, we will not insist more on this aspect.

We devote the following section to a short review of the extension theory in the parabolic setting, which allows to transform a nonlocal problem into a local one. The switch to a local problem is a key point of our discussion, since it allows us to take advantage of some monotonicity formulas of Almgren-Poon type (see   Theorem \ref{Theorem:GeneralizedFrequency}) and ultimately to carry out the blow-up procedure. In light of the equivalence with the problem, we can divert our attention from the notion of the solution to \eqref{eq:IntegralRepresentationFHO} just defined and focus on the local solutions of the extended problem \eqref{eq:ExtensionEquationLocalForward}.

\subsection{The extension theory: from nonlocal to local}\label{subsec:extension} To justify the switch to a local problem, we need to recall some known facts from the extension theory developed independently by Nystr\"om and Sande \cite{NystromSande2016:art}, and Stinga and Torrea \cite{StingaTorrea2015:art} (see  also \cite[Section 2,3,4]{BanGarofalo2017:art}, \cite{GarofaloTralli2019:art} and \cite{CaffSil2007:art}). We give a short summary of their results, adapting them to our notations/framework.

Let $s \in (0,1)$ and let $u = u(x,t)$ be bounded and smooth in $\mathbb{R}^N\times(-\infty,0)$. Setting $a := 1-2s$ (and so $a \in (-1,1)$), we define the function
\begin{equation}\label{eq:DEFOFEXTENDEDVERSIONOFU}
\overline{u}(X,t) := \int_{-\infty}^t \int_{\mathbb{R}^N} u(z,\tau)P_y^a(x-z,t-\tau) dz d\tau,
\end{equation}
where $X := (x,y) \in \mathbb{R}^{N+1}_+$ and the Poisson kernel is defined as (see   \eqref{eq:EXPRESSIONFORGAUSSIANINRNHEATEQUATION})
\begin{equation}\label{eq:POISSONKERNELWITHANOTATION}
P_y^a(x,t) := \frac{1}{2^{1-a}\Gamma(\frac{1-a}{2})} G_N(x,t)  \,\frac{y^{1-a}}{t^{1+\frac{1-a}{2}}} e^{-\frac{y^2}{4t}} \qquad (x,y) \in \mathbb{R}_+^{N+1}, \; t > 0.
\end{equation}
As showed in \cite[Proof of Theorem 1]{NystromSande2016:art} (see also \cite[Theorem 1.7]{StingaTorrea2015:art}), we have that $\overline{u}$ is smooth in $\mathbb{R}^{N+1}_+\times(-\infty,0)$ and, furthermore,
\[
\int_0^{\infty} \int_{\mathbb{R}^{N}}P_y^a(x,t)\,dxdt = 1 \quad \text{ for all } y > 0, \qquad P_y^a(x,t) \to \delta_{(0,0)}(x,t) \quad \text{ as } y \to 0^+,
\]
where $\delta_{(0,0)}$ is the Dirac delta at $(0,0)$ (of course, the above limit is intended in the sense of distributions). Consequently, it follows by standard computations
\begin{equation}\label{eq:DirichletCondition}
\overline{u}(x,y,t) \to u(x,t) \quad \text{ as } y \to 0^+,
\end{equation}
point-wise in $B_1\times(-1,0)$. Furthermore, if $\nabla = \nabla_X$ and $\nabla \cdot = \nabla_X\cdot$ denote the gradient and the divergence operator, respectively, a direct computation shows that %
\[
\partial_t P^a - y^{-a} \nabla \cdot (y^a \nabla P^a) = 0 \quad  \text{ in } \mathbb{R}_+^{N+1} \times (-\infty,0)
\]
and so, it is easily seen that
\begin{equation}\label{eq:EquationExtension}
\partial_t \overline{u} - y^{-a} \nabla \cdot(y^a \nabla \overline{u}) = 0 \quad \text{ in } \mathbb{R}_+^{N+1} \times (-\infty,0),
\end{equation}
in the point-wise sense (see  \cite[Proof of Theorem 1]{NystromSande2016:art} again). Finally, we have
\[
- \frac{\overline{u}(X,t) - u(x,t)}{y^{1-a}} = \frac{1}{2^{1-a}\Gamma(\frac{1-a}{2})} \int_{-\infty}^t \int_{\mathbb{R}^N} [u(x,t) - u(z,\tau)] \frac{G_N(x-z,t-\tau)}{(t-\tau)^{1+s}} e^{-\frac{y^2}{4(t-\tau)}} dzd\tau,
\]
for all $(X,t) \in \mathbb{R}^{N+1}_+ \times (-\infty,0)$. Thus, using the fact that $e^{-\frac{y^2}{4(t-\tau)}} \to 1$ as $y \to 0^+$, we can pass to the limit to obtain
\begin{equation}\label{eq:NeumannDerivative}
-\partial_y^a \overline{u}(X,t) := -\lim_{y \to 0^+} \frac{\overline{u}(X,t) - u(x,t)}{y^{1-a}} = \frac{|\Gamma(-s)|}{2^{2s}\Gamma(s)} H^su(x,t),
\end{equation}
point-wise in $\mathbb{R}^N\times(-\infty,0)$. Notice that the limit can be finite or not.

The main consequence of the above discussion is the possibility  to convert the nonlocal problem \eqref{eq:NONLOCALHEPOTENTIAL} into a \emph{local} one, involving the \emph{extended variable} $y \in \mathbb{R}_+$:
\begin{equation}\label{eq:EXTENDEDVERSIONOFSTINGATORREAOPERATORALLSPACE}
\begin{cases}
\partial_t\overline{u} - y^{-a}\nabla\cdot( y^a \nabla \overline{u}) = 0 \quad &\text{in } \mathbb{B}_1^+\times (-1,0), \\
\overline{u} = u \quad &\text{in } B_1 \times\{0\} \times (-1,0) \\
-\partial_y^a \overline{u} = q(x,t)u  \quad &\text{in } B_1 \times\{0\} \times (-1,0),
\end{cases}
\end{equation}
which is nothing more than \eqref{eq:ExtensionEquationLocalForward} (here we simply stress the relation $\overline{u} = u$ in $B_1 \times\{0\} \times (-1,0)$, given in the first boundary equation). We point out that, due to the presence of the constant $|\Gamma(-s)|/(2^{2s}\Gamma(s))$ in front of $H^su$ in \eqref{eq:NeumannDerivative}, the last equation in \eqref{eq:EXTENDEDVERSIONOFSTINGATORREAOPERATORALLSPACE} holds up to a multiplicative constant which can be reabsorbed in the potential $q$. Further, with respect to
\cite[Section 4]{BanGarofalo2017:art}, both the boundary conditions in \eqref{eq:EXTENDEDVERSIONOFSTINGATORREAOPERATORALLSPACE} are \emph{local}, in the sense that are satisfied only in $B_1 \times\{0\} \times (-1,0)$ instead of $\mathbb{R}^{N+1}_+\times\mathbb{R}$.

From now on, we will work in the extended framework, which allows us to bypass the non-locality of our equation. The main drawback is the presence of an additional, though fictitious variable $y \in \RR_+$ that requires the passage to the (nontrivial) limit \eqref{eq:NeumannDerivative}, in order to recover the properties of the solutions to the original nonlocal problem.
\begin{defn}\label{def:WeakSolutions}(Banerjee and Garofalo \cite[Definition 4.1 and 4.3]{BanGarofalo2017:art})
\\
Let $a \in (-1,1)$. We say that the pair $(u,\overline{u})$ is a weak solution to \eqref{eq:EXTENDEDVERSIONOFSTINGATORREAOPERATORALLSPACE} if $\overline{u} \in L^2(0,1;W^{1,2}(\mathbb{B}_1^+,y^a))$, the second equation is intended in the sense of trace and the following relation is satisfied
\[
\int_{t_1}^{t_2} \int_{\mathbb{B}_1^+} y^a \overline{u} \partial_t\eta - \int_{t_1}^{t_2} \int_{\mathbb{B}_1^+} y^a\nabla \overline{u} \cdot \nabla \eta  = \int_{\mathbb{B}_1^+} y^a \overline{u}(t_2)\eta(t_2) - \int_{\mathbb{B}_1^+} y^a \overline{u}(t_1)\eta(t_1) + \int_{t_1}^{t_2} \int_{B_1\times\{0\}} q u \eta\, dxdt,
\]
for all $\eta \in W^{1,2}(\mathbb{B}_1^+\times(t_1,t_2),y^a)$ with compact support in $(\BB_1^+\cup B_1) \times [t_1,t_2]$ and a.e. $-1 < t_1 < t_2 < 0$, where $dX := dxdy$. Furthermore, we say that $\overline{u}$ is a weak global solution to \eqref{eq:ExtensionEquationLocalBackward} (i.e in the whole space $\mathbb{R}^{N+1}_+ \times (-\infty,0)$) if $\overline{u} \in L_{loc}^2(-\infty,0;W_{loc}^{1,2}(\mathbb{R}^{N+1}_+,y^a)$ and the above integral formulation is satisfied on every set of the type $\mathbb{B}_r^+ \times (t_1,t_2)$, with $r > 0$ and  $t_1 < t_2 < 0$.
\end{defn}
Above we have used the notations
\[
\begin{aligned}
B_r(x_0) &:= \{x \in \RR^N: |x-x_0| < r \}, \\
\BB_r(X_0) &:= \{X = (x,y) \in \RR^{N+1}: |x-x_0|^2 + |y-y_0|^2 < r^2 \}, \\
\BB_r^+(X_0) &:= \BB_r(X_0) \cap \{y > 0\},
\end{aligned}
\]
where $X_0 = (x_0,y_0)$ and $r > 0$. We are now in the position to recall some regularity results proved in \cite{BanGarofalo2017:art}. Their proofs are based on the De Giorgi-Nash-Moser approach. We stress that the class of equations \eqref{eq:EXTENDEDVERSIONOFSTINGATORREAOPERATORALLSPACE} belongs is made by a wider group of problems studied in the 80's by Chiarenza and Serapioni  \cite{ChiarenzaSerapioni1985:art} (see also \cite{GutierrezWheeden1991:art,Ishige1999:art}, and \cite{FabesKenigSerapioni1982:art} for the elliptic setting), since the weight $w(y) = y^a$ belongs to the $\mathcal{A}_2$ Muckenhoupt class (see, for instance,  \cite{Baouendi1967:art,Muckenhoupt1972:art}). %
\begin{thm*}(Banerjee and Garofalo \cite[Corollary 4.6]{BanGarofalo2017:art})
\\
Let $a \in (-1,1)$, $q$ satisfying \eqref{eq:ASSUMPTIONSPOTENTIAL} and let $u$ be a solution to \eqref{eq:NONLOCALHEPOTENTIAL}. Then the function $\overline{u}$ defined in \eqref{eq:DEFOFEXTENDEDVERSIONOFU} is a weak solution to the extended problem \eqref{eq:EXTENDEDVERSIONOFSTINGATORREAOPERATORALLSPACE} in $\mathbb{B}_1^+ \times (-1,0)$.
\end{thm*}
This result guarantees the possibility of redirecting the analysis of solutions $u$ to our original problem \eqref{eq:NONLOCALHEPOTENTIAL} to the study of solutions $\overline{u} = \overline{u}(X,t)$ to problem \eqref{eq:EXTENDEDVERSIONOFSTINGATORREAOPERATORALLSPACE} and, subsequently, transferring the properties of $\overline{u}$ to their traces $u$. A fundamental step is to gain regularity from the distributional definition of solutions.

\begin{thm*}(Banerjee and Garofalo \cite[Theorem 5.1, Lemma 5.5 and 5.6, Remark 5.7]{BanGarofalo2017:art})
\\
Let $a \in (-1,1)$, $q$ satisfying \eqref{eq:ASSUMPTIONSPOTENTIAL} and let $\overline{u}$ be a weak solution to \eqref{eq:EXTENDEDVERSIONOFSTINGATORREAOPERATORALLSPACE}. Then:

(i) H\"{o}lder regularity: there exist $\alpha \in (0,1)$ depending only on $a$ and $N$, such that $\overline{u}$ is locally $\alpha$-H\"{o}lder continuous up to the boundary $y = 0$.

(ii) Higher regularity: there exist $\alpha' \in (0,1)$ depending only on $a$ and $N$, such that $\nabla_x \overline{u}$, $\partial_t\overline{u}$, and $y^a\partial_y\overline{u}$ are locally $\alpha'$-H\"{o}lder continuous up to the boundary $y = 0$.
\end{thm*}
In view of the extension theory and the statements above, we give the following definition.
\begin{defn} Fix $s \in (0,1)$. We say that $u = u(x,t)$ is a solution to \eqref{eq:NONLOCALHEPOTENTIAL} if the pair $(u,\overline{u})$ is a weak solution to \eqref{eq:EXTENDEDVERSIONOFSTINGATORREAOPERATORALLSPACE}.
\end{defn}

\begin{rem}
It is wortwhile stressing the two main differences between the notion of solution in \eqref{eq:IntegralRepresentationFHO} and \eqref{eq:FourierDefHs},\eqref{eq:SemigroupDefHs}, of great importance  in the study of \eqref{eq:NONLOCALHEPOTENTIAL}.   The first one is that in  \eqref{eq:IntegralRepresentationFHO} the function $u = u(x,t)$ on which $H^s$ operates is required to be defined only in $\mathbb{R}^N\times(-\infty,0)$, whereas in \eqref{eq:FourierDefHs},\eqref{eq:SemigroupDefHs} we need $u$ to be defined in the whole $\mathbb{R}^{N+1}$. A second remarkable difference is that equation \eqref{eq:NONLOCALHEPOTENTIAL} needs not  to be verified everywhere in $\mathbb{R}^{N+1}$ and holds only in $B_1\times(-1,0)$.  As we have mentioned before, we will focus on solutions to the auxiliary local problem \eqref{eq:EXTENDEDVERSIONOFSTINGATORREAOPERATORALLSPACE}. This choice is relevant in the light of some classical results concerning the same problem for the perturbed plain parabolic Heat operator (i.e. \eqref{eq:NONLOCALHEPOTENTIAL} in the case $s=1$) \cite{Chen1998:art,HanLin1994:art,Lin1991:art}  in $B_1\times(-1,0)$. In this regard, it is worth recalling two phenomena known for local solutions in the local case. First of all, the lack of the unique continuation principle is known (as in the celebrated Tychonoff's  counterexample); for the same reason, moreover, the solutions may have an infinite order of vanishing; in connection with this fact, the results on the structure of the nodal set require an a priori hypothesis of finiteness of the vanishing order.
We shall keep the same perspective here,  differentiating in a substantial way our work from the one of Banerjee and Garofalo \cite{BanGarofalo2017:art} in which the authors studied global solutions to \eqref{eq:NONLOCALHEPOTENTIAL} and a parabolic unique continuation principle holds (see  Section \ref{subsection:KnownResults}). In our setting, it seems to be very hard to exclude a priori that nontrivial solutions may vanish of infinite order: as we will see, this seems to be linked to the problem of defining nonlocal operators as $(-\Delta)^s$, $\partial_t + (-\Delta)^s$, \ldots \, for functions with polynomial growth and local approximations of functions through solutions to nonlocal equations (see  \cite{DipierroSavinVald2017:art1,DipierroSavinVald2016:art1,DipierroSavinVald2016:art2}). In connection with this, it is important to keep in mind that there are  solutions to the extension problem which are not than global solutions to \eqref{eq:NONLOCALHEPOTENTIAL}, since they do not belong to $dom(H^s)$: the main examples are the homogeneous blow-ups we will obtain in Section \ref{Section:BlowUp1}.
\end{rem}
\subsection{Main results} Let us now state our main results. We present here simplified versions of our main theorems, postponing their precise statements to the next sections, after having introduced some auxiliary  matter. As mentioned above, we study solutions having \emph{finite vanishing order} at their nodal points (we will say ``solutions having finite vanishing order in $B_1\times(-1,0)$''), a crucial assumption in our statements. The precise notion of vanishing order we use is technical and it will postponed to  Section \ref{Section:VanOrdMonFor}.

Below, we provide the classical notion of points having finite (or infinite) vanishing order. To do so, we consider a weak solution $(u,\overline{u})$  to \eqref{eq:EXTENDEDVERSIONOFSTINGATORREAOPERATORALLSPACE}, we fix $p_0 = (x_0,0,t_0) \in B_1\times\{0\}\times(-1,0) \in \overline{u}^{-1}(0)$, and define
\[
\QQ_{r,\delta}^+(p_0) := \BB_r^+(X_0) \times (t_0-r^2,t_0-\delta r^2),
\]
for $\delta \in (0,1)$, $r > 0$. We consider the limit
\[
\limsup_{r\to 0^+}\dfrac 1{r^{\sigma}}\sup_{\QQ_{r,\delta}^+(p_0)}|\overline{u}|,
\]
for $\sigma > 0$. Notice that if the limit above is zero for some $\sigma_1 > 0$, then it is zero for any $\sigma < \sigma_1$, while if it is infinite for some $\sigma_2 > 0$, then it is infinite for any $\sigma > \sigma_2$.

Consequently, for any $\delta \in (0,1)$, it is well-defined the number $\sigma_\delta(p_0) \in (0,\infty]$ such that
\begin{equation}\label{eq:vanishing_order}
\limsup_{r\to 0^+}\dfrac 1{r^{\sigma}}\sup_{\QQ_{r,\delta}^+(p_0)}|\overline{u}| =
\begin{cases}
0              \quad &\text{ if } \sigma < \sigma_{\delta} \\
+\infty        \quad &\text{ if } \sigma > \sigma_{\delta}.
\end{cases}
\end{equation}
Notice that $\sigma_{\delta}(p_0) > 0$ for all $\delta > 0$, thanks to the H\"older regularity of the weak solutions $\overline{u}$. Furthermore, since $\QQ_{r,\delta_2}^+(p_0) \subset \QQ_{r,\delta_1}^+(p_0)$ for any $0 < \delta_1 < \delta_2 < 1$, we have
\[
\sup_{\QQ_{r,\delta_2}^+(p_0)}|\overline{u}| \leq \sup_{\QQ_{r,\delta_1}^+(p_0)}|\overline{u}|
\]
and so, it easily follows that the function
\[
\delta \to \sigma_{\delta}(p_0)
\]
is monotone non-decreasing for $\delta \in (0,1)$.

\begin{defn}\label{Def:ClassicalVanOrd}
Let $(u,\overline{u})$ be a weak solution to \eqref{eq:EXTENDEDVERSIONOFSTINGATORREAOPERATORALLSPACE}, $p_0 = (x_0,0,t_0) \in B_1\times\{0\}\times(-1,0) \in \overline{u}^{-1}(0)$. Let us define
\[
\sigma_\ast(p_0) := \lim_{\delta \to 0^+} \sigma_{\delta}(p_0) = \inf_{\delta \in (0,1)} \sigma_{\delta}(p_0),
\]
where $\sigma_{\delta}(p_0)$ is defined as in \eqref{eq:vanishing_order}. We say that $\overline{u}$ (resp. $u$)  vanishes of infinite order at $p_0$ if $\sigma_\ast(p_0) = +\infty$. Otherwise, the vanishing order of $\overline{u}$ (resp. $u$) at $p_0$ is the number $\sigma_\ast(p_0) \in (0,\infty)$.
\end{defn}
Our first main result is concerned with the \emph{asymptotic behaviour} of solutions near their nodal points. i.e., points
\[
p_0=(x_0,t_0) \in \Gamma(u) := u^{-1}(0).
\]
Nothing can be said on at points having infinite vanishing order. On the other hand, it will turn out that the order of vanishing, as a function of $p_0$, is upper semicontinuous, allowing us to focus our attention on open subsets where the vanishing order is bounded above. The main fact here is that any solution to \eqref{eq:NONLOCALHEPOTENTIAL} behaves (near its nodal points) as a suitable \emph{parabolically homogeneous polynomial}, called \emph{tangent map} of $u$
(see   Definition \ref{DEFBLOWUPSET} and Theorem \ref{CONTINUOUSDEPOFBLOWUPLIMITS}) and the vanishing order is characterized. We recall that $u = u(x,t)$ is parabolically homogeneous of degree $2\kappa \geq 0$ if
\[
u(rx,r^2t) = r^{2\kappa}u(x,t),
\]
for all $r > 0$ and all $(x,t) \in \mathbb{R}^{N+1}$.

\begin{thm}\label{Thm:IntroTangentMap} (Existence and uniqueness of the tangent map at nodal points)
Let $s \in (0,1)$ and $q$ satisfying \eqref{eq:ASSUMPTIONSPOTENTIAL}. Let $u = u(x,t)$ be a solution to \eqref{eq:NONLOCALHEPOTENTIAL} having finite vanishing order at $p_0 \in \Gamma(u)$, and set $u_{p_0}(x,t) := u(x+x_0,t+t_0)$.
Then the following assertions hold true:

(i) The vanishing order $\sigma_{\ast}(p_0)$ is a positve integer, and the function $p_0\to \sigma_{\ast}(p_0)$ is upper semicontinuous.

(ii) Set $\sigma_\ast = \sigma_{\ast}(p_0)=2\kappa$. Then there exists a unique nontrivial parabolically $2\kappa$-homogeneous polynomial $\vartheta = \vartheta(x,t)$ such that
\[
u_{p_0}(x,t) = \vartheta_{p_0}(x,t) + o(\|(x,t)\|^{2\kappa}),
\]
as $\|(x,t)\|^2 := |x|^2 + |t| \to 0^+$.

(iii) If $\Gamma_{2\kappa}(u)$ and $\mathfrak{B}_{2\kappa}(u)$ denote the set of nodal points of order $2\kappa > 0$ and the set of parabolically $2\kappa$-homogeneous polynomials, $\kappa \in \NN/2$, respectively (see   Definitions \ref{DEFREGULARSINGULARSET} and \ref{DEFBLOWUPSET} below), then the map
\[
\Gamma_{2\kappa}(u) \to \mathfrak{B}_{2\kappa}(u): p_0 \to \vartheta_{p_0}
\]
is continuous.
\end{thm}
The proof of the theorem relies on the fine analysis of the blow-up family
\[
u_{p_0,r}(x,t) := \frac{u_{p_0}(r x, r^2t)}{r^{2\kappa}},
\]
where $p_0 \in \Gamma(u)$ and $2\kappa > 0$ is the \emph{vanishing order} of $u$ at $p_0$, together with  the analysis of the spectrum of an auxiliary operator of Ornstein-Uhlenbeck type (\eqref{eq:ORNSTEINUHLENEBCKOPERATOR}), and the validity of some monotonicity formulas of Weiss and Monneau type (see Weiss \cite{Weiss1999a:art,Weiss1999b:art} and Monneau \cite{Monneau2003:art,Monneau2009:art}). Given a solution $u$ and a nodal point $p_0$, if finite, the vanishing order $2\kappa$ will be obtained as (twice) the limit of the above mentioned Almgren-Poon quotient of $u$ at $p_0$ (see  Theorem \ref{THEOREMBLOWUP1New}). A remarkable fact is that the admissible values for $\kappa$ are the positive half-integers, i.e. $\kappa \in \mathbb{N}/2$ ($\kappa > 0$), thus making consistent the notation $\Gamma_{2\kappa}(u)$ introduced in the statement of the theorem above: it is the subset of $\Gamma(u)$ in which the limit of the Almgren-Poon quotient equals $\kappa$ and the vanishing order is $2\kappa$. The half-integers $\kappa$ also appear  as eigenvalues of the associated Ornstein-Uhlenbeck  operator.  Due to the importance of such quantities in our analysis, we have decided to keep separate the notations for the vanishing order $2\kappa$ and the limit of the Almgren-Poon quotient $\kappa$, and consequently, $\kappa$ will always be a positive half-integer.

The existence of tangent maps is thoroughly connected with the compactness  of the blow-up families $\{u_{p_0,r}\}_{r > 0}$ (see Theorem \ref{THEOREMBLOWUP1New} again).
Such compactness will be a consequence of a second important result, of independent interest, concerning \emph{uniform H\"older bounds} for locally uniformly bounded solutions to \eqref{eq:ExtensionEquationLocalForward}, which improves the regularity estimates proved in \cite{BanGarofalo2017:art} and, actually, gives an independent proof of the H\"older regularity of weak solutions
to \eqref{eq:ExtensionEquationLocalForward} (see  Theorem \ref{CONVERGENCELINFINITYLOC} for the precise statement). These H\"older bounds  are obtained through a blow-up technique combined with some new \emph{Liouville type results} (see Theorem \ref{CONVERGENCELINFINITYLOC}, Corollary \ref{Corollary:THMECONVERGENCELINFINITYLOCUNIQUENESS} and Lemma \ref{1LIOUVILLETYPETHEOREMFORWEIGHTEDHE}).

\begin{thm}\label{Theorem:UniformHolderBounds}
Let $s \in (0,1)$, $q$ satisfying \eqref{eq:ASSUMPTIONSPOTENTIAL} and let $u = u(x,t)$ be a bounded solution to \eqref{eq:NONLOCALHEPOTENTIAL} in $B_1\times(-1,0)$. Define $\nu_\ast := \min\{1,1-a\}/2$. Then for any $\nu \in (0,\nu_\ast)$, there exists a constant $C = C_\nu > 0$, such that
\[
\| u \|_{\mathcal{C}^{2\nu,\nu}(B_{1/2}\times(-1/2,0))} \leq C \|u \|_{L^{\infty}(B_1\times(-1,0))}.
\]
\end{thm}
%
%_
%
%
%
Finally, we state our main result concerning the \emph{structure} and the \emph{regularity} of the nodal set of solutions to \eqref{eq:NONLOCALHEPOTENTIAL}. In particular, we prove it possesses remarkable \emph{non-degeneracy} and \emph{stratification} properties (see Theorem \ref{THMDIMENSIONALESTIMATENODALSET} and Theorem \ref{STRUCTUREOFSINGULARSET} for the precise statements).

\begin{thm}\label{Theorem:StructureRegularity} (Structure and regularity of the nodal set)
Let $s \in (0,1)$, $q$ satisfying \eqref{eq:ASSUMPTIONSPOTENTIAL} and let $u = u(x,t)$ be a solution to \eqref{eq:NONLOCALHEPOTENTIAL} having finite vanishing order in $B_1\times(-1,0)$. Then the following assertions hold true:

(i) The nodal set $\Gamma(u)$ has parabolic Hausdorff codimension at least one, i.e.
\[
\dim_{\mathcal{\mathcal{P}}}(\Gamma(u)) \leq N+1,
\]
where $\dim_{\mathcal{\mathcal{P}}}(\cdot)$ denotes the parabolic Hausdorff dimension (see Definition \ref{DEFINPARHAUSDIMENSION}).

(ii) $\Gamma(u)$ can be written as the disjoint union of a regular and singular part, denoted by ${\rm Reg}(u)$ and ${\rm Sing}(u)$, respectively. Furthermore:

$\bullet$ ${\rm Reg}(u) = \Gamma(u)\cap\{|\nabla_x u| > 0\}$ is a locally $\mathcal{C}^1$-manifold of parabolic Hausdorff codimension at least one, i.e.
\[
\dim_{\mathcal{\mathcal{P}}}({\rm Reg}(u)) \leq N+1.
\]

$\bullet$ ${\rm Sing}(u) := \Gamma(u)\setminus{\rm Reg}(u)$ can be written as the countable unions
\[
{\rm Sing}(u) = \bigcup_{\substack{\kappa \in \NN/2 \\ \kappa > 0}} \Gamma_{2\kappa}(u) = \bigcup_{\substack{\kappa \in \NN/2 \\ \kappa > 0}} \bigcup_{d = 0,\ldots, N} \Gamma_{2\kappa}^d(u),
\]
and each set $\Gamma_{2\kappa}^d(u)$ is contained in a countable union of $(d+1)$-dimensional space-like $\mathcal{C}^{1,0}$ manifolds for any $d = 0,\ldots,N-1$, whilst $\Gamma_{2\kappa}^N(u)$ in a countable union of $N$-dimensional time-like $\mathcal{C}^1$ manifolds (see  Definitions \ref{DEFOFDIMENSIONSINGULARPOINTS} and \ref{DEFSPACELIKETIMELIKE}, below).
\end{thm}
The proof of this latter theorem relies on the validity of the first two, combining the blow-up analysis  with a parabolic version of the Federer's reduction principle and of the Whitney's extension theorem  (see Theorems \ref{FEDERERTHEOREM} and\ref{PARABOLICWHITNEYEXTENSION}, respectively).

Let us briefly comment on the above statement. First of all, it tells us that the nodal sets of solutions to \eqref{eq:NONLOCALHEPOTENTIAL} having finite vanishing order is ``thin'', in the sense that it has at least codimension one in $\mathbb{R}^{N+1}$ in the parabolic sense (we recall that the parabolic dimension $\dim_{\mathcal{\mathcal{P}}}(\mathbb{R}^{N+1}) = N + 2$). Secondly, it says that $\Gamma(u)$ can be split in two parts: the first one is a locally regular manifold with at least codimension one, whilst the second one possesses a  stratified structure: the main fact here is that we can control the size of the strata $\Gamma_{2\kappa}^{d}(u)$ (see Definition \ref{DEFOFDIMENSIONSINGULARPOINTS} for the precise definition) for each $\kappa \in \mathbb{N}/2$ and $d \in \{0,\ldots,N-1\}$. This kind of estimates are remarkable since, with respect to the local framework, the singular set has \emph{not} higher codimension (see Section \ref{subsection:KnownResults} below): an explicit counterexample has been given in the elliptic case by Tortone \cite{Tortone2018:book} in his PhD dissertation. He proved the existence of an $s$-harmonic function in the interval $(-1,1)$ satisfying $\Gamma(u) = {\rm Sing}(u) = \{0\}$ (${\rm Reg}(u) = \emptyset$). In this paper, we do not provide a corresponding result in the parabolic setting due to the lack of theoretical tools, but we will see how to recover examples of this kind of functions from the blow-up analysis.

\subsection{Known and related results}\label{subsection:KnownResults} The study of the nodal properties of solutions to parabolic equations has a long history. The first results were obtained in the local framework ($s = 1$) for both the elliptic and the parabolic settings by Han and Lin \cite{HanLin1994:art,Lin1991:art}, and Chen \cite{Chen1998:art}. These authors studied solutions $u = u(x,t)$ to the Heat equation (and more general versions of it), satisfying also a \emph{doubling property} (see  formula (1.5) of \cite{HanLin1994:art}), which is a non-degeneracy assumption needed to exclude solutions of infinite vanishing order (see for instance Jones's example \cite{Jones1977:art}). For this class of functions they proved the following estimates on the Hausdorff dimension of the nodal set of $u$ (see for instance \cite{Chen1998:art}):
\begin{equation}\label{eq:CHENESTIMATEFARAWAYSIGMA}
\begin{aligned}
&\dim_{\mathcal{\mathcal{P}}}(\Gamma(u)) \leq N+1 \\
&\dim_{\mathcal{\mathcal{P}}}({\rm Sing}(u)) \leq N \\
&\dim_{\mathcal{H}}(\Gamma_t(u)) \leq N-1,
\end{aligned}
\end{equation}
where $\Gamma_t(u) := \{x \in \RR^N: u(x,t) = 0\}$, $t > 0$. Taking into account the definition and the properties of the parabolic Hausdorff dimension, it follows that $\Gamma(u)$ is composed by a regular part having at least codimension 1 in $\mathbb{R}^{N+1}$ (as in our setting) and a singular part ${\rm Sing}(u)$ of higher codimension (as mentioned above this bound is false in the nonlocal setting). At the same time, Poon \cite{Poon1996:art} proved a \emph{parabolic unique continuation principle} for the Heat equation and some more general parabolic equations. To our knowledge, these are the first results concerning the parabolic framework. Later, various types of parabolic unique continuation properties where investigated with Carleman's estimates and blow-up analysis techniques (\cite{EscauriazaFernandez2003:art,EscauriazaFernandezVess2006:art,FelliPrimo1949:art} and the references therein).

The \emph{nonlocal elliptic} framework was studied by Sire, Terracini and Tortone in \cite{SireTerTor2018:art} (see also \cite{SireTerVita2019:art,Tortone2018:book,Vita2018:book}) where the authors proved the elliptic counterpart of our results. We mention also Fall and Felli \cite{FallFelli2014:art} and Yu \cite{Yu2017:art}, where the authors proved some unique continuation properties in the nonlocal elliptic setting.

Our results can be viewed as the natural follow-up of the work of Banerjee and Garofalo \cite{BanGarofalo2017:art}, where the authors showed a parabolic unique continuation principle for \emph{global} solutions $u = u(x,t)$ to \eqref{eq:NONLOCALHEPOTENTIAL}. Namely, they showed that if $u$ vanishes of infinite order backward in time at some point $(x_0,t_0) \in \mathbb{R}^{N+1}$, in the sense of Definition \ref{Def:ClassicalVanOrd}, then $u = 0$ in $\mathbb{R}^N\times(-\infty,t_0]$ (we recall that $Q_r(x_0,t_0) = B_r(x_0)\times(t_0-r^2,t_0)$). This significant result relies on the validity of an Almgren-Poon monotonicity formula (that the authors proved as well) and the fact that the solutions they consider are \emph{global}, i.e. they belong to the space $\text{dom}(H^s)$ defined in \eqref{eq:NatDomHs}. In our setting, a unique continuation principle is not expected to hold since we consider \emph{local} solutions without any additional \emph{doubling properties requirements}  (\cite{HanLin1994:art,Jones1977:art,Lin1991:art}): for this reason we restrict our analysis to the nodal points in which the vanishing order is \emph{finite} (see Definition \ref{def:VanishingOrder}).

It is worth to stress that a big deal of work has been devoted to the associated obstacle problem
\begin{equation}\label{eq:NonlocParObProb}
\min \{H^su,u - \psi\} = 0  \quad \text{ in } \mathbb{R}^{N+1},
\end{equation}
where $\psi$ can be though as a smooth given function (the obstacle). The analysis of the regular part of the free boundary of solutions as been firstly investigated by Athanasopoulos et al. \cite{AthCaffaMilakis2016:art} and, later, by Banerjee et al. \cite{BanGarDanPetr2019:art} with different techniques, while the singular part was studied by Danielli et al. \cite{DainelliGarofaloPetTo2017:book} when $s = 1/2$ and generalized to the range $s \in (0,1)$ by Banerjee et al. in \cite{BanGarDanPetr2018:art}. Their approach, as ours, is based on the extension theory as developed in  \cite{CaffSil2007:art,NystromSande2016:art,StingaTorrea2009:art,StingaTorrea2015:art} and on the validity of some monotonicity formulas (see \cite{BanGarofalo2017:art,BanGarDanPetr2019:art}) of Almgren-Poon, Monneau and Weiss type.

Finally, we want to mention the works of Caffarelli and Figalli \cite{CafFigalli2013:art} and Barrios et al. \cite{BarriosFigalliRosOton2016:art}, where the authors proved the optimal regularity and the properties of the regular part of the free boundary of solutions to \eqref{eq:NonlocParObProb}, where the operator $H^s$ is replaced by $\partial_t + (-\Delta)^s$: in this setting no monotonicity formulas are known and the analysis of the singular set seems to be a very challenging open problem.

\subsection{Structure of the paper} The paper is organized as follows.

Section \ref{Section:NotPrel} contains the notational background of the work: we introduce the \emph{fundamental solution} to the extended problem in the whole space and a class of Gaussian type spaces.

In Section \ref{Section:VanOrdMonFor} we give the definition of \emph{solutions in strips} and present a different notion of \emph{vanishing order} we need to prove an almost-monotonicity formula of Almgren-Poon type for solutions vanishing of \emph{finite} order. Even though the proof of the monotonicity formula combines well-known techniques employed in \cite{BanGarofalo2017:art,BanGarDanPetr2018:art,DainelliGarofaloPetTo2017:book} with our definition of vanishing order, it has a different nature, since the constants appearing in our formulas depends on the solution itself but not on the blow-up parameter. This will suffice to carry out the blow-up analysis.

In Section \ref{Section:SpectralAnalysis} we develop a spectral analysis for a class of Ornstein-Uhlenbeck type problems and, as a main consequence, we prove the validity of a class of Gaussian-Poincar\'e type inequalities. The description of the spectrum of these eigenvalue problems is related to the characterization of parabolically homogeneous solutions to \eqref{eq:EXTENDEDVERSIONOFSTINGATORREAOPERATORALLSPACE} (posed in the whole space with $q = 0$) and the characterization of the admissible frequencies of the Almgren-Poon quotient.

In Section \ref{Section:BlowUp1} we prove Theorem \ref{THEOREMBLOWUP1New}: we show that any blow-up sequence converges to a homogeneous solution in Gaussian type energy spaces. We classify the blow-ups, the admissible frequencies and the vanishing order. Notice that, compared to \cite{BanGarofalo2017:art,BanGarDanPetr2018:art,DainelliGarofaloPetTo2017:book}, our approach is different: we do not need a priori Gaussian bounds and, using the results proved in Section \ref{Section:SpectralAnalysis}, we prove the uniqueness of the blow-up limit.

In Section \ref{Section:BlowUp2} we show a locally uniform H\"older bound for locally bounded solutions to \eqref{eq:EXTENDEDVERSIONOFSTINGATORREAOPERATORALLSPACE}  (Theorem \ref{CONVERGENCELINFINITYLOC}). As an immediate consequence, Theorem \ref{Theorem:UniformHolderBounds} follows and, further, we obtain that the convergence of the blow-up sequences holds not only in energy spaces but also locally uniformly (actually locally in suitable H\"older spaces). The proof is based on an independent blow-up procedure and the validity of some new Liouville type results for parabolic equations that we prove as well.

In Section \ref{SECTIONNODALSET} we prove Theorem \ref{THMDIMENSIONALESTIMATENODALSET}, Theorem \ref{CONTINUOUSDEPOFBLOWUPLIMITS} and Theorem \ref{STRUCTUREOFSINGULARSET}, from which we deduce the statements of Theorem \ref{Thm:IntroTangentMap} and Theorem \ref{Theorem:StructureRegularity} as main consequences. As mentioned above, we employ the Federer's reduction principle and the Whitney's extension theorem, as well as Weiss and Monneau type monotonicity formulas (see  for instance with \cite[Section 12]{DainelliGarofaloPetTo2017:book}).

Section \ref{Appendix} is a complementary section. We present the blow-up analysis for a class of equations related to \eqref{eq:EXTENDEDVERSIONOFSTINGATORREAOPERATORALLSPACE}: the techniques we employ are very similar to the ones of Sections \ref{Section:BlowUp1} - \ref{Section:BlowUp2}.

\section{Notations and preliminaries: the fundamental solution}\label{Section:NotPrel}

Let us consider the problem
\begin{equation}\label{eq:DEFINITIONPARABOLICEQUATIONPOSITVETIMES}
\begin{cases}
\partial_t\mathcal{G} - y^{-a}\nabla \cdot \left(y^a\nabla \mathcal{G} \right) = 0 \quad &\text{in } \RR_+^{N+1}\times(0,\infty) \\
-\partial_y^a \mathcal{G} = 0 \quad &\text{in } \RR^N\times\{0\}\times(0,\infty),
\end{cases}
\end{equation}
which is the \emph{global} version of \eqref{eq:EXTENDEDVERSIONOFSTINGATORREAOPERATORALLSPACE}, with $q = 0$. If $G_N = G_N(x,t)$ denotes the fundamental solution to the Heat equation (see  \eqref{eq:EXPRESSIONFORGAUSSIANINRNHEATEQUATION}) and
\[
G_{a+1}(y,t) =  \frac{1}{2^a\Gamma(\frac{1+a}{2})} \frac{1}{t^{\frac{1+a}{2}}} e^{-\frac{y^2}{4t}},
\]
we define
\begin{equation}\label{eq:SEPARATIONOFVARIABLESSOLUTION}
\mathcal{G}_a(X,t) = G_N(x,t) G_{a+1}(y,t) = C_{N,a} \frac{1}{t^{\frac{N+a+1}{2}}} e^{-\frac{|X|^2}{4t}},
\end{equation}
where $|X|^2 = |x|^2 + y^2$ and
\[
C_{N,a} := \frac{1}{2^a \Gamma(\frac{1+a}{2})(4\pi)^{N/2}}.
\]
It is straightforward to verify that $\mathcal{G}_a$ is a smooth solution to \eqref{eq:DEFINITIONPARABOLICEQUATIONPOSITVETIMES}, with the following key scaling property
\[
\mathcal{G}_a(rX,r^2t) = r^{-(N+a+1)} \mathcal{G}_a(X,t),
\]
where $rX = (rx,ry)$, for all $r > 0$, and it satisfies also
\begin{equation}\label{eq:GradientPropFundSol}
\nabla_X\mathcal{G}_a(X,t) = - \frac{X}{2t} \, \mathcal{G}_a(X,t).
\end{equation}
We could have proceeded by looking for solutions to \eqref{eq:DEFINITIONPARABOLICEQUATIONPOSITVETIMES} with the \emph{separate variable} form \eqref{eq:SEPARATIONOFVARIABLESSOLUTION}. This procedure easily leads to the expressions for $G_N$ and $G_{a+1}$, and it is based on the decomposition of the operator
\[
y^{-a}\nabla \cdot \left(y^a\nabla \cdot \right) = \Delta_x + \partial_{yy} + \frac{a}{y}\partial_y,
\]
as the sum of the generator of the Heat and the Bessel semigroup, respectively. Similarily, the expression $\mathcal{G}_a$ could have been worked out by looking for solutions to \eqref{eq:DEFINITIONPARABOLICEQUATIONPOSITVETIMES} in \emph{self-similar} form, i.e. solutions satisfying the above scaling property. Finally, we mention that the expression of $\mathcal{G}_a$ was obtain firstly in \cite[Theorem 1.8, Remark 1.9]{StingaTorrea2015:art}, as the function realizing the following identity
\[
-y^{-a} \partial_y \mathcal{G}_{-a}(x,y,t) = P_y^a(x,t) \quad \text{in } \RR_+^{N+1}\times(0,\infty),
\]
where $P_y^a$ is the Poisson kernel defined in \eqref{eq:POISSONKERNELWITHANOTATION} (see also \cite[formula (2.3)]{CaffSil2007:art} in the elliptic case). Our choice of the normalization constant $C_{N,a} > 0$ differs from the one in \cite{StingaTorrea2015:art} and makes the family of functions
\[
X \to y^a \mathcal{G}_a(X,t),
\]
a family of probability densities for any $t > 0$. We thus define the family of probability measures on $\RR_+^{N+1}$
\begin{equation}\label{eq:EXTENDEDGAUSSIANMEASURE}
d\mu_t = d\mu_t(X) := y^a \mathcal{G}_a(X,t)dX,
\end{equation}
shortened, for $t = 1$, to $d\mu := d\mu_1$. For $p_0 = (x_0,t_0) \in B_1\times(0,1)$, we set
\begin{equation}\label{eq:GaussianMeasures}
\begin{aligned}
d\mu_t^{p_0}(X) &:= y^a \mathcal{G}_a(x-x_0,y,t-t_0) \,dxdy \\
d\mu_t^{p_0}(x,0) &:= \mathcal{G}_a(x-x_0,0,t-t_0) \,dx,
\end{aligned}
\end{equation}
where $X \in \mathbb{R}^{N+1}_+$ and $t > t_0$, abbreviating $d\mu^{p_0} = d\mu_1^{p_0}$. In some cases (see Section \ref{Section:SpectralAnalysis})
, we will use the notations
\[
d\mu_x := G_N(x,1) dx = \frac{1}{(4\pi)^{N/2}} e^{-\frac{|x|^2}{4}} dx, \qquad d\mu_y := \frac{1}{2^a \Gamma(\frac{1+a}{2})} y^a e^{-\frac{y^2}{4}} dy,
\]
to denote the marginals of the measure $d\mu$. We conclude the presentation of the fundamental solution by stressing that, since $\partial_y^a\mathcal{G}_a = 0$ in $\RR^N\times\{0\}\times(0,\infty)$, it follows that the function
\begin{equation}\label{eq:FUNDAMENTALSOLUTIONTONONLOCALHEATEQ1}
\widetilde{\mathcal{G}}_a(X,t) := \frac{1}{2} \mathcal{G}_a(X,t),
\end{equation}
is a smooth solution to equation
\begin{equation}\label{eq:DEFINITIONPARABOLICEQUATION1POSITIVETIMES}
\partial_t\mathcal{G} - |y|^{-a}\nabla\cdot( |y|^a \nabla \mathcal{G}) = 0 \quad \text{in } \RR^{N+1}\times(0,\infty).
\end{equation}
We will call it fundamental solution to equation \eqref{eq:DEFINITIONPARABOLICEQUATION1POSITIVETIMES} and, with some abuse of notation, we will denote it again with $\mathcal{G}$. Of course, it satisfies the same scaling properties of the fundamental solution in the half-space $\RR_+^{N+1}$. We will set again
\begin{equation}\label{eq:EXTENDEDGAUSSIANMEASURE1}
d\mu_t = d\mu_t(X) := |y|^a \mathcal{G}_a(X,t)dX,
\end{equation}
which is now a probability measure on $\RR^{N+1}$, for all $t > 0$. Again we will shorten with $d\mu$ instead of $d\mu_1$. Also in this case we will indicate with
\[
d\mu_y := \frac{1}{2^{1+a} \Gamma(\frac{1+a}{2})} |y|^a e^{-\frac{y^2}{4}} dy,
\]
the marginal of the measure $d\mu$ with respect to to  $y > 0$. We have decided not to distinguish between the two fundamental solutions and related probability measures not to exceed in notations. It will be always clear from the context in which framework we will work.

The fundamental solution will play a main role what follows. Here we employ it to define a class of Gaussian spaces which will appear later. For any $t > 0$, we define
\[
L^2(\RR_+^{N+1},d\mu_t) := \left\{V = V(X) \text{ measurable: } \int_{\RR_+^{N+1}} V^2 d\mu_t < +\infty \right\},
\]
where $d\mu_t$ is as in \eqref{eq:EXTENDEDGAUSSIANMEASURE}, endowed with the natural $L^2$ type norm. We will often simplify the notation to
\[
L_{\mu_t}^2 := L^2(\RR_+^{N+1},d\mu_t) \qquad \text{ and } \qquad L_{\mu}^2 := L^2(\RR_+^{N+1},d\mu).
\]
Then, for any $t > 0$, we introduce the family of weighted Sobolev norms
\[
\|V\|_{H^1(\RR_+^{N+1},d\mu_t)}^2 := \int_{\RR_+^{N+1}} V^2 d\mu_t + t \int_{\RR_+^{N+1}} |\nabla V|^2 d\mu_t,
\]
and the Hilbert spaces
\[
H^1(\RR_+^{N+1},d\mu_t) \qquad \text{ and } \qquad H_0^1(\RR_+^{N+1},d\mu_t),
\]
obtained as the closure of the spaces $\mathcal{C}_c^{\infty}(\overline{\RR_+^{N+1}})$ and $\mathcal{C}_c^{\infty}(\RR_+^{N+1})$ with respect to the norm $\|\cdot\|_{H^1(\RR_+^{N+1},d\mu_t)}$ defined above, respectively. Again, we will shorten
\[
\begin{aligned}
H_{\mu_t}^1 &:= H^1(\RR_+^{N+1},d\mu_t), \quad H_{0,\mu_t}^1 := H_0^1(\RR_+^{N+1},d\mu_t) \\
H_{\mu}^1 &:= H^1(\RR_+^{N+1},d\mu), \quad H_{0,\mu}^1 := H_0^1(\RR_+^{N+1},d\mu).
\end{aligned}
\]
Finally, we will need some space-time $L^2$-Sobolev type spaces:
\[
\begin{aligned}
L^2(0,T;L_{\mu_t}^2) &:= \left\{ V = V(X,t) \text{ measurable}: \int_0^T \|V(t)\|_{L_{\mu_t}^2}^2 dt < +\infty \right\} \\
L^2(0,T;H_{\mu_t}^1) &:= \left\{ V = V(X,t) \text{ measurable}: \int_0^T \|V(t)\|_{H_{\mu_t}^1}^2 dt < +\infty \right\} \\
L^2(0,T;H_{0,\mu_t}^1) &:= \left\{ V = V(X,t) \text{ measurable}: \int_0^T \|V(t)\|_{H_{0,\mu_t}^1}^2 dt < +\infty \right\},
\end{aligned}
\]
which are again Hilbert spaces with the natural scalar product and associated norm.

In the same way, we consider the family of spaces
\[
L^2(\RR^{N+1},d\mu_t) := \left\{V = V(X) \text{ measurable: } \int_{\RR^{N+1}} V^2 d\mu_t < +\infty \right\},
\]
where $d\mu_t$ is defined in \eqref{eq:EXTENDEDGAUSSIANMEASURE1}, for all $t > 0$. Again we will use the simplified notations
\[
L_{\mu_t}^2 := L^2(\RR^{N+1},d\mu_t) \qquad \text{ and } \qquad L_{\mu}^2 := L^2(\RR^{N+1},d\mu).
\]
As before, considering the family of weighted Sobolev norms
\[
\|V\|_{H^1(\RR^{N+1},d\mu_t)}^2 := \int_{\RR^{N+1}} V^2 d\mu_t + t \int_{\RR^{N+1}} |\nabla V|^2 d\mu_t,
\]
we can define the family of spaces $H^1(\RR^{N+1},d\mu_t)$ obtained as the closure of the space $\mathcal{C}_c^{\infty}(\RR^{N+1})$ with respect to the above norm. Also in this case we will shorten $H_{\mu_t}^1 := H^1(\RR^{N+1},d\mu_t)$ and $H_{\mu}^1 := H^1(\RR^{N+1},d\mu)$. Finally, the Hilbert spaces $L^2(0,T;L_{\mu_t}^2)$ and $L^2(0,T;H_{\mu_t}^1)$ are defined as above substituting $\RR_+^{N+1}$ with $\RR^{N+1}$.

\section{Solutions in strips, vanishing order and monotonicity formulas}\label{Section:VanOrdMonFor}
\subsection{Solutions in strips} For expository reasons, we introduce the function
\[
\overline{v}(X,t) := \overline{u}(X,-t) \quad \text{ in } \mathbb{R}^{N+1}_+\times(0,1)
\]
which satisfies (in the weak sense) the \emph{backward problem} of \eqref{eq:EXTENDEDVERSIONOFSTINGATORREAOPERATORALLSPACE}:
\begin{equation}\label{eq:ExtensionEquationLocalBackward}
\begin{cases}
\partial_t\overline{v} + y^{-a}\nabla\cdot( y^a \nabla \overline{v}) = 0 \quad &\text{in } \mathbb{B}_1^+\times (0,1), \\
\overline{v} = v \quad &\text{in } B_1 \times\{0\} \times (0,1) \\
-\partial_y^a \overline{v} = q(x,t)v  \quad &\text{in } B_1 \times\{0\} \times (0,1),
\end{cases}
\end{equation}
with $\overline{v}(x,0,t) = u(x,-t) := v(x,t)$. This change of variables does affect the local properties of solutions to \eqref{eq:NONLOCALHEPOTENTIAL}, but allows us to simplify the treatise.

Once problem \eqref{eq:NONLOCALHEPOTENTIAL} is localized by switching to  \eqref{eq:ExtensionEquationLocalBackward}, it is convenient to further localize the function $\overline{v}$, introducing a new one which is defined in the whole space $\mathbb{R}^{N+1}_+\times(0,1)$ and has compact support in space. These solutions are called \emph{solutions in strips} (see  \cite{BanGarDanPetr2018:art}) and are simply obtained by multiplying $\overline{v}$ by a suitable cut-off.

So let us consider a weak solution $\overline{v}$ to \eqref{eq:ExtensionEquationLocalBackward} (see Definition \ref{def:WeakSolutions}). Following \cite[Section 3]{BanGarDanPetr2018:art}, we consider a smooth cut-off function $\zeta = \zeta(X)$ satisfying
\begin{equation}\label{eq:PropCutOff}
\supp(\zeta) \subset \mathbb{B}_1, \qquad \zeta = 1 \quad \text{in } \mathbb{B}_{1/2}, \qquad 0 \leq \zeta \leq 1  \quad \text{in } \mathbb{R}^{N+1},
\end{equation}
with $\zeta(\cdot,y) = \zeta(\cdot,-y)$. We also assume that $\zeta_y = O(|y|)$, $y^{-a} \nabla \cdot (y^a \nabla\zeta) = O(1)$, and $\overline{v}_y\zeta_y = O(1)$ as $y \to 0$ (these last requirements are guaranteed thanks to the regularity of $\zeta$ and $\overline{v}$). We next define the function
\[
W := \zeta(X) \overline{v},
\]
which is well-defined in $\mathbb{R}^{N+1}_+ \times (0,1)$ and has compact support in $(\BB_1^+\cup B_1) \times (0,1)$. Further, setting
\begin{equation}\label{eq:DefRightHS}
F := 2 \nabla \zeta \cdot \nabla \overline{v} + \zeta y^{-a} \nabla \cdot (y^a \nabla\zeta) \overline{v},
\end{equation}
it is easily seen that $F \in L^{\infty}(\mathbb{R}^{N+1}_+\times(0,1))$, and $W$ satisfies
\begin{equation}\label{eq:TruncatedExtensionEquationLocalBackward}
\begin{cases}
\partial_t W + y^{-a} \nabla \cdot(y^a \nabla W) = F \quad &\text{ in } \mathbb{R}_+^{N+1} \times (0,1) \\
-\partial_y^a W = q(x,t)w \quad &\text{ in } \mathbb{R}^N \times \{0\} \times (0,1),
\end{cases}
\end{equation}
in the weak sense (here $w(x,t) := W(x,0,t) = \zeta(x,0) v(x,t)$). Notice that, with respect to \eqref{eq:ExtensionEquationLocalBackward}, the Neumann boundary conditions do not change since $\zeta$ is symmetric with respect to $y$. Clearly, the function $W$ inherits all the regularity properties of $\overline{v}$ (compare with the above mentioned regularity results obtained in \cite{BanGarofalo2017:art}): in particular, using the fact that $\partial_t W$ is a well-defined and H\"older continuous function, we can integrate by parts w.r.t to time to find that $W$ satisfies
\[
\begin{aligned}
\int_{t_1}^{t_2} \int_{\mathbb{R}^{N+1}_+} y^a \partial_t W \eta - \int_{t_1}^{t_2} \int_{\mathbb{R}^{N+1}_+} y^a\nabla W \cdot \nabla \eta = \int_{t_1}^{t_2} \int_{\mathbb{R}^{N+1}_+} y^a F \eta + \int_{t_1}^{t_2} \int_{\mathbb{R}^N\times\{0\}} q w \eta\, dxdt,
\end{aligned}
\]
for all $\eta \in L_{loc}^2(0,1;W^{1,2}(\mathbb{B}_1^+,y^a)))$ with compact support in $(\BB_1^+\cup B_1) \times [t_1,t_2]$ and a.e. $0 < t_1 < t_2 < 1$. Consequently, using the regularity of $W$ again it follows
\begin{equation}\label{eq:StrongSolutionW}
\begin{aligned}
\int_{\mathbb{R}^{N+1}_+} y^a \partial_t W \eta \,dX- \int_{\mathbb{R}^{N+1}_+} y^a\nabla W \cdot \nabla \eta \,dX = \int_{\mathbb{R}^{N+1}_+} y^a F \eta \,dX - \int_{\mathbb{R}^N\times\{0\}} q w \eta\, dx,
\end{aligned}
\end{equation}
for all $\eta \in L_{loc}^2(0,1;W^{1,2}(\mathbb{B}_1^+,y^a)))$ with compact support in $(\BB_1^+\cup B_1)$, and all $t \in (0,1)$.

An equivalent formulation is obtained by testing \eqref{eq:StrongSolutionW} with $\eta \mathcal{G}_a$, where $\mathcal{G}_a$ is the fundamental solution to \eqref{eq:DEFINITIONPARABOLICEQUATIONPOSITVETIMES}. Recalling the property in \eqref{eq:GradientPropFundSol}, it is not difficult to obtain
\begin{equation}\label{eq:StrongSolutionW1}
\begin{aligned}
\int_{\mathbb{R}^{N+1}_+} ( 2t\partial_t W + X \cdot \nabla W  ) \eta \, &y^a\mathcal{G}_a(X,t) dX - 2t\int_{\mathbb{R}^{N+1}_+} \nabla W \cdot \nabla \eta \, y^a\mathcal{G}_a(X,t) dX \\
&= 2t\int_{\mathbb{R}^{N+1}_+} F \eta \, y^a\mathcal{G}_a(X,t) dX - 2t\int_{\mathbb{R}^N\times\{0\}} q w \eta \,\mathcal{G}_a(x,0,t) dx,
\end{aligned}
\end{equation}
for all $\eta \in L_{loc}^2(0,1;W^{1,2}(\mathbb{B}_1^+,y^a)))$ with compact support in $(\BB_1^+\cup B_1) \times (0,1)$, and all $t \in (0,1)$. From now on, we assume that a weak solution $W$ to \eqref{eq:TruncatedExtensionEquationLocalBackward} satisfies \eqref{eq:StrongSolutionW} or \eqref{eq:StrongSolutionW1}, depending on the context we work. Since $W = \overline{v}$ in $\mathbb{B}_{1/2}^+\times(0,1)$, we will subsequently transfer the information on $W$ to the function $\overline{v}$.

Finally, we introduce the function
\begin{equation}\label{eq:RescaledOfW}
\widetilde{W}(X,t) = W(\sqrt{t}X,t),
\end{equation}
which is just the solution $W$ to \eqref{eq:StrongSolutionW}-\eqref{eq:StrongSolutionW1}, computed at the parabolic variables $(\sqrt{t}X,t)$. Passing to the function $\widetilde{W}$ will be important both in the blow-up and in the classification of the blow-up limits. An easy computation shows it is a weak solution to
\begin{equation}\label{eq:RescaledTruncatedExtensionEquationLocalBackward}
\begin{cases}
t\partial_t \widetilde{W} + \mathcal{O}_a\widetilde{W} = t\widetilde{F} \quad &\text{ in } \mathbb{R}_+^{N+1} \times (0,1) \\
-\partial_y^a \widetilde{W} = t^{\frac{1-a}{2}} \widetilde{q}(x,t)\widetilde{w} \quad &\text{ in } \mathbb{R}^N \times \{0\} \times (0,1),
\end{cases}
\end{equation}
where we have set $\widetilde{w}(x,t) = \widetilde{W}(x,0,t) = w(\sqrt{t}x,t)$, $\widetilde{q}(x,t) = q(\sqrt{t}x,t)$, $\widetilde{F}(x,y,t) = F(\sqrt{t}X,t)$, and
\begin{equation}\label{eq:ORNSTEINUHLENEBCKOPERATOR}
\mathcal{O}_a\widetilde{W} = \frac{1}{y^a\mathcal{G}_a(X,1)} \nabla \cdot \left( y^a\mathcal{G}_a(X,1) \nabla \widetilde{W} \right)
\end{equation}
is an Ornstein-Uhlenbeck type operator. Such operator plays an important role in the blow-up classification and we devote the entire Section \ref{Section:SpectralAnalysis} to the analysis of its spectral properties.

\subsection{Almgren-Poon and Weiss functions} In what follows, we introduce some relevant quantities needed in the definition of the vanishing order and in the blow-up procedure. To do so, we basically follow the notations in \cite{BanGarofalo2017:art,BanGarDanPetr2018:art}.

Let $W$ be a weak solution to \eqref{eq:TruncatedExtensionEquationLocalBackward} and $p_0 = (x_0,t_0) \in B_{1/2}\times(0,1)$. We will also write $p_0 = (x_0,0,t_0) \in B_{1/2}\times\{0\}\times(0,1)$ when it is useful to imagine it in the extended space, and $O := (0,0,0)$. Keeping in mind \eqref{eq:GaussianMeasures}, we consider the functions (see also  \cite{BanGarofalo2017:art,BanGarDanPetr2018:art})
\begin{equation}\label{eq:TimeEnergies}
\begin{aligned}
h(W,p_0,t) &:= \int_{\mathbb{R}^{N+1}_+} W^2 d\mu_t^{p_0} \\
d_0(W,p_0,t) &:= (t-t_0)\int_{\mathbb{R}^{N+1}_+} |\nabla W|^2 d\mu_t^{p_0} \\
d(W,p_0,t) &:= d_0(W,p_0,t) - (t-t_0) \int_{\mathbb{R}^N\times\{0\}} q w^2 \, d\mu_t^{p_0}(x,0) \\
i(W,p_0,t) &:= \frac{1}{2} \int_{\mathbb{R}^{N+1}_+} W Z_{p_0}W d\mu_t^{p_0},
\end{aligned}
\end{equation}
where we have set
\begin{equation}\label{eq:TimeHomogeneityFactors}
\begin{aligned}
Z_{p_0} W &:= 2(t-t_0)\partial_t W + (x-x_0,y) \cdot \nabla W \\
Z W &:= 2t\partial_t W + X \cdot \nabla W.
\end{aligned}
\end{equation}
Notice that all the quantities are well-defined for $t \in (t_0,1)$ thanks to the regularity of $W$ (see  \cite[Section 5]{BanGarofalo2017:art}). Further, defining
\[
\begin{aligned}
W_{p_0}(x,y,t) &:= W(x+x_0,y,t+t_0) \\
w_{p_0}(x,t) &:= w(x+x_0,t+t_0),
\end{aligned}
\]
it is not difficult to verify that $h(W,p_0,t_0 + t) = h(W_{p_0},O,t)$, $d_0(W,p_0,t_0 + t) = d_0(W_{p_0},O,t)$, $i(W,p_0,t_0 + t) = i(W_{p_0},O,t)$ as well as
\[
d(W,p_0,t+t_0) := d_0(W,p_0,t) - t \int_{\mathbb{R}^N\times\{0\}} q_{p_0} w_{p_0}^2 \, d\mu_t(x,0).
\]
Moreover, testing the equation of $W_{p_0}$ (i.e. \eqref{eq:StrongSolutionW}) with $ty^a W_{p_0} \mathcal{G}_a$, integrating by parts and using \eqref{eq:GradientPropFundSol}, we easily find that $d$ and $i$ are linked by the following relation
\[
i(W,p_0,t) = d(W,p_0,t) + (t-t_0) \int_{\mathbb{R}^{N+1}_+} W F \,d\mu_t^{p_0},
\]
for all $t \in (t_0,1)$ (see  \cite[Lemma 6.1]{BanGarofalo2017:art} and \cite[Lemma 4.1]{BanGarDanPetr2018:art}). All these relations are very useful to simplify the exposition: in the proofs, we will always reduce to the case $p_0 = O$ and the general case is obtained by translation (it is enough to consider $W_{p_0}$ instead of $W$). Finally, we mention that we will always abbreviate $h(W,O,t) = h(W,t)$, $d(W,O,t) = d(W,t)$ and so on.

Following \cite{BanGarofalo2017:art,BanGarDanPetr2018:art,BanGarDanPetr2019:art}, we will prevalently work with the following averaged versions of the functions $h$, $d_0$, $d$ and $i$:
\begin{equation}\label{eq:AveragedEnergies}
\begin{aligned}
H(W,p_0,r) &:= \frac{1}{r^2} \int_{t_0}^{t_0+r^2} h(W,p_0,t) \, dt \\
D_0(W,p_0,r) &:= \frac{1}{r^2} \int_{t_0}^{t_0+r^2} d_0(W,p_0,t) \, dt \\
D(W,p_0,r) &:= \frac{1}{r^2} \int_{t_0}^{t_0+r^2} d(W,p_0,t) \, dt \\
I(W,p_0,r) &:= \frac{1}{r^2} \int_{t_0}^{t_0+r^2} i(W,p_0,t) \, dt,
\end{aligned}
\end{equation}
for $r \in (0,1)$. Notice that we have the relation
\[
I(W,p_0,r) = D(W,p_0,r) + \frac{1}{r^2} \int_{t_0}^{t_0+r^2} (t-t_0) \int_{\mathbb{R}^{N+1}_+} W F \, d\mu_t ^{p_0} dt
\]
and, similarly to the non-averaged case, $H(W,p_0,r) = H(W_{p_0},O,r)$, $D_0(W,p_0,r) = D_0(W_{p_0},O,r)$, $D(W,p_0,r) = D(W_{p_0},O,r)$, and $I(W,p_0,r) = I(W_{p_0},O,r)$. Also in this case, we will employ the abbreviations: $H(W,O,r) = H(W,r)$, $D(W,O,r) = D(W,r)$ and so on.

\begin{defn} (Almgren-Poon quotients)
We define the Almgren-Poon quotients of the function $W$ at $p_0 \in B_{1/2}\times\{0\}\times(0,1)$ as:
\begin{equation}\label{eq:AlmgrenPoonQuotients}
N_0(W,p_0,r) := \frac{D_0(W,p_0,r)}{H(W,p_0,r)}, \qquad N_D(W,p_0,r) := \frac{D(W,p_0,r)}{H(W,p_0,r)}, \qquad  N_I(W,p_0,r) := \frac{I(W,p_0,r)}{H(W,p_0,r)},
\end{equation}
for all $r \in (0,1)$ such that $H(W,p_0,r) > 0$.
\end{defn}
\begin{defn} (Weiss energy)
We define the Weiss energy (see  \cite[Theorem 10.1]{BanGarDanPetr2018:art}) of the function $W$ at $p_0 \in B_{1/2}\times\{0\}\times(0,1)$ as:
\begin{equation}\label{eq:AveragedWeiss}
\mathcal{W}_{\sigma}(W,p_0,r) := \frac{1}{r^{2\sigma}} \left[ D(W,p_0,r) - \frac{\sigma}{2} H(W,p_0,r) \right], \quad r \in (0,1),
\end{equation}
for any fixed $\sigma > 0$. Finally, we introduce the quantity
\begin{equation}\label{eq:AveragedTau}
\mathcal{T}_{\sigma}(W,p_0,r) := \frac{H(W,p_0,r)}{r^{2\sigma}}, \quad r \in (0,1),
\end{equation}
for any fixed $\sigma > 0$.
\end{defn}
As above we will abbreviate $\mathcal{W}_{\sigma}(W,O,r) = \mathcal{W}_{\sigma}(W,r)$, $\mathcal{T}_{\sigma}(W,O,r) = \mathcal{T}_{\sigma}(W,r)$, $N_0(W,O,r) = N_0(W,r)$ and so on.
\subsection{Definition of vanishing order} The precise formulation of the vanishing order needs some introductory/technical material that we present in a series of lemmas.
\begin{lem}\label{Lemma:LimitTauSigmaSmall}
	Let $W$ be a weak solution to \eqref{eq:TruncatedExtensionEquationLocalBackward}, with $q$ satisfying \eqref{eq:ASSUMPTIONSPOTENTIAL} and $p_0 \in B_{1/2}\times\{0\}\times(0,1)$. Then there exists $\sigma_0 > 0$ such that
	\[
	\lim_{r \to 0^+} \mathcal{T}_{\sigma}(W,p_0,r) := \mathcal{T}_{\sigma}(W,p_0,0^+) = 0,
	\]
	for any fixed $\si \in (0,\si_0)$.
\end{lem}
\emph{Proof.} Assume for simplicity $p_0 = O$. To show our claim it is enough to verify that
\begin{equation}\label{eq:HeightHolderBound}
H(W,r) \leq C r^{2\sigma_0},
\end{equation}
for some $C > 0$, $\sigma_0 > 0$ as $r \to 0^+$. So, recalling that $W$ is supported in $\mathbb{B}_1^+$, it follows
\[
\begin{aligned}
H(W,r) &= \frac{1}{r^2}\int_0^{r^2} \int_{\mathbb{B}_1^+} W^2(X,t) \,d\mu_t dt \leq \frac{C}{r^2}\int_0^{r^2} \int_{\mathbb{B}_1^+} \left( |X|^2 + t \right)^{\sigma_0} \,d\mu_t dt \\
& = \int_0^1 \int_{\mathbb{B}_{1/r}^+} \left( |rX|^2 + r^2t \right)^{\sigma_0} \,d\mu_t dt \leq C r^{2\sigma_0} \int_0^1 \int_{\mathbb{R}^{N+1}_+} \left( |X|^2 + t \right)^{\sigma_0} \,d\mu_t dt \\
& = C r^{2\sigma_0} \int_0^1 t^{\sigma_0} \left[\int_{\mathbb{R}^{N+1}_+} \left( |X|^2 + 1 \right)^{\sigma_0} \,d\mu \right] dt \leq C r^{2\sigma_0},
\end{aligned}
\]
where, in the first inequality, we have used that $W$ is (parabolically) $\sigma_0$-H\"older continuous in $\BB_1^+$, for some $\sigma_0 > 0$. $\Box$
\begin{lem}\label{Lemma:PrelDefVanOrder}
	Let $W$ be a weak solution to \eqref{eq:TruncatedExtensionEquationLocalBackward}, with $q$ satisfying \eqref{eq:ASSUMPTIONSPOTENTIAL} and $p_0 \in B_{1/2}\times\{0\}\times(0,1)$. Then if there exists $\sigma_1 > 0$ such that
	\[
	\liminf_{r \to 0^+} \mathcal{T}_{\sigma_1}(W,p_0,r) > 0,
	\]
	it must be
	\[
	\lim_{r \to 0^+} \mathcal{T}_{\sigma}(W,p_0,r) = + \infty,
	\]
	for any $\si > \si_1$.
\end{lem}
\emph{Proof.} Assume $p_0 = O$. Since $r \in (0,1)$, we immediately see that the function
\[
\sigma \to \mathcal{T}_{\sigma}(W,r)
\]
is strictly increasing on $(0,+\infty)$, uniformly w.r.t. $r \in (0,1)$. Consequently,
\[
\lim_{r \to 0^+} \frac{H(W,r)}{r^{2\si}} > \liminf_{r \to 0^+} \frac{H(W,r)}{r^{2\si_1}} \lim_{r \to 0^+} \frac{1}{r^{2(\si-\si_1)}} = + \infty,
\]
since $\si > \si_1$, and our claim follows. $\Box$
\begin{lem}\label{Lemma:PrelDefVanOrderWeiss}
	Let $W$ be a weak solution to \eqref{eq:TruncatedExtensionEquationLocalBackward}, with $q$ satisfying \eqref{eq:ASSUMPTIONSPOTENTIAL} and $p_0 \in B_{1/2}\times\{0\}\times(0,1)$. Then if there exists $\sigma_1 > 0$ such that
	\[
	\limsup_{r \to 0^+} \mathcal{W}_{\sigma_1}(W,p_0,r) < 0,
	\]
	it must be
	\[
	\lim_{r \to 0^+} \mathcal{W}_{\sigma}(W,p_0,r) = \mathcal{W}_{\sigma}(W,p_0,0^+) = - \infty,
	\]
	for any $\si > \si_1$.
\end{lem}
\emph{Proof.} Assume $p_0 = O$. Thanks to our assumption, we can assume the existence of $r_0 \in (0,1)$, such that $\mathcal{W}_{\sigma_1}(W,r) < 0$ for all $r \in (0,r_0)$. Further, it easy to check that
\[
\mathcal{W}_{\sigma}(W,r) < \frac{1}{r^{2(\sigma-\sigma_1)}} \mathcal{W}_{\sigma_1}(W,r)
\]
for all $\sigma > \sigma_1$ and all $r \in (0,r_0)$. Consequently,
\[
\lim_{r \to 0^+} \mathcal{W}_{\sigma}(W,r) < \limsup_{r \to 0^+}\mathcal{W}_{\sigma_1}(W,r)\lim_{r \to 0^+}\frac{1}{r^{2(\sigma-\sigma_1)}} = -\infty,
\]
and our claim follows. $\Box$

\bigskip

The following lemma establishes the connections between the limits $\mathcal{T}_{\sigma}(W,p_0,0^+)$ and $\mathcal{W}_{\sigma}(W,p_0,0^+)$, and it is crucial in the definition of the vanishing order.
\begin{lem}\label{Lemma:PropWeissTau}
	Let $W$ be a weak solution to \eqref{eq:TruncatedExtensionEquationLocalBackward}, with $q$ satisfying \eqref{eq:ASSUMPTIONSPOTENTIAL} and $p_0 \in \Gamma(W)$. Let us define
	\[
	\overline{\sigma}_{\ast} := \inf \left\{\sigma > 0: \limsup_{r \to 0^+} \mathcal{W}_{\sigma}(W,p_0,r) < 0 \right\} \in[0,+\infty].
	\]
	Then the following statements hold true:

	(i) $\overline{\sigma}_{\ast} > 0$. Furthermore, if $\overline{\sigma}_{\ast} \in (0,+\infty)$, then
	\[
	\overline{\sigma}_{\ast} = \inf \left\{\sigma > 0: \mathcal{W}_{\sigma}(W,p_0,0^+) = -\infty \right\}.
	\]

	(ii) If $\overline{\sigma}_{\ast} \in (0,+\infty)$, then
	\[
	\liminf_{r \to 0^+} \mathcal{W}_{\sigma}(W,p_0,r) = 0
	\]
	for all $\sigma \in (0,\overline{\sigma}_{\ast}]$.

	(iii) If $\overline{\sigma}_{\ast} \in (0,+\infty)$, then
	\[
	\overline{\sigma}_{\ast} = \inf \{\sigma > 0: \mathcal{T}_{\sigma}(W,p_0,0^+) = +\infty \} = \sup \{\sigma > 0: \mathcal{T}_{\sigma}(W,p_0,0^+) = 0 \}.
	\]
\end{lem}
\emph{Proof.} As usual, we fix $p_0 = O$. In view of Lemma \ref{Lemma:PrelDefVanOrderWeiss}, we obtain that if $\overline{\sigma}_{\ast} = 0$, then $\mathcal{W}_{\sigma}(W,0^+) = -\infty$ for all $\sigma > 0$. On the other hand, from Lemma \ref{Lemma:LimitTauSigmaSmall}, we have
\[
\mathcal{T}_{\sigma}(W,0^+) = 0
\]
if $\sigma > 0$ is small enough. Now, let us re-write $\mathcal{W}_{\sigma}(W,r)$ as
\begin{equation}\label{eq:SecondFormWeiss}
\mathcal{W}_{\sigma}(W,r) = \mathcal{T}_{\sigma}(W,r) \left[ N_I(W,r) - \frac{\sigma}{2} \right] - \frac{1}{r^{2\si +2}} \int_0^{r^2} t \int_{\mathbb{R}^{N+1}_+} FW \, d\mu_t dt,
\end{equation}
where we have used \eqref{eq:AveragedEnergies} and notice that $N_I$ is bounded from below (Lemma \ref{Lemma:BoundBelowNI}). If the second term in the r.h.s. is uniformly bounded, we get a contradiction with the assumption $\mathcal{W}_{\sigma}(W,0^+) = -\infty$ for all $\sigma > 0$ and the fact that $\mathcal{T}_{\sigma}(W,0^+) = 0$ for small $\sigma > 0$. Now, recalling the definition of $F$ (see  \eqref{eq:DefRightHS}), we immediately see that $F\not\equiv 0$ only  in $\mathbb{B}_1^+ \setminus\mathbb{B}_{1/2}^+$ only and, furthermore,
\[
|F| \leq C,
\]
for some constant $C > 0$ depending on the solution to \eqref{eq:ExtensionEquationLocalBackward}. Consequently, since $1/4 \leq |x|^2 + y^2 \leq 1$ and the weight $y^a$ is integrable near $y = 0$, we deduce
\[
\begin{aligned}
\int_0^{r^2} t^2 \int_{\mathbb{R}^{N+1}_+} F^2 d\mu_t dt &= \int_0^{r^2} t^2 \int_{\mathbb{B}_1^+ \setminus\mathbb{B}_{1/2}^+} F^2 d\mu_t dt \leq  C r^2 \int_0^{r^2} \int_{\mathbb{B}_1^+ \setminus\mathbb{B}_{1/2}^+}  y^a t^{-\frac{N+a+1}{2}} e^{-\frac{|x|^2+y^2}{4t}}dxdy dt \\
& \leq C r^2 \int_0^{r^2} t^{-\frac{N+a+1}{2}} e^{-\frac{1}{16t}} dt,
\end{aligned}
\]
for some new constant $C > 0$ depending on $N$, $a$, and the solution to \eqref{eq:ExtensionEquationLocalBackward}. Now, noticing that the function
\[
t \to t^{-\frac{N+a+1}{2}} e^{-\frac{1}{16t}},
\]
is decreasing for $t \in (0,t_{N,a})$, for some $t_{N,a} \in (0,1)$ depending only on $N$ and $a$, we obtain
\begin{equation}\label{eq:ExpBoundIntF}
\int_0^{r^2} t^2 \int_{\mathbb{R}^{N+1}_+} F^2 d\mu_t dt  \leq C r^{\sigma_0} e^{-\frac{1}{16r^2}},
\end{equation}
for all $r \in (0,r_0)$, for some $r_0 \in (0,1)$ and $\sigma_0 \in \mathbb{R}$, depending only on $N$ and $a$. Consequently, using that $H(W,r)$ is bounded and the H\"older inequality it follows
\begin{equation}\label{eq:ExpBoundIntFW}
\left| \int_0^{r^2} t \int_{\mathbb{R}^{N+1}_+} FW \, d\mu_t dt \right| \leq C r^{\sigma_0} e^{-\frac{1}{32r^2}},
\end{equation}
for all $r \in (0,r_0)$ and some new $\sigma_0 \in \mathbb{R}$. For what established above, it follows $\overline{\sigma}_{\ast} > 0$. The fact that $\overline{\sigma}_{\ast} = \inf \{\sigma > 0: \mathcal{W}_{\sigma}(W,p_0,0^+) = -\infty \}$ follows by Lemma \ref{Lemma:PrelDefVanOrderWeiss}.

Now let us prove part (ii). Assume $\overline{\sigma}_{\ast} \in (0,+\infty)$ and assume there exists $\sigma \in (0,\overline{\sigma}_{\ast}]$ such that $\liminf_{r \to 0^+} \mathcal{W}_{\sigma}(W,r) > 0$. Then $\mathcal{W}_{\sigma}(W,r) \geq l_0$, for some $l_0 > 0$ and small $r$'s. On the other hand, standard computations give
\begin{equation}\label{eq:DerivativeTer}
\frac{d}{dr} \mathcal{T}_{\sigma}(W,r) = \frac{4}{r} \mathcal{W}_{\sigma}(W,r) + \frac{4}{r^{2\sigma+3}} \int_0^{r^2} t \int_{\mathbb{R}^{N+1}_+} F W \, d\mu_t dt.
\end{equation}
From \eqref{eq:ExpBoundIntFW}, we know that the second term in the r.h.s. of \eqref{eq:DerivativeTer} goes to zero (exponentially fast) as $r \to 0^+$, while, for some small $r_0$, we have
\[
\frac{d}{dr} \mathcal{T}_{\sigma}(W,r) \geq \frac{4l_0}{r}  + \frac{4}{r^{2\sigma+3}} \int_0^{r^2} t \int_{\mathbb{R}^{N+1}_+} F W \, d\mu_t dt,
\]
for all $r \in (0,r_0)$. Integrating between $r$ and $r_0$, rearranging terms and using \eqref{eq:ExpBoundIntFW}, we pass to the limit as $r \to 0$ to deduce $\mathcal{T}_{\sigma}(W,0^+) = -\infty$, in contradiction with the non negativity of $\mathcal{T}_\sigma$.

Proof of part (iii). Using that $\mathcal{W}_{\sigma}(W,0^+) = -\infty$ for all $\sigma > \overline{\sigma}_{\ast}$  and \eqref{eq:SecondFormWeiss}, it is easy to see that $\mathcal{T}_{\sigma}(W,0^+) = +\infty$ for all $\sigma > \overline{\sigma}_{\ast}$ (we use again that $N_I$ is bounded from below). Consequently, in view of Lemma \ref{Lemma:LimitTauSigmaSmall} and Lemma \ref{Lemma:PrelDefVanOrder} (and the strict monotonicity of the function $\sigma \to \mathcal{T}_{\sigma}(W,r)$), the number
\[
\widetilde{\si}_{\ast} = \inf \{\sigma > 0: \mathcal{T}_{\sigma}(W,p_0,0^+) = +\infty \} = \sup \{\sigma > 0: \mathcal{T}_{\sigma}(W,p_0,0^+) = 0 \},
\]
is well-defined and satisfies $0 < \widetilde{\si}_{\ast} \leq \overline{\sigma}_{\ast}$. If by contradiction $\widetilde{\si}_{\ast} < \overline{\sigma}_{\ast}$, then for any $\si \in (\widetilde{\si}_{\ast},\overline{\sigma}_{\ast})$, we have
\[
0 = \liminf_{r \to 0^+} \mathcal{W}_{\si}(W,r) = \lim_{r \to 0^+} \mathcal{T}_{\si}(W,r) \liminf_{r \to 0^+} \left[N_D(W,r) - \frac{\si}{2}\right].
\]
Consequently, since $\mathcal{T}_{\si}(W,0^+) = +\infty$ by definition ($\si > \widetilde{\si}_{\ast}$), it must be
\[
\liminf_{r \to 0^+} N_D(W,r) = \frac{\si}{2},
\]
for all $\si \in (\widetilde{\si}_{\ast},\overline{\sigma}_{\ast})$, which is impossible since $N_I$ does not depend on $\si$. $\Box$

\begin{lem}\label{Lemma:ExistenceLimitNegPartWeiss}
	Let $W$ be a weak solution to \eqref{eq:TruncatedExtensionEquationLocalBackward}, with $q$ satisfying \eqref{eq:ASSUMPTIONSPOTENTIAL} and $p_0 \in B_{1/2}\times\{0\}\times(0,1)$. Then if there exists $\si_0 > 0$, such that
	\begin{equation}\label{eq:boundsuT}
	\limsup_{r \to 0^+} \mathcal{T}_{\si_0}(W,p_0,r) < +\infty,
	\end{equation}
	then both limits
	\begin{equation}\label{eq:boundsuWeT}
	\lim_{r \to 0^+} \mathcal{W}_{\sigma}^{-}(W,p_0,r) \quad \text{ and } \quad \lim_{r \to 0^+} \mathcal{T}_{\si}(W,p_0,r)
	\end{equation}
	exist for all $\si \in (\si_0,\si_0 + (1-a)/2)$. The symbol $\mathcal{W}_{\sigma}^{-}$ denotes the negative part of $\mathcal{W}_{\sigma}$.
\end{lem}
\emph{Proof.} Assume $p_0 = O$. Standard computations show that
\[
\begin{aligned}
r^{2\si + 3} \frac{d}{dr} \mathcal{W}_{\si}(W,r) &= \int_0^{r^2}\int_{\mathbb{R}_+^{N+1}} \left( ZW - \si W - tF \right)^2 d\mu_t dt - \int_0^{r^2} t^2 \int_{\mathbb{R}^{N+1}_+} F^2 d\mu_t dt \\
& \quad - \int_0^{r^2} t \int_{\mathbb{R}^N\times\{0\}} [(1-a)q + 2t\partial_t q + x \cdot \nabla q] w^2 \,d\mu_t(x,0)dt \\
& \geq -C r^{\overline{\si}} e^{-\frac{1}{16r^2}} - C\int_0^{r^2} t \int_{\mathbb{R}^N\times\{0\}} w^2 \,d\mu_t(x,0)dt,
\end{aligned}
\]
for some suitable constants $\overline{\si} \in \RR$ and $C > 0$, thanks to \eqref{eq:ExpBoundIntF} and the assumptions \eqref{eq:ASSUMPTIONSPOTENTIAL} on the potential $q$. Now, from the definition of the Weiss function it immediately follows that
\[
\mathcal{W}_\si(W,r) < 0 \quad \text{ if and only if } \quad  D(W,r) < \frac{\si}{2} H(W,r).
\]
Consequently, applying Lemma \ref{Lemma:TraceIneqL2Norm}, we obtain
\[
\begin{aligned}
\frac{1}{r^2} \int_0^{r^2} t \int_{\mathbb{R}^N\times\{0\}} w^2 \,d\mu_t(x,0)dt &\leq \frac{C_0}{r^2} \int_0^{r^2}  t^{\frac{1-a}{2}} \left[ d(W,t) + h(W,t) \right]dt \\
& \leq \frac{C_0}{r^2} \int_{(0,r^2)\cap\{d < 0\} } t^{\frac{1-a}{2}} \left[ d(W,t) + h(W,t) \right]dt \\
& \quad + \frac{C_0}{r^2} \int_{(0,r^2)\cap\{d \geq 0\}} t^{\frac{1-a}{2}} \left[ d(W,t) + h(W,t) \right]dt \\
& \leq  C_0r^{1-a} [D(W,r) + H(W,r)] \leq C_0 r^{1-a} H(W,r),
\end{aligned}
\]
for a new constant $C_0 > 0$ depending on $\si$, whenever $\mathcal{W}_\si(W,r) < 0$. Thus, passing to the negative part of $\mathcal{W}_\si$, it follows
\begin{equation}\label{eq:DiffIneqNegPartWeiss}
\frac{d}{dr} \mathcal{W}_{\si}^-(W,r) \leq C r^{\overline{\si}} e^{-\frac{1}{16r^2}} + C \frac{H(W,r)}{r^{2\si + a}} \leq \frac{C}{r^{2(\si - \si_0) + a}},
\end{equation}
for all $r\in(0,r_0)$ and for some suitable $r_0 \in (0,1)$, thanks to our assumption $\mathcal{T}_{\si_0}(W,r)  < +\infty$, for small $r$'s. Now, since the last term in the above chain of inequalities is integrable at $r = 0$ if and only if $\si < \si_0 + (1-a)/2$, we can integrate between any $0 < r_1 < r_2 < r_0$, and deduce that the function
\[
r \to \mathcal{W}_\si^-(W,r) - C' r^{1-2(\si - \si_0) - a}
\]
is non-increasing in $(0,r_0)$, where $C' = C/[1-2(\si - \si_0) - a]$. Thus $\mathcal{W}_\si^-(W,r)$ has a limit as $r \to 0^+$ and the first part of our claim follows.

To complete the proof, we assume first that $\mathcal{W}_\si^-(W,0^+) > 0$, i.e. $\mathcal{W}_\si(W,0^+) < 0$. In this case, it must be $\mathcal{W}_\si(W,r) < -\varepsilon$ for small $r$'s and a suitable $\varepsilon > 0$, and so, using \eqref{eq:ExpBoundIntFW} and \eqref{eq:DerivativeTer}, we deduce
\[
\frac{d}{dr} \mathcal{T}_{\sigma}(W,r) < -\frac{4C\varepsilon}{r},
\]
for some $C > 0$ and all small $r$'s. It thus follows that $r \to \mathcal{T}_\si(W,r)$ is decreasing and so it has a limit as $r \to 0^+$.

On the other hand, if $\mathcal{W}_\si^-(W,0^+) = 0$, we can integrate \eqref{eq:DiffIneqNegPartWeiss} to obtain
\[
\mathcal{W}_\si(W,r) \geq - C' r^{1-2(\si - \si_0) - a},
\]
and so, using \eqref{eq:DerivativeTer} again, we find
\[
\frac{d}{dr} \mathcal{T}_{\sigma}(W,r) \geq - \frac{C}{r^{2(\si-\si_0)+a}} - Cr^{\overline{\si}} e^{-\frac{1}{32r^2}},
\]
for some new constants $C > 0$ and $\overline{\si} \in \RR$, and so, exactly as above, we conclude that  $\mathcal{T}_\si(W,r)$ has a limit as $r \to 0^+$. $\Box$
\begin{cor}\label{Cor:ExistenceLimitUnBound}
	Let $W$ be a weak solution to \eqref{eq:TruncatedExtensionEquationLocalBackward}, with $q$ satisfying \eqref{eq:ASSUMPTIONSPOTENTIAL} and $p_0 \in B_{1/2}\times\{0\}\times(0,1)$. Assume that there exist $\si_0, r_0, C_0,  \delta_0 > 0$, such that
	\begin{equation}\label{eq:UnBoundonT}
	\mathcal{T}_{\si_0}(W,p,r) \leq C_0,	
	\end{equation}
	for all $r \in (0,r_0)$ and all $p \in Q_{\delta_0}$. Then both limits
	\[
	\lim_{r \to 0^+} \mathcal{W}_{\sigma}^{-}(W,p,r) \quad \text{ and } \quad \lim_{r \to 0^+} \mathcal{T}_{\si}(W,p,r)
	\]
	exist for all $\si \in (\si_0,\si_0 + (1-a)/2)$.

	Furthermore, there exists $\delta >0$ such that for any $\si \in (\si_0,\si_0 + (1-a)/2)$ and for any $p \in Q_{\delta}$ such that $\mathcal{W}_{\sigma}^{-}(W,p,0^+) = 0$, there exist $r_1, C_1 > 0$ independent of $p$ such that
	\[
	\mathcal{T}_{\si}(W,p,r) \leq C_1,
	\]
	for all $r \in (0,r_1)$.
\end{cor}
\emph{Proof.} The first part of the proof works exactly as in Lemma \ref{Lemma:ExistenceLimitNegPartWeiss}, thanks to \eqref{eq:UnBoundonT}. To show the second part, we first notice that
\[
\begin{aligned}
r^{2\si + 3} \frac{d}{dr} \mathcal{W}_{\si}(W,p,r) \geq - \int_0^{r^2} t^2 \int_{\mathbb{R}^{N+1}_+} F_p^2 d\mu_t dt  - \int_0^{r^2} t \int_{\mathbb{R}^N\times\{0\}} w_p^2 \,d\mu_t(x,0)dt
\end{aligned}
\]
thanks to \eqref{eq:ASSUMPTIONSPOTENTIAL} and the translation properties mentioned above. Now, it is not difficult to see that there exists $\delta \in (0,\delta_0)$ and $\overline{r}_0 \in (0,r_0)$ such that
\[
\int_0^{r^2} t^2 \int_{\mathbb{R}^{N+1}_+} F_p^2 d\mu_t \leq C r^{\overline{\si}} e^{-\frac{1}{32r^2}}
\]
for all $p \in Q_\delta$, $r \in (0,\overline{r}_0)$ and some constants $C,\overline{\sigma}$ independent of $p$ and $r$. Then, proceeding as before and using that \eqref{eq:UnBoundonT} is uniform w.r.t. $p \in Q_{\delta}$, we find that the constant $C > 0$ appearing in \eqref{eq:DiffIneqNegPartWeiss} does not depend on $p \in Q_{\delta_0}$, and so, taking eventually $\overline{r}_0$ smaller,
\begin{equation}\label{eq:DiffIneqNegPartWeissBis}
\frac{d}{dr} \mathcal{W}_{\si}^-(W,p,r) \leq \frac{C}{r^{2(\si - \si_0) + a}},
\end{equation}
for some constant $C$ independent of $p \in Q_\delta$ and all $r \in (0,\overline{r}_0)$. Thus, whenever $\mathcal{W}_{\si}^-(W,p,0^+)$, we integrate \eqref{eq:DiffIneqNegPartWeissBis}, to obtain
\[
\mathcal{W}_\si(W,p,r) \geq - C' r^{1-2(\si - \si_0) - a},
\]
for some new constant $C$ independent of $p \in Q_\delta$ and all $r \in (0,\overline{r}_0)$. Now, taking $r_1 \in (0,\overline{r}_0)$ small enough (independent of $p$) and using \eqref{eq:DerivativeTer}, it follows
\[
\frac{d}{dr} \mathcal{T}_{\sigma}(W,p,r) \geq - \frac{C}{r^{2(\si-\si_0)+a}},
\]
for all $r \in (0,r_1)$. Consequently, integrating between $r$ and $r_1$, we deduce
\[
\mathcal{T}_{\sigma}(W,p,r) \leq \mathcal{T}_{\sigma}(W,p,r_1) + Cr_1^{1 - 2(\sigma - \sigma_0) - a} \leq \sup_{p \in Q_\delta} \mathcal{T}_{\sigma}(W,p,r_1) + Cr_1^{1 - 2(\sigma - \sigma_0) - a} \leq C,
\]
for all $r \in (0,r_1)$ and some new constant $C > 0$, since $\mathcal{T}_{\sigma}(W,p,r_1)$ is bounded uniformly  in $p \in Q_\delta$. $\Box$

\bigskip

Let us summarize what we have proved so far. In Lemma \ref{Lemma:LimitTauSigmaSmall} we have shown the existence of $\si_0 > 0$, such that the limit
\[
\lim_{r \to 0^+} \mathcal{T}_{\sigma}(W,p_0,r) = 0,
\]
for every $\si \in (0,\si_0]$.  Now, in view of Lemma \ref{Lemma:ExistenceLimitNegPartWeiss}, there are two possible scenarios:
\begin{enumerate}
	\item There exists $\si \in (\si_0,\si_0 + (1-a)/2)$ such that $\mathcal{T}_{\sigma}(W,p_0,0^+) > 0$. In this case, thanks to Lemma \ref{Lemma:PrelDefVanOrder}, we can define
	\[
	\si_\ast := \inf \{\sigma > 0: \mathcal{T}_{\sigma}(W,p_0,0^+) > 0 \} = \sup \{\sigma > 0: \mathcal{T}_{\sigma}(W,p_0,0^+) = +\infty \}.
	\]
	\item $\mathcal{T}_{\sigma}(W,p_0,0^+) = 0$ for all $\si \in (\si_0,\si_0 + (1-a)/2)$. In this case, we can apply Lemma \ref{Lemma:ExistenceLimitNegPartWeiss} again and check whether $\mathcal{T}_{\sigma}(W,p_0,0^+)$ is positive or zero, for $\si \in (\si_1,\si_1 + (1-a)/2)$, where $\si_1 := \si_0 + (1-a)/2$ and so on.
\end{enumerate}
Iterating this procedure, it follows that either $\mathcal{T}_{\sigma}(W,p_0,0^+) = 0$ for all $\si > 0$, or there is $\si \in (0,+\infty)$ such that $\mathcal{T}_{\sigma}(W,p_0,0^+) > 0$. So, we can give the following definition.
\begin{defn} (Vanishing order)\label{def:VanishingOrder}
	Let $W$ be a solution to \eqref{eq:TruncatedExtensionEquationLocalBackward} with $q$ satisfying \eqref{eq:ASSUMPTIONSPOTENTIAL} and $p_0 \in \Gamma(W)$. We say that $W$ vanishes of infinite order at $p_0$ if
	\[
	\lim_{r \to 0^+} \mathcal{T}_{\sigma}(W,p_0,r) = 0 \quad  \text{ for all } \sigma > 0.
	\]
	We say that $W$ vanishes of order $\sigma_{\ast} > 0$ at $p_0$ if
	\begin{equation}\label{eq:DefVanishingOrder}
	\sigma_{\ast} := \inf \{\sigma > 0: \mathcal{T}_{\sigma}(W,p_0,0^+) = +\infty \} = \sup \{\sigma > 0: \mathcal{T}_{\sigma}(W,p_0,0^+) = 0 \}
	\end{equation}
	is finite.
\end{defn}
\begin{rem} The definition of vanishing order does not depend on the cut-off function $\zeta$ defined in \eqref{eq:PropCutOff}. Indeed, since $W = \overline{v}\zeta$ and $\zeta = 1$ in $\BB_{1/2}^+$ and $\supp \zeta \subset \BB_1^+$, it follows
\[
\int_0^{r^2} \int_{\mathbb{R}^{N+1}_+} W^2 d\mu_t dt = \int_0^{r^2} \int_{\BB_{1/2}^+} \overline{v}^2 d\mu_t dt + \int_0^{r^2} \int_{\BB_1^+ \setminus \BB_{1/2}^+} (\overline{v}\zeta)^2 d\mu_t dt
\]
Following the computations that have led us to \eqref{eq:ExpBoundIntF}, we obtain that
\[
\int_0^{r^2} \int_{\BB_1^+ \setminus \BB_{1/2}^+} (\overline{v}\zeta)^2 d\mu_t dt \leq C r^{\sigma_0} e^{-\frac{1}{16r^2}},
\]
for some $\si_0 \in \RR$ and some constant $C > 0$ depending only on $\overline{v}$. Consequently, for any $\si > 0$,
\[
\lim_{r \to 0^+} \mathcal{T}_{\sigma}(W,O,r) = \lim_{r \to 0^+} \frac{1}{r^{2\si + 2}} \int_0^{r^2} \int_{\BB_{1/2}^+} \overline{v}^2 d\mu_t dt.
\]
In particular, given two cut-off functions $\zeta_1$ and $\zeta_2$ (see  \eqref{eq:PropCutOff}), supported in $\mathbb{B}_{r_1}^+$ and $\mathbb{B}_{r_2}^+$, respectively ($r_1 \not= r_2$), and denoting with $W_1 = \overline{v}\zeta_1$ and $W_2 = \overline{v}\zeta_2$, we can repeat the computation above to deduce
\[
\lim_{r \to 0^+} \mathcal{T}_{\sigma}(W_1,O,r) = \lim_{r \to 0^+} \mathcal{T}_{\sigma}(W_2,O,r),
\]
which implies that the vanishing order of the solutions at $O$ does not depend on the cut-off $\zeta$.
\end{rem}
\begin{rem}
In the next subsection we will prove that if $\si_{\ast} \in (0,+\infty)$, then $N_D$ is bounded (compare, for instance, with Theorem \ref{Theorem:GeneralizedFrequency} and Corollary \ref{Lemma:LimitAlmgrenPoon}). As a consequence, in view of Lemma \ref{Lemma:PropWeissTau}, we obtain that if $\sigma_{\ast} \in (0,+\infty)$, then
\[
\sigma_{\ast} = \inf \{\sigma > 0: \mathcal{W}_{\sigma}(W,p_0,0^+) = -\infty \} = \sup \{\sigma > 0: \mathcal{W}_{\sigma}(W,p_0,0^+) = 0 \}.
\]
Indeed, if $N_D$ is bounded and $\si \in (0,\si_{\ast})$, then it is enough to pass to the limit as $r \to 0^+$ in the definition of $\mathcal{W}_{\si}$ and using that $\mathcal{T}_{\si}(W,0^+) = 0$ (see   \eqref{eq:SecondFormWeiss} and \ref{Lemma:PropWeissTau} part (i) and (iii)).
\end{rem}
\begin{rem}
Notice that, as in Definition \ref{Def:ClassicalVanOrd}, we have employed the notation $\sigma_\ast = \sigma_\ast(p_0)$ to denote the vanishing order of $W$ at $p_0$ even if, a priori, the two notions are different. In the next lemma we show that they are in fact equivalent. In particular, we have that $\sigma_\ast$ is the vanishing order of $W$ at $p_0$, as given in Definition \ref{Def:ClassicalVanOrd}, if and only if it is the vanishing order of $W$ at $p_0$ w.r.t. Definition \ref{def:VanishingOrder}. In doing this, we use some results we will show in Theorem \ref{THEOREMBLOWUP1New} and Theorem \ref{CONTINUOUSDEPOFBLOWUPLIMITS} that depend only on the notion of vanishing order given in Definition \ref{def:VanishingOrder}. In some sense, as a nontrivial consequence of our analysis, we recover the equivalence between the two definitions of vanishing order.
\end{rem}
\begin{lem}\label{Lemma:EquivalenceVanOrder}
Let $W$ be a solution to \eqref{eq:TruncatedExtensionEquationLocalBackward} with $q$ satisfying \eqref{eq:ASSUMPTIONSPOTENTIAL} and $p_0 \in \Gamma(W)$. Then the notions of vanishing order of $W$ at $p_0$ given in Definition \ref{Def:ClassicalVanOrd} and Definition \ref{def:VanishingOrder} are equivalent.
\end{lem}
\emph{Proof.} Fix $p_0 = O$. First, we show that if $W$ vanishes of infinite order at $O$ according with Definition  \ref{def:VanishingOrder}, then it vanishes of infinite order at $O$ in the sense of Definition\ref{Def:ClassicalVanOrd}. So, we assume that for any $\sigma > 0$, it holds
\[
\lim_{r \to 0^+} \frac{H(W,r)}{r^{2\sigma}} = 0.
\]
In view of the De Giorgi-Nash-Moser type estimate proved in \cite[Formula 5.7]{BanGarofalo2017:art}, there exists a universal constant $C> 0$ such that
\[
\| W \|_{L^{\infty}(\BB_r^+ \times (\delta r^2,r^2))}^2 \leq \frac{C}{r^{N + a + 3}} \int_{\delta r^2}^{4r^2} \int_{\BB_{2r}^+} y^a W^2 dXdt,
\]
for any fixed $\delta \in (0,1)$ and all $r \in (0,1/4)$ (notice that $\zeta = 1$ in $\BB_{1/2}^+$ and so $W$ and $\overline{v}$ coincide in $\BB_{1/2}^+ \times (0,1)$). On the other hand, for any $\delta \in (0,1)$, the fundamental solution satisfies
\[
\G_a(X,t) \geq \frac{C_{N,a,\delta}}{r^{N+a+1}} \quad \text{ in }  \BB_{2r}^+ \times (\delta r^2,4r^2),
\]
and so
\[
\frac{1}{r^{N + a + 3}} \int_{\delta r^2}^{4r^2} \int_{\BB_{2r}^+} y^a W^2 dXdt \leq \frac{C_{N,a,\delta}}{r^2} \int_{\delta r^2}^{4r^2} \int_{\BB_{2r}^+} y^a W^2 \G_a(X,t) dXdt \leq C_{N,a,\delta} H(W,r),
\]
for any fixed $\delta \in (0,1)$ and $r \in (0,1/4)$. Consequently, dividing by $r^{2\sigma}$, combining the two inequalities and passing to the limit as $r \to 0^+$, it follows
\[
\lim_{r \to 0^+} r^{-\sigma} \| W \|_{L^{\infty}(\BB_r^+ \times (\delta r^2,r^2))}  = 0,
\]
for any fixed $\delta \in (0,1)$ and any fixed $\sigma > 0$, which is exactly what we wanted to show.

\smallskip

To complete the proof, it is enough to show that if $W$ vanishes of order $\sigma_\ast \in (0,\infty)$ at $O$ in the sense of Definition \ref{def:VanishingOrder}, then it vanishes of the same order at $O$ also according with Definition\ref{Def:ClassicalVanOrd}.

This is an immediate consequence of Theorem \ref{THEOREMBLOWUP1New} and Theorem \ref{CONTINUOUSDEPOFBLOWUPLIMITS}, where we prove that if $W$ vanishes of order $\sigma_\ast \in (0,\infty)$ at $O$ w.r.t. Definition \ref{def:VanishingOrder}, then $\sigma_\ast$ is a positive integer (Theorem \ref{THEOREMBLOWUP1New}) and, furtheremore,
\[
W(X,t) = \Theta(X,t) + o(\|(X,t)\|^{\sigma_\ast}).
\]
as $\|(X,t)\|^2 := |X|^2 + |t| \to 0^+$ (Theorem \ref{CONTINUOUSDEPOFBLOWUPLIMITS}), where $\Theta$ is a parabolically homogeneous polynomial of degree $\sigma_\ast$, meaning that
\[
\sup_{(X,t) \in \BB_r^+ \times (0,r^2)} |W(X,t) - \Theta(X,t)| = o(r^{\sigma_\ast}),
\]
as $r \to 0^+$. Consequently, since
\[
\frac{1}{r^\sigma} \sup_{\BB_r^+ \times (0,r^2)} |\Theta| = \frac{1}{r^{\sigma-\sigma_\ast}} \sup_{\BB_1^+ \times (0,1)} |\Theta| \to
\begin{cases}
0           \quad &\text{ if } 0 < \sigma < \sigma_\ast \\
+ \infty    \quad &\text{ if } \sigma > \sigma_\ast,
\end{cases}
\]
as $r \to 0^+$, the same holds for $W$ and the proof is completed. $\Box$

\begin{lem}\label{Lemma:UpperSemicontinuity}
Let $W$ be a solution to \eqref{eq:TruncatedExtensionEquationLocalBackward} with $q$ satisfying \eqref{eq:ASSUMPTIONSPOTENTIAL} and $p_0 \in \Gamma(W)$. Then there exists a neighbourhood $Q$ of $p_0$, such that the function
\[
p \to \sigma_{\ast}(p)
\]
is upper semi-continuous on $\Gamma(W) \cap Q$.
\end{lem}

\emph{Proof.} After a possible translation, we can assume that $p_0 = O \in \Gamma(W)$ and that the order of vanishing $\si_{\ast} = \si_{\ast}(O)$ is finite. We fix $h \in (0,(1-a)/2)$ and proceed with an iterative argument.

\emph{Step 0:} By uniform H\"older regularity (Lemma \ref{Lemma:LimitTauSigmaSmall}), we know that there exist $\si_0, r_0, C_0,  \delta_0 > 0$ such that \eqref{eq:UnBoundonT} holds, i.e.
\[
\mathcal{T}_{\si_0}(W,p,r) \leq C_0,	
\]	
for all $r \in (0,r_0)$ and all $p \in Q_{\delta_0}$ (notice that we can take $\delta_0 = 1/2$). Then, by Corollary \ref{Cor:ExistenceLimitUnBound}, there are two possibilities:

\smallskip

\noindent (1) $\si_{\ast}(O) \in [\si_0,\si_0 + h)$ and so $\mathcal{T}_{\si}(W,O,0^+) = +\infty$ for all $\si \in (\si_\ast(O),\si_0 + h)$.

\noindent (2) $\si_{\ast}(O) \geq \si_0 + h$ and so $\mathcal{T}_{\si}(W,O,0^+) = +\infty$ for all $\si \in (\si_0,\si_0 + h)$.

\smallskip

Assume that (1) holds and fix arbitrarily $\si \in (\si_\ast(O),\si_0 + h)$. Then there are two possible scenarios:

\smallskip

\noindent (1.1) There exists $\delta > 0$ such that $\mathcal{W}_{\si}^-(W,p,0^+) > 0$ for all $p \in Q_\delta$.

\noindent (1.2) There exist a sequence $p_j \to O$ as $j \to + \infty$ such that $\mathcal{W}_{\si}^-(W,p_j,0^+) = 0$.

\smallskip

Assume that (1.1) holds. Then $\mathcal{W}_{\si}(W,p,0^+) < 0$ for all $p \in Q_\delta$ which implies $\mathcal{T}_{\si}(W,p,0^+) = +\infty$ for all $p \in Q_\delta$ in view of Lemma \ref{Lemma:PropWeissTau}. Using the arbitrariness of $\si \in (\si_\ast(O),\si_0 + h)$, it follows that $\si_\ast(p) \leq \si_\ast(O)$, i.e. $p \to \si_\ast(p)$ is upper semi-continuous.

Now we show that (1.2) cannot be valid. If $\mathcal{W}_{\si}^-(W,p_j,0^+) = 0$ for all $j$, then there exist $r_1, C_1 > 0$ independent of $j$ such that
\begin{equation}\label{eq:UnBoundpjStep0}
\mathcal{T}_{\si}(W,p_j,r) \leq C_1,
\end{equation}
for all $r \in (0,r_1)$, by Corollary \ref{Cor:ExistenceLimitUnBound}. Since $\mathcal{T}_{\si}(W,O,0^+) = +\infty$, we can assume
\[
\mathcal{T}_{\si}(W,O,r) > 2C_1,
\]
for all $r \in (0,r_1]$, taking eventually $r_1$ smaller (independent of $j$). Furthermore, using the continuity of the function $p \to \mathcal{T}_{\si}(W,p,r)$ for any fixed $r \in (0,1)$ and $\si > 0$, it follows that there exists $\delta > 0$ such that
\[
\mathcal{T}_{\si}(W,p,r_1) > C_1,
\]
for all $p \in Q_\delta$. Consequently, taking $p = p_j$ and passing to the limit as $j \to +\infty$, we find a contradiction with \eqref{eq:UnBoundpjStep0}.

Summing up, either $\si_{\ast}(O) \in [\si_0,\si_0 + h)$ and the function $p \to \si_\ast(p)$ is upper-semicontinuous at $p = O$, or $\si_{\ast}(O) \geq \si_0 + h$.

\emph{Step 1: Iteration.} Since $\si_\ast(O) \in (0,+\infty)$, there exists $k \in \NN$ such that $\si_\ast(O) \in [\si_0 + kh,\si_0 + (k+1)h)$. Thus, to conclude the proof, it is enough to repeat the argument of \emph{Step 0}. Since, $\si_\ast(O) \in [\si_0 + kh,\si_0 + (k+1)h)$, given any $\si \in (\si_\ast(O),\si_0 + (k+1)h)$, there are two options: either (1.1) or (1.2).

If (1.1) holds then it easily follows that $p \to \si_\ast(p)$ is upper-semicontinuous at $p = O$ by the arbitrariness of $\si \in (\si_\ast(O),\si_0 + (k+1)h)$.

If we assume (1.2), then $\mathcal{W}_{\si}^-(W,p_j,0^+) = 0$ for all $j$ and some sequence $p_j \to O$, as $j \to +\infty$. Thus, applying Corollary \ref{Cor:ExistenceLimitUnBound}, we obtain a bound of the form \eqref{eq:UnBoundpjStep0} (uniform in $j$) and we deduce a contradiction with the fact that $\mathcal{T}_{\si}(W,O,0^+) = +\infty$, for all $\si > \si_\ast(O)$. This ends the proof of the lemma. $\Box$
\subsection{Almost monotonicity of the Almgren-Poon quotient} The main goal of this section is to show an almost-monotonicity of the Almgren-Poon quotient at nodal points having finite vanishing order. The proof will be obtained in a sequence of steps and, even though it is more general, it mainly combines the ideas and the techniques of \cite{BanGarDanPetr2018:art} and the definition of vanishing order given in Definition \ref{def:VanishingOrder}.
\begin{thm}\label{Theorem:GeneralizedFrequency} (compare with \cite[Theorem 6.3]{DainelliGarofaloPetTo2017:book} and \cite[Theorem 4.8]{BanGarDanPetr2018:art})
Let $W$ be a weak solution to \eqref{eq:TruncatedExtensionEquationLocalBackward}, with $q$ satisfying \eqref{eq:ASSUMPTIONSPOTENTIAL}. Assume that $W$ vanishes of finite order at $p_0 \in B_{1/2}\times\{0\}\times(0,1)$. Then there exist $r_0 \in (0,1)$ and $C > 0$ depending on $N$, $s$, $K$ in \eqref{eq:ASSUMPTIONSPOTENTIAL} and on the solution to \eqref{eq:ExtensionEquationLocalBackward}, such that the function
\[
r \to  \Phi_a(W,p_0,r) := e^{Cr^{1-a}} N_I(W,p_0,r) + e^{Cr^{1-a}} - 1
\]
is monotone nondecreasing in $r \in (0,r_0)$. In particular, the limit
\[
\Phi_a(W,p_0,0^+) := \lim_{r \to 0^+} \Phi_a(W,p_0,r)
\]
exists and it is finite.
\end{thm}
\begin{rem} The main difference with respect to \cite[Theorem 6.3]{DainelliGarofaloPetTo2017:book} and \cite[Theorem 4.8]{BanGarDanPetr2018:art} is that in Theorem \ref{Theorem:GeneralizedFrequency} the constants $r_0, C > 0$ depend on the solution itself. Even if it is unusual, we will see that this dependence does not create problems in the blow-up procedure. In particular, it will turn out that the same monotonicity formula holds for any blow-up family, with constants not depending on the blow-up parameter, but only on the original solution to problem \eqref{eq:ExtensionEquationLocalBackward}.
\end{rem}

In the following we repeatedly use the notations in \eqref{eq:GaussianMeasures}:
\[
\begin{aligned}
d\mu_t^{p_0}(X) &:= y^a \mathcal{G}_a(x-x_0,y,t-t_0)dX \\
d\mu_t^{p_0}(x,0) &:= \mathcal{G}_a(x-x_0,0,t-t_0)dx,
\end{aligned}
\]
where $X \in \mathbb{R}^{N+1}_+$ and $t > t_0$.

\begin{lem}\label{Lemma:HeightsVariations} (compare with  \cite[Lemma 6.5]{BanGarofalo2017:art} and \cite[Lemma 4.3]{BanGarDanPetr2018:art})
Let $W$ be a weak solution to \eqref{eq:TruncatedExtensionEquationLocalBackward}, with $q$ satisfying \eqref{eq:ASSUMPTIONSPOTENTIAL}. Then
\[
h'(W,t) = \frac{2}{t} i(W,t) = \frac{2}{t} \left[ d(W,t) +  t\int_{\mathbb{R}^{N+1}_+}  F W \,d\mu_t \right]
\]
for any $t \in (0,1)$, and
\[
H'(W,r) = \frac{4}{r} I(W,r) = \frac{4}{r} \left[ D(W,r) + \frac{1}{r^2} \int_0^{r^2} t \int_{\mathbb{R}^{N+1}_+} F W \, d\mu_t dt \right],
\]
for any $r \in (0,1)$.
\end{lem}
\emph{Proof.} Even if our setting is more general, the proof is very similar to the ones given in \cite[Lemma 6.5]{BanGarofalo2017:art} and \cite[Lemma 4.3]{BanGarDanPetr2018:art}. Here we recall only the basic ideas. First of all, for any $\delta \in (0,1)$, we consider the approximations of $H(W,r)$ and $D(W,r)$, given by
\[
H_{\delta}(W,r) = \frac{1}{r^2} \int_{\delta r^2}^{r^2} h(W,t) \, dt = \int_{\delta}^1 h(W,r^2t) \, dt,
\]
and
\[
D_{\delta}(W,r) = \frac{1}{r^2} \int_{\delta r^2}^{r^2} d(W,t) \, dt = \int_{\delta}^1 d(W,r^2t) \, dt,
\]
respectively. Using the regularity properties of $W$ and its derivatives, the fact that $\mathcal{G}_a$ is bounded for any $t \in (\delta,1)$, we can apply the Lebesgue dominated convergence theorem to obtain
\[
h'(W,t) = 2 \int_{\mathbb{R}^{N+1}_+} y^a W \partial_tW \mathcal{G}_a  + \int_{\mathbb{R}^{N+1}_+} y^a W^2 \partial_t \mathcal{G}_a,
\]
for any $t \in (\delta,1)$. Furthermore, using the equations of $W$ and $\mathcal{G}_a$, and integrating by parts, it follows
\[
h'(W,t) = \frac{2}{t} \left[ d(W,t) +  t\int_{\mathbb{R}^{N+1}_+}  F W \,d\mu_t \right] = \frac{2}{t} i(W,t),
\]
for any $t \in (\delta,1)$. We recall that the integrations by parts can be justified by integrating on the region $\{(x,y) \in \mathbb{R}^{N+1}_+ :y > \varepsilon\}$ and then passing to the limit as $\varepsilon \to 0^+$, thanks to the regularity properties of $W$. Consequently,
\[
H'_{\delta}(W,r) = 2r \int_{\delta}^1 t h'(W,r^2t) dt = \frac{4}{r} D_{\delta}(W,r) + \frac{4}{r^3} \int_{\delta r^2}^{r^2} t \int_{\mathbb{R}^{N+1}_+} F W \,d\mu_t dt,
\]
for any $\delta \in (0,1)$ and any $r \in (0,1)$. Now, using the regularity of $W$ and $F$ (see  \cite[Lemma 6.5]{BanGarofalo2017:art}), it follows that for any $\delta_0 \in (0,1)$
\[
D_{\delta}(W,r) \to D(W,r) \quad \text{ uniformly in } r \in [\delta_0,1),
\]
and
\[
\int_{\delta r^2}^{r^2} t \int_{\mathbb{R}^{N+1}_+} F W \,d\mu_t dt \to \int_0^{r^2} t \int_{\mathbb{R}^{N+1}_+} F W \,d\mu_t dt \quad \text{ uniformly in } r \in [\delta_0,1),
\]
as $\delta \to 0^+$, and thus $H'_{\delta}(W,r) \to H'(W,r)$ as $\delta \to 0^+$ uniformly in $r \in [\delta_0,1)$. Using the arbitrariness of $\delta_0 \in (0,1)$, our statement follows. $\Box$
\begin{lem}\label{Lemma:BoundBelowNI} (Banerjee and Garofalo \cite[Lemma 4.5]{BanGarDanPetr2018:art})
Let $W$ be a weak solution to \eqref{eq:TruncatedExtensionEquationLocalBackward}, with $q$ satisfying \eqref{eq:ASSUMPTIONSPOTENTIAL}. Then for all $r \in (0,1)$ such that $H(W,r) > 0$, we have
\[
N_I(W,r) + \frac{1}{2} \geq 0.
\]
\end{lem}
\emph{Proof.} Assume $H(W,r) > 0$. Then
\[
2r\int_{\mathbb{R}^{N+1}_+\times\{t = r^2\}} W^2 \,d\mu_t = \frac{d}{dr} \left[r^2 H(W,r)\right] = 2rH(W,r) + r^2H'(W,r) = 2rH(W,r) \left[ 1 + 2N_I(W,r) \right],
\]
and the thesis follows since the l.h.s. is nonnegative. $\Box$

\begin{lem}\label{Lemma:DerivativeDd} (compare with \cite[Lemma 6.6]{BanGarofalo2017:art} and \cite[Lemma 4.6]{BanGarDanPetr2018:art})
Let $W $ be a weak solution to \eqref{eq:TruncatedExtensionEquationLocalBackward}, with $q$ satisfying \eqref{eq:ASSUMPTIONSPOTENTIAL}. Then
\[
\begin{aligned}
d'(W,t) &=  \frac{1}{2t} \int_{\mathbb{R}^{N+1}_+} (Z W)^2 \,d\mu_t - \int_{\mathbb{R}^{N+1}_+} Z W F \,d\mu_t \\
&- \frac{1-a}{2} \int_{\mathbb{R}^N\times\{0\}} q w^2 \,d\mu_t(x,0) - \frac{1}{2} \int_{\mathbb{R}^N\times\{0\}} [2t\partial_t q + x \cdot \nabla q] w^2 \,d\mu_t(x,0),
\end{aligned}
\]
for all $t \in (0,1)$, and
\[
\begin{aligned}
D'(W,r) &= \frac{1}{r^3} \int_0^{r^2} \int_{\mathbb{R}^{N+1}_+} (Z W)^2 \,d\mu_tdt  - \frac{2}{r^3} \int_0^{r^2} t \int_{\mathbb{R}^{N+1}_+} Z W F \,d\mu_tdt \\
& - \frac{1-a}{r^3} \int_0^{r^2} t \int_{\mathbb{R}^N\times\{0\}} q w^2 \,d\mu_t(x,0)dt - \frac{1}{r^3} \int_0^{r^2} t \int_{\mathbb{R}^N\times\{0\}} [2t\partial_t q + x \cdot \nabla q] w^2 \,d\mu_t(x,0)dt,
\end{aligned}
\]
for all $r \in (0,1)$.
\end{lem}
\emph{Proof.} As in Lemma \ref{Lemma:HeightsVariations}, we present the main ideas without focusing on the more technical parts which has been treated in \cite{BanGarofalo2017:art,BanGarDanPetr2018:art}. Compared with these papers, we work with the re-scaled version $\widetilde{W}$ of $W$ defined in \eqref{eq:RescaledOfW}, which has the advantage to make the computations more intuitive. Thus, since
\[
d(W,t) = d_0(\widetilde{W},1) + t^{\frac{1-a}{2}} \int_{\mathbb{R}^N\times\{0\}} \widetilde{q} \widetilde{w}^2 \mathcal{G}_a(x,0,1) \, dx,
\]
we compute
\[
\begin{aligned}
d'(W,t) &=  \frac{d}{dt} \int_{\mathbb{R}^{N+1}_+} |\nabla \widetilde{W}|^2 d\mu - \frac{1-a}{2} t^{-\frac{1+a}{2}} \int_{\mathbb{R}^N\times\{0\}} \widetilde{q} \widetilde{w}^2 \,d\mu(x,0) - t^{\frac{1-a}{2}} \int_{\mathbb{R}^N\times\{0\}} \partial_t \widetilde{q} \widetilde{w}^2 \,d\mu(x,0) \\
&- 2 t^{\frac{1-a}{2}} \int_{\mathbb{R}^N\times\{0\}} \widetilde{q} \widetilde{w} \partial_t\widetilde{w} \,d\mu(x,0).
\end{aligned}
\]
Using equation \eqref{eq:RescaledTruncatedExtensionEquationLocalBackward} (satisfied by $\widetilde{W}$) with test $\partial_t\widetilde{W}$, it follows
\[
\begin{aligned}
\frac{d}{dt} \int_{\mathbb{R}^{N+1}_+} |\nabla \widetilde{W}|^2 d\mu &= 2\int_{\mathbb{R}^{N+1}_+} \nabla \widetilde{W} \nabla (\partial_t\widetilde{W}) \,d\mu \\
&= 2t \int_{\mathbb{R}^{N+1}_+} (\partial_t\widetilde{W})^2 \,d\mu -2t\int_{\mathbb{R}^{N+1}_+} \partial_t\widetilde{W} \widetilde{F} \,d\mu + 2 t^{\frac{1-a}{2}} \int_{\mathbb{R}^N\times\{0\}} \widetilde{q} \widetilde{w}\partial_t\widetilde{w} \,d\mu(x,0).
\end{aligned}
\]
and thus, substituting above, it follows
\[
\begin{aligned}
d'(W,t) &=  2t \int_{\mathbb{R}^{N+1}_+} (\partial_t\widetilde{W})^2 \,d\mu -2t\int_{\mathbb{R}^{N+1}_+} \partial_t\widetilde{W} \widetilde{F} \,d\mu \\
&- \frac{1-a}{2} t^{-\frac{1+a}{2}} \int_{\mathbb{R}^N\times\{0\}} \widetilde{q} \widetilde{w}^2 \,d\mu(x,0) - t^{\frac{1-a}{2}} \int_{\mathbb{R}^N\times\{0\}} \partial_t \widetilde{q} \widetilde{w}^2 \,d\mu(x,0).
\end{aligned}
\]
Testing with $\partial_t \widetilde{W}$ is admissible thanks to the regularity properties of $\widetilde{W}$ (see  \cite[Lemma 5.5, Lemma 5.6]{BanGarofalo2017:art}) and the integrations by parts are done first on $\{(x,y) \in \mathbb{R}^{N+1}_+ :y > \varepsilon\}$ and then extended to $\mathbb{R}^{N+1}_+$, passing to the limit as $\varepsilon \to 0^+$. Equivalently,
\[
\begin{aligned}
d'(W,t) &=  \frac{1}{2t} \int_{\mathbb{R}^{N+1}_+} (Z W)^2 \,d\mu_t - \int_{\mathbb{R}^{N+1}_+} Z W F \,d\mu_t \\
&- \frac{1-a}{2} \int_{\mathbb{R}^N\times\{0\}} q w^2 \,d\mu_t(x,0) - \frac{1}{2} \int_{\mathbb{R}^N\times\{0\}} [2t\partial_t q + x \cdot \nabla q] w^2 \,d\mu_t(x,0),
\end{aligned}
\]
where we recall that $Z W := 2t\partial_t W + (x,y) \cdot \nabla W$. Now, we fix $\delta \in (0,1)$ and we proceed as in Lemma \ref{Lemma:HeightsVariations}. So, differentiating under the integral sign, we obtain
\[
\begin{aligned}
D_{\delta}'(W,r) & = 2r \int_{\delta}^1 t d'(W,r^2t) \, dt = \frac{1}{r^3} \int_{\delta r^2}^{r^2} \int_{\mathbb{R}^{N+1}_+} (Z W)^2 \,d\mu_tdt  - \frac{2}{r^3} \int_{\delta r^2}^{r^2} t \int_{\mathbb{R}^{N+1}_+} Z W F \,d\mu_tdt \\
&- \frac{1-a}{r^3} \int_{\delta r^2}^{r^2} t \int_{\mathbb{R}^N\times\{0\}} q w^2 \,d\mu_t(x,0) - \frac{1}{r^3} \int_{\delta r^2}^{r^2} t \int_{\mathbb{R}^N\times\{0\}} [2t\partial_t q + x \cdot \nabla q] w^2 \,d\mu_t(x,0),
\end{aligned}
\]
and, arguing as in \cite[Lemma 6.6]{BanGarofalo2017:art} and \cite[Lemma 4.6]{BanGarDanPetr2018:art}, for any fixed $\delta_0 \in (0,1)$, it follows
\[
\begin{aligned}
D_{\delta}'(W,r) &\to \frac{1}{r^3} \int_0^{r^2} \int_{\mathbb{R}^{N+1}_+} (Z W)^2 \,d\mu_tdt  - \frac{2}{r^3} \int_0^{r^2} t \int_{\mathbb{R}^{N+1}_+} Z W F \,d\mu_tdt \\
& - \frac{1-a}{r^3} \int_0^{r^2} t \int_{\mathbb{R}^N\times\{0\}} q w^2 \,d\mu_t(x,0) - \frac{1}{r^3} \int_0^{r^2} t \int_{\mathbb{R}^N\times\{0\}} [2t\partial_t q + x \cdot \nabla q] w^2 \,d\mu_t(x,0),
\end{aligned}
\]
uniformly in $r \in [\delta_0,1)$, as $\delta \to 0^+$. The arbitrariness of $\delta_0 \in (0,1)$ gives us the thesis. $\Box$
\begin{lem}\label{Lemma:DerivativeI} (compare with \cite[Lemma 4.7]{BanGarDanPetr2018:art})
Let $W$ be a weak solution to \eqref{eq:TruncatedExtensionEquationLocalBackward}, with $q$ satisfying \eqref{eq:ASSUMPTIONSPOTENTIAL}. Then
\[
\begin{aligned}
I'(W,r) &= D'(W,r) - \frac{2}{r^3} \int_0^{r^2} t \int_{\mathbb{R}^{N+1}_+} W F \,d\mu_tdt + 2r \int_{\mathbb{R}^{N+1}_+\times\{t = r^2 \}} W F \,d\mu_t,
\end{aligned}
\]
for all $r \in (0,1)$.
\end{lem}
\emph{Proof.} It is enough to differentiate $I(W,r) = D(W,r) + r^{-2} \int_0^{r^2} t \int_{\mathbb{R}^{N+1}_+} F W \, d\mu_tdt$. $\Box$
\begin{lem}\label{Lemma:TraceIneqL2Norm} (Banerjee and Garofalo \cite[Lemma 6.7]{BanGarofalo2017:art})
Let $W$ be a weak solution to \eqref{eq:TruncatedExtensionEquationLocalBackward}, with $q$ satisfying \eqref{eq:ASSUMPTIONSPOTENTIAL}. Then there exists $t_0 \in (0,1)$ and $C_0 > 0$ depending only on $N$, $s$, $K$ in \eqref{eq:ASSUMPTIONSPOTENTIAL} such that
\[
\int_{\mathbb{R}^N\times\{0\}} w^2 d\mu_t(x,0) \leq C_0 t^{-\frac{1+a}{2}} \left[ d(W,t) + h(W,t) \right],
0\]
for all $t \in (0,t_0)$. Furthermore, if $s \in [1/2,1)$, one also has
\[
\left|\int_{\mathbb{R}^N\times\{0\}} (x\cdot \nabla_x q) w^2 d\mu_t(x,0) \right|\leq C_0 t^{-\frac{1+a}{2}} \left[ d(W,t) + h(W,t) \right],
0\]
for all $t \in (0,t_0)$.
\end{lem}
\emph{Proof.} The proof is based on two types of trace inequalities and works exactly as in \cite{BanGarofalo2017:art}. $\Box$
\begin{lem}\label{Lemma:DerivativeD} (compare with \cite[Lemma 6.6]{BanGarofalo2017:art})
Let $W$ be a weak solution to \eqref{eq:TruncatedExtensionEquationLocalBackward}, with $q$ satisfying \eqref{eq:ASSUMPTIONSPOTENTIAL}. Then there exists $t_0,r_0 \in (0,1)$ and $C_0 > 0$ depending only on $N$, $s$, $K$ in \eqref{eq:ASSUMPTIONSPOTENTIAL} such that
\[
d'(W,t) \geq  \frac{1}{2t} \int_{\mathbb{R}^{N+1}_+} (Z W)^2 \,d\mu_t - \int_{\mathbb{R}^{N+1}_+} Z W F \,d\mu_t - C t^{-\frac{1+a}{2}} \left[ d(W,t) + h(W,t) \right],
\]
for all $t \in (0,t_0)$, and
\[
D'(W,r) \geq \frac{1}{r^3} \int_0^{r^2} \int_{\mathbb{R}^{N+1}_+} (Z W)^2 \,d\mu_tdt  - \frac{2}{r^3} \int_0^{r^2} t \int_{\mathbb{R}^{N+1}_+} Z W F \,d\mu_tdt - C_0 r^{-a} \left[H(W,r) + D(W,r)\right],
\]
for all $r \in (0,r_0)$.
\end{lem}
\emph{Proof.} Assume first $s \in (0,1/2)$. Thanks to our assumptions on the potential $q = q(x,t)$, we easily see that
\[
d'(W,t) \geq  \frac{1}{2t} \int_{\mathbb{R}^{N+1}_+} (Z W)^2 \,d\mu_t - \int_{\mathbb{R}^{N+1}_+} Z W F \,d\mu_t - C \int_{\mathbb{R}^N\times\{0\}} w^2 \,d\mu_t(x,0)
\]
for some suitable constant $C > 0$ and all $t \in (0,1)$. A direct application of Lemma \ref{Lemma:TraceIneqL2Norm} proves our first claim. Now, combing the inequality in our first claim and exploiting the fact that $D'(W,r) = 2r \int_0^1 t d'(W,r^2t) \, dt$, we obtain
\[
D'(W,r) \geq \frac{1}{r^3} \int_0^{r^2} \int_{\mathbb{R}^{N+1}_+} (Z W)^2 \,d\mu_tdt  - \frac{2}{r^3} \int_0^{r^2} t \int_{\mathbb{R}^{N+1}_+} Z W F \,d\mu_tdt - C_0 r^{-a} \left[H(W,r) + D(W,r)\right],
\]
for any $r \in (0,r_0)$, with $r_0 = \sqrt{t_0}$ and a suitable $C_0 > 0$, which is exactly the second inequality in our statement. When $s \in [1/2,1)$ the proof is similar and it immediately follows by employing also the second inequality in Lemma \ref{Lemma:TraceIneqL2Norm}. $\Box$
\begin{lem}\label{Lemma:EstimateFBelowIPrime}
Let $W $ be a weak solution to \eqref{eq:TruncatedExtensionEquationLocalBackward}, with $q$ satisfying \eqref{eq:ASSUMPTIONSPOTENTIAL}. Then there exists $r_0 \in (0,1)$ and $C_0 > 0$ depending only on $N$, $s$, $K$ in \eqref{eq:ASSUMPTIONSPOTENTIAL} such that
\[
\begin{aligned}
I'(W,r) &\geq \frac{1}{r^3} \int_0^{r^2} \int_{\mathbb{R}^{N+1}_+} (Z W)^2 \,d\mu_tdt  - \frac{2}{r^3} \int_0^{r^2} t \int_{\mathbb{R}^{N+1}_+} Z W F \,d\mu_tdt \\
& + \left( \frac{C_0}{r^{2+a}} -\frac{2}{r^3} \right) \int_0^{r^2} t \int_{\mathbb{R}^{N+1}_+} W F \,d\mu_tdt \\
&+ 2r \int_{\mathbb{R}^{N+1}_+\times\{t = r^2\}} W F \,d\mu_t  - C_0 r^{-a} \left[H(W,r) + I(W,r)\right],
\end{aligned}
\]
for all $r \in (0,r_0)$.
\end{lem}
\emph{Proof.} It is enough to combine Lemma \ref{Lemma:DerivativeI} and Lemma \ref{Lemma:DerivativeD}. $\Box$
\begin{lem}\label{Lemma:DerivativeNI}
Let $W$ be a weak solution to \eqref{eq:TruncatedExtensionEquationLocalBackward}, with $q$  satisfying \eqref{eq:ASSUMPTIONSPOTENTIAL}. Assume that $W$ vanishes of finite order $\sigma_{\ast} > 0$ at $O$. Then there exists $r_0 \in (0,1)$ and $C > 0$ depending on $N$, $s$, $K$ in \eqref{eq:ASSUMPTIONSPOTENTIAL} and on the solution to \eqref{eq:ExtensionEquationLocalBackward}, such that
\[
\frac{d}{dr}[1 + N_I(W,r)] \geq - C [1 + N_I(W,r)] r^{-a}.
\]
for all $r \in (0,r_0)$.
\end{lem}
\emph{Proof.} First we notice that the assumption $\sigma_{\ast} \in (0,\infty)$ implies the existence of a small $\varepsilon > 0$ such that
\begin{equation}\label{eq:BoundsHGoodPoints}
C^{-1}r^{2\sigma_\ast + \varepsilon} \leq H(W,r) \leq Cr^{2\sigma_\ast - \varepsilon}
\end{equation}
for all $r \in (0,r_0)$, for some $r_0,C > 0$ depending on the solution to problem \eqref{eq:ExtensionEquationLocalBackward}. In particular, it guarantees $H(W,r) > 0$ and so the function $N_I$ is well-defined.

In the first part of the proof, we follow \cite[Theorem 4.8]{BanGarDanPetr2018:art}, while in the second, we use the information in \eqref{eq:BoundsHGoodPoints} to show our claim. Notice that compared with \cite{BanGarDanPetr2018:art}, there extra terms coming from the inequality in Lemma \ref{Lemma:EstimateFBelowIPrime}. Applying Lemma \ref{Lemma:EstimateFBelowIPrime} and abbreviating $H(W,r) = H(r)$, $I(W,r) = I(r)$ and so on, we compute
\[
\begin{aligned}
r^5 H(r)^2 N_I'(r) & = r^5 \left[ I'(r) H(r) - I(r) H'(r) \right] = r^5 H(r) \left[ I'(r) - \frac{4I^2(r)}{rH(r)} \right] \\
& \geq r^5 H(r) \bigg\{  \frac{1}{r^3} \int_0^{r^2} \int_{\mathbb{R}^{N+1}_+} (Z W)^2 \,d\mu_tdt  - \frac{2}{r^3} \int_0^{r^2} t \int_{\mathbb{R}^{N+1}_+} Z W F \,d\mu_tdt \\
& + \left( \frac{C_0}{r^{2+a}} -\frac{2}{r^3} \right) \int_0^{r^2} t \int_{\mathbb{R}^{N+1}_+} W F \,d\mu_tdt + 2r \int_{\mathbb{R}^{N+1}_+\times\{t = r^2\}} W F \,d\mu_t \\
& - C_0 r^{-a} \left[H(r) + I(r)\right] - \frac{4I(r)^2}{r H(r)}  \bigg\} \\
& = r^5 H(r) \bigg\{  \frac{1}{r^3} \int_0^{r^2} \int_{\mathbb{R}^{N+1}_+} (Z W - tF)^2 \,d\mu_tdt  - \frac{1}{r^3} \int_0^{r^2} t^2 \int_{\mathbb{R}^{N+1}_+} F^2 \,d\mu_tdt \\
& + \left( \frac{C_0}{r^{2+a}} -\frac{2}{r^3} \right) \int_0^{r^2} t \int_{\mathbb{R}^{N+1}_+} W F \,d\mu_tdt + 2r \int_{\mathbb{R}^{N+1}_+\times\{t = r^2\}} W F \,d\mu_t \\
& - C_0 r^{-a} \left[H(r) + I(r)\right] - \frac{4I(r)^2}{r H(r)}  \bigg\},
\end{aligned}
\]
for all $r \in (0,r_0)$. Now, we have
\[
\begin{aligned}
r^5 H(r)^2 N_I'(r) & = r^2 H(r) \int_0^{r^2} \int_{\mathbb{R}^{N+1}_+} (Z W - tF)^2 \,d\mu_tdt  - r^2 H(r) \int_0^{r^2} t^2 \int_{\mathbb{R}^{N+1}_+} F^2 \,d\mu_tdt \\
& + \left( C_0 r^{3-a} - 2r^2 \right) H(r) \int_0^{r^2} t \int_{\mathbb{R}^{N+1}_+} W F \,d\mu_tdt + 2r^6 H(r) \int_{\mathbb{R}^{N+1}_+\times\{t = r^2\}} W F \,d\mu_t \\
& - C_0 r^{5-a} H^2(W,r) \left[1 + N_I(W,r)\right] - \left( \int_{\mathbb{R}^{N+1}_+} W Z W \,d\mu_tdt \right)^2,
\end{aligned}
\]
and so, by Cauchy-Schwartz inequality, we obtain
\[
\begin{aligned}
r^5 H(r)^2 N_I'(r) & \geq  \left[ \int_0^{r^2} \int_{\mathbb{R}^{N+1}_+} W(Z W - tF) \,d\mu_tdt \right]^2  - r^2 H(r) \int_0^{r^2} t^2 \int_{\mathbb{R}^{N+1}_+} F^2 \,d\mu_tdt \\
& + \left( C_0 r^{3-a} - 2r^2 \right) H(r) \int_0^{r^2} t \int_{\mathbb{R}^{N+1}_+} W F \,d\mu_tdt + 2r^6 H(r) \int_{\mathbb{R}^{N+1}_+\times\{t = r^2\}} W F \,d\mu_t \\
& - C_0 r^{5-a} H^2(W,r) \left[1 + N_I(W,r)\right] - \left( \int_{\mathbb{R}^{N+1}_+} W Z W \,d\mu_tdt \right)^2 \\
& \geq  \left( \int_0^{r^2} t \int_{\mathbb{R}^{N+1}_+} WF \,d\mu_tdt \right)^2 - 4r^2 I(r) \int_0^{r^2} t \int_{\mathbb{R}^{N+1}_+} WF \,d\mu_tdt \\
& - r^2 H(r) \int_0^{r^2} t^2 \int_{\mathbb{R}^{N+1}_+} F^2 \,d\mu_tdt + \left( C_0 r^{3-a} - 2r^2 \right) H(r) \int_0^{r^2} t \int_{\mathbb{R}^{N+1}_+} W F \,d\mu_tdt \\
& + 2r^6 H(r) \int_{\mathbb{R}^{N+1}_+\times\{t = r^2\}} W F \,d\mu_t - C_0 r^{5-a} H^2(W,r) \left[1 + N_I(r)\right],
\end{aligned}
\]
for all $r \in (0,r_0)$. Thus, dividing by $r^5 H(r)^2$, it follows
\[
\begin{aligned}
N_I'(r) & \geq - \frac{4}{r^3} N_I(r) \frac{\int_0^{r^2} t \int_{\mathbb{R}^{N+1}_+} WF \,d\mu_tdt}{H(r)} - \frac{1}{r^3}\frac{\int_0^{r^2} t^2 \int_{\mathbb{R}^{N+1}_+} F^2 \,d\mu_tdt}{H(r)} \\  %
& + \frac{C_0 r^{1-a}-2}{r^3} \frac{\int_0^{r^2} t \int_{\mathbb{R}^{N+1}_+} W F \,d\mu_tdt}{H(r)} + 2r \frac{\int_{\mathbb{R}^{N+1}_+\times\{t = r^2\}} W F \,d\mu_t}{H(r)} - C_0 r^{-a} \left[1 + N_I(r)\right] \\
& = - \frac{4}{r^3} \left[ N_I(r)  + \frac{1}{2} \right] \frac{\int_0^{r^2} t \int_{\mathbb{R}^{N+1}_+} WF \,d\mu_tdt}{H(r)} + C_0 r^{-2-a} \frac{\int_0^{r^2} t \int_{\mathbb{R}^{N+1}_+} WF \,d\mu_tdt}{H(r)} \\
&- \frac{1}{r^3}\frac{\int_0^{r^2} t^2 \int_{\mathbb{R}^{N+1}_+} F^2 \,d\mu_tdt}{H(r)} + 2r \frac{\int_{\mathbb{R}^{N+1}_+\times\{t = r^2\}} W F \,d\mu_t}{H(r)} - C_0 r^{-a} \left[1 + N_I(r)\right],
\end{aligned}
\]
for all $r \in (0,r_0)$. Using the Cauchy-Schwartz inequality again and recalling that $N_I(r) + 1/2 \geq 0$, we obtain
\[
\begin{aligned}
N_I'(r) & \geq - \frac{4}{r^2} \left[ \frac{1}{2} + N_I(r) \right] \frac{\int_0^{r^2} t^2 \int_{\mathbb{R}^{N+1}_+} F^2 \,d\mu_tdt}{H^{1/2}(r)} - \frac{1}{r^3}\frac{\int_0^{r^2} t^2 \int_{\mathbb{R}^{N+1}_+} F^2 \,d\mu_tdt}{H(r)} \\
& - 2r \frac{\left( \int_{\mathbb{R}^{N+1}_+\times\{t = r^2\}} W^2 \,d\mu_t \right)^{1/2} \left( \int_{\mathbb{R}^{N+1}_+\times\{t = r^2\}} F^2 \,d\mu_t \right)^{1/2}}{H(r)} \\
& - C_0 r^{-1-a} \frac{ \left( \int_0^{r^2} t^2 \int_{\mathbb{R}^{N+1}_+} F^2 \,d\mu_tdt  \right)^{1/2}}{H^{1/2}(r)} - C_0 r^{-a} \left[1 + N_I(r)\right],
\end{aligned}
\]
and so
\[
\begin{aligned}
N_I'(r) & \geq - \frac{4}{r^2} \left[ 1 + N_I(r) \right] \frac{\int_0^{r^2} t^2 \int_{\mathbb{R}^{N+1}_+} F^2 \,d\mu_tdt}{H^{1/2}(r)} - \frac{1}{r^3}\frac{\int_0^{r^2} t^2 \int_{\mathbb{R}^{N+1}_+} F^2 \,d\mu_tdt}{H(r)} \\
& - 2 r \frac{ \left[ 1 + N_I(r) \right] \left( \int_{\mathbb{R}^{N+1}_+\times\{t = r^2\}} F^2 \,d\mu_t \right)^{1/2}}{H^{1/2}(r)} \\
& - C_0 r^{-1-a} \frac{ \left( \int_0^{r^2} t^2 \int_{\mathbb{R}^{N+1}_+} F^2 \,d\mu_tdt  \right)^{1/2}}{H^{1/2}(r)} - C_0 r^{-a} \left[1 + N_I(r)\right] \\
& \geq - \left[ 1 + N_I(r) \right] \bigg\{ \frac{4}{r^2} \frac{\int_0^{r^2} t^2 \int_{\mathbb{R}^{N+1}_+} F^2 \,d\mu_tdt}{H^{1/2}(r)} +  2r \frac{ \Big(\int_{\mathbb{R}^{N+1}_+\times\{t = r^2\}} F^2 \,d\mu_t \Big)^{1/2}}{H^{1/2}(r)} + C_0r^{-a} \bigg\} \\
& - \frac{1}{r^3}\frac{\int_0^{r^2} t^2 \int_{\mathbb{R}^{N+1}_+} F^2 \,d\mu_tdt}{H(r)} - C_0 r^{-1-a} \frac{ \left( \int_0^{r^2} t^2 \int_{\mathbb{R}^{N+1}_+} F^2 \,d\mu_tdt  \right)^{1/2}}{H^{1/2}(r)},
\end{aligned}
\]
for all $r \in (0,r_0)$, where we have used the fact that $\int_{\mathbb{R}^{N+1}_+\times\{t = r^2\}} W^2 \,d\mu_t = H(r)[1 + 2N_I(r) ]$ (see  Lemma \ref{Lemma:BoundBelowNI}) and the fact that $(1 + 2N)^{1/2} \leq 1 + N$.

At this point, proceeding as in \eqref{eq:ExpBoundIntF}, we deduce
\[
\int_{\mathbb{R}^{N+1}_+\times\{t = r^2\}} F^2 \,d\mu_t \leq C r^{\sigma_0} e^{-\frac{1}{16r^2}},   \qquad   \int_0^{r^2} t^2 \int_{\mathbb{R}^{N+1}_+} F^2 \,d\mu_tdt \leq C r^{\sigma_0} e^{-\frac{1}{16r^2}},
\]
for all $r \in (0,r_0)$, for some $r_0 \in (0,1)$ and $\sigma_0 \in \mathbb{R}$, depending only on $N$ and $a$ and $C > 0$ depending on the solution to \eqref{eq:ExtensionEquationLocalBackward}. So, choosing a suitable $\sigma_0 \in \mathbb{R}$ (let us say, negative enough), we obtain
\[
\begin{aligned}
\frac{d}{dr}[ 1 + N_I(r) ] & \geq - C\left[ 1 + N_I(r) \right] \bigg\{ \frac{r^{\sigma_0}}{H^{1/2}(r)}e^{-\frac{1}{16r^2}}  + r^{-a} \bigg\} \\
& - C\frac{r^{\sigma_0}}{H(r)} e^{-\frac{1}{16r^2}} - C \frac{r^{\sigma_0}}{H^{1/2}(r)}e^{-\frac{1}{32r^2}}
\end{aligned}
\]
for all $r \in (0,r_0)$ and, taking eventually $r_0 \in (0,1)$ smaller (depending on the solution to \eqref{eq:ExtensionEquationLocalBackward}) we can use \eqref{eq:BoundsHGoodPoints} and Lemma \ref{Lemma:BoundBelowNI} to deduce
\[
\begin{aligned}
\frac{d}{dr}[ 1 + N_I(r) ] & \geq - C\left[ 1 + N_I(r) \right] \bigg[ r^{\sigma_0} e^{-\frac{1}{16r^2}}  + r^{-a} \bigg] \\
& - Cr^{\sigma_0} e^{-\frac{1}{16r^2}} - C r^{\sigma_0}e^{-\frac{1}{32r^2}} \\
&\geq - C\left[ 1 + N_I(r) \right]  r^{-a}
\end{aligned}
\]
for all $r \in (0,r_0)$, for some new $r_0,C > 0$ depending on the solution to \eqref{eq:ExtensionEquationLocalBackward}. $\Box$

\bigskip

\emph{Proof of Theorem \ref{Theorem:GeneralizedFrequency}.} It is enough to compute
\[
\begin{aligned}
\frac{d}{dr} \Phi_a(W,O,r) &= \frac{d}{dr} \left\{  e^{Cr^{1-a}} \left[ N_I(r) + 1 \right] - 1 \right\} \\
& = r^{-a}e^{Cr^{1-a}} \left\{ C(1-a) \left[N_I(r) + 1 \right] + r^a \frac{d}{dr} \left[N_I(r) + 1 \right] \right\} \\
& \geq r^{-a}e^{Cr^{1-a}} \left[N_I(r) + 1 \right] \left[ C(1-a)  - C_1 \right],
\end{aligned}
\]
for all $r \in (0,r_0)$, thanks to Lemma \ref{Lemma:DerivativeNI}. Choosing $C \geq C_1/(1-a)$, the proof is completed. $\Box$
\begin{cor}\label{Lemma:LimitAlmgrenPoon} (see  \cite[Corollary 6.12]{BanGarofalo2017:art} and \cite[Lemma 6.1]{BanGarDanPetr2018:art})
Let $W $ be a weak solution to \eqref{eq:TruncatedExtensionEquationLocalBackward}, with $q$ satisfying \eqref{eq:ASSUMPTIONSPOTENTIAL}. Assume that $W$ vanishes of finite order $\sigma_{\ast} > 0$ at $p_0$.

Then
\[
\Phi_a(W,p_0,0^+) = \lim_{r \to 0^+} N_I(W,p_0,r) = \lim_{r \to 0^+} N_D(W,p_0,r) = \lim_{r \to 0^+} N_0(W,p_0,r) = \kappa,
\]
for some finite limit $\kappa \geq 0$.
\end{cor}
\emph{Proof.} Let us take $p_0 = O$. First, since the limit of $\Phi_a(W,O,r)$ as $r \to 0^+$ exists and it is finite, it must equal the limit of $N_I(W,r)$ by definition. On the other hand,
\[
N_I(W,r) := N_D(W,r) + \frac{r^{-2}\int_0^{r^2} t \int_{\mathbb{R}^{N+1}_+} WF \,d\mu_tdt}{H(r)},
\]
and thus, since
\[
\frac{\left|r^{-2}\int_0^{r^2} t \int_{\mathbb{R}^{N+1}_+} WF \,d\mu_tdt\right|}{H(r)} \leq \frac{1}{r} \frac{\left( \int_0^{r^2} t^2 \int_{\mathbb{R}^{N+1}_+} F^2 \,d\mu_tdt \right)^{1/2}}{H^{1/2}(W,r)} \leq C r^{\sigma_0} e^{-\frac{1}{32r^2}} \to 0,
\]
as $r \to 0$ thanks to \eqref{eq:ExpBoundIntF}, we obtain that $N_I(W,r)$ and $N_D(W,r)$ have the same finite limit as $r \to 0$ (recall that both $N_I(W,r)$ and $N_D(W,r)$ are bounded from below). Similar,
\[
N_D(W,r) := N_0(W,r) - \frac{r^{-2}\int_0^{r^2} t \int_{\mathbb{R}^N \times \{0\}} qw^2 \,d\mu_t(x,0)dt}{H(r)},
\]
and, by Lemma \ref{Lemma:TraceIneqL2Norm}, it follows
\[
\int_{\mathbb{R}^N\times\{0\}} w^2 d\mu_t(x,0) \leq C_0 t^{-\frac{1+a}{2}} \left[ d(W,t) + h(W,t) \right], \quad t \in (0,t_0),
\]
for some constants $C_0 > 0$ and $t_0 \in (0,1)$. Consequently, since
\[
\begin{aligned}
\frac{\left| r^{-2}\int_0^{r^2} t \int_{\mathbb{R}^N \times \{0\}} qw^2 \,d\mu_t(x,0)dt\right|}{H(r)} &\leq K \frac{\left| r^{-2}\int_0^{r^2} t \int_{\mathbb{R}^N \times \{0\}} w^2 \,d\mu_t(x,0)dt\right|}{H(r)} \\
&\leq C_0 K r^{1-a} \left[N_D(W,r) + 1 \right] \to 0,
\end{aligned}
\]
as $r \to 0$, also $N_D(W,r)$ and $N_0(W,r)$ have the same limit as $r \to 0$ (here we are exploiting the fact that $N_D(W,r)$ is bounded from below). Finally, since $N_0(W,r) \geq 0$, it follows $\kappa \geq 0$. $\Box$

\begin{cor}
Let $W$ be a weak solution to \eqref{eq:TruncatedExtensionEquationLocalBackward}, with $q$ satisfying \eqref{eq:ASSUMPTIONSPOTENTIAL}. Assume that $W$ vanishes of finite order at some $p_0 \in B_{1/2}\times\{0\}\times(0,1)$. Then there is a neighbourhood $Q$ of $p_0$ such that for any $p \in Q$, $r_0 \in (0,1)$ and $C > 0$ depending on $N$, $s$, $K$ in \eqref{eq:ASSUMPTIONSPOTENTIAL}, the solution to \eqref{eq:ExtensionEquationLocalBackward} and $p$, such that the function
\[
r \to  \Phi_a(W,p,r) := e^{Cr^{1-a}} N_I(W,p,r) + e^{Cr^{1-a}} - 1
\]
is monotone nondecreasing in $r \in (0,r_0)$. Furthermore, for any $p \in Q$
\[
\Phi_a(W,p,0^+) = \lim_{r \to 0^+} N_I(W,p,r) = \lim_{r \to 0^+} N_D(W,p,r) = \lim_{r \to 0^+} N_0(W,p,r) = \kappa,
\]
for some finite limit $\kappa \geq 0$, depending on $p$.
\end{cor}
\emph{Proof.} The proof follows by a translation argument and Lemma \ref{Lemma:UpperSemicontinuity}. $\Box$

\begin{cor}\label{Corollary:NonDegeneracy}
Let $W$ be a weak solution to \eqref{eq:TruncatedExtensionEquationLocalBackward}, with $q$ satisfying \eqref{eq:ASSUMPTIONSPOTENTIAL}. Assume that $W$ vanishes of finite order at some $p_0 \in B_{1/2}\times\{0\}\times(0,1)$. Then there exist $r_0 \in (0,1)$ and $C > 0$ depending on $N$, $s$, $K$ in \eqref{eq:ASSUMPTIONSPOTENTIAL}, the solution to \eqref{eq:ExtensionEquationLocalBackward} and $p_0$, such that
\[
\frac{H(W,p_0,r_2)}{H(W,p_0,r_1)} \leq \left(\frac{r_2}{r_1}\right)^C,
\]
for all $0 < r_1 < r_2 < r_0$.
\end{cor}
\emph{Proof.} The proof is standard. Take $p_0 = O$ and notice that since
\[
\lim_{r \to 0^+} N_I(W,r) = \lim_{r \to 0^+} \frac{r H'(W,r)}{4H(W,r)} = \kappa \geq 0,
\]
there exists $\kappa_0 > \kappa$ and $r_0 \in (0,1)$ such that
\[
\frac{H'(W,r)}{H(W,r)} \leq \frac{4\kappa_0}{r},
\]
for all $r \in (0,r_0)$. Taking $C = 4\kappa_0$ and integrating between $r_1$ and $r_2$ we obtain the thesis. $\Box$

\subsection{Monotonicity of the Almgren-Poon quotient for global solutions} We end the section by recalling that if $U$ is a global solutions to \eqref{eq:ExtensionEquationLocalBackward} (i.e. defined in the whole space) with $q = 0$ and having finite vanishing order at $O$, our methods can be adapted to prove the pure monotonicity of the Almgren-Poon quotient (see  \cite{BanGarofalo2017:art}).
\begin{thm}\label{Theorem:MonAlmgrenPoon} (see   Banerjee and Garofalo \cite[Theorem 1.3, Theorem 8.3]{BanGarofalo2017:art})
Let $U$ be a global solution to \eqref{eq:ExtensionEquationLocalBackward} with $q = 0$, having finite vanishing order at $O$. Then there exists a universal $r_0 > 0$, depending only on $N$ and $a$ such that, under the assumption $H(U,r) > 0$ for $r \in (0,r_0)$, then the function
\[
r \to N_0(U,r)
\]
is monotone nondecreasing in $(0,r_0)$. Furthermore, $N_0(U,r) = \kappa$ for all $r \in (0,r_0)$ if and only if $U$ is parabolically homogeneous of degree $2\kappa$.

Similar, there exists a universal $t_0 > 0$, depending only on $N$ and $a$ such that, under the assumption $h(U,t) > 0$ for $t \in (0,t_0)$, then the function
\[
t \to n_0(U,t) := \frac{d_0(U,t)}{h(U,t)}
\]
is monotone nondecreasing in $(0,t_0)$. Furthermore, $n_0(U,t) = \kappa$ for all $t \in (0,t_0)$ if and only if $U$ is parabolically homogeneous of degree $2\kappa$.
\end{thm}
\begin{rem}\label{Remark:HomogeneousStuff} As we have recalled in the introduction, $U = U(X,t)$ is parabolically homogeneous of degree $2\kappa$ if $U(rX,r^2t) = r^{2\kappa}U(X,t)$ for all $r > 0$. If $U$ is $C^1$, it is not difficult to show that this condition is equivalent to say that
\[
Z U = 2\kappa U,
\]
where $Z U := 2t\partial_t U + X \cdot \nabla U$,  see  \eqref{eq:TimeHomogeneityFactors}. Now, if $U$ is a global solution to \eqref{eq:ExtensionEquationLocalBackward} with $q = 0$, its re-scaled version
\[
\widetilde{U}(X,t) := U(\sqrt{t}X,t)
\]
is a solution to problem \eqref{eq:RescaledTruncatedExtensionEquationLocalBackward} (with $\widetilde{F} = \widetilde{q} = 0$). Consequently, noticing that $2t\partial_t\widetilde{U} = ZU$, it follows that $U$ is parabolically $2\kappa$-homogeneous if and only if
\[
\begin{cases}
- \mathcal{O}_a\widetilde{U} = \kappa\widetilde{U} \quad &\text{ in } \mathbb{R}_+^{N+1} \\
-\partial_y^a \widetilde{U} = 0 \quad &\text{ in } \mathbb{R}^N \times \{0\},
\end{cases}
\]
for all $t \in (0,1)$, in the weak sense, which means that, at each time level, $\widetilde{U}$ is a solution to an Ornstein-Uhlenbeck eigenvalue problem. With the aim of classifying all global homogeneous solutions to \eqref{eq:ExtensionEquationLocalBackward} with $q = 0$, we devote the next section to the study of the above eigenvalue problem.
\end{rem}
\section{Spectral Analysis and Applications}\label{Section:SpectralAnalysis}
\subsection{Spectral Theorem} As mentioned at the end of the previous section, the monotonicity of the Almgren-Poon quotient along global solutions to \eqref{eq:ExtensionEquationLocalBackward} (with $q = 0$) naturally leads to the spectral analysis of the following eigenvalue problem
\begin{equation}\label{eq:ORNSTUHLENFOREXTENSIONNEUMANN}
\begin{cases}
-\mathcal{O}_aV = \kappa V \quad &\text{in } \RR^{N+1}_+ \\
-\partial_y^a V = 0 \quad &\text{in } \RR^N\times\{0\},
\end{cases}
\end{equation}
where $\mathcal{O}_a$ is defined in \eqref{eq:ORNSTEINUHLENEBCKOPERATOR}. Together with the solutions to problem \eqref{eq:ORNSTUHLENFOREXTENSIONNEUMANN}, we will study also the boundary value problem
\begin{equation}\label{eq:ORNSTUHLENFOREXTENSION}
\begin{cases}
-\mathcal{O}_aV = \kappa V \quad &\text{in } \RR^{N+1}_+ \\
V = 0 \quad &\text{in } \RR^N\times\{0\},
\end{cases}
\end{equation}
which is \eqref{eq:ORNSTUHLENFOREXTENSIONNEUMANN} with homogeneous Dirichlet boundary conditions, and the problem
\begin{equation}\label{eq:ORNSTUHLENFOREXTENSION1}
-\frac{1}{|y|^a \mathcal{G}_a} \nabla\cdot(|y|^a \mathcal{G}_a \nabla V) = \kappa V \quad \text{in } \RR^{N+1}.
\end{equation}
We will see that solutions to \eqref{eq:ORNSTUHLENFOREXTENSION1} which are even w.r.t. $y \in \mathbb{R}$ correspond to solutions to \eqref{eq:ORNSTUHLENFOREXTENSIONNEUMANN}, while odd solutions to \eqref{eq:ORNSTUHLENFOREXTENSION1} (w.r.t. to $y \in \mathbb{R}$) correspond to solutions to \eqref{eq:ORNSTUHLENFOREXTENSION}. Thus, the study of \eqref{eq:ORNSTUHLENFOREXTENSION1} will give us all the information to describe both the spectrum of \eqref{eq:ORNSTUHLENFOREXTENSIONNEUMANN} and \eqref{eq:ORNSTUHLENFOREXTENSION}, and viceversa.

\begin{defn}\label{DEFINITIONWEAKEIGENFUNCTIONS}
We proceed with three different definitions. We refer to Section \ref{Section:NotPrel} for the definitions of the Gaussian spaces.

\noindent $\bullet$ A nontrivial function $V \in H^1(\RR_+^{N+1},d\mu)$ is said to be a weak eigenfunction to problem \eqref{eq:ORNSTUHLENFOREXTENSIONNEUMANN} with eigenvalue $\kappa \in \RR$ if
\begin{equation}\label{eq:DEFINITIONWEAKEIGENFUNCTIONNEUMANN}
\int_{\RR_+^{N+1}} \nabla V \cdot \nabla \eta \,d\mu = \kappa \int_{\RR_+^{N+1}} V \eta \,d\mu, \quad \text{ for all } \eta \in H^1(\RR_+^{N+1},d\mu).
\end{equation}

\noindent $\bullet$ A nontrivial function $V \in H_0^1(\RR_+^{N+1},d\mu)$ is said to be a weak eigenfunction to problem \eqref{eq:ORNSTUHLENFOREXTENSION} with eigenvalue $\kappa \in \RR$ if
\begin{equation}\label{eq:WEAKDEFINITIONOFORNSTEINUHLENBECK}
\int_{\RR_+^{N+1}} \nabla V \cdot \nabla \eta \,d\mu = \kappa \int_{\RR_+^{N+1}} V \eta \,d\mu, \quad \text{ for all } \eta \in H_0^1(\RR_+^{N+1},d\mu).
\end{equation}

\noindent $\bullet$ A nontrivial function $V \in H^1(\RR^{N+1},d\mu)$ is said to be a weak eigenfunction to problem \eqref{eq:ORNSTUHLENFOREXTENSION1} with eigenvalue $\kappa \in \RR$ if
\begin{equation}\label{eq:WEAKFORMULATIONEIGENVALUEPROBLEM1}
\int_{\RR^{N+1}} \nabla V \cdot \nabla \eta \,d\mu = \kappa \int_{\RR^{N+1}} V \eta \,d\mu, \quad \text{ for all } \eta \in H^1(\RR^{N+1},d\mu).
\end{equation}
\end{defn}
Before moving forward, we point out that using the fact that $\nabla \mathcal{G}_a(X,1) = -\frac{X}{2} \, \mathcal{G}_a(X,1)$, a simple integration by parts shows that classical eigenfunctions to \eqref{eq:ORNSTUHLENFOREXTENSIONNEUMANN}, \eqref{eq:ORNSTUHLENFOREXTENSION} and \eqref{eq:ORNSTUHLENFOREXTENSION1} (belonging to $H^1_{\mu}$) are a weak eigenfunctions to \eqref{eq:ORNSTUHLENFOREXTENSIONNEUMANN}, \eqref{eq:ORNSTUHLENFOREXTENSION} and \eqref{eq:ORNSTUHLENFOREXTENSION1}, respectively. Moreover, it is easily seen that smooth weak eigenfunctions are classical eigenfunctions. Notice that in both \eqref{eq:ORNSTUHLENFOREXTENSIONNEUMANN}/\eqref{eq:ORNSTUHLENFOREXTENSION} and \eqref{eq:ORNSTUHLENFOREXTENSION1} we obtain the classical Ornstein-Uhlenbeck eigenvalue problem by taking $a = 0$.

In the following theorem we completely characterize the spectrum of the above problems.
%
%
%%%%%%%%%%%%%%%%%%%%%%%%%%%%%%%%%%%%%%%%%%%%%%%%%%%%%%%%%%
%
%
\begin{thm}\label{SPECTRALTHEOREMEXTENDEDORNUHL1}
Fix $a \in (-1,1)$. Then the following three assertions hold:

\smallskip

(i) The set of eigenvalues of the homogeneous Neumann problem \eqref{eq:ORNSTUHLENFOREXTENSIONNEUMANN} is
\[
\{\widetilde{\kappa}_{n,m}\}_{n,m \in \NN}, \qquad \text{ where } \qquad \widetilde{\kappa}_{n,m} := \frac{n}{2} + m, \quad n,m \in \NN,
\]
with finite geometric multiplicity. For any $n_0,m_0 \in \NN$, the eigenspaces are given by
\[
\mathcal{V}_{n_0,m_0} = \text{span}\left\{V_{\alpha,m}(x,y) = H_{\alpha}(x)L_{(\frac{a-1}{2}),m}(y^2/4) : (\alpha,m) \in \widetilde{J}_0 \right\},
\]
where
\[
\widetilde{J}_0 := \left\{ (\alpha,m) \in \ZZ_{\geq 0}^N\times\NN : |\alpha| = n \in \NN \text{ and } \widetilde{\kappa}_{n,m} = \widetilde{\kappa}_{n_0,m_0} \right\},
\]
while $H_{\alpha}(\cdot)$ is a N-dimensional Hermite polynomial of order $|\alpha|$, while $L_{(\frac{a-1}{2}),m}(\cdot)$ is the $m^{th}$ Laguerre polynomial of order $(a-1)/2$. Furthermore, the set of eigenfunctions $\{V_{\alpha,m}\}_{\alpha,m}$ is an orthogonal basis of $L^2(\RR_+^{N+1},d\mu)$.

\smallskip

(ii) The set of eigenvalues of the homogeneous Dirichlet problem \eqref{eq:ORNSTUHLENFOREXTENSION} is
\[
\{\widehat{\kappa}_{n,m}\}_{n,m \in \NN}, \qquad \text{ where } \qquad \widehat{\kappa}_{n,m} := \frac{n}{2} + m + \frac{1-a}{2}, \quad n,m \in \NN,
\]
with finite geometric multiplicity. For all $n_0,m_0 \in \NN$, the eigenspaces are given by
\[
\mathcal{V}_{n_0,m_0} = \text{span}\left\{V_{\alpha,m}(x,y) = H_{\alpha}(x)\,y^{1-a}L_{(\frac{1-a}{2}),m}(y^2/4) : (\alpha,n) \in J_0 \right\},
\]
where
\[
\widehat{J}_0 := \left\{ (\alpha,m) \in \ZZ_{\geq 0}^N\times\NN : |\alpha| = n \in \NN \text{ and } \widehat{\kappa}_{n,m} = \widehat{\kappa}_{n_0,m_0} \right\},
\]
while now $L_{(\frac{1-a}{2}),m}(\cdot)$ is the $m^{th}$ Laguerre polynomial of order $(1-a)/2$. Again, the set of eigenfunctions $\{V_{\alpha,m}\}_{\alpha,m}$ is an orthogonal basis of $L^2(\RR_+^{N+1},d\mu)$.

\smallskip

(iii) The set of eigenvalues of problem \eqref{eq:ORNSTUHLENFOREXTENSION1} is
\[
\{\kappa_{n,m}\}_{n,m \in \NN} = \{\widehat{\kappa}_{n,m}\}_{n,m \in \NN}\cup\{\widetilde{\kappa}_{n,m}\}_{n,m \in \NN},
\]
with finite geometric multiplicity ($\widetilde{\kappa}_{n,m}$ and $\widehat{\kappa}_{n,m}$ are defined in part (i) and (ii), respectively). For any $n_0,m_0 \in \NN$, the eigenspaces corresponding to $\widehat{\kappa}_{n_0,m_0}$ and $\widetilde{\kappa}_{n_0,m_0}$ are given by
\[
\begin{aligned}
\widetilde{\mathcal{V}}_{n_0,m_0} &= \text{span}\left\{\widetilde{V}_{\alpha,m}(x,y) = H_{\alpha}(x) L_{(\frac{a-1}{2}),m}(y^2/4) : (\alpha,m) \in \widetilde{J}_0\right\}, \\
\widehat{\mathcal{V}}_{n_0,m_0} &= \text{span}\left\{\widehat{V}_{\alpha,m}(x,y) = H_{\alpha}(x) \,y|y|^{-a}L_{(\frac{1-a}{2}),m}(y^2/4) : (\alpha,m) \in \widehat{J}_0 \right\},
\end{aligned}
\]
respectively, where $\widetilde{J}_0$ and $\widehat{J}_0$ are defined in part (i) and (ii), respectively. Finally, similarly to the previous cases, the set $\{\widetilde{V}_{\alpha,m}(x,y)\}_{(\alpha,m)}\cup\{\widehat{V}_{\alpha,m}(x,y)\}_{(\alpha,m)}$ is an orthogonal basis of $L^2(\RR^{N+1},d\mu)$.
\end{thm}
%
%
%%%%%%%%%%%%%%%%%%%%%%%%%%%%%%%%%%%%%%%%%%%%%%%%%%%%%%%%%%%%%%
%
%
The proof is based on a separation of variables approach, together with some known results about Hermite and Laguerre polynomials (see  also with \cite{FelliPrimo1949:art} from which we borrow some ideas). The above theorem is crucial in the blow-up classification and we anticipate that, as a corollary, we will obtain some Gaussian-Poincar\'e type inequalities with \emph{optimal constants}. Similar inequalities were known since a long time (see  \cite{BatheRoberto2003:art,Muckenhoupt1972:art}) but, to the best of our knowledge, our result is new and it will play an important role in the proofs of some Liouville type theorems.

The spectral analysis of problems \eqref{eq:ORNSTUHLENFOREXTENSIONNEUMANN}, \eqref{eq:ORNSTUHLENFOREXTENSION}, and \eqref{eq:ORNSTUHLENFOREXTENSION1}, is based on the following one-dimensional eigenvalue problem (where we set $\psi' = d\psi/dy$ for simplicity)
\begin{equation}\label{eq:HALFLINEORNUHLEN}
-y^{-a}\left(y^a \psi'\right)' + \frac{y}{2} \, \psi' = \sigma \psi, \quad y > 0,
\end{equation}
which is an equivalent way to write
\begin{equation}\label{eq:HALFLINEORNUHLENSECONDVERSION}
-\frac{1}{y^a G_{a+1}} \left(y^a G_{a+1}\, \psi' \right)' = \sigma \psi, \quad y > 0,
\end{equation}
where $G_{a+1} = G_{a+1}(y) = \mathcal{G}_a(0,y,1)$,  see  Section \ref{Section:NotPrel}). Miming the steps followed before, we consider the following eigenvalue problems (corresponding to \eqref{eq:ORNSTUHLENFOREXTENSIONNEUMANN}, \eqref{eq:ORNSTUHLENFOREXTENSION} and \eqref{eq:ORNSTUHLENFOREXTENSION1}):
\begin{equation}\label{eq:ONEDIMORNSUHLENNEUMANN}
\begin{cases}
-y^{-a}\left(y^a \psi'\right)' + \frac{y}{2} \, \psi' = \sigma \psi, \quad &\text{for } y > 0,  \\
-\partial_y^a\psi = 0           \quad  &\text{for } y = 0,
\end{cases}
\end{equation}
\begin{equation}\label{eq:ONEDIMORNSUHLEN}
\begin{cases}
-y^{-a}\left(y^a \psi'\right)' + \frac{y}{2} \, \psi' = \sigma \psi, \quad &\text{for } y > 0,  \\
\psi = 0           \quad  &\text{for } y = 0,
\end{cases}
\end{equation}
and
\begin{equation}\label{eq:ONEDIMORNSUHLEN1}
-|y|^{-a}\left(|y|^a \psi'\right)' + \frac{y}{2} \, \psi' = \sigma \psi, \quad y \not= 0.
\end{equation}
\begin{defn}\label{DEFINITIONWEAKEIGENFUNCTIONSONED}
Again we proceed in three cases.

\noindent $\bullet$ A nontrivial function $\psi \in H^1(\RR_+,d\mu_y)$ is said to be a weak eigenfunction to problem \eqref{eq:ONEDIMORNSUHLENNEUMANN} with eigenvalue $\sigma \in \RR$ if
\begin{equation}\label{eq:ONEDDEFWEAKEIGENFUNCTIONNEUMANN}
\int_{\RR_+} \psi' \eta' \,d\mu_y = \sigma \int_{\RR_+} \psi \eta \,d\mu_y, \quad \text{ for all } \eta \in H^1(\RR_+,d\mu_y).
\end{equation}

\noindent $\bullet$ A nontrivial function $\psi \in H_0^1(\RR_+,d\mu_y)$ is said to be a weak eigenfunction to problem \eqref{eq:ONEDIMORNSUHLEN} with eigenvalue $\sigma \in \RR$ if
\begin{equation}\label{eq:ONEDWEAKDEFORNSTEINUHLENBECK}
\int_{\RR_+} \psi' \eta' \,d\mu_y = \sigma \int_{\RR_+} \psi \eta \,d\mu_y, \quad \text{ for all } \eta \in H_0^1(\RR_+,d\mu_y).
\end{equation}

\noindent $\bullet$ A nontrivial function $\psi \in H^1(\RR,d\mu_y)$ is said to be a weak eigenfunction to problem \eqref{eq:ONEDIMORNSUHLEN1} with eigenvalue $\sigma \in \RR$ if
\begin{equation}\label{eq:ONEDWEAKEIGENVALUEPROBLEM1}
\int_{\RR} \psi' \eta' \,d\mu_y = \sigma \int_{\RR} \psi \eta \,d\mu_y, \quad \text{ for all } \eta \in H^1(\RR,d\mu_y).
\end{equation}
\end{defn}
Notice that in the previous definition we employ the probability measure $d\mu_y := G_{a+1}(y)dy$ in problems \eqref{eq:ONEDDEFWEAKEIGENFUNCTIONNEUMANN} and \eqref{eq:ONEDWEAKDEFORNSTEINUHLENBECK}, while its even extension in \eqref{eq:ONEDWEAKEIGENVALUEPROBLEM1} (see  Section \ref{Section:NotPrel} again). We begin with the following lemma.

\begin{lem}\label{PROPEIGENVALUESONEDIMNEUMANNDIRICHLET}
Fix $a \in (-1,1)$. Then the following two assertions hold:

(i) The eigenvalues and the weak eigenfunctions to the homogeneous Neumann problem \eqref{eq:ONEDIMORNSUHLENNEUMANN} are
\[
\widetilde{\sigma}_m = m \in \NN, \qquad \widetilde{\psi}_m(y) = \widetilde{A}_m L_{(\frac{a-1}{2}),m}\left(\frac{y^2}{4}\right), \quad y > 0,
\]
where $\widetilde{A}_m \in \RR$ is arbitrary and $L_{(\frac{a-1}{2}),m}(\cdot)$ is the $m^{th}$ Laguerre polynomial of order $(a-1)/2$.

(ii) The eigenvalues and the weak eigenfunctions to the homogeneous Dirichlet problem \eqref{eq:ONEDIMORNSUHLEN} are
\[
\widehat{\sigma}_m = \frac{1-a}{2} + m \in \NN, \qquad \widehat{\psi}_m(y) = \widetilde{A}_m y^{1-a} L_{(\frac{1-a}{2}),m}\left(\frac{y^2}{4}\right), \quad y > 0,
\]
where $\widetilde{A}_m \in \RR$ is arbitrary and $L_{(\frac{1-a}{2}),m}(\cdot)$ is the $m^{th}$ Laguerre polynomial of order $(1-a)/2$. Finally, both sets $\{\widetilde{\psi}_m\}_{m \in \NN}$ and $\{\widehat{\psi}_m\}_{m \in \NN}$ are orthogonal basis of $L^2(\RR_+,d\mu_y)$.
\end{lem}
\emph{Proof.} We carry out a detailed analysis of the spectrum of equation \eqref{eq:HALFLINEORNUHLEN}. In order to work in the framework as general as possible, we do not impose boundary conditions at $y = 0$ but only the following integrability conditions (which are required in our definition)
\begin{equation}\label{eq:INTEGRABILITYCONDITIONSONPSI}
\int_0^{\infty} \psi^2(y) \, y^a e^{-\frac{y^2}{4}} \,dy < +\infty, \qquad \int_0^{\infty} \left|\psi'(y)\right|^2 \, y^ae^{-\frac{y^2}{4}} \,dy < +\infty.
\end{equation}
Our procedure will naturally distinguish the solutions with Dirichlet and/or Neumann boundary conditions.

So, let us set $\psi(y) = \zeta (y^2/4)$. It is easily seen that the equation for $\zeta = \zeta(r)$, $r = y^2/4$ is:
\begin{equation}\label{eq:KUMMERCONFHYPEREQUATIONEIGENVALUE11}
r\frac{d^2\zeta}{dr^2} + \left(1 + \frac{a-1}{2} - r\right)\frac{d\zeta}{dr} + \sigma \zeta = 0, \quad r > 0,
\end{equation}
which can be seen as a Kummer Confluent Hypergeometric type equation, with
\[
b_1 = -\sigma \qquad \text{ and } \qquad b_2 = 1 + \frac{a-1}{2} = \frac{1+a}{2},
\]
and/or a Laguerre equation with
\[
\alpha = \frac{a-1}{2} > -1.
\]
A detail report about these topics is given in \cite{AbramowitzStegun1972:art}, see also the classical reference \cite{Szego1939:art}. We know that all solutions are given by
\begin{equation}\label{eq:ALLSOLUTIONSTOKUMMERCONFLUENT}
\zeta(r) = A_1 M\left(-\sigma,\frac{1+a}{2}, r\right) + A_2 \widetilde{M}\left(-\sigma,\frac{1+a}{2}, r\right), \quad r > 0,
\end{equation}
where $A_1, A_2 \in \RR$, and $M(\cdot,\cdot,\cdot)$ and $\widetilde{M}(\cdot,\cdot,\cdot)$ are the Kummer and the Tricomi functions, as explained in the appendix mentioned above. We divide the analysis in three cases:

\noindent $\bullet$ \emph{Case $\sigma \in \NN$.}  As explained in \cite{Szego1939:art}, when $\sigma = m \in \NN = \{0,1,\ldots \}$, then to each $\sigma = m$ it corresponds a unique solution (up to multiplicative constants) to \eqref{eq:KUMMERCONFHYPEREQUATIONEIGENVALUE11} given by the $m^{th}$ Laguerre polynomial of order $(a-1)/2$:
\[
\zeta_m (r) = L_{(\frac{a-1}{2}),m}(r), \quad r > 0, \quad m \in \NN,
\]
and so, for any $\sigma = m \in \NN$, we obtain the solutions to \eqref{eq:HALFLINEORNUHLEN}, given by
\begin{equation}\label{eq:SOLUTIONTOKUMMERFORMUEQUALTOM}
\widetilde{\psi}_m (y) = \widetilde{A}_m \zeta_m\left(\frac{y^2}{4}\right) = \widetilde{A}_m L_{(\frac{a-1}{2}),m}\left(\frac{y^2}{4}\right), \quad y > 0,
\end{equation}
where $\widetilde{A}_m$ is an arbitrary real constant. Note that since $\{L_{(\frac{a-1}{2}),m}\}_{m\in\NN}$ is an orthogonal basis for the space
\[
L^2(\RR_+,d\nu) \qquad \text{ where } \qquad d\nu(r) = r^{\frac{a-1}{2}} e^{-r} dr,
\]
it is immediate to see that $\{\widetilde{\psi}_m\}_{m \in \NN}$ is an orthogonal basis for $L^2(\RR_+,d\mu_y)$.
\\
Similarly, one could have used formula \eqref{eq:ALLSOLUTIONSTOKUMMERCONFLUENT} for the explicit expression of all solutions to \eqref{eq:KUMMERCONFHYPEREQUATIONEIGENVALUE11} and note that for $\sigma = m \in \NN$, we have
\[
M\left(-m, \frac{1+a}{2}, r\right) := \widetilde{Q}_m(r) = \sum_{j=0}^m \frac{(-m)_j}{(\frac{1+a}{2})_j} \frac{r^j}{j!}, \quad r > 0,
\]
while (see  formula 13.1.3 of \cite{AbramowitzStegun1972:art})
\[
\begin{aligned}
\widetilde{M}\left(-m,\frac{1+a}{2},r\right) & = \frac{\pi}{\sin \left(\frac{1+a}{2}\pi\right)} \left[ \frac{M(-m,\frac{1+a}{2},r)}{\Gamma(\frac{1-a}{2}-m)\Gamma(\frac{1+a}{2})} - r^{\frac{1-a}{2}}\frac{M(\frac{1-a}{2}-m,1+\frac{1-a}{2},r)}{\Gamma(-m)\Gamma(1+\frac{1-a}{2})}\right] \\
& = \frac{\Gamma(\frac{1-a}{2})}{\Gamma(\frac{1-a}{2}-m)} \, \widetilde{Q}_m(r) , \quad r > 0,
\end{aligned}
\]
where we have used the well-known properties of the Gamma function
\[
\frac{1}{\Gamma(-z)} = 0, \quad z \in \NN  \qquad \text{ and } \qquad \frac{\pi}{\sin \left(z\pi\right)} = \Gamma(1-z)\Gamma(z), \quad z \not\in \ZZ.
\]
Consequently, it follows an equivalent expression for $\widetilde{\psi}_m = \widetilde{\psi}_m(y)$ given by
\[
\widetilde{\psi}_m(y) = \widetilde{A}_m \widetilde{Q}_m\left(\frac{y^2}{4}\right) = \widetilde{A}_m \sum_{j=0}^m \frac{(-m)_j}{4^j j! (\frac{1+a}{2})_j} y^{2j}, \quad y > 0,
\]
which is equivalent to the expression found in \eqref{eq:SOLUTIONTOKUMMERFORMUEQUALTOM} in view of formula 22.5.54 of \cite{AbramowitzStegun1972:art}. In what follows we will adopt the notation used in \eqref{eq:SOLUTIONTOKUMMERFORMUEQUALTOM}, but always keeping in mind that
\begin{equation}\label{eq:LaguerreExpression}
L_{(\frac{a-1}{2}),m}\left(\frac{y^2}{4}\right) = \sum_{j=0}^m \frac{(-m)_j}{4^j j! (\frac{1+a}{2})_j} y^{2j}, \quad y > 0,
\end{equation}
up to a multiplicative constant.

\noindent $\bullet$ \emph{Case $\sigma \not\in \NN$ but $\sigma - (1-a)/2 \in \NN$.} If $\sigma = (1-a)/2 + m$, $m \in \NN$, it follows
\[
\begin{aligned}
\widetilde{M}\left(-\sigma,\frac{1+a}{2},r\right) & = \frac{\pi}{\sin \left(\frac{1+a}{2}\pi\right)}
\left[ \frac{M(-\sigma,\frac{1+a}{2},r)}{\Gamma(-m)\Gamma(\frac{1+a}{2})} - r^{\frac{1-a}{2}}\frac{M(-m,1+\frac{1-a}{2},r)}{\Gamma(-\sigma)\Gamma(1+\frac{1-a}{2})}\right] \\
& = -\frac{2\Gamma(\frac{1+a}{2})}{(1-a)\Gamma(-\frac{1-a}{2}-m)} \, r^{\frac{1-a}{2}} \widetilde{P}_m(r), \quad r > 0,
\end{aligned}
\]
where we have employed again the properties of the Gamma function and we have defined
\[
\widetilde{P}_m(r) = M\left(-m,1+\frac{1-a}{2},r\right) =  \sum_{j=0}^m \frac{(-m)_j}{(1+ \frac{1-a}{2})_j} \frac{r^j}{j!}, \quad r > 0.
\]
On the other hand,
\[
M\left(-\sigma, \frac{1+a}{2}, r\right) =  \sum_{j=0}^{\infty} \frac{(-\sigma)_j}{(\frac{1+a}{2})_j}. \frac{r^j}{j!} \sim \frac{\Gamma(\frac{1+a}{2})}{\Gamma(-\sigma)}\, e^r r^{-\sigma - \frac{1+a}{2}}, \quad \text{ for } r \sim +\infty,
\]
since $\sigma \not\in \NN$ (see  formula 13.1.4 of \cite{AbramowitzStegun1972:art}). From the above asymptotic expansion for $r \sim +\infty$, it is not difficult to see that the function
\[
y \to M\left(-\sigma, \frac{1+a}{2}, \frac{y^2}{4}\right)
\]
does not satisfies the first bound in \eqref{eq:INTEGRABILITYCONDITIONSONPSI} and so, we have to take $A_1 = 0$. Consequently, for $\sigma = (1-a)/2 + m$, $m \in \NN$, we deduce that
\[
\zeta_{a,m}(r) = \widetilde{A}_m \, r^{\frac{1-a}{2}} \widetilde{P}_m (r) = \widetilde{A}_m \sum_{j=0}^m \frac{(-m)_j}{j! (1+\frac{1-a}{2})_j} r^{\frac{1-a}{2} + j}, \quad r > 0,
\]
and so, coming back to the variable $y$ (we recall that $r = y^2/4$),
\[
\widehat{\psi}_m(y) = \widetilde{A}_m \, y^{1-a} \widetilde{P}_m \left(\frac{y^2}{4}\right) = \widetilde{A}_m \sum_{j=0}^m \frac{(-m)_j}{4^j j! (1+\frac{1-a}{2})_j} y^{1+2j-a}, \quad y > 0.
\]
Exactly as before, we can recall formulas 13.5.6, 13.5.7, 13.5.8 of \cite{AbramowitzStegun1972:art}, to deduce $L_{(\frac{1-a}{2}),m}(r) = M(-m,1+\frac{1-a}{2},r) = \widetilde{P}_m(r)$ (up to a multiplicative constant), and so
\begin{equation}\label{eq:EXPRESSIONFORPHIMALAGUERREPOL}
\widehat{\psi}_m(y) = \widetilde{A}_m \, y^{1-a} L_{(\frac{1-a}{2}),m}\left(\frac{y^2}{4}\right), \quad y > 0.
\end{equation}
Again, since $\{L_{(\frac{1-a}{2}),m}\}_{m\in\NN}$ is an orthogonal basis for the space (note the difference in the measure $\nu$ with respect to the previous case)
\[
L^2(\RR_+,d\nu) \qquad \text{ where now } \qquad d\nu(r) = r^{\frac{1-a}{2}} e^{-r} dr,
\]
it follows again that $\{\widehat{\psi}_m\}_{m \in \NN}$ is an orthogonal basis for $L^2(\RR_+,d\mu_y)$.

\noindent $\bullet$ \emph{Case $\sigma \not\in \NN$ and $\sigma - (1-a)/2 \not\in \NN$.} Proceeding as before, we obtain
\[
\zeta(r) = \left( A_1 + \frac{\Gamma(\frac{1-a}{2})}{\Gamma(\frac{1-a}{2}-\sigma)} A_2 \right) \sum_{j = 0}^{\infty} \frac{(-\sigma)_j}{(\frac{1+a}{2})_j} \frac{r^j}{j!}
-\frac{2\Gamma(\frac{1+a}{2})}{(1-a)\Gamma(-\sigma)} A_2 \, r^{\frac{1-a}{2}}\sum_{j = 0}^{\infty} \frac{(\frac{1-a}{2}-\sigma)_j}{(1 + \frac{1-a}{2})_j} \frac{r^j}{j!},
\]
and so, since the coefficient in front of the first series as to be zero (for the same reason of the above case), we deduce
\[
\begin{aligned}
\zeta(r) & = \frac{2\Gamma(\frac{1+a}{2}) \Gamma(\frac{1-a}{2}-\sigma) }{(1-a)\Gamma(\frac{1-a}{2})\Gamma(-\sigma)} A_2 \, r^{\frac{1-a}{2}}\sum_{j = 0}^{\infty} \frac{(\frac{1-a}{2}-\sigma)_j}{(1 + \frac{1-a}{2})_j} \frac{r^j}{j!} = \overline{A}_2 r^{\frac{1-a}{2}} M\left(\frac{1-a}{2}-\sigma,1 + \frac{1-a}{2},r\right) \\
& \sim \frac{\Gamma(1 + \frac{1-a}{2})}{\Gamma(\frac{1-a}{2}-\sigma)}\overline{A}_2 \, e^r r^{-1 - \sigma}, \quad \text{for } r \sim +\infty.
\end{aligned}
\]
Exactly as before, the above expansion for $r \sim +\infty$ tells us that $\psi(y) = \zeta(y^2/4)$ does not satisfies the first bound in \eqref{eq:INTEGRABILITYCONDITIONSONPSI} and so we to take $\overline{A_2} = 0$, i.e. $A_2 = A_1 = 0$.

\noindent $\bullet$ \emph{Conclusions.} From the analysis carried out, it follows that the set of eigenvalues for equation \eqref{eq:HALFLINEORNUHLEN} is given by
\[
\{\widetilde{\sigma}_m\}_{m \in \NN} \cup \{\widehat{\sigma}_m\}_{m \in \NN} \qquad \text{ where } \qquad \widetilde{\sigma}_m = m, \qquad \widehat{\sigma}_m = \frac{1-a}{2} + m,
\]
and to the eigenvalue $\widetilde{\sigma}_m = m$ it corresponds (up to multiplicative constants) the eigenfunction
\[
\widetilde{\psi}_m (y) = \widetilde{A}_m L_{(\frac{a-1}{2}),m}\left(\frac{y^2}{4}\right), \quad y > 0,
\]
where $L_{(\frac{a-1}{2}),m}(\cdot)$ it the $m^{th}$ Laguerre polynomial of order $(a-1)/2$, whilst the eigenvalue $\widehat{\sigma}_m = (1-a)/2 + m$ possesses (up to multiplicative constants) the eigenfunction
\[
\widehat{\psi}_m (y) = y^{1-a}L_{(\frac{1-a}{2}),m}\left(\frac{y^2}{4}\right), \quad y > 0.
\]
We point out that the existence of two distinct orthogonal basis of eigenfunctions of $L^2(\RR_+,d\mu_y)$ comes from the fact that $\{\widetilde{\psi}_m\}_{m\in\NN}$ and $\{\widehat{\psi}_m\}_{m \in \NN}$ are solutions to different eigenvalue problems. Indeed, for all $m \in \NN$, we have
\[
\widetilde{\psi}_m(0) \not= 0 \qquad \text{ and } \qquad \partial_y^a\widetilde{\psi}_m = 0,
\]
while
\[
\widehat{\psi}_m(0) = 0 \qquad \text{ and } \qquad \partial_y^a\widehat{\psi}_m \not= 0,
\]
which mean that the functions $\widetilde{\psi}_m$ correspond to the eigenfunctions to equation \eqref{eq:HALFLINEORNUHLEN} with homogeneous Neumann boundary condition at $y = 0$, i.e., to problem \eqref{eq:ONEDIMORNSUHLENNEUMANN}, while $\widehat{\psi}_m$ are the eigenfunctions to equation \eqref{eq:HALFLINEORNUHLEN} with homogeneous Dirichlet boundary condition at $y = 0$, i.e., to problem \eqref{eq:ONEDIMORNSUHLEN}. Note that it is not hard to verify that $\widetilde{\psi}_m$ and $\widehat{\psi}_m$ are weak eigenfunctions to \eqref{eq:ONEDIMORNSUHLENNEUMANN} and \eqref{eq:ONEDIMORNSUHLEN}, respectively, i.e. they satisfy \eqref{eq:ONEDDEFWEAKEIGENFUNCTIONNEUMANN} and \eqref{eq:ONEDWEAKDEFORNSTEINUHLENBECK}, respectively. This easily follows by an integration by parts and noting that $\widetilde{\psi}_m$ and $\widehat{\psi}_m$ satisfy both bounds in \eqref{eq:INTEGRABILITYCONDITIONSONPSI}. The proof is completed. $\Box$
\begin{lem}\label{PROPEIGENVALUESONENDIM}
Fix $a \in (-1,1)$. Then the set of eigenvalues of problem \eqref{eq:ONEDIMORNSUHLEN1} is
\[
\{\widetilde{\sigma}_m\}_{m\in\NN} \cup \{\widehat{\sigma}_m\}_{m\in\NN},
\]
where $\widetilde{\sigma}_m$ and $\widehat{\sigma}_m$ correspond to the eigenvalues of the Neumann problem \eqref{eq:ONEDIMORNSUHLENNEUMANN} and of the Dirichlet on \eqref{eq:ONEDIMORNSUHLEN}, respectively, and are defined in Lemma \ref{PROPEIGENVALUESONEDIMNEUMANNDIRICHLET}. Moreover, the eigenfunctions corresponding to $\widetilde{\sigma}_m$ are
\[
\widetilde{\psi}_m(y) = \widetilde{A}_m L_{(\frac{a-1}{2}),m}\left(\frac{y^2}{4}\right), \quad y \in \RR,
\]
where $\widetilde{A}_m \in \RR$ is arbitrary and $L_{(\frac{a-1}{2}),m}(\cdot)$ is the $m^{th}$ Laguerre polynomial of order $(a-1)/2$, while the eigenfunctions corresponding to $\widehat{\sigma}_m$ are
\[
\widehat{\psi}_m(y) = \widetilde{A}_m \, y |y|^{-a} L_{(\frac{1-a}{2}),m}\left(\frac{y^2}{4}\right), \quad y \in \RR,
\]
where $\widetilde{A}_m \in \RR$ is arbitrary and $L_{(\frac{1-a}{2}),m}(\cdot)$ is the $m^{th}$ Laguerre polynomial of order $(1-a)/2$. Finally, the set $\{\widetilde{\psi}_m\}_{m \in \NN} \cup \{\widehat{\psi}_m\}_{m \in \NN}$ is an orthogonal basis of $L^2(\RR,d\mu_y)$.
\end{lem}
\emph{Proof.} Let us consider equation \eqref{eq:ONEDIMORNSUHLEN1}. We first notice that if $\psi = \psi(y)$ is a weak eigenfunction to equation \eqref{eq:ONEDIMORNSUHLEN1}, then $\widetilde{\psi}(y) = \psi(-y)$ is a weak eigenfunction too, and so are all their linear combinations. Consequently, it suffices to consider solutions $\psi^+ = \psi^+(y)$ to equation \eqref{eq:HALFLINEORNUHLEN} (already analyzed before) and then using their reflections (with respect to $y = 0$) to obtain all solutions to \eqref{eq:ONEDIMORNSUHLEN1} defined for all $y \in \RR$.

\noindent From the analysis carried out in the proof of Lemma \ref{PROPEIGENVALUESONEDIMNEUMANNDIRICHLET}, we that the only admissible values for $\sigma$ are
\[
\widetilde{\sigma}_m := m \quad \text{ with eigenfunctions } \quad \widetilde{\psi}_m^+(y) = \widetilde{A}_m L_{(\frac{a-1}{2}),m}\left(\frac{y^2}{4}\right), \quad m \in \NN,
\]
defined for all $y > 0$, satisfying $\partial_y^a \psi_{m}^+ = 0$ and both bounds in \eqref{eq:INTEGRABILITYCONDITIONSONPSI}, and
\[
\widehat{\sigma}_m := \frac{1-a}{2} + m \quad \text{ with eigenfunctions } \quad \widehat{\psi}_m^+(y) = \widetilde{A}_m y^{1-a} L_{(\frac{1-a}{2}),m}\left(\frac{y^2}{4}\right), \quad m \in \NN,
\]
defined for all $y > 0$, satisfying $\widehat{\psi}_m^+(0) = 0$ and both bounds in \eqref{eq:INTEGRABILITYCONDITIONSONPSI}. As always $\widetilde{A}_m \in \RR$ and $m \in \NN$.

We begin by focusing on $\widehat{\psi}_m^+ = \widehat{\psi}_m^+(y)$. Clearly, there are two different ways to build a solution to equation \eqref{eq:ONEDIMORNSUHLEN1} with $\sigma = \widehat{\sigma}_m$, through reflection methods. We can consider both odd and even extensions:
\begin{equation}\label{eq:EXTENSIONOFPHI}
\widehat{\psi}_m^o(y) :=
\begin{cases}
\widehat{\psi}_m^+(y) \quad &\text{if } y \geq 0 \\
-\widehat{\psi}_m^+(|y|) \quad &\text{if } y < 0,
\end{cases}
\qquad \qquad
\widehat{\psi}_m^e(y) :=
\begin{cases}
\widehat{\psi}_m^+(y) \quad &\text{if } y \geq 0 \\
\widehat{\psi}_m^+(|y|) \quad &\text{if } y < 0.
\end{cases}
\end{equation}
Note that both $\widehat{\psi}_m^o$ and $\widehat{\psi}_m^e$ satisfy both bounds in \eqref{eq:INTEGRABILITYCONDITIONSONPSI} (this follows from the analysis done above). However, it is not difficult to see that the even extension (the second formula in \eqref{eq:EXTENSIONOFPHI}) produces a candidate eigenfunction $\widehat{\psi}_m^e = \widehat{\psi}_m^e(y)$ which is not a weak eigenfunction to problem \eqref{eq:ORNSTUHLENFOREXTENSION1}, i.e., it does not satisfy \eqref{eq:WEAKFORMULATIONEIGENVALUEPROBLEM1}. Indeed, using the equation of $\widehat{\psi}_m^e = \widehat{\psi}_m^e(y)$ it is not hard to show that
\begin{equation}\label{eq:EVENEXTENSIONNOTSOLUTION}
\int_{\RR} (\widehat{\psi}_m^e)'\eta ' \,d\mu_y = \widehat{\sigma}_m \int_{\RR} \widehat{\psi}_m^e \eta \,d\mu_y - 2(1-a)\eta(0), \quad \text{ for all } \eta \in H^1(\RR,d\mu_y),
\end{equation}
where we recall that $d\mu_y := |y|^a G_{a+1}(y)dy$. On the other hand, repeating the same procedure for the odd extension $\widehat{\psi}_m^o = \widehat{\psi}_m^o(y)$, we obtain that \eqref{eq:EVENEXTENSIONNOTSOLUTION} holds for $\widehat{\psi}_m^o$ without the extra term $2(1-a)\eta(0)$ (it cancels thanks to the oddness of $\widehat{\psi}_m^o$), and so $\widehat{\psi}_m^o = \widehat{\psi}_m^o(y)$ is an eigenfunction associated to the eigenvalue $\widehat{\sigma}_m$.

Repeating the same procedure for $\widetilde{\psi}_m^+ = \widetilde{\psi}_m^+(y)$ it is easily seen that its even reflection is a weak eigenfunction corresponding to the eigenvalue $\widetilde{\sigma}_m = m$, while its odd extension cannot be a weak eigenfunction since it has a jump discontinuity at $y = 0$.

Consequently, we can conclude that problem \eqref{eq:ONEDIMORNSUHLEN1} has eigenvalues and corresponding weak eigenfunctions defined by
\[
\widetilde{\sigma}_m := m \quad \qquad \widetilde{\psi}_m(y) = \widetilde{A}_m L_{(\frac{a-1}{2}),m}\left(\frac{y^2}{4}\right), \quad m \in \NN,
\]
\[
\widehat{\sigma}_m := \frac{1-a}{2} + m \quad \qquad \widehat{\psi}_m(y) = \widetilde{A}_m \,y|y|^{-a} L_{(\frac{1-a}{2}),m}\left(\frac{y^2}{4}\right), \quad m \in \NN,
\]
defined for all $y \in \RR$. We stress that, using the parity and the oddness of $\widetilde{\psi}_m = \widetilde{\psi}_m(y)$ and $\widehat{\psi}_m = \widehat{\psi}_m(y)$, respectively, we immediately see that (up to multiplicative constants)
\[
\int_{\RR} \widetilde{\psi}_{m}(y)\widetilde{\psi}_{n}(y)d\mu_y = \int_{\RR} \widehat{\psi}_{m}(y)\widehat{\psi}_{n}(y)d\mu_y = \delta_{m,n} \quad \text{ for all } m,n \in \NN,
\]
while
\[
\int_{\RR} \widetilde{\psi}_{m}(y)\widehat{\psi}_n(y)d\mu_y = 0 \quad \text{ for all } m,n \in \NN,
\]
so that it follows that the family $\{\psi_{m}\}_{m\in\NN}\cup\{\widehat{\psi}_m\}_{m\in\NN}$ is an orthogonal set of $L^2(\RR,d\mu_y)$. The fact that the family $\{\psi_{m}\}_{m\in\NN}\cup\{\widehat{\psi}_m\}_{m\in\NN}$ is a basis of $L^2(\RR,d\mu_y)$ follows since both $\{\widetilde{\psi}_m\}_{m\in\NN}$ and $\{\widehat{\psi}_m\}_{m\in\NN}$ are orthogonal basis of $L^2(\RR_+,d\mu_y)$. $\Box$

\bigskip

\emph{Proof of Theorem \ref{SPECTRALTHEOREMEXTENDEDORNUHL1}.} We begin by proving part (i) and (ii). We look for solutions to problems \eqref{eq:ORNSTUHLENFOREXTENSIONNEUMANN} and/or \eqref{eq:ORNSTUHLENFOREXTENSION} in separate variables form $V(x,y) = \varphi(x)\psi(y)$, $x \in \RR^N$, $y > 0$. So, substituting into \eqref{eq:ORNSTUHLENFOREXTENSIONNEUMANN} and/or \eqref{eq:ORNSTUHLENFOREXTENSION}, it is not difficult to see that
\[
\left(-\Delta_x\varphi + \frac{x}{2}\cdot\nabla_x\varphi \right)\psi + \left[-y^{-a}\left(y^a \psi'\right)' + \frac{y}{2} \, \psi' \right]\varphi = \kappa \varphi \psi \quad \text{in } \RR^{N+1}_+,
\]
where, as always, $\psi' = d\psi/dy$. Hence, we recover the eigenvalue $\kappa \in \RR$ as the sum of $\nu \in \RR$ and $\sigma \in \RR$ eigenvalues to
\begin{equation}\label{eq:CLASSICALORNUHLEN}
-\Delta_x\varphi + \frac{x}{2}\cdot\nabla_x\varphi = \nu \varphi, \quad x \in \RR^N
\end{equation}
and equation \eqref{eq:HALFLINEORNUHLEN}:
\[
-y^{-a}\left(y^a \psi'\right)' + \frac{y}{2} \, \psi' = \sigma \psi, \quad y > 0,
\]
respectively.

\smallskip

\emph{Step1: Analysis of \eqref{eq:CLASSICALORNUHLEN}.} Evidently \eqref{eq:CLASSICALORNUHLEN} is the classical Ornstein-Uhlenbeck eigenvalue problem in the all Euclidean space. It possesses the sequence of eigenvalues
\[
\nu_n = \frac{n}{2}, \quad n \in \NN = \{0,1,\ldots \},
\]
and a eigenfunction basis composed by the so-called Hermite polynomials
\[
H_{\alpha}(x) = H_{n_1}(x_1) \dots H_{n_N}(x_N), \quad x = (x_1,\ldots,x_N) \in \RR^N,
\]
where $\alpha = (n_1,\ldots,n_N) \in \ZZ_{\geq 0}^N$ and $H_{n_j}(\cdot)$, $j = 1,\ldots,N$, is the $n_j^{th}$ 1 dimensional Hermite polynomial (see  \cite{Bell2017:art,Szego1939:art}).

\smallskip

\emph{Step2: Analysis of \eqref{eq:HALFLINEORNUHLEN} and conclusion.} The complete analysis has been carried out above. Consequently, part (i) and (ii) of Theorem \ref{SPECTRALTHEOREMEXTENDEDORNUHL1} follow by combining the above mentioned facts about Hermite polynomials and classical Ornstein-Uhlenbeck eigenvalue problem, with the statement of Lemma \ref{PROPEIGENVALUESONEDIMNEUMANNDIRICHLET}. In particular, the fact that $\widetilde{V}_{\alpha,m}(x,y) = H_{\alpha}(x) L_{(\frac{a-1}{2}),m}(y)$ and $\widehat{V}_{\alpha,m}(x,y) = H_{\alpha}(x)\,y^{1-a}L_{(\frac{1-a}{2}),m}(y)$ form an orthogonal basis of $L^2(\RR_+^{N+1},d\mu)$ comes from the fact that both $H_{\alpha}(x)$ and $L_{(\frac{a-1}{2}),m}(y)$ (resp. $y^{1-a}L_{(\frac{1-a}{2}),m}(y)$) are orthogonal basis for $L^2(\RR^N,d\mu_x)$ and $L^2(\RR_+,d\mu_y)$, respectively, since the measure $d\mu = d\mu(x,y)$ defined on $\RR_+^{N+1}$ is obtained the measure product of $d\mu_x$ and $d\mu_y$.

\smallskip

Let us now focus on part (iii). Exactly as before we look for solutions to problem  \eqref{eq:ORNSTUHLENFOREXTENSION1} with form $V(x,y) = \varphi(x)\psi(y)$, $x \in \RR^N$, $y \in \RR$ and we obtain that $V = V(x,y)$ is an eigenfunction with eigenvalue $\kappa = \nu + \sigma$ if $\varphi = \varphi(x)$ satisfies \eqref{eq:CLASSICALORNUHLEN} and $\psi = \psi(y)$ satisfies \eqref{eq:ONEDIMORNSUHLEN1}:
\[
-|y|^{-a}\left(|y|^a \psi'\right)' + \frac{y}{2} \, \psi' = \sigma \psi, \quad y \not= 0.
\]
Since the analysis of the above equation has already been done, we deduce part (iii) of Theorem \ref{SPECTRALTHEOREMEXTENDEDORNUHL1} exactly as before, but in this case we apply Lemma \ref{PROPEIGENVALUESONENDIM} instead of Lemma \ref{PROPEIGENVALUESONEDIMNEUMANNDIRICHLET}. We end the proof by stressing that the family $\{\widetilde{V}_{\alpha,m}(x,y)\}_{(\alpha,m)}\cup\{\widehat{V}_{\alpha,m}(x,y)\}_{(\alpha,m)}$ is an orthogonal basis of $L^2(\RR^{N+1},d\mu)$. This follows exactly as before. $\Box$

%
%%%%%%%%%%%%%%%%%%%%%%%%%%%%%%%%%%%%%%%%%%%%%%%%%%%%%%%%%%%%%%%%%%%%%%%%%%%%%%%%%%%%%%%%%%%%%%%%%%%%%%%%%%%
%
%
%%%%%%%%%%%%%%%%%%%%%%%%%%%%%%%%%%%%%%%%%%%%%%%%%%%%%%%%%%%%%%%%%%%%%%%%%%%%%%%%%%%%%%%%%%%%%%%%%%%%%%%%%%%
%
\subsection{Gaussian-Poincar\'e type inequalities}\label{POINCAREGAUSSIANINEQUALITY}
In this section we study some applications of the spectral analysis carried out above. In particular, we will show three Gaussian-Poincar\'e type inequalities that basically follows from the spectral decomposition of the space $L^2_{\mu}$ in eigenspaces generated by problems \eqref{eq:ORNSTUHLENFOREXTENSIONNEUMANN}, \eqref{eq:ORNSTUHLENFOREXTENSION} and \eqref{eq:ORNSTUHLENFOREXTENSION1}. There is a wide literature about Gaussian-Poincar\'e inequalities: such kind of inequalities were studied in the works of Muckenhoupt \cite{Muckenhoupt1972:art} and Barthe and Roberto \cite{BatheRoberto2003:art} and the references therein. Those we prove in Lemma \ref{GAUSSIANPOICAREINEQUALITYONED} and Theorem \ref{GAUSSIANPOICAREINEQUALITYSEVERALD} are written with \emph{optimal constants} which, to the best of our knowledge, is a new result. This is essentially due to the availability of a orthogonal basis of eigenfunctions constructed in the previous section. These kind of inequalities have a self interest and, moreover, they will play an important role in the proof of some Liouville type theorems, essential tools to get locally uniform convergence of the blow-up sequences. Miming the steps followed in the spectral analysis, we proceed first with the one dimensional case.
%
%
%
%%%%%%%%%%%%%%%%%%%%%%%%%%%%%%%%%%%%%%%%%%%%%%%%%%%%%%%%%%%%%%%%%%%%%%%%%%%%%%%%%%%%%%%%%%%%%%%%%%%%%%%%%
%
%
%
\subsection{One-dimensional Gaussian Poincar\'e inequality}
We have dealt with the eigenvalue problem \eqref{eq:ONEDIMORNSUHLEN1}:
\[
-|y|^{-a}\left(|y|^a \psi'\right)' + \frac{y}{2} \, \psi' = \sigma \psi, \quad y \not= 0,
\]
and we have proved Lemma \ref{PROPEIGENVALUESONENDIM}, characterizing the entire spectrum. In particular, it provides a basis of eigenfunctions $\{\psi_m\}_{m\in\NN} = \{\widetilde{\psi}_m\}_{m \in \NN} \cup \{\widehat{\psi}_m\}_{m \in \NN}$ for the space $L^2(\RR,d\mu_y)$ (as in Lemma \ref{PROPEIGENVALUESONENDIM}), where, as always,
\[
d\mu_y = \frac{1}{2^{1+a} \Gamma(\frac{1+a}{2})} |y|^a e^{-\frac{y^2}{4}} dy.
\]
Before moving on, we highlight a very elementary but crucial fact. Since $\{\psi_m\}_{m\in\NN} \subset H^1(\RR,d\mu_y)$ and
\[
\int_{\RR} \psi_m' \psi' \,d\mu_y = \sigma_m \int_{\RR} \psi_m \psi\,d\mu_y, \quad \text{ for all } m \in \NN,
\]
and for any $\psi \in H^1(\RR,d\mu_y)$, where $\{\sigma_m\}_{m\in\NN} = \{\widetilde{\sigma}_m\}_{m\in\NN}\cup\{\widehat{\sigma}_m\}_{m\in\NN}$ (as in Lemma \ref{PROPEIGENVALUESONENDIM}), we easily deduce
\begin{equation}\label{eq:POINCAREEQUALITYEIGNFUNCTIONS}
\int_{\RR} \psi_m' \psi_n' \,d\mu_y = \sigma_m \int_{\RR} \psi_m \psi_n\,d\mu_y = \sigma_m \delta_{m,n}, \quad \text{ for all } m,n \in \NN,
\end{equation}
up to normalization, thanks to the $L^2(\RR,d\mu_y)$-orthogonality of the eigenfunctions ($\delta_{n,m}$ denotes the Kronecker delta). This means that not only the eigenfunctions $\psi_m = \psi_m(y)$ are orthogonal in $L^2(\RR,d\mu_y)$, but also in $H^1(\RR,d\mu_y)$. This will be important later. We also recall that the space $L^2(\RR,d\mu_y)$ is the closure of the space $\mathcal{C}_c^{\infty}(\RR)$ with respect to the norm
\[
\|\psi\|_{L^2_{\mu_y}}^2 := \int_{\RR} \psi^2 (y) \,d\mu_y,
\]
while $H^1(\RR,d\mu_y)$ is the closure of the same space but with respect to the norm
\[
\|\psi\|_{H^1_{\mu_y}}^2 := \int_{\RR} \psi^2 (y) \,d\mu_y + \int_{\RR} (\psi')^2 (y) \,d\mu_y.
\]
\begin{lem}\label{GAUSSIANPOICAREINEQUALITYONED}
Fix $a \in (-1,1)$. The following three statements hold:

\noindent (i) For any $\psi \in H^1(\RR,d\mu_y)$, it holds
\[
\int_{\RR} \psi^2 \,d\mu_y - \left(\int_{\RR} \psi \,d\mu_y \right)^2 \leq \frac{2}{1-a} \int_{\RR}(\psi')^2 \,d\mu_y.
\]
Furthermore, the equality is attained if and only if $\psi(y) = A$ or $\psi(y) = Ay|y|^{-a}$, $A \in \RR$.

\noindent (ii) For any $\psi \in H_0^1(\RR_+,d\mu_y)$, it holds
\[
\int_{\RR_+} \psi^2 \,d\mu_y \leq \frac{2}{1-a} \int_{\RR_+}(\psi')^2 \,d\mu_y.
\]
Furthermore, the equality is attained if and only if $\psi(y) = Ay^{1-a}$,  $A \in \RR$.

\noindent (iii) For any $\psi \in H^1(\RR_+,d\mu_y)$, it holds
\[
\int_{\RR_+} \psi^2 \,d\mu_y - \left(\int_{\RR_+} \psi \,d\mu_y \right)^2 \leq \int_{\RR_+}(\psi')^2 \,d\mu_y.
\]
Furthermore, the equality is attained if and only if $\psi(y) = A$ or $\psi(y) = A(1 - \frac{1-a}{2} - \frac{y^2}{4})$, $A \in \RR$.
\end{lem}

\emph{Proof.} Let us prove part (i). Take $\psi \in H^1(\RR,d\mu_y)$ with $\overline{\psi} := \int_{\RR} \psi \,d\mu_y$, and consider $\Psi_M := \sum_{m=1}^M c_m \psi_m$ such that
\[
\Psi_M \to \psi - \overline{\psi} \quad \text{ in } L^2(\RR,d\mu_y),
\]
where, of course, $c_m := \int_{\RR}\psi\psi_m \,d\mu_y$. Note that the sum $\Psi_M$ starts from the first nonconstant eigenfunction, since $\psi - \overline{\psi}$ is orthogonal to the eigenspace generated by the constants. Now, using the definition of weak eigenfunction and the orthogonality condition in \eqref{eq:POINCAREEQUALITYEIGNFUNCTIONS}, it is not difficult to see that
\[
\int_{\RR} (\Psi_M')^2 \,d\mu_y = \int_{\RR} \Psi_M'\psi' \,d\mu_y  \quad \text{ for all } M \in \NN_0,
\]
i.e. $\Psi_M'$ and $\Psi_M' - \psi'$ are orthogonal in $L^2(\RR,d\mu_y)$. Consequently, we have
\[
0 \leq \int_{\RR} (\psi' - \Psi_M')^2 d\mu_y = \int_{\RR}(\psi')^2 d\mu_y - 2 \int_{\RR} \Psi_M'\psi' d\mu_y + \int_{\RR}(\Psi_M')^2 d\mu_y = \int_{\RR}(\psi')^2 d\mu_y - \int_{\RR}(\Psi_M')^2 d\mu_y,
\]
and so $\int_{\RR}(\Psi_M')^2 d\mu_y \leq \int_{\RR}(\psi')^2 d\mu_y < +\infty$, uniformly in $M \in \NN_0$ (note that here we have used the fact that $\psi \in H^1(\RR,d\mu_y)$). Finally, from the above bound and employing \eqref{eq:POINCAREEQUALITYEIGNFUNCTIONS} again, we have
\[
\begin{aligned}
\int_{\RR}(\psi')^2 d\mu_y \geq \int_{\RR} (\Psi_M')^2 \,d\mu_y &= \sum_{m=1}^M c_m^2 \int_{\RR} (\psi_m')^2 \,d\mu_y = \sum_{m=1}^M \sigma_m c_m^2 \int_{\RR} \psi_m^2 \,d\mu_y \\
&\geq \underline{\sigma}\sum_{m=1}^M c_m^2 \int_{\RR} \psi_m^2 \,d\mu_y = \underline{\sigma} \int_{\RR} \Psi_M^2 \,d\mu_y,
\end{aligned}
\]
where we have set $\underline{\sigma} := \min_{m\in\NN_0} \{\sigma_m \} = \widehat{\sigma}_0 = (1-a)/2$. Consequently, passing to the limit as $M \to +\infty$, we deduce
\[
\int_{\RR} (\psi')^2 \,d\mu_y \geq \int_{\RR} (\psi - \overline{\psi})^2 \,d\mu_y = \int_{\RR} \psi^2 \,d\mu_y - \left(\int_{\RR} \psi \,d\mu_y \right)^2,
\]
obtaining the desired inequality. The last part of the statement follows from the fact that $\psi(y) = A y|y|^{-a}$ is the first non constant eigenfunction to problem \eqref{eq:ONEDIMORNSUHLEN1} (with eigenvalue $\widehat{\sigma}_0 = (1-a)/2$).

The proof of (ii) and (iii) is very similar to the one of part (i) and we omit it (note that we just have to use the bases of eigenfunctions $\{\widehat{\psi}_m\}_{m\in\NN}$ for (ii) and $\{\widetilde{\psi}_m\}_{m\in\NN}$ for (iii), instead of $\{\psi_m\}_{m\in\NN}$). We just mention that in (ii) we work with the space $H_0^1(\RR,d\mu_y)$ which does not contain constant functions different from the trivial one (this is the reason because we do not need to employ the average of $\psi$ in Gaussian Poincar\'e inequality), while in (iii) it is easily seen that the first nonconstant eigenfunction $\psi(y) = A(1 - \frac{1-a}{2} - \frac{y^2}{4})$, $A \in \RR$ has zero mean (with respect to the measure $d\mu_y$). Notice that, since in this last case we are now working with functions defined for $y>0$, the normalization constant in front of the probability measure $d\mu_y$ is different from the one of Lemma \ref{GAUSSIANPOICAREINEQUALITYONED}. $\Box$

\bigskip

For the general case, we proceed as before, but we consider problem \eqref{eq:ORNSTUHLENFOREXTENSION1}:
\[
-\frac{1}{|y|^a \mathcal{G}_a} \nabla\cdot(|y|^a \mathcal{G}_a \nabla V) = \kappa V \quad \text{in } \RR^{N+1}.
\]
From Theorem \ref{SPECTRALTHEOREMEXTENDEDORNUHL1} we know that it has the set of eigenvalues $\{\kappa_{n,m}\}_{n,m \in \NN} = \{\widehat{\kappa}_{n,m}\}_{n,m \in \NN}\cup\{\widetilde{\kappa}_{n,m}\}_{n,m \in \NN}$ where
\[
\widetilde{\kappa}_{n,m} := \frac{n}{2} + m, \qquad \widehat{\kappa}_{n,m} := \frac{n}{2} + m + \frac{1-a}{2} \quad \text{ for all } n,m \in \NN,
\]
while the set of eigenfunctions is $\{V_{\alpha,m}\}_{(\alpha,m)} = \{\widetilde{V}_{\alpha,m}\}_{(\alpha,m)}\cup\{\widehat{V}_{\alpha,m}\}_{(\alpha,m)}$, where the eigenfunctions can be of two types
\[
\begin{aligned}
\widetilde{V}_{\alpha,m}(x,y) = H_{\alpha}(x)\widetilde{\psi}_m(y), \quad \widehat{V}_{\alpha,m}(x,y) = H_{\alpha}(x) \widehat{\psi}_m(y), \quad  \text{ for all } (\alpha,m) \in \ZZ_{\geq 0}^N\times\NN,
\end{aligned}
\]
where as before $\widetilde{\psi}_m(y) = L_{(\frac{a-1}{2}),m}(y^2/4)$ and $\widehat{\psi}_m(y) =  y|y|^{-a}L_{(\frac{1-a}{2}),m}(y^2/4)$ and $H_{\alpha}(\cdot)$ is a N-dimensional Hermite polynomial of order $|\alpha|$. We recall that similarly to the $1$-dimensional case the set $\{V_{\alpha,m}\}_{(\alpha,m)}$ is an orthogonal basis of $L^2(\RR^{N+1},d\mu)$, where
\[
d\mu(X) = \frac{1}{2^{1+a} \Gamma(\frac{1+a}{2})(4\pi)^{N/2}} |y|^a e^{-\frac{|X|^2}{4}} dxdy.
\]
Since the sets of eigenvalues and eigenfunctions are countable, we can drop one index and denote by $\kappa_j$ the $j^{th}$ eigenvalue with corresponding eigenfunctions $V_j$, $j \in \NN$. In this setting, we have that the first nonzero eigenvalue depends on the parameter $a \in (-1,1)$:
\[
\nu_{\ast} := \min_{j\in\NN} \{\kappa_j \not= 0\} = \frac{1}{2} \min\{1,1-a\}.
\]
This is the unique remarkable difference with respect to the $1$-dimensional case. We have proved the following theorem.
\begin{thm}\label{GAUSSIANPOICAREINEQUALITYSEVERALD}
The following three statements hold:

(i) For any $V \in H^1(\RR^{N+1},d\mu)$, it holds
\begin{equation}\label{eq:GaussianPoincWholeSpace}
\int_{\RR^{N+1}} V^2 \,d\mu - \left(\int_{\RR^{N+1}} V \,d\mu \right)^2 \leq  P_a \int_{\RR^{N+1}} |\nabla V|^2 \,d\mu,
\end{equation}
where $P_a := 1/\nu_{\ast} = 2/\min\{1,1-a\}$. Furthermore, the equality is attained if and only if $V(x,y) = A$ or, depending on $-1 < a < 1$:
\[
V(x,y) =
\begin{cases}
Ax_j         \quad &\text{ if }  a < 0  \quad \text{ for some } j \in \{1,\ldots,N\}\\
Ax_j         \quad &\text{ if }  a = 0  \quad \text{ for some } j \in \{1,\ldots,N+1\}\\
Ay|y|^{-a}   \quad &\text{ if }  a > 0,
\end{cases}
\]
where $A \in \RR$ and we have used the convention $x_{N+1} = y$.

\noindent (ii) For any $V \in H_0^1(\RR^{N+1}_+,d\mu)$, it holds
\begin{equation}\label{eq:GaussianPoincWholeSpace0}
\int_{\RR^{N+1}_+} V^2 \,d\mu \leq \frac{2}{1-a} \int_{\RR^{N+1}_+} |\nabla V|^2 \,d\mu.
\end{equation}
Furthermore, the equality is attained if and only if $V(x,y) = Ay^{1-a}$,  $A \in \RR$.

\noindent (iii) For any $V \in H^1(\RR^{N+1}_+,d\mu)$, it holds
\begin{equation}\label{eq:GaussianPoincWholeSpace+}
\int_{\RR^{N+1}_+} V^2 \,d\mu - \left(\int_{\RR^{N+1}_+} V \,d\mu \right)^2 \leq 2 \int_{\RR^{N+1}_+}|\nabla V|^2 \,d\mu.
\end{equation}
Furthermore, the equality is attained if and only if $V(x,y) = A$ or $V(x,y) = Ax_j$, for some $j \in \{1,\ldots,N\}$ and $A \in \RR$.
\end{thm}
\begin{rem} Note that, since the Gaussian-Poincar\'e constants do not depend on the spacial dimension $N$, they can be extended to infinite dimensional spaces (see  Beckner \cite{Beckner1989:art}).
\end{rem}
\emph{Proof.} We begin by proving part (i), following the ideas of the case $N=1$. Take $V \in H^1(\RR^{N+1},d\mu)$ with $\overline{V} := \int_{\RR^{N+1}} V \,d\mu$, and we approximate it with the sequence $\Psi_J := \sum_{j=1}^J c_j V_j$ such that
\[
\Psi_J \to V - \overline{V} \quad \text{ in } L^2(\RR^{N+1},d\mu), \qquad c_j := \int_{\RR^{N+1}}V_jV \,d\mu.
\]
Proceeding as before, we easily find that
\[
\int_{\RR^{N+1}} |\nabla \Psi_J|^2 \,d\mu = \int_{\RR^{N+1}} \nabla \Psi_J \nabla V \,d\mu  \quad \text{ for all } J \in \NN_0,
\]
and so $\int_{\RR^{N+1}} |\nabla \Psi_J|^2 \,d\mu \leq \int_{\RR^{N+1}} |\nabla V|^2 \,d\mu$. Hence,
\[
\begin{aligned}
\int_{\RR^{N+1}}|\nabla V|^2 d\mu \geq \int_{\RR^{N+1}} |\nabla\Psi_J|^2 \,d\mu &= \sum_{j=1}^J c_j^2 \int_{\RR^{N+1}} |\nabla V_j|^2 \,d\mu = \sum_{j=1}^J \sigma_j c_j^2 \int_{\RR^{N+1}} V_j^2 \,d\mu \\
&\geq \nu_{\ast}\sum_{j=1}^J c_j^2 \int_{\RR^{N+1}} V_j^2 \,d\mu = \nu_{\ast} \int_{\RR^{N+1}} \Psi_J^2 \,d\mu,
\end{aligned}
\]
where now the minimum of the eigenvalues $\nu_{\ast} := \min_{j} \{\kappa_j \not= 0 \} = \min\{1,1-a\}/2$. Passing to the limit as $J \to \infty$ we obtain the inequality of statement (i).

\noindent To prove the second part of the statement, let us firstly fix $a \in (-1,0)$. In this case we have $\nu_{\ast} = 1/2$ and the corresponding eigenfunctions are $V(x,y) = A x_j$ for $j \in \{1,\ldots,N\}$. We thus conclude by the definition of weak eigenfunction and the fact that $V(x,y) = A x_j$ has mean zero for any $j \in \{1,\ldots,N\}$ (with respect to the measure $d\mu$). Similar, if $a \in (0,1)$ it turns out that $\nu_{\ast} = (1-a)/2$ and the corresponding eigenfunctions are $V(x,y) = A y|y|^{-a}$ (also in this case they have mean zero). Finally, if $a = 0$ we have $\nu_{\ast} = 1/2$ and the eigenfunctions are $V(x,y) = A x_j$ for $j \in \{1,\ldots,N\}$ and $V(x,y) = A y$ (note that we get back the classical statement, see  for instance with \cite{Beckner1989:art}). This conclude the proof of part (i).

To prove part (ii) and (iii) we employ the bases of eigenfunctions $\{\widehat{V}_j\}_{j} = \{\widehat{V}_{(\alpha,m)}\}_{(\alpha,m)}$ for (ii) and $\{\widetilde{V}_j\}_{j} = \{\widetilde{V}_{(\alpha,m)}\}_{(\alpha,m)}$ for (iii), instead of $\{V_j\}_j = \{V_{(\alpha,m)}\}_{(\alpha,m)}$. Note that compared with the proof of Lemma \ref{GAUSSIANPOICAREINEQUALITYONED}, the first eigenvalue to the Neumann problem \eqref{eq:ORNSTUHLENFOREXTENSION} is $\nu_{\ast} = 1/2$ with corresponding eigenfunctions $V(x,y) = A x_j$ for $j \in \{1,\ldots,N\}$. $\Box$

\section{Blow-up Analysis: Part I}\label{Section:BlowUp1}
Let $W$ be a weak solution to \eqref{eq:TruncatedExtensionEquationLocalBackward} and $p_0 = (x_0,0,t_0) \in B_{1/2}\times\{0\}\times(0,1)$. We consider the blow-up family
\begin{equation}\label{eq:BLOWUPPLUSTRANSSEQUENCE}
W_{p_0,r}(X,t) := \frac{W_{p_0}^r(X,t)}{\sqrt{H(W,p_0,r)}} := \frac{W_{p_0}(rX,r^2t)}{\sqrt{H(W,p_0,r)}}, \quad r > 0,
\end{equation}
satisfying
\[
H(W_{p_0,r},O,1) = H(W_r,p_0,1) = 1,
\]
where $W_{p_0}(x,y,t) = W(x+x_0,y,t+t_0)$ and $W_{O,r}(X,t) = W_r(X,t)$. For any fixed $r > 0$, $W_r$ is a weak solution to
\begin{equation}\label{eq:TruncatedExtensionEquationBlowUp}
\begin{cases}
\partial_t W_r + y^{-a} \nabla \cdot(y^a \nabla W_r) = F_r \quad &\text{ in } \mathbb{R}_+^{N+1} \times (0,1/r^2) \\
-\partial_y^a W_r = r^{1-a} q^r(x,t)w_r \quad &\text{ in } \mathbb{R}^N \times \{0\} \times (0,1/r^2),
\end{cases}
\end{equation}
where $q^r(x,t) := q(r x,r^2 t)$, while
\[
w_r(x,t) := \frac{w(r x,r^2t)}{\sqrt{H(W,r)}}, \qquad F_r(X,t) := \frac{r^2 F(rX,r^2t) }{\sqrt{H(W,r)}}.
\]
Now, as mentioned above, we apply Lemma \ref{Lemma:DerivativeNI} and Theorem \ref{Theorem:GeneralizedFrequency} to obtain the almost-monotonicity of the Almgren-Poon quotient computed along the blow-up family \eqref{eq:BLOWUPPLUSTRANSSEQUENCE}, \emph{independently} on the blow-up parameter. We stress that this fact is not trivial and it strongly depends on the definition of blow-up sequence and
the scaling of problem \eqref{eq:TruncatedExtensionEquationLocalBackward}.

\begin{cor}\label{Corollary:GeneralizedFrequencyBlowUp}
Let $W$ be a weak solution to \eqref{eq:TruncatedExtensionEquationLocalBackward}, with $q$ satisfying \eqref{eq:ASSUMPTIONSPOTENTIAL} and let $W_\varrho$ be defined as in \eqref{eq:BLOWUPPLUSTRANSSEQUENCE} for all $\varrho \in (0,1)$. Assume that $W$ vanishes of finite order at $O$. If $r_0,C > 0$ are chosen as in Theorem \ref{Theorem:GeneralizedFrequency} (not depending on $\varrho \in (0,1)$), then the function
\[
r \to  \Phi_a(W_\varrho,O,r) := e^{Cr^{1-a}} N_I(W_\varrho,r) + e^{Cr^{1-a}} - 1
\]
is monotone nondecreasing in $r \in (0,r_0)$. Furthermore,
\[
\Phi_a(W_\varrho,O,0^+) = \lim_{r \to 0^+} N_I(W_\varrho,r) = \lim_{r \to 0^+} N_D(W_\varrho,r) = \lim_{r \to 0^+} N_0(W_\varrho,r) = \kappa,
\]
for all $\varrho \in (0,1)$, where $\kappa \geq 0$ is defined as in Corollary \ref{Lemma:LimitAlmgrenPoon}.
\end{cor}
\emph{Proof.} The proof is based on the following scaling relation:
\begin{equation}\label{eq:ScalingPropAlmgrenPoon}
N_I(W_\varrho,r) = N_I(W,r\varrho) = N_I(W_r,\varrho),
\end{equation}
for all $r,\varrho > 0$. Indeed, setting $\rho := r\varrho$ and using \eqref{eq:ScalingPropAlmgrenPoon}, we have
\[
\frac{d}{dr} \left[ 1 + N_I(W_\varrho,r) \right] = \varrho \frac{d}{d\rho} \left[ 1 + N_I(W,\rho), \right]
\]
and so, applying Lemma \ref{Lemma:DerivativeNI}, we deduce the existence of $r_0 \in (0,1)$ and $C_1 > 0$ depending only on the solution to \eqref{eq:ExtensionEquationLocalBackward}, such that
\[
\frac{d}{dr} \left[ 1 + N_I(W_\varrho,r) \right] \geq - C_1 \varrho \left[ 1 + N_I(W,\rho) \right] \rho^{-a} = - C_1 \varrho^{1-a} \left[ 1 + N_I(W_\varrho,r) \right]r^{-a}.
\]
Notice that the assumption $\varrho \in (0,1)$ is needed to guarantee $r\varrho \in (0,r_0)$. Consequently, using the above inequality, the fact that $\varrho \in (0,1)$ and $a \in (-1,1)$, it follows
\[
\begin{aligned}
\frac{d}{dr} \Phi_a(W_\varrho,O,r) & = r^{-a}e^{Cr^{1-a}} \left\{ C(1-a) \left[N_I(W_\varrho,r) + 1 \right] + r^a \frac{d}{dr} \left[N_I(W_\varrho,r) + 1 \right] \right\} \\
& \geq r^{-a}e^{Cr^{1-a}} \left[N_I(W_\varrho,r) + 1 \right] \left\{ C(1-a) - C_1\varrho^{1-a} \right\} \\
& \geq r^{-a}e^{Cr^{1-a}} \left[N_I(W_\varrho,r) + 1 \right] \left\{ C(1-a) - C_1 \right\} \geq 0,
\end{aligned}
\]
if $C \geq C_1/(1-a)$, i.e. $C > 0$ is chosen as in Theorem \ref{Theorem:GeneralizedFrequency}. $\Box$

\bigskip

We are ready to state and prove the first main blow-up classification result.
\begin{thm}\label{THEOREMBLOWUP1New}
Let $a \in (-1,1)$ and let $W$ be a weak solution to \eqref{eq:TruncatedExtensionEquationLocalBackward}, with $q$ satisfying \eqref{eq:ASSUMPTIONSPOTENTIAL}. Assume that $W$ vanishes of finite order at $p_0 = (x_0,0,t_0) \in B_{1/2}\times\{0\}\times(0,1)$. Then there exist $n_0,m_0 \in \NN$ such that the following assertions hold:

(i) The Almgren-Poon quotient satisfies
\[
\lim_{r \to 0^+} N_D(W_{p_0},r) = \kappa := \widetilde{\kappa}_{n_0,m_0},
\]
where the admissible values for $\widetilde{\kappa}_{n,m}$ corresponds to the eigenvalues of problem \eqref{eq:ORNSTUHLENFOREXTENSIONNEUMANN}:
\[
\widetilde{\kappa}_{n,m} := \frac{n}{2} + m,
\]
for any $m,n \in \NN$.

(ii) The vanishing order of $W$ at $p_0$ is $\si_\ast = 2\kk$ and, furthermore,
\begin{equation}\label{eq:LIMITOFTHEQUOTIENTFORSMALLTSTATEMENT1New}
\lim_{r \to 0^+}  \mathcal{T}_{2\kappa}(W_{p_0},r) := \lim_{r \to 0^+}  \frac{H(W_{p_0},r)}{r^{4\kappa}} = L_0,
\end{equation}
for some constant $L_0 > 0$ (depending on $p_0$ and $W$).

(iii) There exists a parabolically homogeneous polynomial of degree $2\widetilde{\kappa}_{n_0,m_0}$ given by the expression
\[
\Theta_{p_0}(X,t) := t^{\widetilde{\kappa}_{n_0,m_0}} \sum_{(\alpha,m) \in \widetilde{J}_0} v_{\alpha,m} \overline{V}_{\alpha,m}\left(\frac{X}{\sqrt{t}}\right),
\]
such that, for all $0 < t_{\ast} < T_{\ast}$, it holds
\[
W_{p_0,r} \to \Theta_{p_0} \quad \text{ in } L^2(0,T_{\ast};H^1(\RR_+^{N+1},d\mu_t))\cap\mathcal{C}^0(t_\ast,T_{\ast};L^2(\RR_+^{N+1},d\mu_t)),
\]
as $r \to 0^+$, i.e.:
\begin{equation}\label{eq:ENERGYCONVERGENCESBLOWUP1New}
\lim_{r \to 0^+} \int_0^{T_{\ast}} \left\| W_{p_0,r} - \Theta_{p_0} \right\|_{H_{\mu_t}^1}^2 dt  = \lim_{r \to 0^+}\sup_{t \in [t_\ast,T_{\ast}]} \left\| W_{p_0,r} - \Theta_{p_0} \right\|_{L_{\mu_t}^2}^2 = 0.
\end{equation}
The $v_{\alpha,m}$'s are suitable constants, the sum is done over the set of indexes
\[
\widetilde{J}_0 := \left\{\right (\alpha,m) \in \ZZ_{\geq 0}^N\times\NN : |\alpha| = n \in \NN \text{ and } \widetilde{\kappa}_{n,m} = \widetilde{\kappa}_{n_0,m_0} \},
\]
and the integration probability measure is defined in \eqref{eq:EXTENDEDGAUSSIANMEASURE}. Moreover,
\[
\overline{V}_{\alpha,m}(X) := \frac{V_{\alpha,m}(X)}{\|V_{\alpha,m}\|_{L_{\mu}^2}}, \qquad \alpha \in \ZZ_{\geq 0}^N, \; m \in \NN,
\]
are the normalized eigenfunctions $V_{\alpha,m}$ to problem \eqref{eq:ORNSTUHLENFOREXTENSIONNEUMANN} corresponding to the eigenvalue $\widetilde{\kappa}_{n,m}$, defined in the statement of Theorem \ref{SPECTRALTHEOREMEXTENDEDORNUHL1} (part (i)).
\end{thm}
\emph{Proof.} Let us fix $T_\ast \geq 1$, $T = \sqrt{T_{\ast}}$ and $p_0 = O$. We consider a weak solution $W$ to \eqref{eq:TruncatedExtensionEquationLocalBackward} and the blow-up family $W_r$ defined in \eqref{eq:BLOWUPPLUSTRANSSEQUENCE}, which satisfy problem \eqref{eq:TruncatedExtensionEquationBlowUp} for any fixed $r \in (0,1)$, in the sense of \eqref{eq:StrongSolutionW1}. Combining Corollary \ref{Corollary:NonDegeneracy}, Theorem \ref{Theorem:GeneralizedFrequency} and Corollary \ref{Corollary:GeneralizedFrequencyBlowUp}, we easily deduce the existence of $r_0\in(0,1)$ and $C > 0$ (depending on $W$ and $T$) such that
\begin{equation}\label{eq:BoundOnNDBlowUp}
N_D(W,rT) \leq C, \qquad rT < 1,
\end{equation}
for all $r \in (0,r_0)$.

Now, we re-scale and consider the function
\[
\widetilde{W}_r(X,t) := W_r(\sqrt{t}X,t),
\]
which is a weak solution to
\begin{equation}\label{eq:RescaledBlowUp}
\begin{cases}
t\partial_t \widetilde{W}_r + \mathcal{O}_a\widetilde{W}_r = t\widetilde{F}_r \quad &\text{ in } \mathbb{R}_+^{N+1} \times (0,1/r^2) \\
-\partial_y^a \widetilde{W}_r = (tr^2)^{\frac{1-a}{2}} \widetilde{q}^r(x,t)\widetilde{w}_r \quad &\text{ in } \mathbb{R}^N \times \{0\} \times (0,1/r^2),
\end{cases}
\end{equation}
for each $r \in (0,r_0)$, in the sense that
\[
\begin{aligned}
t \int_{\mathbb{R}^{N+1}_+} \partial_t \widetilde{W}_r \eta \,d\mu = \int_{\mathbb{R}^{N+1}_+}  \nabla \widetilde{W}_r \cdot \nabla \eta \,d\mu + t \int_{\mathbb{R}^{N+1}_+} \widetilde{F}_r \eta \,d\mu - (tr^2)^{\frac{1-a}{2}} \int_{\mathbb{R}^N\times\{0\}} \widetilde{q}^r \widetilde{w}_r \eta \,\mathcal{G}_a(x,0,1)dx,
\end{aligned}
\]
for all $\eta \in L^2(0,1/r^2;H_{\mu}^1)$ and for all $t \in (0,1/r^2)$,  see  \eqref{eq:RescaledTruncatedExtensionEquationLocalBackward}.

\smallskip

\emph{Step 1: Uniform estimates in Gaussian spaces.} Changing variables twice ($\sqrt{t}X \to X$, $(rX,r^2t) \to (X,t)$), we obtain
\[
\begin{aligned}
\int_0^{T^2} \int_{\mathbb{R}_+^{N+1}} \widetilde{W}_r^2(X,t) d\mu(X) dt &= \frac{1}{H(W,r)} \int_0^{T^2} \int_{\mathbb{R}_+^{N+1}} W^2(rX,r^2t) d\mu_t(X) dt \\
& = T^2 \frac{H(W,rT)}{H(W,r)} \leq T^{2+C} \leq C_{W,T},
\end{aligned}
\]
for all $r \in (0,r_0)$, where $C > 0$ is chosen as in Corollary \ref{Corollary:NonDegeneracy} and $C_{W,T} > 0$ depends on $W$ and $T$. Similar, we find
\[
\begin{aligned}
\int_0^{T^2} \int_{\mathbb{R}_+^{N+1}} |\nabla \widetilde{W}_r(X,t)|^2 d\mu(X) dt &= \frac{r^2}{H(W,r)} \int_0^{T^2} t \int_{\mathbb{R}_+^{N+1}} |\nabla W(rX,r^2t)|^2 d\mu_t(X) dt \\
& = T^2 \frac{D(W,rT)}{H(W,r)} = T^2 N_D(W,rT) \frac{H(W,rT)}{H(W,r)} \leq C_{W,T},
\end{aligned}
\]
for some new $C_{W,T} > 0$, thanks to Corollary \ref{Corollary:NonDegeneracy} and \eqref{eq:BoundOnNDBlowUp}. It thus follows that the family $\{\widetilde{W}_r\}_{r \in (0,r_0)}$ is uniformly bounded (w.r.t. $r \in (0,r_0)$) in $L^2(0,T_{\ast};H^1_{\mu})$.

Now, let us consider the integral formulation of \eqref{eq:RescaledBlowUp}. Thanks to the uniform estimates found above, we easily see that
\begin{equation}\label{eq:FirstTermBlowUp}
\int_0^{T^2} \int_{\mathbb{R}^{N+1}_+}  \nabla \widetilde{W}_r \cdot \nabla \eta \,d\mu dt \leq \|\widetilde{W}_r\|_{L^2(0,T^2;H^1_\mu)} \|\eta\|_{L^2(0,T^2;H^1_\mu)} \leq C_{W,T} \|\eta\|_{L^2(0,T^2;H^1_\mu)},
\end{equation}
for some new $C_{W,T} > 0$, uniformly w.r.t. $r \in (0,r_0)$. Passing to the second term in the r.h.s., we have
\[
\left| \int_0^{T^2} t \int_{\mathbb{R}^{N+1}_+} \widetilde{F}_r \eta \,d\mu dt \right| \leq \|t\widetilde{F}_r\|_{L^2(0,T^2;L^2_\mu)} \|\eta\|_{L^2(0,T^2;L^2_\mu)},
\]
while, using \eqref{eq:ExpBoundIntF} and the fact that $W$ vanishes of finite order at $O$, it follows
\begin{equation}\label{eq:SecondTermBlowUp}
\begin{aligned}
\|t\widetilde{F}_r\|_{L^2(0,T^2;L^2_\mu)}^2 &= \int_0^{T^2} t^2 \int_{\mathbb{R}_+^{N+1}} |\widetilde{F}_r(X,t)|^2 d\mu(X) dt = \frac{1}{r^6 H(W,r)} \int_0^{(rT)^2} t^2 \int_{\mathbb{R}_+^{N+1}} F^2 d\mu_t(X) dt \\
& \leq \frac{C_{W,T}}{H(W,r)} r^{\sigma_0} e^{-\frac{1}{16T^2r^2}} \leq C_{W,T} r^{\sigma_0} e^{-\frac{1}{16T^2r^2}},
\end{aligned}
\end{equation}
for some new $C_{W,T} > 0$ and some $\sigma_0 \in \mathbb{R}$ depending on $W$, for all $r \in (0,r_0)$. Finally,
\begin{equation}\label{eq:ThirdTermBlowUp}
r^{1-a} \left| \int_0^{T^2} t^{\frac{1-a}{2}} \int_{\mathbb{R}^N\times\{0\}} \widetilde{q}^r \widetilde{w}_r \eta \,\mathcal{G}_a(x,0,1)dx \right| \leq C r^{1-a} \|t^{\frac{1-a}{2}} \widetilde{w}_r \|_{L^2(0,T^2;L^2_{\mu_x})} \|\eta\|_{L^2(0,T^2;L^2_{\mu_x})}
\end{equation}
for some $C > 0$ depending on $K > 0$ in \eqref{eq:ASSUMPTIONSPOTENTIAL}, $N$ and $a$. Now, applying Lemma \ref{Lemma:TraceIneqL2Norm} to $w_r$, we obtain
\[
\begin{aligned}
\|t^{\frac{1-a}{2}} \widetilde{w}_r\|_{L^2(0,T^2;L^2_{\mu_x})}^2 & \leq T^{2(1-a)} \int_0^{T^2} \int_{\mathbb{R}_+^{N+1}} \widetilde{w}_r^2 \mathcal{G}_a(x,0,1)dx dt \\
& =  \int_0^{T^2} t^{\frac{1+a}{2}}\int_{\mathbb{R}_+^{N+1}} w_r^2 \mathcal{G}_a(x,0,t)dx dt \leq C_0 T^{2(1-a)} \int_0^{T^2} [d(W_r,t) + h(W_r,t)] dt\\
& = C_0 T^{2(2-a)} [D(W_r,T) + H(W_r,T)] = C_0 T^{2(2-a)} \frac{D(W,rT) + H(W,rT)}{H(W,r)}  \\
& \leq C_0 T^{2(2-a)} N_D(W,rT) \frac{H(W,rT)}{H(W,r)} \leq C_{W,T},
\end{aligned}
\]
for some new $C_{W,T} > 0$, thanks to Corollary \ref{Corollary:NonDegeneracy} and \eqref{eq:BoundOnNDBlowUp}. Now, coming back to the integral formulation \eqref{eq:RescaledBlowUp}, it follows that for each $t \in (0,T_{\ast})$, the function $t \partial_t\widetilde{W}_r$ can be seen as a linear and continuous functional defined on $H^1_{\mu}$ or, equivalently, an element of $H^{-1}_\mu := (H^1_{\mu})'$. Furthermore, integrating between $0$ and $T_{\ast}$, and using \eqref{eq:FirstTermBlowUp}, \eqref{eq:SecondTermBlowUp} and \eqref{eq:ThirdTermBlowUp}, we obtain
\[
\left| \int_0^{T_\ast} t \int_{\mathbb{R}^{N+1}_+} \partial_t \widetilde{W}_r \eta \,d\mu dt \right| \leq C_{W,T} \|\eta\|_{L^2(0,T_\ast;H^1_\mu)},
\]
which, thanks to the Riesz representation theorem and the fact that $H^1_\mu \subset L^2_\mu$, implies that the family $\{t\partial_t\widetilde{W}_r\}_{r \in (0,r_0)}$ is uniformly bounded in $L^2(0,T_\ast;L^2_\mu)$. Consequently, for any fixed $t_\ast \in (0,T_\ast)$, $\{\widetilde{W}_r\}_{r \in (0,r_0)}$ is uniformly bounded in $L^2(0,T_{\ast};H^1_{\mu})$, while $\{\partial_t\widetilde{W}_r\}_{r \in (0,r_0)}$ is uniformly bounded in $L^2(t_\ast,T_\ast;L^2_\mu)$ and so, applying Lemma \ref{SIMONLEMMARECALL}  (\cite[Corollary 8]{Simon1987:art}) with $X = H_{\mu}^1$,  $B = Y = L_{\mu}^2$ and $p = r = 2$, it follows that
\[
\{\widetilde{W}_r\}_{r \in (0,r_0)} \subset \mathcal{C}^0(t_\ast,T_\ast; L^2_\mu)
\]
is relatively compact.

\smallskip

\emph{Step 2: Compactness and properties of the limit.} So, for any $r_n \to 0^+$ and any fixed $0 < t_{\ast} < T_{\ast}$, we can extract a sub-sequence $r_{n_j} \to 0^+$ (that we rename $r_j := r_{n_j}$ by convenience) such that
\begin{equation}\label{eq:C0CONVERGENCEOFWJ}
\widetilde{W}_{r_j} \to \widetilde{\Theta} \quad \text{ in } \mathcal{C}^0\left(t_{\ast},T_{\ast};L_{\mu}^2\right), \quad \text{as } j \to +\infty,
\end{equation}
where $\widetilde{\Theta} \in \bigcap_{t_{\ast} \in (0,T_{\ast})} \mathcal{C}^0 \left(t_{\ast},T_{\ast}; L_{\mu}^2 \right)$. A standard diagonal procedure allows to make the sequence $r_j$ independent of $t_{\ast}$ (we skip that not to weight down the presentation). Further, the uniform bounds found above and the reflexivity of the spaces $L^2(t_{\ast},T_{\ast};H_{\mu}^1)$ and $L^2(t_{\ast},T_{\ast};L_{\mu}^2)$ allows us to assume also
\[
\begin{aligned}
\widetilde{W}_{r_j} \rightharpoonup \widetilde{\Theta} \quad &\text{ weakly in } L^2(t_{\ast},T_{\ast};H_{\mu}^1), \\
\partial_t \widetilde{W}_{r_j} \rightharpoonup \partial_t \widetilde{\Theta} \quad &\text{ weakly in } L^2(t_{\ast},T_{\ast};L_{\mu}^2).
\end{aligned}
\]
Consequently, combining the above limits with \eqref{eq:SecondTermBlowUp} and \eqref{eq:ThirdTermBlowUp}, we can pass to the limit in \eqref{eq:RescaledBlowUp} to deduce
\[
\begin{aligned}
\int_{t_1}^{t_2} t \int_{\mathbb{R}^{N+1}_+} \partial_t \widetilde{\Theta} \eta \,d\mu dt = \int_{t_1}^{t_2} \int_{\mathbb{R}^{N+1}_+}  \nabla \widetilde{\Theta} \cdot \nabla \eta \,d\mu dt,
\end{aligned}
\]
for all $\eta \in L^2(t_1,t_2,;H_{\mu}^1)$ and for all $t_\ast \leq t_1 \leq t_2 \leq T_\ast$, that is $\widetilde{\Theta}$ is a global weak solution to
\[
\begin{cases}
t\partial_t \widetilde{\Theta} + \mathcal{O}_a\widetilde{\Theta} = 0 \quad &\text{ in } \mathbb{R}_+^{N+1} \times (0,\infty) \\
-\partial_y^a \widetilde{\Theta} = 0 \quad &\text{ in } \mathbb{R}^N \times \{0\} \times (0,\infty).
\end{cases}
\]
Now, subtracting the equations of $\widetilde{W}_{r_j}$ and $\widetilde{\Theta}$, testing with $\eta = \widetilde{W}_{r_j} - \widetilde{\Theta}$ and using \eqref{eq:SecondTermBlowUp} and \eqref{eq:ThirdTermBlowUp} again, we obtain
\[
\int_{t_1}^{t_2}t \left\langle \partial_t(\widetilde{W}_{r_j} - \widetilde{\Theta})(t), (\widetilde{W}_{r_j} - \widetilde{\Theta})(t) \right\rangle_{L_{\mu}^2} dt = \int_{t_1}^{t_2} \int_{\RR^{N+1}}|\nabla (\widetilde{W}_{\lambda
_j} - \widetilde{\Theta})(t)|^2 \,d\mu dt + o(1),
\]
for large $j$'s. Integrating by parts in time and using that $2 ( \partial_tU, U )_{L_{\mu}^2} = \partial_t(\|U\|_{L_{\mu}^2}^2)$, we easily find
\[
\begin{aligned}
\int_{t_1}^{t_2} \|(\widetilde{W}_{r_j} - \widetilde{\Theta})(t)\|_{L_{\mu}^2}^2 & + \|\nabla (\widetilde{W}_{r_j} - \widetilde{\Theta})(t)\|_{L_{\mu}^2}^2 \,dt  \\
& = \frac{1}{2}\left[\|(\widetilde{W}_{r_j} - \widetilde{\Theta})(t_1)\|_{L_{\mu}^2}^2 + \|(\widetilde{W}_{r_j} - \widetilde{\Theta})(t_2)\|_{L_{\mu}^2}^2\right] + o(1) \to 0, \quad \text{ as } j \to +\infty,
\end{aligned}
\]
thanks to \eqref{eq:C0CONVERGENCEOFWJ}. This implies
\begin{equation}\label{eq:L2STRONGCONVERGENCEOFWJ}
\widetilde{W}_{r_j} \to \widetilde{\Theta} \quad \text{ in } L^2\left(t_{\ast},T_{\ast};H_{\mu}^1\right), \quad \text{ as } j \to +\infty.
\end{equation}

\smallskip

\emph{Step 3: Improved convergence.} We prove that
\begin{equation}\label{eq:FINALLIMITH1BLOWUP}
\lim_{j \to  0^+} \int_0^{T_{\ast}} \left\| W_{r_j} - \Theta \right\|_{L_{\mu_t}^2}^2 + t\left\| \nabla W_{r_j} - \nabla \Theta \right\|_{L_{\mu_t}^2}^2 dt = 0,
\end{equation}
where
\[
\Theta(X,t) = \widetilde{\Theta}\left(\frac{X}{\sqrt{t}},t\right).
\]
Scaling back to $W_{r} = W_{r}(X,t)$, we deduce that for any fixed $t_\ast \in (0,T_\ast)$, it holds
\[
\lim_{j \to +\infty} \int_{t_{\ast}}^{T_{\ast}}\left\| W_{r_j}(\cdot,\cdot,t) - \Theta(\cdot,\cdot,t) \right\|_{H^1_{\mu_t}}^2 \,dt \to 0, \quad \text{ as } j \to +\infty,
\]
\[
\lim_{j \to +\infty} \sup_{t \in [t_{\ast},T_{\ast}]}\left\| W_{r_j}(\cdot,\cdot,t) - \Theta(\cdot,\cdot,t) \right\|_{L^2_{\mu_t}}^2  \to 0, \quad \text{ as } j \to +\infty,
\]
which are nothing more than \eqref{eq:L2STRONGCONVERGENCEOFWJ} and \eqref{eq:C0CONVERGENCEOFWJ}, respectively. So, it is enough to prove \eqref{eq:FINALLIMITH1BLOWUP} when $T_\ast$ is fixed and small (for instance $T_{\ast} < 1$). Now, let us consider the differences
\[
S_j := W_{r_j} - \Theta,
\]
which satisfy \eqref{eq:TruncatedExtensionEquationBlowUp} for any $j$, as $W_{r_j}$ do. A direct application of Lemma \ref{Lemma:HeightsVariations} implies
\[
\begin{aligned}
\int_0^{T_{\ast}} \int_{\RR_+^{N+1}} |\nabla S_j|^2 \,d\mu_t dt &- r_j^{1-a}\int_0^{T_{\ast}} \int_{\mathbb{R}^N\times\{0\}} q^{r_j} w_{r_j}^2 \, d\mu_t^{p_0}(x,0) + \int_0^{T_{\ast}} t\int_{\mathbb{R}_+^{N+1}} F_{r_j}W_{r_j} \, d\mu_tdt \\
&\leq \frac{1}{2} h(S_j,T_\ast),
\end{aligned}
\]
and so, since the second and the third terms in the r.h.s. go to zero as $j \to +\infty$ (see   \eqref{eq:ExpBoundIntF}, \eqref{eq:SecondTermBlowUp} and \eqref{eq:ThirdTermBlowUp}) and $T_\ast \leq 1$, we obtain
\[
\begin{aligned}
\int_0^{T_{\ast}} t\left\| \nabla W_{r_j} - \nabla \Theta \right\|_{L_{\mu_t}^2}^2 dt &\leq  \int_0^{T_{\ast}} \left\| \nabla W_{r_j} - \nabla \Theta \right\|_{L_{\mu_t}^2}^2 dt  \\
&= \int_0^{T_{\ast}} \int_{\RR_+^{N+1}} |\nabla S_j|^2 \,d\mu_t dt \leq h(S_j,T_\ast) \to 0,
\end{aligned}
\]
as $j \to +\infty$, thanks to \eqref{eq:C0CONVERGENCEOFWJ}. On the other hand, applying the Gaussian-Poincar\'e inequality to $S_j$ (see  Theorem \eqref{GAUSSIANPOICAREINEQUALITYSEVERALD}, part (iii)), we obtain
\[
\begin{aligned}
\int_0^{T_{\ast}} \|W_{r_j}-\Theta\|_{L^2_{\mu_t}}^2 dt = \int_0^{T_{\ast}} \int_{\RR^{N+1}_+} S_j^2 \,d\mu_t dt &\leq 2 \int_0^{T_{\ast}} t \int_{\RR^{N+1}_+} |\nabla S_j|^2 \,d\mu_t dt \\
&\leq 2 \int_0^{T_{\ast}}  \int_{\RR_+^{N+1}} |\nabla S_j|^2 \,d\mu_t dt \to 0,
\end{aligned}
\]
as $j \to +\infty$, obtaining \eqref{eq:L2STRONGCONVERGENCEOFWJ}.

\smallskip

\emph{Step 4: The limit $\Theta$ is a re-scaled eigenfunction.} First, notice that from \eqref{eq:FINALLIMITH1BLOWUP}, we have
\[
N_D(W_{r_j},\varrho) := \frac{D(W_{r_j},\varrho)}{H(W_{r_j},\varrho)} \to \frac{D_0(\Theta,\varrho)}{H(\Theta,\varrho)} = N_0(\Theta,\varrho),
\]
for any fixed $\varrho \in (0,T)$, as $j \to +\infty$. On the other hand, since $N_D(W_{r_j},\varrho) = N_D(W,r_j\varrho)$, it follows
\[
N_D(W_{r_j},\varrho) \to \kappa := \lim_{r \to 0} N_D(W,r),
\]
as $j \to +\infty$, and so
\begin{equation}\label{eq:BLOWUPPOONQUOTIENTCONSTANT}
N_D(\Theta,\varrho) \equiv \kappa,
\end{equation}
for any fixed $\varrho \in (0,T)$. Consequently, $\widetilde{\Theta}$ must be a weak eigenfunction to the Ornstein-Uhlenbeck eigenvalue problem type \eqref{eq:DEFINITIONWEAKEIGENFUNCTIONNEUMANN}, in the sense that the identity
\begin{equation}\label{eq:BLOWUPLIMITRESCALEDEIGENFUNCTION}
\int_{\RR_+^{N+1}} \nabla \widetilde{\Theta} \cdot \nabla \eta(t) \,d\mu = \kappa \int_{\RR_+^{N+1}} \widetilde{\Theta} \, \eta(t) \,d\mu,
\end{equation}
is satisfied for all $t \in (0,T_{\ast})$ and all $\eta(t) \in H_{\mu}^1$. It thus follows that $\kappa$ is an eigenvalue of the Ornstein-Uhlenbeck operator $\mathcal{O}_a$ (see  Remark \ref{Remark:HomogeneousStuff}) and we complete the proof of part (i). In particular, the blow-up limit has the form
\[
\Theta(X,t) = t^{\kappa} V\left(\frac{X}{\sqrt{t}}\right),
\]
where $V = V(X)$ is an eigenfunction to \eqref{eq:DEFINITIONWEAKEIGENFUNCTIONNEUMANN} with eigenvalue $\kappa$.

\smallskip

\emph{Step 5: The vanishing order equals $2\kappa$.} We prove \eqref{eq:LIMITOFTHEQUOTIENTFORSMALLTSTATEMENT1New}:
\[
\lim_{r \to 0^+} \mathcal{T}_{2\kappa}(W,r) = \lim_{r \to 0^+}  r^{-4\kappa} H(W,r) = L_0,
\]
for some $L_0 > 0$. This fact will be essential to show the uniqueness of the blow-up limit and implies that the vanishing order of $W$ at $O$ is $2\kk$. To do so, we first need an improvement of the formula shown in Corollary \ref{Corollary:NonDegeneracy}. From Theorem \ref{Theorem:GeneralizedFrequency} and Corollary \ref{Lemma:LimitAlmgrenPoon}, it follows that
\[
\Phi_a(W,r) = \frac{r}{4}e^{Cr^{1-a}} \frac{H'(W,r)}{H(W,r)} + e^{Cr^{1-a}} - 1 \geq \kappa
\]
for all $r \in (0,r_0)$. Now, rearranging terms, we deduce
\[
\frac{H'(W,r)}{H(W,r)} \geq \frac{4}{r} \left[ (\kappa + 1) e^{-Cr^{1-a}} - 1 \right]
\]
and so, using that $e^z \geq 1 + z$, it follows
\[
\frac{H'(W,r)}{H(W,r)} \geq \frac{4\kappa}{r} - C r^{-a},
\]
for some new constant $C > 0$. Integrating between $r_1$ and $r_2$, we obtain
\[
\frac{H(W,r_1)}{r_1^{4\kappa}} \leq e^{C(r_2^{1-a} - r_1^{1-a})} \frac{H(W,r_2)}{r_2^{4\kappa}}.
\]
In particular, it follows that the function $r \to \mathcal{T}_{2\kk}(W,r)$ is bounded for all $r \in (0,r_0)$. Now, integrating formula \eqref{eq:DerivativeTer} (with $\sigma = 2\kk$) between $r_1$ and $r_2$, it follows
\[
\mathcal{T}_{2\kk}(W,r_2) - \mathcal{T}_{2\kk}(W,r_1) = \int_{r_1}^{r_2} \frac{4}{\varrho} \mathcal{W}_{2\kk}(W,\varrho) d\varrho + \int_{r_1}^{r_2}\frac{4}{\varrho^{4\kk+3}} \int_0^{\varrho^2} t \int_{\mathbb{R}^{N+1}_+} F W \, d\mu_t dt d\varrho,
\]
and so, since the second term in the r.h.s. is integrable at $0$ (see  \eqref{eq:ExpBoundIntF}), we obtain
\[
\lim_{r_1 \to 0^+} \int_{r_1}^{r_2} \frac{4}{\varrho} \mathcal{W}_{2\kk}(W,\varrho) d\varrho = \int_0^{r_2} \frac{4}{\varrho} \mathcal{W}_{2\kk}(W,\varrho) d\varrho < +\infty,
\]
which implies that also $\mathcal{T}_{2\kk}(W,r_1)$ has a (finite) limit as $r_1 \to 0^+$. In particular, the limit in \eqref{eq:LIMITOFTHEQUOTIENTFORSMALLTSTATEMENT1New} exists and it is finite.

Assume by contradiction $r^{-4\kappa} H(W,r) \to 0$ as $r \to 0^+$ and consider the family of eigenfunctions
\[
\left\{V_{\alpha,m} = V_{\alpha,m}(X): \alpha \in \ZZ_{\geq 0}^N \;\text{ with } |\alpha| = n, \; n \in \NN, m \in \NN \right\},
\]
found in Theorem \ref{SPECTRALTHEOREMEXTENDEDORNUHL1}. An orthonormal basis of $L_{\mu}^2$ is given by
\[
\overline{V}_{\alpha,m}(X) = \frac{V_{\alpha,m}(X)}{\|V_{\alpha,m}\|_{L_{\mu}^2}}, \qquad \alpha \in \ZZ_{\geq 0}^N \; \text{ with } |\alpha| = n, \; n \in \NN, m \in \NN.
\]
Now, by scaling, we have that the parabolic re-scaling of $W = W(X,1)$:
\begin{equation}\label{eq:FOURIEREXPANSIONOFWLAMBDA}
W^{r}(X,1) :=  W(r X,r^2)
\end{equation}
belongs to $H_{\mu}^1$ for all $r \in (0,r_0)$ and so we can write
\[
W^{r}(X,1) = \sum_{\alpha,m}  w_{\alpha,m}(r) \overline{V}_{\alpha,m}(X)  \quad \text{ in } L_{\mu}^2,
\]
where the coefficients are given by
\[
w_{\alpha,m}(r) = \int_{\RR^{N+1}_+} W^{r}(X,1) \overline{V}_{\alpha,m}(X)\,d\mu(X).
\]
Using the de l'Hospital's rule, it is not difficult to see that the assumption $r^{-4\kappa} H(W,r) \to 0$ as $r \to 0^+$, implies $r^{-4\kappa} h(W,r^2) \to 0$ as $r \to 0^+$, where we recall that $H$ and $h$ are related by the relation
\[
H(W,r) = \frac{1}{r^2} \int_0^{r^2} h(W,t) dt.
\]
So by the orthonormality of the eigenfunctions $\overline{V}_{\alpha,m}$ and the usual scaling properties, it is easy to see that
\[
\begin{aligned}
h(W,r^2) &= \int_{\RR^{N+1}_+} \left[W^{r}(X,1)\right]^2 \, d\mu(X) \\
& = \int_{\RR^{N+1}_+} \left[\sum_{\alpha,m} w_{\alpha,m}(r) \overline{V}_{\alpha,m}(X)\right]^2 d\mu(X) \geq w_{\alpha,m}^2(r),
\end{aligned}
\]
for any $\alpha \in \ZZ_{\geq 0}^N$, $m \in \NN$, so that the assumption $r^{-4\kappa} H(W,r) \to 0$ as $r \to 0^+$ implies
\begin{equation}\label{eq:FALSEASSUMPTIONINBLOWUP}
r^{-2\kappa} w_{\alpha,m}(r) \to 0 \quad \text{ as } r \to 0^+,
\end{equation}
where we have used the fact that $h(W,r^2) = h(W^r,1)$. Now, writing \eqref{eq:StrongSolutionW1} for $t = r^2$ and changing variable $X \to rX$, we obtain
\begin{equation}\label{eq:EQUATIOFFORULAMBDAUNIQUENESSPOTENTIAL}
\begin{aligned}
\int_{\mathbb{R}^{N+1}_+} \frac{d W^r}{dr} \eta \, d\mu - \frac{2}{r} \int_{\mathbb{R}^{N+1}_+} \nabla W^r \cdot \nabla \eta \, d\mu = 2r \int_{\mathbb{R}^{N+1}_+} F^r \eta \, d\mu - 2r^{-a} \int_{\mathbb{R}^N\times\{0\}} q^r w^r \eta \, d\mu(x,0),
\end{aligned}
\end{equation}
for all $\eta \in H_{\mu}^1$ and all $r \in (0,r_0)$, where $F^r(X) = F^r(rX,r^2)$ and $w^r(x) = w(rx,r^2)$. Further, differentiating w.r.t. to $r$, it follows
\begin{equation}\label{eq:FOURIEREXPANSIONFORDERIVATIVEWRTLAMBDA}
\frac{d}{dr} W^{r}(X,1) = \sum_{\alpha,m} w_{\alpha,m}'(r) \overline{V}_{\alpha,m}(X) \quad \text{ in } L_{\mu}^2,
\end{equation}
where $w_{\alpha,m}'$ denotes the weak derivative of $w_{\alpha,m}$. Furthermore, since $\overline{V}_{\alpha,m}$ are solutions to \eqref{eq:WEAKFORMULATIONEIGENVALUEPROBLEM1} with $\kappa = \kappa_{n,m}$ (defined in the statement of Theorem \ref{SPECTRALTHEOREMEXTENDEDORNUHL1}), we deduce
\[
\begin{aligned}
\int_{\RR^{N+1}} \nabla W^{r} \nabla \eta \,d\mu &= \sum_{\alpha,m}  w_{\alpha,m}(r) \int_{\RR^{N+1}} \nabla \overline{V}_{\alpha,m} \nabla \eta \,d\mu \\
&= \sum_{\alpha,m}  \kappa_{n,m} w_{\alpha,m}(r) \int_{\RR^{N+1}} \overline{V}_{\alpha,m} \eta \,d\mu.
\end{aligned}
\]
Consequently, re-writing \eqref{eq:EQUATIOFFORULAMBDAUNIQUENESSPOTENTIAL} as
\[
\begin{aligned}
\sum_{\alpha,m} \left[w_{\alpha,m}'(r) - \frac{2\kappa_{n,m}}{r} w_{\alpha,m}(r) \right] \int_{\mathbb{R}^{N+1}_+} \overline{V}_{\alpha,m}\eta \, d\mu &= 2r \int_{\mathbb{R}^{N+1}_+} F^r \eta \, d\mu \\
& \quad - 2r^{-a} \sum_{\alpha,m} w_{\alpha,m}(r) \int_{\mathbb{R}^N\times\{0\}} q^r \overline{V}_{\alpha,m}(x,0) \eta \, d\mu(x,0),
\end{aligned}
\]
testing \eqref{eq:EQUATIOFFORULAMBDAUNIQUENESSPOTENTIAL} with $\overline{V}_{\alpha,m}$ and using the orthogonality properties of $\overline{V}_{\alpha,m}$ and their traces, it follows
\[
\left|w_{\alpha,m}'(r) - \frac{2\kappa_{n,m}}{r} w_{\alpha,m}(r) \right| \leq f_{\alpha,m}(r) + Cr^{-a}w_{\alpha,m}(r)
\]
for all $\alpha \in \ZZ_{\geq 0}^N$ with $|\alpha| = n$, $n \in \NN$, $m \in \NN$ and all $r \in (0,r_0)$, where we have set
\[
f_{\alpha,m}(r) := 2r \int_{\mathbb{R}^{N+1}_+} F^r V_{\alpha,m} \, d\mu,
\]
and chosen $C > 0$ depending on $N$, $a$ and the constant $K > 0$ in \eqref{eq:ASSUMPTIONSPOTENTIAL}. Integrating the above differential inequalities between $\varrho$ and $r$, we obtain
\begin{equation}\label{eq:FUNDAMENTALBOUNDSCOEFFPOTENTIAL}
-\int_{\varrho}^r \frac{f_{\alpha,m}(\rho)}{\rho^{2\kappa}} e^{C\rho^{1-a}} d\rho \leq \frac{w_{\alpha,m}(r)}{r^{2\kappa}} e^{Cr^{1-a}} - \frac{w_{\alpha,m}(\varrho)}{\varrho^{2\kappa}} e^{C\varrho^{1-a}} \leq \int_{\varrho}^r \frac{f_{\alpha,m}(\rho)}{\rho^{2\kappa}} e^{C\rho^{1-a}} d\rho.
\end{equation}
Taking the limit as $\varrho \to 0$ and rearranging terms, we deduce
\[
\frac{|w_{\alpha,m}(r)|}{r^{2\kappa}} \leq \int_0^r \frac{f_{\alpha,m}(\rho)}{\rho^{2\kappa}} d\rho.
\]
In particular, using the same methods that led to \eqref{eq:ExpBoundIntF} and \eqref{eq:SecondTermBlowUp} it follows
\[
|w_{\alpha,m}(r)| \leq r^{2\kappa} \int_0^r \frac{f_{\alpha,m}(\rho)}{\rho^{2\kappa}} d\rho  \leq C r^{\sigma_0} e^{-\frac{1}{16r^2}},
\]
for all $r \in (0,r_0)$, all $\alpha \in \ZZ_{\geq 0}^N$ with $|\alpha| = n$, $n \in \NN$, $m \in \NN$ and constants $\sigma_0, C > 0$ depending eventually on $n$ and $m$. Finally, since $W$ vanishes of finite order at $O$, it turns out
\[
\frac{|w_{\alpha,m}(r)|}{\sqrt{H(W,r)}} \leq C r^{\sigma_0} e^{-\frac{1}{16r^2}},
\]
for all $r \in (0,r_0)$, all $n \in \NN$, $m \in \NN$ and for some new constants $\sigma_0, C > 0$. So, if $r_j \to 0^+$ and $\widetilde{\Theta}$ is the blow-up limit of $W$ along the subsequence $r_j$, we obtain
\[
(W_{r_j},\overline{V}_{\alpha,m})_{L^2_{\mu}} = \frac{w_{\alpha,m}(r_j)}{\sqrt{H(W,r_j)}} \to 0,
\]
as $j \to +\infty$ and so, $(\widetilde{\Theta},\overline{V}_{\alpha,m})_{L^2_{\mu}}$ = 0, for all all $\alpha \in \ZZ_{\geq 0}^N$ with $|\alpha| = n$, $n \in \NN$, $m \in \NN$, i.e. $\widetilde{\Theta} = 0$. This is in contradiction with the fact that $H(\widetilde{\Theta},1) = 1$ (see  \eqref{eq:BLOWUPPLUSTRANSSEQUENCE}).

\smallskip

\emph{Step 6: Uniqueness of the blow-up.} Once \eqref{eq:LIMITOFTHEQUOTIENTFORSMALLTSTATEMENT1New} is established and the proof of part (ii) is completed, we show the uniqueness of the blow-up by proving that the limits \eqref{eq:C0CONVERGENCEOFWJ} and \eqref{eq:FINALLIMITH1BLOWUP} do not depend on the sub-sequence $r_j = r_{n_j}$ (this completes the proof of part (iii) and the proof of the theorem). So, let
\[
\Theta(X,t) = t^{\kappa} V\left(\frac{X}{\sqrt{t}}\right),
\]
be the blow-up limit along the sequence $r_j$, where $V = V(X)$ is an eigenfunction to \eqref{eq:DEFINITIONWEAKEIGENFUNCTIONNEUMANN} with eigenvalue $\kappa$. In view of \eqref{eq:LIMITOFTHEQUOTIENTFORSMALLTSTATEMENT1New}, we can re-write \eqref{eq:FINALLIMITH1BLOWUP} and \eqref{eq:C0CONVERGENCEOFWJ} as
\begin{equation}\label{eq:L2CONVERGENCEFORRESCALEDWLAMBDAJ}
\int_0^{T_{\ast}}
\left\| r_j^{-2\kappa} W(r_jX,r_j^2t) - t^{\kappa} \sum_{(\alpha,m) \in \widetilde{J}_0} v_{\alpha,m} \overline{V}_{\alpha,m}\left(\frac{X}{\sqrt{t}}\right) \right\|_{H^1(\RR_+^{N+1},d\mu_t)}^2 dt \to 0,
\end{equation}
\begin{equation}\label{eq:LINFINITYCONVERGENCEFORRESCALEDWLAMBDAJ}
\sup_{t \in [t_{\ast},T_{\ast}]}\left\| r_j^{-2\kappa} U(r_jX,r_j^2t) - t^{\kappa} \sum_{(\alpha,m) \in \widetilde{J}_0} v_{\alpha,m} \overline{V}_{\alpha,m}\left(\frac{X}{\sqrt{t}}\right) \right\|_{L^2(\RR_+^{N+1},d\mu_t)}^2  \to 0,
\end{equation}
as $j \to +\infty$, respectively. Above, we have written $V$ as the sum
\[
V(X) = \sum_{(\alpha,m) \in \widetilde{J}_0} v_{\alpha,m} \overline{V}_{\alpha,m}(X) \quad \text{ in } L_{\mu}^2,
\]
and implicitly substituted $v_{\alpha,m} \to L_0 v_{\alpha,m}$, where $L_0 > 0$ is the limit in \eqref{eq:LIMITOFTHEQUOTIENTFORSMALLTSTATEMENT1New}. In particular, choosing $t = 1$ (recall that we are assuming $T_{\ast} \geq 1$ for simplicity), we deduce
\begin{equation}\label{eq:POINTWISECONVERGENCEFORTAU1}
r_j^{-2\kappa} U(r_jX,r_j^2) \to \sum_{(\alpha,m) \in \widetilde{J}_0} v_{\alpha,m} \overline{V}_{\alpha,m}(X) \quad \text{ in } L_{\mu}^2.
\end{equation}
So, in other to prove that \eqref{eq:L2CONVERGENCEFORRESCALEDWLAMBDAJ} and \eqref{eq:LINFINITYCONVERGENCEFORRESCALEDWLAMBDAJ} hold for any subsequence $r_{n_j}\to 0^+$, or, equivalently, hold for $r \to 0^+$, we are left to show that the coefficients $v_{\alpha,m}$ are independent from $r_j$. Thus, proceeding as in the first part of this step, we expand
\[
W(r X,r^2) = \sum_{\alpha,m}  w_{\alpha,m}(r) \overline{V}_{\alpha,m}(X),
\]
in series of normalized eigenfunctions. Multiplying \eqref{eq:POINTWISECONVERGENCEFORTAU1} by $y^a \overline{V}_{\alpha,m}(X) \mathcal{G}_a(X,1)$ and integrating on $\RR_+^{N+1}$, we obtain
\[
r_j^{-2\kappa} w_{\alpha,m}(r_j) \to v_{\alpha,m}
\]
as $j \to +\infty$, where we have used again the orthogonality of the eigenfunctions $\overline{V}_{\alpha,m} = \overline{V}_{\alpha,m}(X)$. Now, re-writing \eqref{eq:FUNDAMENTALBOUNDSCOEFFPOTENTIAL} taking $\varrho = r_j$ and passing to the limit as $j \to + \infty$, it follows
\[
v_{\alpha,m} -\int_0^r \frac{f_{\alpha,m}(\rho)}{\rho^{2\kappa}} e^{C\rho^{1-a}} d\rho \leq \frac{w_{\alpha,m}(r)}{r^{2\kappa}} e^{Cr^{1-a}} \leq v_{\alpha,m} + \int_0^r \frac{f_{\alpha,m}(\rho)}{\rho^{2\kappa}} e^{C\rho^{1-a}} d\rho,
\]
for all $r \in (0,r_0)$. Taking the limit as $r \to 0^+$ and using that $\int_0^r \rho^{-2\kappa} f_{\alpha,m}(\rho) e^{C\rho^{1-a}} d\rho \to 0$, we deduce that $r^{-2\kappa}w_{\alpha,m}(r) \to v_{\alpha,m}$ as $r \to 0^+$ and so $v_{\alpha,m}$ does not depend on the subsequence $r_j$. This complete the proof of the theorem. $\Box$

\begin{ex} Out of clarity, we complete the section with some concrete examples of blow-up profiles. Let us put ourselves in the easiest case when the spacial dimension is $N=1$ and denote by
\[
\widetilde{\Theta}_{\alpha,m}(x,y,t) = t^{\widetilde{\kappa}_{n,m}} V_{\alpha,m}\left(\frac{x}{\sqrt{t}},\frac{y}{\sqrt{t}}\right) = t^{\widetilde{\kappa}_{n,m}} H_n\left(\frac{x}{\sqrt{t}}\right)L_{(\frac{a-1}{2}),m}\left(\frac{y^2}{4t}\right),
\]
the blow-up profile corresponding to $\widetilde{\kappa}_{n,m} = \frac{n}{2}+m$ (i.e. to the eigenfunction $V_{\alpha,m} = V_{\alpha,m}(x,y)$). Then we have:
\[
\begin{aligned}
   &\widetilde{\kappa}_{0,0} = 0           \qquad\qquad\,          \widetilde{\Theta}_{0,0}(x,y,t) = 1 \\
   &\widetilde{\kappa}_{1,0} = \frac{1}{2}    \qquad\qquad         \widetilde{\Theta}_{1,0}(x,y,t) = x \\
   &\widetilde{\kappa}_{2,0} = 1         \qquad\qquad\,         \widetilde{\Theta}_{2,0}(x,y,t) = x^2-2t \\
   &\widetilde{\kappa}_{0,1} = 1          \qquad\qquad\,           \widetilde{\Theta}_{0,1}(x,y,t) = \left(\frac{1+a}{2}\right)t - \frac{y^2}{4} \\
   &\widetilde{\kappa}_{3,0} = \frac{3}{2}    \qquad\qquad           \widetilde{\Theta}_{3,0}(x,y,t) = x(x^2 - 6t) \\
   &\widetilde{\kappa}_{1,1} = \frac{3}{2}     \qquad\qquad         \widetilde{\Theta}_{1,1}(x,y,t) = x\left[\left( \frac{1+a}{2}\right)t - \frac{y^2}{4}\right] \\
   &\widetilde{\kappa}_{4,0} = 2    \qquad\qquad \,          \widetilde{\Theta}_{4,0}(x,y,t) = x^4 - 12x^2t + 12t^2 \\
   &\widetilde{\kappa}_{2,1} = 2    \qquad\qquad \,          \widetilde{\Theta}_{2,1}(x,y,t) = (x^2-2t)\left[\left( \frac{1+a}{2}\right)t - \frac{y^2}{4}\right] \\
   &\widetilde{\kappa}_{0,2} = 2    \qquad\qquad \,          \widetilde{\Theta}_{0,2}(x,y,t) = \frac{1}{8}\left[(1+a)(3+a)t^2 - (3+a)y^2t + \frac{y^4}{4}\right],
\end{aligned}
\]
and so on.
\end{ex}
\begin{lem}\label{SIMONLEMMARECALL}(Simon \cite[Corollary 8]{Simon1987:art})
	Let $X \subset B \subseteq Y$ be Banach spaces satisfying the following two properties:
	\begin{itemize}
		\item $X$ is compactly embedded in $Y$.
		\item There exist $0 < \theta < 1$ and $C > 0$ such that $\|V\|_B \leq C\|V\|_X^{1-\theta}\|V\|_Y^{\theta}$, for all $V \in X \cap Y$.
	\end{itemize}
	Fix $-\infty < t_1 < t_2 < +\infty$ and $1\leq p,r \leq \infty$. Let $\mathcal{F}$ be a bounded family in $L^p(t_1,t_2;X)$, with $\partial\mathcal{F}/\partial t$ bounded in $L^r(t_1,t_2;Y)$. Then the following two assertions hold true:
	\begin{itemize}
		\item If $\theta(1-1/r) \leq (1-\theta)/p$ then $\mathcal{F}$ is relatively compact in $L^p(t_1,t_2;B)$, for all $p < p_{\ast}$, where $1/p_{\ast} = (1-\theta)/p - \theta(1-1/r)$.
		\item If $\theta(1-1/r) > (1-\theta)/p$ then $\mathcal{F}$ is relatively compact in $\mathcal{C}^0(t_1,t_2;B)$,
	\end{itemize}
	where we recall that $\mathcal{C}^0(t_1,t_2;B) := \{ V:[t_1,t_2] \to B \text{ continuous}: \|V\|_{\mathcal{C}^0(t_1,t_2;B)} < \infty\}$, where
	\[
	\|V\|_{\mathcal{C}^0(t_1,t_2;B)} := \max_{t \in [t_1,t_2]}\|V(t)\|_B.
	\]
\end{lem}
We end the section with the following lemma: it is necessary in the application of Lemma \ref{SIMONLEMMARECALL} both in the proof of Theorem \ref{THEOREMBLOWUP1New} and Theorem \ref{THEOREMBLOWUP1}.
\begin{lem}
The space $H_{\mu}^1 = H^1(\RR^{N+1},d\mu)$ is compactly embedded in $L_{\mu}^2 = L^2(\RR^{N+1},d\mu)$.
\end{lem}
\emph{Proof.} Let $\{U_j\}_j$ be a bounded sequence in $H_{\mu}^1$. Up to subsequences, we can assume the existence of a function $U \in H_{\mu}^1$ such that $U_j \rightharpoonup U$ (weakly) in $H_{\mu}^1$ and we must prove that
\[
U_j \to U \quad \text{ (strongly) in } L_{\mu}^2,
\]
as $j \to +\infty$. Let us set $V_j := U_j - U$. For any $A > 0$, we can write
\begin{equation}\label{eq:FIRSTSPLITTINGL2NORMCOMPACTIMMERSION}
\int_{\RR^{N+1}} V_j^2 \,d\mu(x,y) = \int_{\BB_A} V_j^2 \,d\mu(x,y) + \int_{\BB_A^c} V_j^2 \,d\mu(x,y),
\end{equation}
where, following \cite{BanGarofalo2017:art}, we have defined
\[
\BB_A := \{(x,y) \in \RR^{N+1}: |x|^2 + y^2 < A^2\} \quad \text{ with } \quad \BB_A^c = \RR^{N+1}\setminus\BB_A.
\]
Now, since the measures $|y|^a\,dxdy$ and $|y|^a \mathcal{G}_a(x,y,1)\,dxdy$ are equivalent on $\BB_A$, we obtain that the first term in the r.h.s. of \eqref{eq:FIRSTSPLITTINGL2NORMCOMPACTIMMERSION} goes to zero as $j \to +\infty$, thanks to the well-known $H^1(\BB_A,|y|^a) \hookrightarrow L^2(\BB_A,|y|^a)$ compact immersions type (see  for instance with \cite{FabesKenigSerapioni1982:art}). Consequently, we are left to prove that also the second integral in the r.h.s. of \eqref{eq:FIRSTSPLITTINGL2NORMCOMPACTIMMERSION} is converging to zero as $j \to +\infty$.

To do so, we repeat the procedure carried out in \cite{BanGarofalo2017:art} (see  Lemma 7.4) with some modifications due to our slightly different framework. Before moving forward, we recall the definition of the fundamental solution and how it splits into the product of two Gaussians:
\[
\mathcal{G}_a(x,y,t) = G_N(x,t) G_{a+1}(y,t) = \frac{1}{(4\pi t)^{\frac{N}{2}}} e^{-\frac{|x|^2}{4t}} \frac{1}{2^{1+a}\Gamma(\frac{1+a}{2})}  \frac{1}{t^{\frac{1+a}{2}}} e^{-\frac{y^2}{4t}},
\]
and that the following two log-Sobolev type inequalities hold (see  Lemma 7.7  of \cite{DainelliGarofaloPetTo2017:book} or formulas (7.16)-(7.17) of \cite{BanGarofalo2017:art}):
\begin{equation}\label{eq:NPLUS1LOGSOBOLEVFORMULA}
\log\left( \frac{1}{c_a \int_{|f| > 0} \mathcal{G}_a(\cdot,1)} \right) \int_{\RR^{N+1}} f^2 \mathcal{G}_a(\cdot,1) \leq 2 \int_{\RR^{N+1}} |\nabla f|^2 \mathcal{G}_a(\cdot,1),
\end{equation}
for all $f \in H^1(\RR^{N+1},d\mu)$, and
\begin{equation}\label{eq:NLOGSOBOLEVFORMULA}
\log\left( \frac{1}{\int_{|f| > 0} G_N(\cdot,1)} \right) \int_{\RR^N} f^2 G_N(\cdot,1) \leq 2 \int_{\RR^N} |\nabla f|^2 G_N(\cdot,1),
\end{equation}
for all $f \in H^1(\RR^N,d\mu_x)$. Note that, with respect to \cite{DainelliGarofaloPetTo2017:book} and/or \cite{BanGarofalo2017:art}, in \eqref{eq:NPLUS1LOGSOBOLEVFORMULA} we have introduced the quantity $c_a := 2^{1+a}\Gamma(\frac{1+a}{2})/\sqrt{4\pi}$ due to the different normalization constant of the fundamental solution (see  also with \cite{BanGarDanPetr2018:art}).

Now, fix $\varepsilon > 0$, and take $A > 4$ large enough such that
\begin{equation}\label{eq:ASSUMPTIONSONACOMPACTEMBEDDING}
c_a\int_{\BB_{A/2}^c} \mathcal{G}_a(x,y,1)\,dxdy \leq e^{-1/\varepsilon}, \qquad \int_{B_{A/2}^c} G_N(x,1)\,dx \leq e^{-1/\varepsilon},
\end{equation}
where $B_A := \{x \in \RR^N: |x| < A\}$ with $B_A^c = \RR^N\setminus B_A$. As explained in \cite{BanGarofalo2017:art}, \eqref{eq:NPLUS1LOGSOBOLEVFORMULA} cannot be directly applied to estimate the integral
\[
\int_{\BB_A^c} V_j^2 \,d\mu(x,y) = \int_{\BB_A^c} V_j^2 \,|y|^a \mathcal{G}_a(x,y,1)dxdy,
\]
due to some regularity issues. So, we write
\begin{equation}\label{eq:SPLITTINGREGINSL2NORM}
\int_{\BB_A^c} |y|^a V_j^2 \, \mathcal{G}_a(x,y,1)dxdy = \int_{\widetilde{R}_A} |y|^a V_j^2 \, \mathcal{G}_a(x,y,1)dxdy + \int_{R_A} |y|^a V_j^2 \, \mathcal{G}_a(x,y,1)dxdy,
\end{equation}
where
\[
\widetilde{R}_A := \BB_A^c \cap (\RR^N\times\{|y|\leq A/2\}) \quad \text{ and } \quad R_A := \BB_A^c \cap (\RR^N\times\{|y| > A/2\}).
\]
Let us start with estimating the first integral in the r.h.s. of \eqref{eq:SPLITTINGREGINSL2NORM}. For any $(x,y) \in \widetilde{R}_A$, we have both $|x|^2 + y^2 \geq A^2$ and $y^2 \leq A^2/4$, so that $|x| \geq (\sqrt{3}/2)A \geq A/2$. Consequently,
\[
\int_{\widetilde{R}_A} |y|^a V_j^2(x,y) \, \mathcal{G}_a(x,y,1)dxdy \leq \int_{-A/2}^{A/2} |y|^a G_{a+1}(y,1) \left(\int_{\{|x| \geq (\sqrt{3}/2)A\}} V_j^2(x,y) G_N(x,1) \,dx\right)dy.
\]
Now, if $\varphi \in \mathcal{C}_c^{\infty}(\RR^N)$ is a radially decreasing cut-off function satisfying $\varphi = 1$ in $B_{A/2}$ and $\varphi = 0$ outside $B_A$ with $|\nabla \varphi| \leq 1$, we firstly observe that
\[
\int_{\{|x| \geq (\sqrt{3}/2)A\}} V_j^2(x,y) G_N(x,1) \,dx \leq C\int_{\RR^N} V_j^2(x,y)[1-\varphi(x)]^2 G_N(x,1) \,dx,
\]
for some constant $C > 0$ depending only on $\varphi$ and for a.e. $y\in\RR$ (the above inequality can be directly verified using the properties of the cut-off function $\varphi = \varphi(x)$). Secondly, observing that the second inequality in \eqref{eq:ASSUMPTIONSONACOMPACTEMBEDDING}, together with the fact that
\[
\int_{|V_j(1-\varphi)| > 0} G_N(x,1)\,dx \leq \int_{B_{A/2}^c} G_N(x,1)\,dx
\]
implies
\[
\log\left( \frac{1}{\int_{|V_j(1-\varphi)| > 0} G_N(\cdot,1)}\right) \geq \frac{1}{\varepsilon}
\]
we deduce (by applying \eqref{eq:NLOGSOBOLEVFORMULA} with $f = V_j(1-\varphi)$):
\[
\int_{\RR^N} V_j^2(x,y)[1-\varphi(x)]^2 G_N(x,1) \,dx \leq \varepsilon \, C \int_{\RR^N} \left[V_j^2 + |\nabla_x V_j|^2\right] G_N(x,1) \,dx \leq \varepsilon \, C,
\]
where $C > 0$ is constant independent of $j\in\NN$ (this follows from the properties of $\varphi = \varphi(x)$ and since $\{V_j\}_j$ is bounded in $H_{\mu}^1$). Consequently, we obtain
\[
\int_{\widetilde{R}_A} |y|^a V_j^2(x,y) \,\mathcal{G}_a(x,y,1)dxdy \leq \varepsilon \, C\int_{-A/2}^{A/2} |y|^a G_{a+1}(y,1) \,dy \leq \varepsilon \,C,
\]
where $C > 0$ is new constant not depending on $j \in \NN$. Let us focus the second integral in the r.h.s. of \eqref{eq:SPLITTINGREGINSL2NORM}. We introduce the new cut-off functions
\begin{itemize}
  \item $\overline{\varphi} \in \mathcal{C}_c^{\infty}(\RR^{N+1})$ with $\overline{\varphi} = 1$ in $\BB_{A/2}$ and $\overline{\varphi} = 0$ outside $\BB_A$ with $|\nabla_{x,y}\overline{\varphi}| \leq 1$.
  \item $\widetilde{\varphi} \in \mathcal{C}_c^{\infty}(\RR)$ with $\widetilde{\varphi}(y) = 1$ for $|y|\leq A/4$ and $\widetilde{\varphi}(y) = 0$ for $|y| \geq A/2$ with $|\widetilde{\varphi}'| \leq 1$,
\end{itemize}
and we immediately see that
\[
\int_{R_A} |y|^a V_j^2 \,\mathcal{G}_a(x,y,1)dxdy \leq \int_{\RR^{N+1}} |y|^a V_j^2 \, [1 - \overline{\varphi}(x,y)]^2[1 - \widetilde{\varphi}(y)]^2 \, \mathcal{G}_a(x,y,1)dxdy.
\]
Now, exactly as in \cite{BanGarofalo2017:art} (see  formula (7.24)), we set $f = |y|^{a/2}V_j(1-\overline{\varphi})(1-\widetilde{\varphi})$ and we estimate
\[
|\nabla f|^2 \leq C\left[ \, |y|^a (V_j^2 + |\nabla V_j|^2) + |y|^{a-2}V_j^2(1 - \overline{\varphi})^2(1 - \widetilde{\varphi})^2 \,\right] \leq C |y|^a (V_j^2 + |\nabla V_j|^2),
\]
where $C > 0$ independent of $j \in \NN$. The last inequality follows since for any $|y| \geq A/4$, we have
\[
|y|^{a-2}V_j^2(1 - \overline{\varphi})^2(1 - \widetilde{\varphi})^2 \leq |y|^{a-2}V_j^2 \leq |y|^a V_j^2,
\]
where we have used the properties of the support of $\widetilde{\varphi} = \widetilde{\varphi}(y)$ and the fact the we have chosen $A > 4$ from the beginning. Consequently, applying \eqref{eq:NPLUS1LOGSOBOLEVFORMULA} (with $f = |y|^{a/2}V_j(1-\overline{\varphi})(1-\widetilde{\varphi})$), we obtain
\[
\int_{R_A} |y|^a V_j^2 \,\mathcal{G}_a(x,y,1)dxdy \leq \varepsilon\,C \int_{\RR^{N+1}} |y|^a \left[V_j^2 + |\nabla V_j|^2\right] \mathcal{G}_a(x,y,1) \,dxdy \leq \varepsilon \, C,
\]
where we have used the first inequality in \eqref{eq:ASSUMPTIONSONACOMPACTEMBEDDING} and the uniform $H_{\mu}^1$ bound on the $V_j's$. Summing up, from \eqref{eq:SPLITTINGREGINSL2NORM} and the above bounds, we have got
\[
\int_{\BB_A^c} |y|^a V_j^2 \, \mathcal{G}_a(x,y,1)dxdy \leq \varepsilon\, C,
\]
which completes the proof, by the arbitrariness of $\varepsilon > 0$. $\Box$
\begin{cor}
	The spaces $H_{\mu}^1 = H^1(\RR_+^{N+1},d\mu)$ and $H_{0,\mu}^1 = H_0^1(\RR_+^{N+1},d\mu)$  are compactly embedded in $L_{\mu}^2 = L^2(\RR^{N+1}_+,d\mu)$.
\end{cor}
\emph{Proof.} The proof is almost identical to the above one. $\Box$

%
%
%
%
%
%%%%%%%%%%%%%%%%%%%%%%%%%%%%%%%%%%%%%%%%%%%%%%%%%%%%%%%%%%%%%%%%%%%%%%%%%%%%%%%%%%%%%%%%%%%%%%%%%%%%
%
%%%%%%%%%%%%%%%%%%%%%%%%%%%%%%%%%%%%%%%%%%%%%%%%%%%%%%%%%%%%%%%%%%%%%%%%%%%%%%%%%%%%%%%%%%%%%%%%%%%%
%
%
%
%
%
\section{Blow-up Analysis: Part II}\label{Section:BlowUp2}
In this section we prove the theorem below, from which the statement of Theorem \ref{Theorem:UniformHolderBounds} will follow as a corollary. We employ the notation
\[
\mathbb{Q}_r^+(X_0,t_0) := \mathbb{B}_r^+(X_0) \times (t_0,t_0+r^2)
\]
for $(X_0,t_0) = (x_0,0,t_0)$ and $r > 0$. Following the notation of Section \ref{Section:SpectralAnalysis}, we set
\[
\nu_{\ast} = \frac{1}{2}\min\{1,1-a\},
\]
where, as always, $a \in (-1,1)$.
\begin{thm}\label{CONVERGENCELINFINITYLOC}
Let $a \in (-1,1)$ and let $W$ be a weak solution to \eqref{eq:TruncatedExtensionEquationLocalBackward}, with $q$ satisfying \eqref{eq:ASSUMPTIONSPOTENTIAL}. Assume that $W$ is bounded in $\mathbb{Q}_1^+$. Then for any $\nu \in (0,\nu_\ast)$, there exists a constant $C = C(\nu) > 0$ such that
\begin{equation}\label{eq:UNIFORMBOUNDHOLDERNORM1}
\|W\|_{\mathcal{C}^{2\nu,\nu}(\overline{\mathbb{Q}_{1/2}^+})} \leq C \|W\|_{L^{\infty}(\mathbb{Q}_1^+)}.
\end{equation}
In other words, locally bounded solutions to \eqref{eq:TruncatedExtensionEquationLocalBackward} are locally $\nu$-H\"older continuous (in the parabolic sense) for any $\nu \in (0,\nu_\ast)$.
\end{thm}
Assuming that Theorem \ref{CONVERGENCELINFINITYLOC} holds true, we can prove the following corollary. Notice that Theorem \ref{Theorem:UniformHolderBounds} follows directly from it by passing to the trace $W_{p_0,r}(x,0,t) = w_{p_0,r}(x,t) = v_{p_0,r}(x,t) = u_{p_0,r}(x,-t)$ in $\mathbb{B}_{1/2}^+\times\{0\}\times(0,1)$, since the cut-off $\zeta = 1$ in $\mathbb{B}_{1/2}$ (see  \eqref{eq:PropCutOff}).
\begin{cor}\label{Corollary:THMECONVERGENCELINFINITYLOC}
Let $a \in (-1,1)$, $p_0 \in \mathbb{B}_{1/2}^+\times\{0\}\times(0,1)$ and let $W$ be a weak solution to \eqref{eq:TruncatedExtensionEquationLocalBackward}, with $q$ satisfying \eqref{eq:ASSUMPTIONSPOTENTIAL}. Let $\{W_{p_0,r}\}_{r \in (0,1)}$ the blow-up family defined in \eqref{eq:BLOWUPPLUSTRANSSEQUENCE}. Then for any compact set $K \subset \overline{\RR^{N+1}_+}\times(0,\infty)$ and any $\nu \in (0,\nu_\ast)$, there exist  constants $C,p_0 > 0$ depending on $K$ and $\nu$ such that
\begin{equation}\label{eq:UNIFORMBOUNDHOLDERBlowUp}
\|W_{p_0,r}\|_{\mathcal{C}^{2\nu,\nu}(K)} \leq C.
\end{equation}
In other words, the blow-up family $\{W_{p_0,r}\}_{r \in (0,1)}$ is uniformly bounded in $\mathcal{C}^{2\nu,\nu}_{loc}$ for any $\nu \in (0,\nu_\ast)$.
\end{cor}
\emph{Proof.} Assume $p_0 = O$ for simplicity and notice that, since the cut-off $\zeta$ defined in \eqref{eq:PropCutOff} satisfies $\zeta = 1$ in $\mathbb{B}_{1/2}$, the blow-up family satisfies
\begin{equation}\label{eq:TruncatedExtensionEquationBlowUpLocal}
\begin{cases}
\partial_t W_r + y^{-a} \nabla \cdot(y^a \nabla W_r) = 0 \quad &\text{ in } \mathbb{B}^+_{1/(2r)} \times (0,1/r^2) \\
-\partial_y^a W_r = r^{1-a} q^r(x,t)w_r \quad &\text{ in } B_{1/(2r)} \times \{0\} \times (0,1/r^2),
\end{cases}
\end{equation}
for any fixed $r \in (0,1)$. Fix a compact $K \subset \overline{\RR^{N+1}_+}\times(0,\infty)$ and choose $r_0 \in (0,1)$ such that $K \subset \overline{\mathbb{B}^+_{1/(2r)}} \times (0,1/r^2)$ for any $r \in (0,r_0]$. Now, using the $L^2 - L^{\infty}$ estimate proved in \cite[Formula (5.7)]{BanGarofalo2017:art} (or the Harnack inequality \cite[Formula (5.12)]{BanGarofalo2017:art}), we easily deduce that $\|W_r\|_{L^{\infty}(K)} \leq C_0$ for all $r \in (0,r_0)$. Applying Theorem \ref{CONVERGENCELINFINITYLOC} to $W_r$ the thesis follows. $\Box$
\begin{cor}\label{Corollary:THMECONVERGENCELINFINITYLOCUNIQUENESS}
Let $a \in (-1,1)$, $p_0 \in \mathbb{B}_{1/2}^+\times\{0\}\times(0,1)$ and let $W$ be a weak solution to \eqref{eq:TruncatedExtensionEquationLocalBackward}, with $q$ satisfying \eqref{eq:ASSUMPTIONSPOTENTIAL}. Let $\{W_{p_0,r}\}_{r \in (0,1)}$ the blow-up family defined in \eqref{eq:BLOWUPPLUSTRANSSEQUENCE}. Then
for any $\nu \in (0,\nu_\ast)$, it holds
\begin{equation}\label{eq:HOLDERCONVERGENCEBlowUp}
W_{p_0,r} \to \Theta_{p_0} \quad \text{ in } \mathcal{C}^{2\nu,\nu}_{loc},
\end{equation}
as $r \to 0$, where $\Theta_{p_0}$ is the blow-up limit of $W$ at $p_0$ found in Theorem \ref{THEOREMBLOWUP1New} part (iii).
\end{cor}
\emph{Proof.} Take $p_0 = O$. From Corollary \ref{Corollary:THMECONVERGENCELINFINITYLOC} we know that for any $\alpha \in (0,\nu_\ast)$, $W_r$ is uniformly bounded in $\mathcal{C}^{2\alpha,\alpha}_{loc}$ and so, by the Ascoli-Arzel\'a theorem, there exists a locally continuous function $W_0$ such that $W_r \to W_0$ locally uniformly and in $\mathcal{C}^{2\nu,\nu}_{loc}$ as $j \to +\infty$, for any fixed $\nu \in (0,\alpha)$, up to passing to a suitable subsequence. If we show that $W_0 = \Theta$, independently of the subsequence, the thesis follows.

So, assume by contradiction there exists a subsequence $r_j \to 0$ such that $W_{r_j} \not\to \Theta$ locally uniformly and fix $0 < t_1 < t_2$. Applying Theorem \ref{THEOREMBLOWUP1New} and Corollary \ref{Corollary:THMECONVERGENCELINFINITYLOC} to $W_{r_j}$, we deduce that, up to passing to a suitable subsequence (that we relabel $r_j$), $W_{r_j} \to \Theta$ in $\mathcal{C}(t_1,t_2;L^2_{\mu_t})$ and $W_{r_j} \to W_0$ locally uniformly, for some continuous function $W_0 \not= \Theta$. Consequently, for any compact set $K \subset \RR_+^{N+1}$, it follows
\[
\sup_{t \in (t_1,t_2)} \int_K |\Theta - W_0|^2 d\mu_t \leq \sup_{t \in (t_1,t_2)} \int_K |W_{r_j} - W_0|^2 d\mu_t + \sup_{t \in (t_1,t_2)} \|W_{r_j} - \Theta\|_{L^2_{\mu_t}}^2,
\]
and so, passing to the limit as $j \to +\infty$ and using the Lebesgue dominated convergence theorem, the two terms in the r.h.s. go to zero. This implies that $\Theta = W_0$ on $K \times(t_1,t_2)$, in contradiction with our assumption. $\Box$
%
%
%%%%%%%%%%%%%%%%%%%%%%%%%%%%%%%%%%%%%%%%%%%%%%%%%%%%%%%%%%%%%%%%%%%%%%%%%%%%%%%%%%%%%%%%%%%%%%%%%%%%%%%%%%%
%
%
%%%%%%%%%%%%%%%%%%%%%%%%%%%%%%%%%%%%%%%%%%%%%%%%%%%%%%%%%%%%%%%%%%%%%%%%%%%%%%%%%%%%%%%%%%%%%%%%%%%%%%%%%%%
%
\subsection{Liouville type theorems}\label{SECTIONLIOUVILLETHEOREMS} The proof of Theorem \ref{CONVERGENCELINFINITYLOC} requires some Liouville type results that we present in the next paragraphs. We prove two types of Liouville theorems: the first one will be obtained as an easy application of the monotonicity of the Almgren-Poon quotient, the Gaussian Poincar\'e inequalities proved in Subsection \ref{POINCAREGAUSSIANINEQUALITY}, and Theorem \ref{THEOREMBLOWUP1}, while the second requires a different monotonicity formula of Alt-Caffarelli-Friedman type that we prove below.

So we consider global weak solutions to problem \eqref{eq:TruncatedExtensionEquationLocalBackward}, with $q = F = 0$, i.e. solutions to
\begin{equation}\label{eq:DEFINITIONPARABOLICEQUATION}
\begin{cases}
\partial_t U + y^{-a}\nabla\cdot( y^a \nabla U) = 0 \quad &\text{in } \RR_+^{N+1}\times(0,1) \\
-\partial_y^a U = 0 \quad &\text{in } \RR^N\times\{0\}\times(0,1).
\end{cases}
\end{equation}
Following the approach of Section \ref{Section:SpectralAnalysis}, we study also solutions
to the homogeneous Dirichlet problem:
\begin{equation}\label{eq:DEFINITIONPARABOLICEQUATION2}
\begin{cases}
\partial_t U + y^{-a}\nabla\cdot( y^a \nabla U) = 0 \quad &\text{in } \RR_+^{N+1}\times(0,1) \\
U = 0 \quad &\text{in } \RR^N\times\{0\}\times(0,1),
\end{cases}
\end{equation}
and the one posed in the whole $\RR^{N+1}$:
\begin{equation}\label{eq:DEFINITIONPARABOLICEQUATION1}
\partial_t U + |y|^{-a}\nabla\cdot( |y|^a \nabla U) = 0 \quad \text{in } \RR^{N+1}\times(0,1).
\end{equation}
This is very natural since the techniques we use to treat solutions to \eqref{eq:DEFINITIONPARABOLICEQUATION} work also for the other two problems with minor changes. The definitions of weak solutions to \eqref{eq:DEFINITIONPARABOLICEQUATION2} and \eqref{eq:DEFINITIONPARABOLICEQUATION1} are almost identical to the one for \eqref{eq:DEFINITIONPARABOLICEQUATION} and we do not enter into the details. In the same way, we can give a notion af vanishing order for weak solutions to \eqref{eq:DEFINITIONPARABOLICEQUATION2} and \eqref{eq:DEFINITIONPARABOLICEQUATION1} in the spirit of Definition \ref{def:VanishingOrder}. So we limit ourselves to state the analogue of Theorem \ref{Theorem:MonAlmgrenPoon} for solutions to \eqref{eq:DEFINITIONPARABOLICEQUATION2} and \eqref{eq:DEFINITIONPARABOLICEQUATION1}.
\begin{thm}\label{Theorem:MonotonicityNewPoon}
The following two statements hold true:

(i) Let $U$ be a weak solution to \eqref{eq:DEFINITIONPARABOLICEQUATION2} having finite vanishing order at $O$. Then the statement of Theorem \ref{Theorem:MonAlmgrenPoon} is satisfied by $U$ for all $t,r > 0$.

(ii) Let $U$ be a weak solution to \eqref{eq:DEFINITIONPARABOLICEQUATION1} having finite vanishing order at $O$. Define
\[
\begin{aligned}
h(U,t) &:= \int_{\mathbb{R}^{N+1}} U^2 |y|^a \mathcal{G}_a(X,t)\,dX \\
d_0(U,t) &:= t\int_{\mathbb{R}^{N+1}} |\nabla U|^2 |y|^a \mathcal{G}_a(X,t)\,dX \\
n_0(U,t) &:= \frac{d_0(U,t)}{h(U,t)}
\end{aligned}
\qquad \qquad
\begin{aligned}
H(U,r) &:= \frac{1}{r^2} \int_0^{r^2} h(U,t) \, dt \\
D_0(U,r) &:= \frac{1}{r^2} \int_0^{r^2} d_0(U,t) \, dt \\
N_0(U,r) &:= \frac{D_0(U,r)}{H(U,r)},
\end{aligned}
\]
for all $t,r > 0$. Then, with these definitions, the statement of Theorem \ref{Theorem:MonAlmgrenPoon} is satisfied by $U$ for all $t,r > 0$.
\end{thm}
\begin{rem}\label{Rem:GrowthCond}
Notice that we are implicitly assuming that weak solutions $U$ to \eqref{eq:DEFINITIONPARABOLICEQUATION1} satisfy some suitable growth conditions when $|X|$ is large. If not, in view of the Jones' result \cite{Jones1977:art}, we cannot expect that $U$ vanishes of finite order (i.e. $N_0(U,r)$ is bounded). In other words, Theorem \ref{Theorem:MonotonicityNewPoon} cannot hold for weak solutions growing too fast at infinity.

In the following lemma, we will consider solutions $U$ having at most ``polynomial growth'', that is
\[
|U(X,t)| \leq C \sqrt{1 + d^{2n}(X,t)}
\]
for some $C >0$ and some $n \in \NN$, where $d = d(X,t)$ denotes the ``parabolic distance'':
\[
d(X,t) := \sqrt{|X|^2 + t}, \quad \text{ for all } (X,t) \in \RR^{N+1}\times(0,\infty).
\]
It is not a distance in the usual sense, since the triangular inequality does not hold, but it is parabolically homogeneous of degree one: $d(rX,r^2t) = r d(X,t)$, for all $r > 0$.
\end{rem}
As mentioned above, the second result of Liouville type is based on an Alt-Caffarelli-Friedman (\cite{AltCafFried1984:art}) type monotonicity formula (see  for instance with \cite[Lemma 5.4]{AthCaffaMilakis2016:art}, \cite[Theorem 1.1.4]{CafKen1998:art}, \cite[Theorem 12.11]{CafSal2005:book}). It turns out to be an easy consequence of the Gaussian-Poincar\'e type inequalities proved in Theorem \ref{GAUSSIANPOICAREINEQUALITYSEVERALD}.

\begin{lem}\label{ALTCAFFARELLIFRIEDMANLEMMAONESPIECES}(Alt-Caffarelli-Friedman type monotonicity formula)
Let $a \in (-1,1)$. Then the following three statements hold:

(i) Let $U$ be a weak solution to problem \eqref{eq:DEFINITIONPARABOLICEQUATION}, having at most polynomial growth. Then the function
\begin{equation}\label{eq:DEFFUNCTIONALJ}
t \to J(U,t) := \frac{1}{t} \int_0^t \int_{\RR_+^{N+1}} |\nabla U|^2(X,\tau) \,d\mu_{\tau}(X) d\tau
\end{equation}
is nondecreasing for all $t \in (0,1)$ and it is constant if and only if $U(X,t) = A$ or $U(X,t) = Ax_j         $ for some $j \in \{1,\ldots,N\}$, $A \in \RR$.

(ii) Let $U$ be a weak solution to problem \eqref{eq:DEFINITIONPARABOLICEQUATION2}, having at most polynomial growth. Then the function
\begin{equation}\label{eq:DEFFUNCTIONALJ2}
t \to J(U,t) := \frac{1}{t^{1-a}} \int_0^t \int_{\RR_+^{N+1}} |\nabla U|^2(X,\tau) \,d\mu_{\tau}(X) d\tau
\end{equation}
is nondecreasing for all $t \in (0,1)$ and it is constant if and only if $U(X,t) = Ay^{1-a}$, $A \in \RR$.

(iii) Let $U$ be a weak solution to equation \eqref{eq:DEFINITIONPARABOLICEQUATION1}, having at most polynomial growth. Then the function
\begin{equation}\label{eq:DEFFUNCTIONALJ1}
t \to J(U,t) := \frac{1}{t^{2\nu_{\ast}}} \int_0^t \int_{\RR^{N+1}} |\nabla U|^2(X,\tau) \,d\mu_{\tau}(X) d\tau
\end{equation}
is nondecreasing for all $t \in (0,1)$ and it is constant if and only if $U(X,t) = A$ or, depending on $a \in (-1,1)$:
\[
U(X,t) =
\begin{cases}
Ax_j         \quad &\text{ if }  a \in (-1,0)  \quad \text{ for some } j \in \{1,\ldots,N\}\\
Ax_j         \quad &\text{ if }  a = 0  \;\quad\qquad  \text{ for some } j \in \{1,\ldots,N+1\}\\
Ay|y|^{-a}   \quad &\text{ if }  a \in (0,1),
\end{cases}
\]
where $A \in \RR$ and we have used the convention $x_{N+1} = y$.
\end{lem}
\emph{Proof.} We begin by proving assertion (i). First of all, we have
\begin{equation}\label{eq:DERIVATIVEOFJLEMMA}
J'(U,t) = -\frac{1}{t^2} \int_0^t \int_{\RR_+^{N+1}} |\nabla U|^2 \,d\mu_{\tau}(X) d\tau + \frac{1}{t^2} \int_{\RR_+^{N+1}} |\nabla U|^2 \,d\mu_{\tau}(X),
\end{equation}
for all $t \in (0,1)$. Furthermore, testing the equation of $U$ with $\eta = U\mathcal{G}_a$ (see  \eqref{eq:StrongSolutionW1} with $q = F = 0$), we obtain
\[
\begin{aligned}
\int_0^t \int_{\RR_+^{N+1}} |\nabla U|^2 \,d\mu_{\tau}(X) d\tau &= \int_0^t \int_{\RR_+^{N+1}} y^a \nabla U \cdot \nabla U \,\mathcal{G}_a(X,\tau)\,dXd\tau \\
& = \int_0^t \int_{\RR_+^{N+1}}  y^a U \partial_{\tau}U \mathcal{G}_a\,dXd\tau - \int_0^t \int_{\RR_+^{N+1}}  y^a U \nabla U \cdot \nabla\mathcal{G}_a\,dXd\tau \\
& = \frac{1}{2} \int_0^t \int_{\RR_+^{N+1}}  y^a \partial_{\tau}\left(U^2\right) \mathcal{G}_a\,dXd\tau + \frac{1}{2}\int_0^t \int_{\RR_+^{N+1}} y^a  U^2 \partial_{\tau} \mathcal{G}_a\,dXd\tau \\
& = \frac{1}{2} \int_{\RR_+^{N+1}} \int_0^t y^a\partial_{\tau}\left(U^2 \mathcal{G}_a\right)\,dXd\tau \\
& \leq \frac{1}{2} \int_{\RR_+^{N+1}} U^2(X,t) \,d\mu_t(X) < +\infty,
\end{aligned}
\]
where we have used the equations of $U$ and $\mathcal{G}_a$, and the identity
\[
2\int_{\RR_+^{N+1}}  y^a U \nabla U \cdot \nabla\mathcal{G}_a\,dX = - \int_{\RR_+^{N+1}} U^2 \nabla\cdot(y^a \nabla \mathcal{G}_a) \,dX,
\]
for all $t \in (0,1)$. Notice that $\eta = U\mathcal{G}_a$ is an admissible test, since it decays exponentially fast at infinity (we are crucially using the assumption of polynomial growth). Consequently, combining \eqref{eq:DERIVATIVEOFJLEMMA} and the above estimate, it suffices to prove
\[
\frac{1}{t^2} \left\{ t \int_{\RR_+^{N+1}} |\nabla U|^2 \,d\mu_t(X) - \frac{1}{2} \int_{\RR_+^{N+1}} U^2 \,d\mu_t(X) \right\} \geq 0,
\]
which can be easily re-written as
\begin{equation}\label{eq:GAUSSIANPOINCAREDEPENDINGONT0}
\frac{t \int_{\RR_+^{N+1}} |\nabla U|^2 \,d\mu_t(X)}{\int_{\RR_+^{N+1}} U^2 \,d\mu_t(X)} \geq \frac{1}{2},
\end{equation}
and, passing to the re-normalized variables $\widetilde{U}(X,t) = U(\sqrt{t}X,t)$, as
\begin{equation}\label{eq:GAUSSIANPOINCAREDEPENDINGONT}
\frac{\int_{\RR_+^{N+1}} |\nabla \widetilde{U}|^2 \,d\mu(X)}{\int_{\RR_+^{N+1}} \widetilde{U}^2 \,d\mu(X)} \geq \frac{1}{2}.
\end{equation}
Now, noticing that $\widetilde{U}(t) \in H^1_{\mu}$ for all $t \in (0,1)$ by definition and that
\begin{equation}\label{eq:ZeroMean}
t \to \int_{\mathbb{R}_+^{N+1}} \widetilde{U}(X,t) d\mu(X)
\end{equation}
is constant, we can assume $\int_{\mathbb{R}_+^{N+1}} \widetilde{U}(t) d\mu(X) = 0$ (if not, we can always translate de solution). Consequently, \eqref{eq:GAUSSIANPOINCAREDEPENDINGONT} follows by the Gaussian-Poincar\'e inequality proved in Theorem \ref{GAUSSIANPOICAREINEQUALITYSEVERALD}, part (i). Notice that \eqref{eq:ZeroMean} formally follows by testing \eqref{eq:RescaledTruncatedExtensionEquationLocalBackward} with $\eta = 1$. This formal procedure can be rigorously justified by a standard approximation method and we omit it.

The proofs of parts (ii) and (iii) are very similar: it is enough to repeat the above procedure using a suitable Gaussian-Poincar\'e inequality (see  part (ii) and (iii) of Theorem \ref{GAUSSIANPOICAREINEQUALITYSEVERALD}).

It remains to establish on which solutions $U$ the functional $t \to J(U,t)$ is constant. To do so, we consider the functional defined in \eqref{eq:DEFFUNCTIONALJ1} which, from this point of view, is more general (the other two cases follows as immediate consequences). So, we first fix $a \in (-1,0)$ (so that $\nu_{\ast} = 1/2$) and notice that Theorem \ref{GAUSSIANPOICAREINEQUALITYSEVERALD} implies that the equality is attained in (see  \eqref{eq:GaussianPoincWholeSpace})
\[
\frac{\int_{\RR^{N+1}} |\nabla \widetilde{U}|^2 \,d\mu(X)}{\int_{\RR^{N+1}} \widetilde{U}^2 \,d\mu(X)} \geq \nu_{\ast},
\]
if and only if
\[
\widetilde{U}(X,t) = A(t)x_j \qquad \text{ i.e. } \qquad U(X,t) = \frac{A(t)}{\sqrt{t}}\,x_j,
\]
for some $j \in \{1,\ldots,N\}$ and function $A = A(t)$. Since such $U$ must be a solution to \eqref{eq:DEFINITIONPARABOLICEQUATION1} and it satisfies $\nabla\cdot(|y|^a\nabla U) = 0$, we deduce that $A(t) = \sqrt{t}$ (up to multiplicative constants), the statement in the case $a \in (-1,0)$ follows. In the case $a \in (0,1)$ ($\nu_{\ast} = (1-a)/2$), we proceed in the same way, noting that the equality is attained if and only if
\[
\widetilde{U}(X,t) = A(t)\,y|y|^{-a} \qquad \text{ i.e. } \qquad U(X,t) = t^{-\frac{1-a}{2}}A(t)\, y|y|^{-a}.
\]
Consequently, noticing that $y|y|^{-a}$ is a stationary solution to \eqref{eq:DEFINITIONPARABOLICEQUATION1} and taking $A(t) = t^{\frac{1-a}{2}}$ (up to multiplicative constants), we conclude the proof of the case $a \in (0,1)$, too. The case $a = 0$ follows exactly in the same way. $\Box$

\bigskip

We are ready to prove our Liouville type results. In this context, it is important to stress that we focus on the \emph{global} quantitative behaviour of solutions, instead of the local one (this is easily seen, for instance, in the interest of the asymptotic behaviour of the quantity $h(U,t)$ as $t \to \infty$, instead of $t \to 0^+$).

\begin{lem}\label{1LIOUVILLETYPETHEOREMFORWEIGHTEDHE}
Let $a \in (-1,1)$. Then the following three statements hold.

(i) Let $U$ be a weak solution to problem \eqref{eq:DEFINITIONPARABOLICEQUATION} and assume it satisfies the bound
\begin{equation}\label{eq:POINTWISEBOUNDFORLIOUVILLETHEOREMFORWEIGHTEDHE}
|U(X,t)| \leq C \sqrt{1 + d^{2\gamma}(X,t)},
\end{equation}
for some $C > 0$ and some exponent
\[
\gamma \in (0,1).
\]
Then $U$ is constant in $\RR_+^{N+1}\times(0,\infty)$.

(ii) Let $U$ be a weak solution to problem \eqref{eq:DEFINITIONPARABOLICEQUATION2} and assume it satisfies the bound in \eqref{eq:POINTWISEBOUNDFORLIOUVILLETHEOREMFORWEIGHTEDHE} for some $C > 0$ and some exponent
\[
\gamma \in (0,1-a).
\]
Then $U$ is identically zero in $\RR_+^{N+1}\times(0,\infty)$.

(iii) Let $U$ be a weak solution to equation \eqref{eq:DEFINITIONPARABOLICEQUATION1} and assume it satisfies the bound in \eqref{eq:POINTWISEBOUNDFORLIOUVILLETHEOREMFORWEIGHTEDHE}, for some $C > 0$ and some exponent
\[
\gamma \in (0,2\nu_{\ast}).
\]
Then $U$ is constant in $\RR^{N+1}\times(0,\infty)$.
\end{lem}
\emph{Proof.} We begin by proving part (i). First of all, we notice that for any $\gamma > 0$, if $U$ satisfies the point-wise bound \eqref{eq:POINTWISEBOUNDFORLIOUVILLETHEOREMFORWEIGHTEDHE}, it holds
\[
\begin{aligned}
h(U,t) & = \int_{\RR_+^{N+1}} U^2(X,t) \,d\mu_t(X) \\
&\leq C \int_{\RR_+^{N+1}} \left[(|X|^2 + t)^{\gamma} + 1 \right] y^a \mathcal{G}_a(X,t) \,dX \\
& = C + t^{\gamma} \int_{\RR_+^{N+1}} (|X|^2 + 1)^{\gamma} y^a \mathcal{G}_a(X,1) \,dX \leq C (1 + t^{\gamma}),
\end{aligned}
\]
for a new constant $C > 0$ and all $t > 0$. Consequently,
\begin{equation}\label{eq:BOUNDFROMABOVEFROTLARGELIOUVILLE}
H(U,r) \leq C(1 + r^{2\gamma}),
\end{equation}
for some new $C > 0$ and all $r > 0$. Now, set as always $\kappa := \lim_{r \to 0^+} N_0(U,r)$ and assume by contradiction that $U$ is non constant. Let us consider first the case $\kappa > 0$. Under this assumption, it must be
\begin{equation}\label{eq:BOUNDFROMBELOWFORKAPPA}
\kappa \geq \frac{1}{2},
\end{equation}
since $\kappa > 0$ has to be an eigenvalue of problem \eqref{eq:ORNSTUHLENFOREXTENSIONNEUMANN} (see  Theorem \ref{THEOREMBLOWUP1New}) and $1/2$ is the first nontrivial eigenvalue (see  Theorem \ref{SPECTRALTHEOREMEXTENDEDORNUHL1} part (i)). On the other hand,
\[
\frac{rH'(U,r)}{4H(U,r)} = N_0(U,r) \geq \kappa,
\]
all $r > 0$, since $r \to N_0(U,r)$ is nondecreasing (see  Theorem \ref{Theorem:MonAlmgrenPoon}).
Therefore, integrating
\[
\frac{H'(U,r)}{H(U,r)} \geq \frac{4\kappa}{r},
\]
it follows that the function $r \to r^{-4\kappa} H(U,r)$ is nondecreasing and so, recalling \eqref{eq:LIMITOFTHEQUOTIENTFORSMALLTSTATEMENT1New}, it follows that
\[
\frac{H(U,r)}{r^{4\kappa}} \geq L_0,
\]
for some $L_0 > 0$ and all $r > 0$. Recalling that $4\kappa \geq 2 > 2\gamma$, we obtain a contradiction with \eqref{eq:BOUNDFROMABOVEFROTLARGELIOUVILLE} taking the limit as $r \to +\infty$.

Now, assume that $\kappa = 0$. In this case, the bound in \eqref{eq:BOUNDFROMABOVEFROTLARGELIOUVILLE} is still true but we do not have a useful bound from below leading us to the desired contradiction. So, we pass to the re-scaling $\widetilde{U}(X,t) = U(\sqrt{t}X,t)$ and we recall that $\widetilde{U}(t) \in H^1(\RR_+^{N+1},d\mu)$ for all $t > 0$. Consequently, we can apply the Gaussian Poincar\'e inequality (see  Theorem \ref{GAUSSIANPOICAREINEQUALITYSEVERALD}, part (i)), to deduce
\begin{equation}\label{eq:GAUSSIANPOINCARETIMELIOUVILLE}
h(U,t) = \int_{\RR_+^{N+1}} \widetilde{U}^2(t) \,d\mu \leq 2 \int_{\RR_+^{N+1}} |\nabla \widetilde{U}|^2(t) \,d\mu = d_0(U,t),
\end{equation}
for all $t > 0$. Notice that we are implicitly assuming $\int_{\RR_+^{N+1}} \widetilde{U}(t) \,d\mu = 0$ for all $t > 0$. The general case is obtained by recalling that the function
\[
t \to \int_{\RR_+^{N+1}} \widetilde{U}(t) \,d\mu
\]
is constant and since
\[
\widetilde{U}(X,t) - \int_{\RR_+^{N+1}} \widetilde{U}(t) \,d\mu,
\]
is still a solution to \eqref{eq:DEFINITIONPARABOLICEQUATION} with zero mean for all $t > 0$. Consequently, integrating \eqref{eq:GAUSSIANPOINCARETIMELIOUVILLE} from $0$ to $r^2$, we obtain
\[
\frac{1}{2} \leq N_0(U,r),
\]
for all $r > 0$ and, since the l.h.s. is strictly positive, while the r.h.s. goes to zero as $r \to 0^+$, we obtain a contradiction (note that we infer $\widetilde{U}\equiv 0$ since we have assumed that $\widetilde{U}$ has zero mean for any $t > 0$).

For what concerns part (ii) and (iii), it is easily seen that the above proof works also in these different settings with straightforward modifications. The most significative is that we have to employ the ``right'' Gaussian Poincar\'e inequality proved in Theorem \ref{GAUSSIANPOICAREINEQUALITYSEVERALD}, depending on the problem (Dirichlet or Neumann) and applying Theorem \ref{Theorem:MonotonicityNewPoon} instead of Theorem \ref{Theorem:MonAlmgrenPoon}. $\Box$

\bigskip

We now prove a second Liouville type result of different nature. In the first one, we impose a critical growth condition one the function $U$, whilst, in the second, a sort of decaying property of the norm of the gradient.

\begin{lem}\label{2LIOUVILLETYPETHEOREMFORWEIGHTEDHE}
Let $a \in (-1,1)$. Then the following three statements hold.

(i) Let $U$ be a weak solution to problem \eqref{eq:DEFINITIONPARABOLICEQUATION} and assume it satisfies the bound
\begin{equation}\label{eq:2POINTWISEBOUNDFORLIOUVILLETHEOREMFORWEIGHTEDHE}
|\nabla U(x,y,t)| \leq C \, d^{\gamma}(x,y,t),
\end{equation}
for some $C > 0$ and some exponent
\[
\gamma < 0.
\]
Then $U$ is constant in $\RR_+^{N+1}\times(0,\infty)$.

(ii) Let $U$ be a weak solution to problem \eqref{eq:DEFINITIONPARABOLICEQUATION2} and assume it satisfies the bound in \eqref{eq:2POINTWISEBOUNDFORLIOUVILLETHEOREMFORWEIGHTEDHE} for some $C > 0$ and some exponent
\[
\gamma < -a.
\]
Then $U$ is identically zero in $\RR_+^{N+1}\times(0,\infty)$.

(iii) Let $U$ be a weak solution to equation \eqref{eq:DEFINITIONPARABOLICEQUATION1} in and assume it satisfies the bound in \eqref{eq:2POINTWISEBOUNDFORLIOUVILLETHEOREMFORWEIGHTEDHE} for some $C > 0$ and some exponent
\[
\gamma < \min\{0,-a\}.
\]
Then $U$ is constant in $\RR^{N+1}\times(0,\infty)$.
\end{lem}
\emph{Proof.} As before, we prove part (i). Part (ii) and (iii) follow similarly. Compared with the previous proof, this one is based on the monotonicity of the functional $t \to J(U,t)$ defined in \eqref{eq:DEFFUNCTIONALJ1} instead of the Almgren-Poon quotient. So, if \eqref{eq:2POINTWISEBOUNDFORLIOUVILLETHEOREMFORWEIGHTEDHE} holds true, then
\[
\begin{aligned}
J(t,U) &= \frac{1}{t} \int_0^t \int_{\RR^{N+1}} |\nabla U|^2(x,y,\tau) \,d\mu^{\tau}(x,y) d\tau \\
&\leq \frac{C}{t} \int_0^t \int_{\RR^{N+1}} (|x|^2 + |y|^2 + \tau)^{\gamma} \,d\mu^{\tau}(x,y) d\tau \\
& \leq \frac{C}{t} \int_0^t \tau^{\gamma} d\tau = Ct^{\gamma} \to 0,
\end{aligned}
\]
as $t \to +\infty$, thanks to the assumption $\gamma < 0$. From Lemma \ref{ALTCAFFARELLIFRIEDMANLEMMAONESPIECES} it thus follows that $J(t,U) = 0$ for any $t > 0$, and so, since the measure $d\mu_t = d\mu_t(x,y)$ is nonnegative, it follows $|\nabla U| = 0$ a.e. in $\RR_+^{N+1}\times(0,\infty)$, i.e. $U$ is constant. $\Box$

\bigskip

%
%
%
%
%
%\color{blue}
\emph{Proof of Theorem \ref{CONVERGENCELINFINITYLOC}.} Let $a \in (-1,1)$ and let $W$ be a weak solution to \eqref{eq:TruncatedExtensionEquationLocalBackward}. Up to modifying the cut-off function $\zeta$ in \eqref{eq:PropCutOff}, in such a way that $\supp(\zeta) \subset \mathbb{B}_2$, and $\zeta = 1$ in $\mathbb{B}_1$, we have that the support of the function $F$ defined in \eqref{eq:DefRightHS} is contained in $(\mathbb{B}_2^+\setminus\mathbb{B}_1^+)\times(0,1)$, and so $W$ satisfies
\begin{equation}\label{eq:TEELBackward}
\begin{cases}
\partial_t W + y^{-a} \nabla \cdot(y^a \nabla W) = 0 \quad &\text{ in } \mathbb{B}_1^+ \times (0,1) \\
-\partial_y^a W = q(x,t)w \quad &\text{ in } B_1 \times \{0\} \times (0,1),
\end{cases}
\end{equation}
in the weak sense, where $w(x,t) := W(x,0,t)$. We have to prove that $W$ is bounded in $\mathcal{C}^{2\nu,\nu}(\overline{\mathbb{Q}_{1/2}^+})$ for any fixed $\nu \in (0,\nu_\ast)$: since we are assuming that $\|W\|_{L^{\infty}(\mathbb{Q}_1^+)}$ is bounded, it is enough to show that
\begin{equation}\label{eq:UNIFORMHOLDERBOUNDBLOWUPSEQUENCE}
[W]_{\mathcal{C}^{2\nu,\nu}(\overline{\mathbb{Q}_{1/2}^+})} := \sup_{(X_1,t_1) \not= (X_2,t_2) \in \overline{\mathbb{Q}_{1/2}^+}} \frac{|W(X_1,t_1) - W(X_2,t_2)|}{(|X_1-X_2|^2 + |t_1-t_2|)^{\nu}} \leq C,
\end{equation}
for some constant $C > 0$ depending on $\nu$. Actually, we will show
\begin{equation}\label{eq:UNIFORMHOLDERBOUNDBLOWUPSEQUENCETEST}
[\eta W]_{\mathcal{C}^{2\nu,\nu}(\QQ_1^+)} \leq C,
\end{equation}
where $\QQ_1^+ := \BB_1^+\times(0,1)$ and $\eta = \eta(X,t)$ is a smooth function satisfying
\[
\begin{cases}
\eta(X,t) = 1 \quad          &\text{ for } (X,t) \in \QQ_{1/2}^+ \\
0 < \eta(X,t) \leq 1 \quad   &\text{ for } (X,t) \in \QQ_1^+ \setminus \QQ_{1/2}^+ \\
\eta(X,t) = 0 \quad          &\text{ for } (X,t) \in \partial \QQ_1^+,
\end{cases}
\]
where $\partial \QQ_1^+ := [\partial \BB_1^+ \times (0,1)] \cup [\BB_1^+\times\{1\}]$. From the definition of $\eta$, \eqref{eq:UNIFORMHOLDERBOUNDBLOWUPSEQUENCE} easily follows from \eqref{eq:UNIFORMHOLDERBOUNDBLOWUPSEQUENCETEST}.

\smallskip

\emph{Step 0: Approximation.} Given a weak solution $W$ to \eqref{eq:TruncatedExtensionEquationLocalBackward}, we consider a family of mollifiers $\{\varrho_j\}_{j \in \mathbb{N}}$ in the  $(x,t) \in \mathbb{R}^{N+1}$ and define
\[
\overline{v}_j(X,t) := (\varrho_j \star W)(X,t) = \int_{\mathbb{R}^{N+1}} \varrho_j(z,\tau)W(x-z,y,t-\tau)dzd\tau.
\]
Thus, $\overline{v}_j \in L^{\infty}(\mathbb{Q}_1^+)$ uniformly w.r.t. $j \in \mathbb{N}$, it is smooth w.r.t. to $x$ and $t$, and satisfies
\[
\begin{cases}
\partial_t \overline{v}_j + y^{-a} \nabla \cdot (y^a \nabla \overline{v}_j) = 0 \quad &\text{ in } \mathbb{B}_1^+\times (0,1) \\
-\partial_y^{a} \overline{v}_j = f_j \quad &\text{ in } B_1 \times \{0\} \times (0,1),
\end{cases}
\]
in the weak sense, where $f_j := \varrho_j \star f$, $f(x,t) := q(x,t) W(x,0,t)$. Further, notice that $\overline{v}_j \to W$ in $L^\infty(\mathbb{Q}_1^+)$ as $j \to +\infty$.

By the properties of $\varrho_j$ and the fact that $\overline{v}_j$ is bounded, the trace $v_j(x,t) := \overline{v}_j(x,0,t)$ is bounded in $\mathbb{R}^N\times\{0\}\times\mathbb{R}$, for each $j \in \mathbb{N}$. Consequently, thanks to \eqref{eq:DEFOFEXTENDEDVERSIONOFU},
\[
\overline{u}_j(X,t) := \int_{-\infty}^{-t} \int_{\mathbb{R}^N} v_j(z,\tau)P_y^a(x-z,-t-\tau) dz d\tau,
\]
is smooth in $\mathbb{R}^{N+1}_+ \times (0,\infty)$ with $\overline{u}_j(x,0,t) = v_j(x,t)$ (in the sense of traces), and satisfies
\[
\partial_t \overline{u}_j + y^{-a} \nabla \cdot (y^a \nabla \overline{u}_j) = 0 \quad \text{ in } \mathbb{R}^{N+1}_+ \times (0,1).
\]
Furthermore, from the explicit expression of the Poisson kernel $P_y^a$ (see \eqref{eq:POISSONKERNELWITHANOTATION}), it follows $\overline{u}_j \in \mathcal{C}_{loc}^{2\nu,\nu}(\overline{\mathbb{R}^{N+1}_+} \times (0,1))$ for all $\nu \in (0,\nu_\ast)$ and all $j \in \mathbb{N}$. Now, the function $\overline{w}_j := \overline{u}_j - \overline{v}_j$, satisfies
\[
\begin{cases}
\partial_t \overline{w}_j + y^{-a} \nabla \cdot (y^a \nabla \overline{w}_j) = 0 \quad &\text{ in } \mathbb{B}_1^+\times (0,1) \\
\overline{w}_j = 0 \quad &\text{ in } \mathbb{R}^N \times \{0\} \times (0,1),
\end{cases}
\]
and so, up to an odd reflection with respect to  $y$, it satisfies
\begin{equation}\label{eq:ODDref}
\partial_t \overline{w}_j + |y|^{-a} \nabla \cdot (|y|^a \nabla \overline{w}_j) = 0 \quad \text{ in } \mathbb{B}_1 \times (0,1).
\end{equation}
Now, we have $\partial_t \overline{w}_j \in L^{\infty}(\mathbb{B}_1\times(0,1))$ and so, we can see it as a r.h.s. for \eqref{eq:ODDref} and apply the elliptic regularity estimates in \cite{SireTerVita2020:art} to deduce that $\overline{w}_j$ is $\nu$-H\"older continuous w.r.t. $y$ for all $\nu \in (0,\nu_\ast)$. Finally, since both $\overline{w}_j$ and $\overline{u}_j$ are smooth w.r.t. to $x$ and $t$, it follows that $\overline{v}_j \in \mathcal{C}^{2\nu,\nu}(\overline{\mathbb{Q}_{1/2}^+})$, for all $j \in \mathbb{N}$, $\nu \in (0,\nu_\ast)$.

We proceed by showing that \eqref{eq:UNIFORMHOLDERBOUNDBLOWUPSEQUENCETEST} holds for each $\overline{v}_j$, for some constant $C > 0$ independent of $j \in \mathbb{N}$. Once \eqref{eq:UNIFORMHOLDERBOUNDBLOWUPSEQUENCETEST} is established for $\overline{v}_j$, we deduce it for $W$, since $\overline{v}_j \to W$ in $L^\infty(\mathbb{Q}_1^+)$ and so, up to passing to a subsequence, $\overline{v}_j \to W$ in $\mathcal{C}^{2\nu,\nu}(\overline{\mathbb{Q}_{1/2}^+})$, for each fixed $\nu \in (0,\nu_\ast)$.
\normalcolor

\smallskip

\emph{Step 1: Main definitions.} Assume by contradiction that \eqref{eq:UNIFORMHOLDERBOUNDBLOWUPSEQUENCETEST} does not hold, i.e. there exists $\nu \in (0,\nu_\ast)$, a sequence $\{\overline{v}_j\}_{j\in\mathbb{N}}$ and two sequences $(X_{1,j},t_{1,j}),(X_{2,j},t_{2,j}) \in \QQ_1^+$ such that
\begin{equation}\label{eq:ACHIVEMENTPOINTSBLOWUP}
[\eta \overline{v}_j]_{\mathcal{C}^{2\nu,\nu}(\overline{\mathbb{Q}_1^+})} := \frac{|\eta(X_{1,j},t_{1,j}) \overline{v}_j(X_{1,j},t_{1,j}) - \eta(X_{2,j},t_{2,j}) \overline{v}_j(X_{2,j},t_{2,j})|}{r_j^{2\nu}} := L_j \to +\infty,
\end{equation}
as $j \to +\infty$, where we have defined
\[
r_j := (|X_{1,j}-X_{2,j}|^2 + |t_{1,j}-t_{2,j}|)^{\frac{1}{2}},
\]
for all $j \in \NN$. From \eqref{eq:ACHIVEMENTPOINTSBLOWUP}, we immediately deduce the bound
\[
L_j \leq \frac{\|\overline{v}_j\|_{L^{\infty}(\QQ_1^+)}}{r_j^{2\nu}} [|\eta(X_{1,j},t_{1,j})| + |\eta(X_{2,j},t_{2,j})|],
\]
and so, since we know that $\{\overline{v}_j\}_{j\in\mathbb{N}}$ is uniformly bounded in $\QQ_1^+$, we obtain that $r_j \to 0$ as $j \to +\infty$. Furthermore, from the same bound on $L_j$ and the smoothness of $\eta = \eta(X,t)$, it is easily seen that
\[
\frac{\text{dist}((X_{1,j},t_{1,j}),\partial \QQ_1^+)}{r_j} + \frac{\text{dist}((X_{2,j},t_{2,j}),\partial \QQ_1^+)}{r_j} \geq \frac{L_j r_j^{2\nu - 1}}{L\|\overline{v}_j\|_{L^{\infty}(\QQ_1^+)}},
\]
where $L > 0$ is taken such that $|\eta(X_1,t_1) - \eta(X_2,t_2)| \leq L (|X_1-X_2|^2 + |t_1-t_2|)^{\frac{1}{2}}$. Consequently, since $\nu \in (0,\nu_\ast)$ and $\nu_\ast \leq 1/2$, we obtain
\begin{equation}\label{eq:SUMOFDISTANCESIMPORTANTSEQUENCE}
\frac{\text{dist}((X_{1,j},t_{1,j}),\partial \QQ_1^+)}{r_j} + \frac{\text{dist}((X_{2,j},t_{2,j}),\partial \QQ_1^+)}{r_j} \to + \infty,
\end{equation}
as $j \to + \infty$.

\emph{Step2: Auxiliary sequences.} Following the ideas of \cite[Section 6]{TerVerZil2017:art} and \cite[Section 4]{TerVerZil2016:art}, the remaining part of the proof is based on the analysis of two different sequences:
\begin{equation}\label{eq:DEFINITIONAUXILIARYSEQUENCES}
\begin{aligned}
W_j(X,t) &:= \eta(\widehat{X}_j,\widehat{t}_j)\frac{\overline{v}_j(\widehat{X}_j + r_jX,\widehat{t}_j + r_j^2t)}{L_jr_j^{2\nu}}, \\
\overline{W}_j(X,t) &:= \frac{(\eta \overline{v}_j)(\widehat{X}_j + r_jX,\widehat{t}_j + r_j^2t)}{L_jr_j^{2\nu}},
\end{aligned}
\end{equation}
where $\widehat{P}_j = (\widehat{X}_j,\widehat{t}_j) \in \QQ_1^+$ will be chosen later and
\[
(X,t) \in \widehat{\mathbb{Q}}_j := (\QQ_1^+ - \widehat{P}_j)/r_j = \BB_{1/r_j}(\widehat{X}_j)\times \{y > -\widehat{y}_j/r_j\} \times I_j,
\]
for all $j \in \NN$, where we have set $I_j := (-\widehat{t}_j/r_j^2,(1-\widehat{t}_j)/r_j^2)$. The definitions of the new sequences in \eqref{eq:DEFINITIONAUXILIARYSEQUENCES} are motivated by the following two facts.

The first one, is that the H\"{o}lder semi-norm of order $\nu$ of $\overline{W}_j = \overline{W}_j(X,t)$ satisfies
\[
\left[ \overline{W}_j \right]_{\mathcal{C}^{2\nu,\nu}(\widehat{\QQ}_j)} = 1,
\]
for all $j \in \NN$. This easily follows by the definition of $\overline{W}_j = \overline{W}_j(X,t)$, $L_j$ and $r_j$. On the other hand, the first sequence satisfies the problem (see  \eqref{eq:TruncatedExtensionEquationBlowUpLocal})
\begin{equation}\label{eq:FIRSTEQUATIONOFWJ}
\begin{cases}
\partial_t W_j + \mathcal{L}_a^jW_j = 0 \quad &\text{ in } \mathbb{B}_{1/r_j}(\widehat{X}_j) \cap \{y > -\widehat{y}_j/r_j\} \times I_j \\
-\partial_y^{a,j} W_j = r_j^{1-a} \widehat{q}_j(x,t) w_j \quad &\text{ in } B_{1/r_j} \times \{-\widehat{y}_j/r_j\} \times I_j,
\end{cases}
\end{equation}
in the weak sense, where $w_j(x,t) =W_j(x,0,t)$, $\widehat{q}_j(x,t) = q(\widehat{x}_j + r_jx,-\widehat{t}_j - r_j^2t)$,
\[
\mathcal{L}_a^jW = \left(\widehat{y}_jr_j^{-1} + y\right)^{-a} \nabla\cdot \left( \left(\widehat{y}_jr_j^{-1} + y\right)^a \nabla W_j \right)
\]
and
\[
-\partial_y^{a,j} W_j = \lim_{y \to -(\widehat{y}_j/r_j)^+} \left(\frac{\widehat{y}_j}{r_j} + y\right)^a \partial_y W_j.
\]
for all $j \in \NN$. Further, notice that $W_j(O) = \overline{W}_j(O)$ by definition. In other words, the sequence $\overline{W}_j$ ``preserves'' the $\nu$-H\"older norm, while the sequence ``preserves'' the equation.

A second important feature of the two sequences is that they are asymptotically equivalent on compact sets of $\RR_+^{N+1}\times\RR$. Indeed, if $K \subset \RR_+^{N+1}\times\RR$ is compact and $(X,t) \in K$, we have
\[
\begin{aligned}
|W_j(X,t) - \overline{W}_j(X,t)| & \leq \frac{\|\eta \overline{v}_j \|_{L^{\infty}(\QQ_1^+)}}{r_j^{2\nu}L_j} |\eta(\widehat{X}_j + r_jX,\widehat{t}_j + r_j^2t) - \eta(\widehat{X}_j,\widehat{t}_j)| \\
& \leq L \frac{\|\eta  \overline{v}_j \|_{L^{\infty}(\QQ_1^+)}}{r_j^{2\nu-1}L_j}(|X|^2 + |t|)^{\frac{1}{2}} \to 0,
\end{aligned}
\]
as $j \to +\infty$ (recall that $r_j \to 0$, $L_j \to +\infty$, and $2\nu < 1$). In particular, it follows
\begin{equation}\label{eq:CONVERGENCEWJOVWJFORLARGEJ}
\|W_j - \overline{W}_j\|_{L^{\infty}(K)} \to 0,
\end{equation}
as $j \to +\infty$. Furthermore, from the bound above and recalling that $\overline{W}_j$ is $\nu$-H\"older continuous uniformly w.r.t. $j \in \mathbb{N}$ with $\overline{W}_j(O) = W_j(O)$, it is easily seen that
\[
\begin{aligned}
|W_j(X,t) - W_j(O)| &\leq |W_j(X,t) - \overline{W}_j(X,t)| + |\overline{W}_j(X,t) - \overline{W}_j(O)| \\
& \leq C \left(\frac{1}{r_j^{2\nu-1}L_j}(|X|^2 + |t|)^{\frac{1}{2}} + (|X|^2 + |t|)^{\nu} \right),
\end{aligned}
\]
from which we deduce the existence of a constant $C > 0$ (depending on the compact set $K \subset \RR_+^{N+1}\times\RR$) such that
\begin{equation}\label{eq:UNIFORMLINFINITYBOUNDONZJ}
\sup_{(X,t) \in K } |W_j(X,t) - W_j(O)| \leq C.
\end{equation}
These last two properties will be crucial in the next step.

\smallskip

\emph{Step3: Asymptotic behaviour of $(X_{1,j},t_{1,j})$ and $(X_{2,j},t_{2,j})$.} We now show that the sequences $(X_{1,j},t_{1,j})$ and $(X_{2,j},t_{2,j})$ approach the set
\[
\Sigma = \{(x,y,t) \in \RR^{N+1}\times\mathbb{R}: y = 0 \}
\]
as $j \to +\infty$. More precisely, we prove the existence of a constant $C > 0$ such that
\begin{equation}\label{eq:CONVERGENCEOFPOINTSCHARACTERMANIFOLD}
\frac{\text{dist}((X_{1,j},t_{1,j}), \QQ_1^+ \cap \Sigma)}{r_j} + \frac{\text{dist}((X_{2,j},t_{2,j}), \QQ_1^+ \cap \Sigma)}{r_j} \leq C,
\end{equation}
for $j \in \NN$ large enough. Arguing by contradiction, we assume that
\begin{equation}\label{eq:ABSURDUMASSUMPTIONS3}
\frac{\text{dist}((X_{1,j},t_{1,j}), \QQ_1^+ \cap \Sigma)}{r_j} + \frac{\text{dist}((X_{2,j},t_{2,j}), \QQ_1^+ \cap \Sigma)}{r_j} \to +\infty,
\end{equation}
as $j \to +\infty$. Let us take $(\widehat{X}_j,\widehat{t}_j) = (X_{1,j},t_{1,j})$ in the definition of $W_j$ and $\overline{W}_j$. Thus, thanks to \eqref{eq:ABSURDUMASSUMPTIONS3} we obtain that
\[
\bigcup_{j\in\mathbb{N}} \mathbb{B}_{1/r_j}(\widehat{X}_j) \cap \{y > -\widehat{y}_j/r_j\} \times I_j
\]
equals one of the following: $\RR^{N+1}\times\RR$ or $\RR^{N+1}\times(0,+\infty)$. From now on, we assume to be in the case $\RR^{N+1}\times\RR$ (the other can be treated in a very similar way). Now, consider the new sequences
\begin{equation}\label{eq:DEFINITIONOFZJ}
\begin{aligned}
U_j(X,t) &:= W_j(X,t) - W_j(O), \\
\overline{U}_j(X,t) &:= \overline{W}_j(X,t) - \overline{W}_j(O),
\end{aligned}
\end{equation}
and let $j \in \NN$ be large enough such that $K \subset \widehat{\QQ}_j$, where $K \subset \RR^{N+1}\times\RR$ is a fixed compact set. Since the sequence $\{\overline{U}_j\}_{j\in\NN}$ is uniformly bounded in $K$ (see  \eqref{eq:CONVERGENCEWJOVWJFORLARGEJ} and \eqref{eq:UNIFORMLINFINITYBOUNDONZJ}) with uniformly bounded $\nu$-H\"{o}lder semi-norm, we can apply the Ascoli-Arzel\`a theorem to deduce the existence of a continuous function $U \in \mathcal{C}(K)$, uniform limit of $\overline{U}_j$ as $j \to +\infty$. Notice that, using a standard diagonal procedure, we obtain that $\overline{U}_j \to U$ uniformly on compact sets of $\RR_+^{N+1}\times\RR$ and, by \eqref{eq:CONVERGENCEWJOVWJFORLARGEJ} and \eqref{eq:UNIFORMLINFINITYBOUNDONZJ}, we obtain $U_j \to U$ uniformly on compact sets of $\RR_+^{N+1}\times\RR$, too. Furthermore, from the definition of $\overline{U}_j$, it follows
\[
|\overline{U}_j(X_1,t_1) - \overline{U}_j(X_2,t_2)| = |\overline{W}_j(X_1,t_1) - \overline{W}_j(X_2,t_2)| \leq (|X_1-X_2|^2 + |t_1-t_2|)^{\nu},
\]
for any choice $(X_1,t_1),(X_2,t_2) \in \RR_+^{N+1}\times\RR$ and $j \in \NN$ large enough. Consequently, taking $j \to +\infty$ in the above inequality and using the arbitrariness of $(X_1,t_1)$ and $(X_2,t_2)$, we obtain $U \in \mathcal{C}^{2\nu,\nu}(\RR_+^{N+1}\times\RR)$. In particular, it satisfies
\begin{equation}\label{eq:HOLDERGLOBALBOUNDONZ}
|U(X,t)| \leq C(1 + (|X|^2 + |t|)^{\nu}] \leq C(1+ d^{2\nu}(X,t))
\end{equation}
for some $C > 0$ and all $(X,t) \in \RR_+^{N+1}\times\RR$.

Now, using \eqref{eq:FIRSTEQUATIONOFWJ}, we find that $U_j$ is a weak solution to
\begin{equation}\label{eq:FIRSTEQUATIONOFUJ}
\begin{cases}
\partial_t U_j + \mathcal{L}_a^jU_j = 0 \quad &\text{ in } \mathbb{B}_{1/r_j}(\widehat{X}_j) \cap \{y > -\widehat{y}_j/r_j\} \times I_j \\
-\partial_y^{a,j} U_j = r_j^{1-a} \widehat{q}_j(x,t) (u_j + W_j(O)) \quad &\text{ in } B_{1/r_j} \times \{-\widehat{y}_j/r_j\} \times I_j,
\end{cases}
\end{equation}
for each $j \in \mathbb{N}$, where $u_j(x,t) = U_j(x,0,t)$, in the sense of Definition \ref{def:WeakSolutions}. So, let us fix $t_1 < t_2$. Testing the equation with $\eta \in C_c^{\infty}(\RR^{N+1}\times(t_1,t_2))$, taking $j$ large enough so that
\[
\text{supp}(\eta) \subset \bigcup_{j\in\mathbb{N}} \mathbb{B}_{1/r_j}(\widehat{X}_j) \cap \{y > -\widehat{y}_j/r_j\} \times (t_1,t_2)
\]
and integrating by parts in space-time, we deduce that
\[
\int_{t_1}^{t_2} \int_{\{y > -\widehat{y}_j/r_j\}} \left[ (1 + r_j\widehat{y}_j^{-1}y)^a\partial_t\eta - \nabla\cdot\left((1 + r_j\widehat{y}_j^{-1}y)^a\nabla \eta \right)\right] U_j \,dXdt = 0,
\]
for all $j \in \NN$. Hence, exploiting the fact that $(1 + r_j\widehat{y}_j^{-1}y)^a \to 1$ on any compact set of $\RR^{N+1}\times\RR$ ($r_j\widehat{y}_j^{-1} \to 0$ thanks to \eqref{eq:ABSURDUMASSUMPTIONS3}) and using the uniform convergence of the sequence $U_j$, we pass to the limit as $j \to +\infty$ in the above relation to conclude
\[
\int_{\RR} \int_{\RR^{N+1}}\left(\partial_t\eta - \Delta_{x,y} \eta \right) U \,dXdt = 0,
\]
for any $\eta \in \mathcal{C}_0^{\infty}(\RR^{N+1}\times\RR)$. In particular, this implies that $U$ is a backward caloric function in $\RR^{N+1}\times\RR$, in the very weak sense. Consequently, from the classical regularity theory of the Heat Equation, the bound in \eqref{eq:HOLDERGLOBALBOUNDONZ}, and the Liouville type theorem proved in Lemma \ref{1LIOUVILLETYPETHEOREMFORWEIGHTEDHE} (part (i) with $a = 0$, see also the classical paper \cite{Hirschman1952:note}), we immediately  deduce that $U$ must be constant.

At this point, since we have chosen $(\widehat{X}_j,\widehat{t}_j) = (X_{1,j},t_{1,j})$, it follows
\[
\frac{(X_{2,j},t_{2,j}) - (\widehat{X}_j,\widehat{t}_j)}{r_j} = \frac{(X_{2,j},t_{2,j}) - (X_{1,j},t_{1,j})}{(|X_{1,j}-X_{2,j}|^2 + |t_{1,j}-t_{2,j}|)^{\frac{1}{2}}} \to (X_2,t_2),
\]
for some $(X_2,t_2) \in \QQ_1$, up to subsequences, and so, using the uniform convergence and the uniform H\"{o}lder bound on $\overline{W}_j = \overline{W}_j(X,t)$, we obtain
\begin{equation}\label{eq:WISNONCONSTANT}
\begin{aligned}
1 &= \left|\overline{W}_j\left(\frac{X_{2,j} - \widehat{X}_j}{r_j},\frac{t_{2,j} - \widehat{t}_j}{r_j^2}\right) - \overline{W}_j\left(\frac{X_{1,j} - \widehat{X}_j}{r_j},\frac{t_{1,j} - \widehat{t}_j}{r_j^2}\right) \right| \\
& = \left|\overline{U}_j\left(\frac{X_{2,j} - \widehat{X}_j}{r_j},\frac{t_{2,j} - \widehat{t}_j}{r_j^2}\right) - \overline{U}_j\left(\frac{X_{1,j} - \widehat{X}_j}{r_j},\frac{t_{1,j} - \widehat{t}_j}{r_j^2}\right) \right| \\
& = \left|\overline{U}_j\left(\frac{X_{2,j} - \widehat{X}_j}{r_j},\frac{t_{2,j} - \widehat{t}_j}{r_j^2}\right) - \overline{U}_j(O) \right| \to |U(X_2,t_2) - U(O)| = 1,
\end{aligned}
\end{equation}
as $j \to +\infty$, i.e. $U$ is nonconstant, in contradiction with the Liouville theorem. This conclude the proof of \eqref{eq:CONVERGENCEOFPOINTSCHARACTERMANIFOLD}.

\emph{Step4: Final part of the proof of \eqref{eq:UNIFORMHOLDERBOUNDBLOWUPSEQUENCE}.} In view of \eqref{eq:CONVERGENCEOFPOINTSCHARACTERMANIFOLD} and the arbitrariness of $\widehat{P}_j = (\widehat{X}_j,\widehat{t}_j) \in \QQ_1$ in the definition of $W_j$ and $\overline{W}_j$ (see  \eqref{eq:DEFINITIONAUXILIARYSEQUENCES}), we can choose $\widehat{P}_j = (x_{1,j},0,t_{1,j})$ (see  \cite{TerVerZil2017:art,TerVerZil2016:art}). Consequently, $U_j$ satisfies
\[
\begin{aligned}
\int_{t_1}^{t_2} &\int_{\RR_+^{N+1}} \left[y^a\partial_t\eta - \nabla\cdot\left( y^a\nabla \eta \right)\right] U_j \,dXdt \\
& +\int_{t_1}^{t_2} \int_{\RR^N\times\{0\}} U_j \partial_y^a \eta \,dxdt =   r_j^{1-a} \int_{t_1}^{t_2} \int_{\RR^N\times\{0\}} \widehat{q}_j [u_j + W_j(O)]\eta \,dXdt
\end{aligned}
\]
for all test functions $\eta$ with compact support in $\overline{\RR_+^{N+1}} \times (t_1,t_2)$ and all $j \in \NN$. Now, noticing that $W_j(O)$ is uniformly bounded (it is enough to use the equation of $W_j$ and the $L^2 \to L^{\infty}$ bounds proved in \cite[Formula (5.7)]{BanGarofalo2017:art}), we can pass to the limit as $j\to+\infty$ to obtain
\[
\int_{t_1}^{t_2} \int_{\RR^{N+1}_+} \left(y^a\partial_t\eta - \nabla\cdot\left(y^a\nabla \eta \right)\right] U \,dXdt + \int_{t_1}^{t_2} \int_{\RR^N\times\{0\}} U \partial_y^a\eta \,dxdt = 0
\]
for all test functions $\eta$ with compact support in $\overline{\RR^{N+1}_+} \times (t_1,t_2)$. Integrating by parts with respect to  space, we deduce that $U$ is a global weak solution to \eqref{eq:DEFINITIONPARABOLICEQUATION} and so, since $U$ satisfies \eqref{eq:HOLDERGLOBALBOUNDONZ} and $\nu \in (0,1/2)$, it follows that $U$ is constant thanks to Lemma \ref{1LIOUVILLETYPETHEOREMFORWEIGHTEDHE} part (i). Since this contradicts \eqref{eq:WISNONCONSTANT}, we complete the proof of the theorem. $\Box$

%
%
%
%
%
%
%
%
%%%%%%%%%%%%%%%%%%%%%%%%%%%%%%%%%%%%%%%%%%%%%%%%%%%%%%%%%%%%%%%%%%%%%%%%%%%%%%%%%%%%%%%%%%%%%%%%%%%%%%%%%%%
%
%
%%%%%%%%%%%%%%%%%%%%%%%%%%%%%%%%%%%%%%%%%%%%%%%%%%%%%%%%%%%%%%%%%%%%%%%%%%%%%%%%%%%%%%%%%%%%%%%%%%%%%%%%%%%
%
\section{Nodal properties of solutions}\label{SECTIONNODALSET}
In this final part of the work, we employ the theory developed in the previous sections to obtain information about the Hausdorff dimension, regularity and structure of the nodal set of solutions to equation \eqref{eq:NONLOCALHEPOTENTIAL}. In other words, we will prove Theorem \ref{Theorem:StructureRegularity}. As we have mentioned in the introduction, the decisive tools turn out to be the blow-up classification of Section \ref{Section:BlowUp1} and Section \ref{Section:BlowUp2} and classical theorems as Federer's Reduction Principle (see  Theorem \ref{FEDERERTHEOREM})
and Whitney's Extension Theorem (see  Theorem \ref{PARABOLICWHITNEYEXTENSION}),
together with some independent technical results.  To do so, we must introduce some additional notations and definitions.

\subsection{Additional notations and background} First, we introduce the sets
\[
\begin{aligned}
\widetilde{\mathcal{K}} &:= \{\widetilde{\kappa}_{n,m}\}_{n,m \in \NN} = \left\{\frac{n}{2} + m \right\}_{n,m \in \NN}, \\
\Sigma &:= \{(x,y,t)\in \RR^{N+1}\times\RR: y = 0\}.
\end{aligned}
\]
\begin{defn}\label{DEFREGULARSINGULARSET}
	Let $W$ be a weak solution to \eqref{eq:TruncatedExtensionEquationLocalBackward}, with $q$ satisfying \eqref{eq:ASSUMPTIONSPOTENTIAL}. For any positive $\kappa \in \widetilde{\mathcal{K}}$, we define
	\[
	\begin{aligned}
	\Sigma(W) &:= \Gamma(W) \cap \overline{\QQ_{1/2}^+} \cap \Sigma, \\
	\Gamma_{2\kappa}(W) &:= \left\{p_0 \in \Sigma(W): \lim_{r \to 0^+} N_D(W_{p_0},r) = \kappa \right\}.
	\end{aligned}
	\]
	Furthermore, we set
	\[
	\begin{aligned}
	{\rm Reg}(W) &:= \Gamma_{1}(W), \\
	{\rm Sing}(W) &:= \Gamma(W) \setminus {\rm Reg}(W).
	\end{aligned}
	\]
	Finally, if $u$ is a solution to \eqref{eq:NONLOCALHEPOTENTIAL}, we define $\Gamma(u) := \Gamma(W)$, $\Gamma_{2\kappa}(u) := \Gamma_{2\kappa}(W)$ and
	\begin{equation}\label{eq:DEFINITIONREGULARSINGULARNONLOCALEQ}
	\begin{aligned}
	{\rm Reg}(u) &:= \Gamma_{1}(u), \\
	{\rm Sing}(u) &:= \Gamma(u) \setminus {\rm Reg}(u).
	\end{aligned}
	\end{equation}
	The sets ${\rm Reg}(u)$ and ${\rm Sing}(u)$ are called the regular and the singular part of $\Gamma(u)$, respectively.
\end{defn}
As a preliminary result, we prove that for any $\Gamma_{2\kappa}(W)$ is $F_{\sigma}$ for any singular frequency $\kappa$. This fact will turn out to be important in the proof of Theorem \ref{STRUCTUREOFSINGULARSET}.
\begin{lem}\label{LEMMAGAMMAKISFSIGMA}
	Let $W$ be a weak solution to \eqref{eq:TruncatedExtensionEquationLocalBackward}, with $q$ satisfying \eqref{eq:ASSUMPTIONSPOTENTIAL}. Then the following two assertions hold:
	
	(i) The set ${\rm Reg}(W)$ is relatively open in $\Sigma(W)$.
	
	(ii) For any $\kappa = 1,\frac{3}{2},2,\ldots$, the set $\Gamma_{2\kappa}(W)$ is a union of countably many closed sets.
	
\end{lem}
\emph{Proof.} Part (i) directly follows from the upper semi-continuity of the map
\[
p_0 \to N_D(W,p_0,0^+) := \lim_{r \to 0} N_D(W,p_0,r) = \lim_{r \to 0} \Phi_a(W,p_0,r),
\]
and the gap condition given by Theorem \ref{THEOREMBLOWUP1New}, which yields
\[
{\rm Reg}(W) = \left\{p_0 \in \Sigma(W): N_D(W,p_0,0^+) = \frac{1}{2}\right\} = \left\{p_0 \in \Sigma(W): N_D(W,p_0,0^+) < 1 \right\}.
\]
Turning to part (ii) we notice that, thanks to formula \eqref{eq:LIMITOFTHEQUOTIENTFORSMALLTSTATEMENT1New}, we have
\begin{equation}\label{eq:SETEJNOPENREGULARCLOSEDSINGULAR}
\Gamma_{2\kappa}(W) = \bigcup_{j=1}^{\infty} E_j  \qquad E_j := \left\{p_0 \in \Gamma_{2\kappa}(W): \frac{1}{j} r^{4\kappa} \leq H(W,p_0,r) \leq j r^{4\kappa}, r \in (0,r_j) \right\},
\end{equation}
for some sequence $r_j \in (0,1)$. Consequently, it is enough to show that all the $E_j$'s are closed sets, i.e. $\overline{E_j} = E_j$ for all $j$. If $p_0 \in \overline{E_j}$, then it is instantly seen that the inequalities in \eqref{eq:SETEJNOPENREGULARCLOSEDSINGULAR} hold by continuity of the function $p_0 \to H(W,p_0,r)$. It is left to show that $p_0 \in \Gamma_{2\kappa}(U)$. Since the function $p_0 \to N(W,p_0,0^+)$ is upper semi-continuous, we infer $N(p_0,t_0^+,U) \geq \kappa$. Finally, if by contradiction $N(W,p_0,0^+) = \overline{\kappa} > \kappa$, combining formula \eqref{eq:LIMITOFTHEQUOTIENTFORSMALLTSTATEMENT1New} with $\kappa = \overline{\kappa}$ and \eqref{eq:SETEJNOPENREGULARCLOSEDSINGULAR}, we easily obtain
\[
\frac{1}{j}r^{4\kappa} \leq H(W,p_0,r) < j r^{4\overline{\kappa}}, \text{ as } r \to 0,
\]
for some $j$, which is contradiction. Consequently, $p_0 \in \Gamma_{2\kappa}(U)$ and so $p_0 \in E_j$. $\Box$

\begin{defn}\label{DEFBLOWUPSET}
	Let $W$ be a weak solution to \eqref{eq:TruncatedExtensionEquationLocalBackward}, with $q$ satisfying \eqref{eq:ASSUMPTIONSPOTENTIAL}, $\kappa \in \widetilde{\mathcal{K}}$ and $p_0 \in \Gamma_{2\kappa}(W)$. We denote with $\mathfrak{B}_{2\kappa}(W)$ the set of all admissible blow-up limits of $W$ at $p_0$.
	Furthermore, the tangent map of $W$ at $p_0$ is the unique $\Theta_{p_0} \in \mathfrak{B}_{2\kappa}(W)$ such that
	\[
	W_{p_0,r} \to \Theta_{p_0} \quad \text{ locally uniformly},
	\]
	as $r \to 0^+$.
\end{defn}
We are now ready to turn to estimate the Hausdorff dimension of the nodal set of solutions to \eqref{eq:NONLOCALHEPOTENTIAL}. The analysis will be performed by applying the Federer's reduction principle (see  for instance with Simon \cite[Appendix A]{Simon1983:book}, or Lin \cite[Section 2]{Lin1991:art}). The parabolic version of it is less employed in literature and so, out of completeness, we review the main concepts in the next paragraphs. We follow the work of Chen \cite[Section 8]{Chen1998:art}, adapting the notation to our framework. We begin with the definition of Hausdorff and \emph{parabolic} Hausdorff dimension.
\begin{defn}\label{DEFINPARHAUSDIMENSION}(Parabolic Hausdorff dimension, \cite[Definition 8.1]{Chen1998:art})
\\
For any $E \subset \RR^N\times\RR$, any real number $d \geq 0$, and $0 < \delta <+ \infty$, we define
\[
\mathcal{P}_{\delta}^d(E) := \inf \left\{\sum_{j=1}^{\infty} r_j^d \; : \; E \subset \bigcup_{j=1}^{\infty} Q_{r_j}(x_j,t_j) \; \text{ with } \; 0< r_j < \delta \right\},
\]
where
\[
Q_r(x,t) := \left\{(z,\tau) \in \RR^N\times\RR: |x-z| < r, |t-\tau| < r^2 \right\}.
\]
Then we define the $d$-dimensional cylindrical Hausdorff measure by
\[
\mathcal{P}^d(E) := \lim_{\delta \to 0^+} \mathcal{P}_{\delta}^d(E) = \sup_{\delta > 0} \mathcal{P}_{\delta}^d(E).
\]
We call parabolic Hausdorff dimension of $E$, the number
\[
\dim_{\mathcal{P}}(E) := \inf \left\{d \geq 0 \; : \; \mathcal{P}^d(E) = 0  \right\} = \sup \left\{d \geq 0 \; : \; \mathcal{P}^d(E) = +\infty  \right\}.
\]
\end{defn}
Before moving forward, we recall that the above definition is just a parabolic version of the spherical Hausdorff measure (Hausdorff dimension, resp.). Indeed, if $E \subset \RR^N$, $d \geq 0$, and $0 < \delta \leq \infty$, we define
\[
\mathcal{H}_{\delta}^d(E) := \inf \left\{\sum_{j=1}^{\infty} r_j^d \; : \; E \subset \bigcup_{j=1}^{\infty} B_{r_j}(x_j) \; \text{ with } \; 0 < r_j < \delta \right\},
\]
where now $B_r(x)$ is the ball of radius $r>0$, centered at $x$. Consequently, the $d$-dimensional spherical Hausdorff measure can be defined as
\[
\mathcal{H}^d(E) := \lim_{\delta \to 0^+} \mathcal{H}_{\delta}^d(E) = \sup_{\delta > 0} \mathcal{H}_{\delta}^d(E).
\]
The Hausdorff dimension of $E$ is the number
\[
\dim_{\mathcal{H}}(E) := \inf \left\{d \geq 0 \; : \; \mathcal{H}^d(E) = 0  \right\} = \sup \left\{d \geq 0 \; : \; \mathcal{H}^d(E) = +\infty  \right\}.
\]
\begin{lem} (\cite[Lemma 8.2]{Chen1998:art}) The following two assertions hold:

\smallskip

(i) For any linear subspace $E \subset \RR^N$ and any $-\infty \leq a < b \leq +\infty$, it holds
\[
\dim_{\mathcal{P}}(E\times(a,b)) = \dim_{\mathcal{H}}(E) + 2.
\]
In particular, $\dim_{\mathcal{P}}(\RR^N\times\RR) = N + 2$.

\smallskip

(ii) For any set $E \subset \RR^N\times\RR$, any point $p_0 = (x_0,t_0) \in \RR^N\times\RR$, and $\lambda >0$, define
\[
E_{p_0,\lambda} = \frac{E - p_0}{\lambda} := \left\{ (x,t) \in \RR^N\times\RR : (x_0 + \lambda x,t_0 + \lambda^2t) \in E \right\}.
\]
Then
\[
\mathcal{P}_{\delta}^d(E_{p_0,\lambda}) = \lambda^{-d} \mathcal{P}_{\delta}^d(E),
\]
for any $d \geq 0$ and $\delta > 0$.
\end{lem}
The above lemma clarifies the relation between the Hausdorff dimension and the little bit more ambiguous parabolic Hausdorff dimension. In particular, the fact that $\dim_{\mathcal{P}}(\RR^N\times\RR) = N + 2$ gives relevance to Theorem \ref{Theorem:StructureRegularity}.
\begin{defn}\label{DEFINITIONSELFSIMILARFAMILY}(Locally asymptotically self-similar family,  see  \cite[Definition 8.3]{Chen1998:art})
\\
Let $\mathcal{F} \subset L_{loc}^{\infty}(\RR^N\times\RR)$ be a family of functions and consider a map
\[
\mathcal{S}:\mathcal{F} \to \mathcal{C} := \{C \subset \RR^N\times\RR : C \text{ is closed}\}.
\]
Moreover, define the blow-up family
\[
u_{p_0,r,\varrho}(x,t) = \frac{u(x_0 + r x,t_0+r^2t)}{\varrho},
\]
for any $p_0 = (x_0,t_0) \in \RR^N\times\RR$ and $r,\varrho > 0$. We say that the pair $(\mathcal{F},\mathcal{S})$ is a locally asymptotically self-similar family if it satisfies the properties (A1), (A2) and (A3) below.

\smallskip

(A1) (Closure under re-scaling, translation and normalization) For any $Q_{r}(x_0,t_0) \subset Q_1(0,0)$, $\varrho > 0$ and $u \in \mathcal{F}$, it holds
\[
u_{p_0,r,\varrho} \in \mathcal{F}.
\]

(A2) (Convergence of the normalized blow-up sequence) For any $p_0 \in Q_1(0,0)$, $u \in \mathcal{F}$ and $r_j \to 0^+$, there exist a number $\kappa \in \RR$ and a function $\vartheta_{p_0} \in \mathcal{F}$ parabolically $2\kappa$-homogeneous such that as $j \to +\infty$, up to extracting  a subsequence of $r_j$, it holds
\[
u_{p_0,r_j,r_j^{2\kappa}} \to \vartheta_{p_0} \quad \text{ locally uniformly in } \, \RR^N\times\RR.
\]
Moreover, if $u_{p_0,r_j,r_j^{2\kappa}} \to \vartheta_{p_0}$ and $u_{p_0,r_j,r_j^{2k}} \to \theta_{p_0}$, then $\kappa = k$ and $\vartheta = \theta$.

\smallskip

(A3) (Singular Set assumptions) The map $\mathcal{S}: \mathcal{F} \to \mathcal{C}$ satisfies the following properties:

\noindent (i) For any $Q_{r}(x_0,t_0) \subset Q_1(0,0)$, $\varrho > 0$ and $u \in \mathcal{F}$, it holds
\[
\mathcal{S}(u_{p_0,r,\varrho}) =  \left(\mathcal{S}(u)\right)_{p_0,r}.
\]

\noindent (ii) For any $p_0 \in Q_1(0,0)$, $u,\vartheta_{p_0} \in \mathcal{F}$, $\kappa \in \RR$ and $r_j \to 0^+$ such that $u_{p_0,r_j,r_j^{2\kappa}} \to \vartheta_{p_0}$  uniformly on compact sets of $\RR^N\times\RR$, the following continuity property holds: for any $\varepsilon > 0$, there exists $j_{\varepsilon} > 0$ such that
\[
\mathcal{S}(u_{p_0,r_j,r_j^{2\kappa}}) \cap Q_1(0,0) \subseteq \{p \in \RR^N\times\RR: \dist(p,\mathcal{S}(\vartheta_{p_0})) < \varepsilon \}, \quad \text{ for all } j \geq j_{\varepsilon}.
\]
\noindent (iii) If $u \in \mathcal{F}$ and $\kappa \in \RR$ are such that $u_{p_0,r,r^{2\kappa}} = u$ for all $(x_0,t_0) \in \RR^N\times\RR$ and $r > 0$, then $\mathcal{S}(u) = \emptyset$.
\end{defn}
We are now ready to state the Federer's reduction principle.
\begin{thm}\label{FEDERERTHEOREM} (Federer's reduction principle, Chen \cite[Theorem 8.5]{Chen1998:art})
\\
Let $(\mathcal{F},\mathcal{S})$ be a locally asymptotically self-similar family and assume there exists at least one $u \in \mathcal{F}$ such that $\mathcal{S}(u)\cap Q_1(0,0) \not= \emptyset$. Then:

\smallskip

(i) There exists an integer $0 \leq d \leq N+1$ such that for all $u \in \mathcal{F}$ it holds
\[
\dim_{\mathcal{P}}\left[\mathcal{S}(u) \cap Q_1(0,0)\right] \leq d.
\]

(ii) There exist $u \in \mathcal{F}$, $\kappa \in \RR$ and a linear subspace of $\RR^N\times\RR$ such that
\[
E_{(0,0),\lambda} = E, \quad \text{ for any } \lambda >0, \qquad \mathcal{S}(u) = E, \qquad \dim_{\mathcal{P}} E = d,
\]
and
\[
u_{p_0,\lambda,\lambda^{2\kappa}} = u \quad \text{ for any } \lambda > 0 \; \text{ and } \; p_0 \in E.
\]
\end{thm}
Finally, we recall the statement of the parabolic Whitney's extension theorem. For a proof of this result we refer to Danielli, Garofalo, Petrosyan and To \cite{DainelliGarofaloPetTo2017:book}.
\begin{thm}\label{PARABOLICWHITNEYEXTENSION} (Parabolic Whitney's extension theorem, \cite[Theorem B.1]{DainelliGarofaloPetTo2017:book})
\\
Let $E$ be a compact subset of $\RR^N\times\RR$, $f: E \to \RR$ a continuous function, and $\{f_{\alpha,j}\}_{|\alpha|+2j \leq 2m}$ with $f_{0,0} = f$ and $m \in \NN$, $\alpha \in \ZZ_{\geq0}^N$ a multi-index. Assume that there exist a family of moduli of continuity $\{\omega_{\alpha,j}\}_{|\alpha|+2j \leq 2m}$ such that
\[
f_{\alpha,j}(x,t) = \sum_{|\beta| + 2i \leq 2m-|\alpha|-2j} \frac{f_{\alpha+\beta,j+i}(x_0,t_0)}{\beta!i!}(x-x_0)^{\beta}(t-t_0)^i + R_{\alpha,j}(x,t;x_0,t_0)
\]
and
\[
|R_{\alpha,j}(x,t;x_0,t_0)| \leq \omega_{\alpha,j}(\|(x-x_0,t-t_0)\|) \|(x-x_0,t-t_0)\|^{2m-|\alpha|-2j}.
\]
Then there exists a function $F \in \mathcal{C}^{2m,m}(\RR^N\times\RR)$ such that $F = f$ on $E$ and $\partial_x^{\alpha}\partial_t^jF = f_{\alpha,j}$ on $E$, for $|\alpha|+2j \leq 2m$.
\end{thm}
\subsection{Proof of Theorem \ref{Theorem:StructureRegularity}} Theorem \ref{Theorem:StructureRegularity} will be proved in a series of steps. We begin by proving that the nodal set of solutions to \eqref{eq:NONLOCALHEPOTENTIAL} has at least codimension one with respect to  the parabolic Hausdorff measure and the regular part is locally smooth.
\begin{thm}\label{THMDIMENSIONALESTIMATENODALSET}
	Let $s \in (0,1)$, $q$ satisfying \eqref{eq:ASSUMPTIONSPOTENTIAL} and let $u$ be a solution to \eqref{eq:NONLOCALHEPOTENTIAL} having finite vanishing order in $B_1\times(-1,0)$. Then:
	\[
	\dim_{\mathcal{\mathcal{P}}}(\Gamma(u)) \leq N+1.
	\]
	Furthermore, $p_0 \in {\rm Reg}(u)$ if and only if $|\nabla_x u(p_0)| \not = 0$. In particular, ${\rm Reg}(u)$ is a locally $\mathcal{C}^1$-manifold of parabolic Hausdorff codimension at least one.	
\end{thm}
\emph{Proof.} We begin by proving the dimensional bound on $\Gamma(u)$: it is obtained combining Federer's theorem and the blow-up classification of Section \ref{Section:BlowUp1} and Section \ref{Section:BlowUp2}.  Let $W$ be a weak solution to \eqref{eq:TruncatedExtensionEquationLocalBackward}, with $q$ satisfying \eqref{eq:ASSUMPTIONSPOTENTIAL} and let $w(x,t) := W(x,0,t)$.
\[
\mathcal{F} := \{ w(x,t) := W(x,0,t) \text{ where $W$ is a weak solution to  \eqref{eq:TruncatedExtensionEquationLocalBackward} } \} \subset L_{loc}^{\infty}(\RR^N\times\RR).
\]
Note that the family $\mathcal{F}$ satisfies assumption $(A_1)$ of Definition \ref{DEFINITIONSELFSIMILARFAMILY} since $W$ does. Now from Corollary \ref{Corollary:THMECONVERGENCELINFINITYLOCUNIQUENESS}, we have
\[
W_{p_0,\lambda_j} \to \Theta_{p_0} \quad \text{ locally uniformly in } \, \overline{\RR_+^{N+1}}\times\RR_+,
\]
as $j \to +\infty$, where $\Theta_{p_0} \in \mathfrak{B}_{2\kappa}(U)$ is the unique (parabolically $2\kappa$-homogeneous) tangent map of $W$ at $p_0 \in \Gamma(W)$ (see  Definitions \ref{DEFBLOWUPSET} and \ref{DEFBLOWUPSET}) and $\lambda_j \to 0$. Consequently,
\[
w_{p_0,\lambda_j} \to \vartheta_{p_0} \quad \text{ locally uniformly in } \, \RR^N\times\RR_+,
\]
as $j \to +\infty$, where $\vartheta_{p_0}(x,t) := \Theta_{p_0}(x,0,t)$ and $w_{p_0,\lambda_j}(x,t) = W_{p_0,\lambda_j}(x,0,t)$. Since $\Theta_{p_0}$ is parabolically $2\kappa$-homogeneous for some $\kappa \geq 0$, we deduce that assumption $(A_2)$ Definition \ref{DEFINITIONSELFSIMILARFAMILY} is satisfied too, by taking
\[
w_{p_0,\lambda_j,\varrho_j}(x,t) = w_{p_0,\lambda_j}(x,t) \qquad (\text{with } \; \varrho_j := r_j^{2\kappa}).
\]
Now, let us consider the map
\[
\mathcal{S}: w \to \mathcal{S}(w) := \Gamma(w) \in \mathcal{C},
\]
since $\Gamma(W)$ is closed by the continuity of $W$. Moreover, by the uniform convergence of $W_{p_0,\lambda}$ towards their tangent maps $\Theta_{p_0}$ (see  Corollary \ref{Corollary:THMECONVERGENCELINFINITYLOCUNIQUENESS}), we easily see that also the assumption $(A_3)$ (of Definition \ref{DEFINITIONSELFSIMILARFAMILY}) is satisfied and so, $(\mathcal{F},\mathcal{S})$ is a locally asymptotically self-similar family. Hence, in view of Federer's Reduction Principle, we conclude the proof since $d \leq N+1$.

\smallskip

Now let us assume $O \in {\rm Reg}(W)$. From the blow up classification of Corollary \ref{Corollary:THMECONVERGENCELINFINITYLOCUNIQUENESS} and the spectral Theorem \ref{SPECTRALTHEOREMEXTENDEDORNUHL1}, we know that as $\lambda \to 0$
\[
W_{\lambda} \to \Theta = \sum_{j=1}^N v_j x_j
\]
locally uniformly in $\overline{\RR_+^{N+1}}\times(0,\infty)$, where $v_i\not=0$ at least for some $i \in {1,\ldots,N}$. Consequently, denoting with $(e_j)_{i=1}^N$ the standard basis of $\RR^N$ and recalling that $O \in \Gamma(W)$, it follows
\[
\lim_{\lambda \to 0^+} W_r(e_i,0,0) = \lim_{\lambda \to 0^+} \frac{W(\lambda e_i,0,0) - W(O)}{\lambda} = \partial_{x_i} W(O).
\]
On the other hand, since $\Theta(O) = 0$, it follows
\[
\begin{aligned}
\lim_{\lambda \to 0^+} W_r(e_i,0,0) & = \lim_{\lambda \to 0^+} \frac{W(\lambda e_i,0,0) - \Theta(\lambda e_i,0,0)}{\lambda} + \lim_{\lambda \to 0^+} \frac{\Theta(\lambda e_i,0,0)- \Theta(O)}{\lambda} \\
& = \lim_{\lambda \to 0^+} W_r(e_i,0,0) - \Theta(e_i) + \partial_{x_i} \Theta(O)
= \partial_{x_i} \Theta(O) = v_i,
\end{aligned}
\]
where we have used the explicit expression of $\Theta$. Consequently, $\partial_{x_i} W(O) = v_i \not= 0$.  To show the converse, we assume $\partial_{x_i} W(O) = v_i \not= 0$ for some $i$ and $O \not\in {\rm Reg}(W)$. So there must be $\kappa > 1/2$ such that
\[
W_r(X,t) = \frac{W(rX,r^2t)}{r^{2\kappa}} \to \Theta(X,t)
\]
locally uniformly as $r \to 0^+$. However,
\[
W_r(e_i,0,0) = r^{1-2\kappa} \frac{W(re_i,0,0) - W(O)}{r} \to +\infty
\]
as $r \to 0^+$, since $v_i \not= 0$. This concludes the proof of the first part of the claim. For what concerns the second part, it is enough to apply the implicit function theorem at a generic point $p_0 \in {\rm Reg}(W)$, using the fact that $|\nabla_x W(p_0)| \not= 0$. $\Box$
\begin{rem}
Theorem \ref{THMDIMENSIONALESTIMATENODALSET} shows that the first bound in \eqref{eq:CHENESTIMATEFARAWAYSIGMA} holds also for the nonlocal framework. As anticipated, the main differences occur in the analysis of the singular set. Indeed, taking $N = 1$, we notice that if $W$ is a weak solution to problem \eqref{eq:TruncatedExtensionEquationLocalBackward}, $O \in \Gamma(W)$ and the limit of the Almgren-Poon quotient is $\kappa = 1$, then its blow up limit must be a linear combination of suitable rescaled eigenfunctions (see  the end of Section \ref{Section:BlowUp1}):
\[
\begin{aligned}
\Theta^{A,B}(x,y,t) & = A\widetilde{\Theta}_{2,0}(x,y,t) + B\widetilde{\Theta}_{0,1}(x,y,t) = A(x^2-2t) + B\left[\left(\frac{1+a}{2}\right)t - \frac{y^2}{4}\right] \\
& = Ax^2 + \left[\frac{1+a}{2}B - 2A\right]t - \frac{B}{4}y^2.
\end{aligned}
\]
In particular, taking $A = 1$, $B = 4/(1+a)$ and $A = 0$, $B = 4/(1+a)$, we obtain
\[
\Theta^{1,\frac{4}{1+a}}(x,y,t) = x^2 - \frac{y^2}{1+a}, \qquad \Theta^{0,\frac{4}{1+a}}(x,y,t) = 2t - \frac{y^2}{1+a},
\]
with traces on $\RR^N\times\RR$
\[
\vartheta^{1,\frac{4}{1+a}}(x,t) := \Theta^{1,\frac{4}{1+a}}(x,0,t) = x^2, \qquad \vartheta^{0,\frac{4}{1+a}}(x,t) := \Theta^{0,\frac{4}{1+a}}(x,0,t) = 2t,
\]
respectively. Consequently, since $\dim_{\mathcal{P}}({\rm Sing}(\vartheta^{1,\frac{4}{1+a}})) =\dim_{\mathcal{P}}({\rm Sing}(\vartheta^{0,\frac{4}{1+a}})) = 2 \not= 1$, we understand that the non-locality of the operator $H^s$ plays a central role the description of the singular set properties. Recalling that
\[
\mathcal{Z}_t(u) := \{x \in \RR^N : u(x,t) = 0\},
\]
we note that also $\dim_{\mathcal{H}}(\mathcal{Z}_0(\vartheta^{0,\frac{4}{1+a}})) = 1 \not= 0$ and so all the relations in \eqref{eq:CHENESTIMATEFARAWAYSIGMA} seems not hold. We stress that an explicit counterexample to inequalities \eqref{eq:CHENESTIMATEFARAWAYSIGMA} has been given for the elliptic-nonlocal setting by Tortone in his Ph.D. dissertation \cite{Tortone2018:book}.
\end{rem}
To prove Theorem \ref{CONTINUOUSDEPOFBLOWUPLIMITS}, we need two monotonicity formulas of Weiss and Monneau type that show below (see  Weiss \cite{Weiss1999a:art,Weiss1999b:art} and Monneau \cite{Monneau2003:art,Monneau2009:art}). We follow the techniques in \cite{DainelliGarofaloPetTo2017:book} adapting them to our framework.
\begin{thm}\label{Lemma:WeissFormula} (Weiss type monotonicity formula)
	Let $W$ be a weak solution to \eqref{eq:TruncatedExtensionEquationLocalBackward}, with $q$ satisfying \eqref{eq:ASSUMPTIONSPOTENTIAL}. Let $\kappa \in \widetilde{\mathcal{K}}$ with $\kappa > 0$, $p_0 \in \Gamma_{2\kappa}(W)$ and $\Theta_{p_0}  \in \mathfrak{B}_{2\kappa}(W)$. Consider the Weiss type functional
	\[
	\mathcal{W}_{2\kappa}(W,p_0,r) := \frac{1}{r^{4\kappa}} [D(W,p_0,r) - \kappa H(W,p_0,r) ]
	\]
	for all $r \in (0,1)$. Then, there exist constants $C' > 0$ and $r_0 \in (0,1)$ such that
	\[
	\frac{d}{dr} \mathcal{W}_{2\kappa}(W,p_0,r) \geq \frac{1}{r^{4\kappa + 3}} \int_0^{r^2}\int_{\mathbb{R}_+^{N+1}} \left( ZW - 2\kk W - tF \right)^2 d\mu_t dt - C'r^{-a},
	\]
	for all $r \in (0,r_0)$. In particular, the function
	\[
	r \to \mathcal{W}_{4\kappa}(W,p_0,r) + Cr^{1-a},
	\]
	is monotone nondecreasing in $r \in (0,r_0)$, for a suitable new constant $C > 0$ depending on $C'$, and
	\[
	\mathcal{W}_{2\kappa}(W,p_0,0^+) := \lim_{r \to 0^+} \mathcal{W}_{2\kappa}(W,p_0,r) = 0.
	\]
\end{thm}
\emph{Proof.} Take for simplcity $p_0 = O$. Using Lemma \ref{Lemma:HeightsVariations} and Lemma \ref{Lemma:DerivativeDd}, we easily compute
\[
\begin{aligned}
r^{4\kappa + 3} \frac{d}{dr} \mathcal{W}_{2\kappa}(W,r) &= \int_0^{r^2}\int_{\mathbb{R}_+^{N+1}} \left( ZW - 2\kk W - tF \right)^2 d\mu_t dt - \int_0^{r^2} t^2 \int_{\mathbb{R}^{N+1}_+} F^2 d\mu_t dt \\
& \quad - \int_0^{r^2} t \int_{\mathbb{R}^N\times\{0\}} [(1-a)q + 2t\partial_t q + x \cdot \nabla q] w^2 \,d\mu_t(x,0)dt.
\end{aligned}
\]
Since the first term in the r.h.s. is nonnegative and the second one satisfies
\[
\int_0^{r^2} t^2 \int_{\mathbb{R}^{N+1}_+} F^2 d\mu_t dt  \leq C r^{\sigma_0} e^{-\frac{1}{16r^2}},
\]
for some $\sigma_0 \in \RR$, $r_0 \in (0,1)$, $C > 0$ and all $r \in (0,r_0)$, it is enough to find a suitable bound for the third term in the r.h.s. Notice that, since the quantity $(1-a)q + 2t\partial_t q + x \cdot \nabla q$ is uniformly bounded in $B_1\times(0,1)$ by \eqref{eq:ASSUMPTIONSPOTENTIAL}, it is enough to focus on $\int_0^{r^2} t \int_{\mathbb{R}^N\times\{0\}} w^2 \,d\mu_t(x,0)dt$.

Now, from Lemma \ref{Lemma:TraceIneqL2Norm} and the definition of $d(W,t)$, we infer
\[
\int_{\mathbb{R}^N\times\{0\}} w^2 \,d\mu_t(x,0) \leq C_0 t^{-\frac{1+a}{2}} \left[ d_0(W,t) + h(W,t) \right] + C_0Kt^{\frac{1-a}{2}} \int_{\mathbb{R}^N\times\{0\}} w^2 \,d\mu_t(x,0)
\]
for some $t_0 \in (0,1)$ and $C_0 > 0$ depending only on $N$, $a$, $K$ in \eqref{eq:ASSUMPTIONSPOTENTIAL}, and all $t \in (0,t_0)$. Consequently, taking eventually $t_0 > 0$ smaller (and depending only on $N$, $a$, $K$), we obtain
\[
t\int_{\mathbb{R}^N\times\{0\}} w^2 \,d\mu_t(x,0) \leq C_0 t^{\frac{1-a}{2}} \left[ d_0(W,t) + h(W,t) \right],
\]
for all $t \in (0,t_0)$ and a new constant $C_0 > 0$. Consequently, there exists $r_0 \in (0,1)$ such that
\[
\frac{1}{r^2} \int_{0}^{r^2} t\int_{\mathbb{R}^N\times\{0\}} w^2 \,d\mu_t(x,0) \leq C_0 r^{1-a} \left[ D_0(W,r) + H(W,r) \right] = C_0 r^{1-a} H(W,r) \left[ N_0(W,r) + 1 \right],
\]
and so
\[
\frac{1}{r^{4\kappa + 3}} \int_{0}^{r^2} t\int_{\mathbb{R}^N\times\{0\}} w^2 \,d\mu_t(x,0) \leq C_0 r^{-a} \frac{H(W,r)}{r^{4\kappa}} \left[ N_0(W,r) + 1 \right] \leq C_1 r^{-a},
\]
for all $r \in (0,r_0)$, up to taking $r_0 > 0$ smaller, in view of Theorem \ref{THEOREMBLOWUP1New} (part (i) and (ii)). Employing such estimate and taking eventually $r_0 > 0$ smaller, we obtain
\[
\frac{d}{dr} \mathcal{W}_{2\kappa}(W,r) \geq - C'r^{-a},
\]
for all $r \in (0,r_0)$ and for any fixed $C' > C_1$, as desired. In particular, the function $r \to \mathcal{W}_{4\kappa}(W,r) + Cr^{1-a}$ is nondecreasing for $C > 0$ large enough ($r \in (0,r_0)$) and the limit of $\mathcal{W}_{4\kappa}(W,r)$ as $r \to 0^+$ exists. Furthermore, writing
\[
\mathcal{W}_{2\kappa}(W,r) := \frac{H(W,r)}{r^{4\kappa}} [N_D(W,r) - \kappa ]
\]
it follows that such limit is zero thanks to Theorem \ref{THEOREMBLOWUP1New} (part (i) and (ii)). $\Box$

\begin{thm}\label{Lemma:MonneauFormula} (Monneau type monotonicity formula)
	Let $W$ be a weak solution to \eqref{eq:TruncatedExtensionEquationLocalBackward}, with $q$ satisfying \eqref{eq:ASSUMPTIONSPOTENTIAL}. Let $\kappa \in \widetilde{\mathcal{K}}$ with $\kappa > 0$, $p_0 \in \Gamma_{2\kappa}(W)$ and $\Theta_{p_0}  \in \mathfrak{B}_{2\kappa}(W)$. Define the Monneau type functional
	\[
	\mathcal{M}_{2\kappa}(W,\Theta,p_0,r) := \frac{1}{r^{4\kappa}} H(W-\Theta,p_0,r)
	\]
	for all $r \in (0,1)$. Then, there exist constants $C'' > 0$ and $r_0 \in (0,1)$ such that
	\[
	\frac{d}{dr} \mathcal{M}_{2\kappa}(W,\Theta,p_0,r) \geq - C'' r^{-a},
	\]
	for all $r \in (0,r_0)$. In particular, the function
	\[
	r \to \mathcal{M}_{2\kappa}(W,\Theta,p_0,r) + Cr^{1-a},
	\]
	is monotone nondecreasing in $r \in (0,r_0)$, for a suitable new constant $C > 0$ depending on $C''$.
\end{thm}
\emph{Proof.} Fix $p_0 = O$ and set $V := W - \Theta$. Repeating the computations of Lemma \ref{Lemma:HeightsVariations}, taking into account that $V$ is a global weak solution to
\[
\begin{cases}
\partial_t V + y^{-a} \nabla \cdot(y^a \nabla V) = F \quad &\text{ in } \RR_+^{N+1} \times (0,1) \\
-\partial_y^a V = q(v + \vartheta) \quad &\text{ in } \RR^N\times\{0\} \times (0,1),
\end{cases}
\]
it easily follows
\[
H'(V,r) =  \frac{4}{r} D(W,r) + \frac{4}{r^3} \int_{0}^{r^2} t \int_{\RR^{N+1}_+} FV d\mu_t dt - \frac{4}{r^3} \int_0^{r^2} t \int_{\RR^N\times\{0\}} qv\vartheta d\mu_t(x,0) dt.
\]
Consequently, using the definition of the Weiss functional, we obtain
\[
\begin{aligned}
\frac{d}{dr} \mathcal{M}_{2\kappa}(W,\Theta,r) &= \frac{H'(V,r)}{r^{4\kappa}} - \frac{4\kappa}{r^{4\kappa+1}} H(V,r) \\
& = \frac{4}{r} \mathcal{W}_{2\kappa}(V,r) + \frac{4}{r^{4\kappa + 3}} \int_{0}^{r^2} t \int_{\RR^{N+1}_+} FV d\mu_t dt - \frac{4}{r^{4\kappa+3}} \int_0^{r^2} t \int_{\RR^N\times\{0\}} qv\vartheta d\mu_t(x,0) dt.
\end{aligned}
\]

\emph{Step 1: Estimates.} Using that $V = W - \Theta$, the fact that $W$ has compact support and that $\Theta$ is a homogeneous polynomial it is not difficult to see that, as in the previous proof, it holds
\[
\Big| \int_0^{r^2} t \int_{\mathbb{R}^{N+1}_+} FV d\mu_t dt \Big| \leq C r^{\sigma_0} e^{-\frac{1}{16r^2}},
\]
for some $\sigma_0 \in \RR$, $r_0 \in (0,1)$, $C > 0$ and all $r \in (0,r_0)$. Consequently, we can focus on the third term. Thanks to the fact that $q$ is bounded and the H\"older inequality, we have
\[
\bigg[ \int_0^{r^2} t \int_{\RR^N\times\{0\}} qv\vartheta d\mu_t(x,0) dt \bigg]^2 \leq K^2 \int_0^{r^2} \int_{\RR^N\times\{0\}} \vartheta^2 d\mu_t(x,0) dt \int_0^{r^2} t^2 \int_{\RR^N\times\{0\}} v^2 d\mu_t(x,0) dt.
\]
Using that $\theta$ is parabolically $2\kappa$-homogeneous, it follows
\begin{equation}\label{eq:EstthMonneau}
\int_0^{r^2} \int_{\RR^N\times\{0\}} \vartheta^2 d\mu_t(x,0) dt = \int_0^{r^2} t^{2\kappa - \frac{1+a}{2}}  \int_{\RR^N\times\{0\}} \vartheta(x,1)^2 d\mu(x,0) dt \leq C r^{4\kappa + 1 - a},
\end{equation}
while, following the ideas of the proof of Theorem \ref{Lemma:WeissFormula}, we obtain
\[
t^2\int_{\mathbb{R}^N\times\{0\}} v^2 \,d\mu_t(x,0) \leq C_0 t^{1+\frac{1-a}{2}} \left[ d_0(V,t) + h(V,t) \right],
\]
and so
\[
\frac{1}{r^2} \int_0^{r^2} t^2 \int_{\mathbb{R}^N\times\{0\}} v^2 \,d\mu_t(x,0)dt \leq C_0 r^{3-a} \left[ D_0(V,r) + H(V,r) \right] = C_0 r^{3-a} H(V,r) \left[ N_0(V,r) + 1 \right].
\]
Now, since $\Theta$ is parabolically $2\kappa$-homogeneous, it is not difficult to see that
\[
H(\Theta,r) = r^{4\kappa},
\]
for all $r > 0$, up to a nonzero multiplicative constant. The above relation and Theorem  \ref{THEOREMBLOWUP1New} (part (ii)) imply
\[
\begin{aligned}
H(V,r) &= H(W,r) + H(\Theta,r) - \frac{2}{r^2} \int_0^{r^2} \int_{\mathbb{R}_+^{N+1}} W \Theta \,d\mu_tdt \\
& \leq H(W,r) + H(\Theta,r) + 2H^{1/2}(W,r)H^{1/2}(\Theta,r) \leq Cr^{4\kappa},
\end{aligned}
\]
for some suitable constants $r_0,C > 0$ and all $r \in (0,r_0)$. In particular,
\[
\int_0^{r^2} t^2 \int_{\mathbb{R}^N\times\{0\}} v^2 \,d\mu_t(x,0)dt \leq C r^{4\kappa + 5-a},
\]
for some new $C > 0$ and all $r \in (0,r_0)$. Putting all these bounds together, it follows
\begin{equation}\label{eq:EstvthMonneau}
\left| \int_0^{r^2} t \int_{\RR^N\times\{0\}} qv\vartheta d\mu_t(x,0) dt \right| \leq Cr^{4\kappa + 3 - a},
\end{equation}
for all $r \in (0,r_0)$, which implies
\begin{equation}\label{eq:EstMonneauDer1}
\frac{d}{dr} \mathcal{M}_{2\kappa}(W,\Theta,r) \geq \frac{4}{r} \mathcal{W}_{2\kappa}(V,r) - Cr^{-a},
\end{equation}
for all $r \in (0,r_0)$ and a new $C > 0$.

\emph{Step 2: An additional property of the Weiss function.} In this step we show that
\begin{equation}\label{eq:WeissPropWV}
\mathcal{W}_{2\kappa}(V,r) \geq \mathcal{W}_{2\kappa}(W,r) - Cr^{1-a}
\end{equation}
for all $r \in (0,r_0)$ and a new $C > 0$. From the explicit formulations of $H(W,r)$ and $D(W,r)$, we obtain
\[
\begin{aligned}
\mathcal{W}_{2\kappa}(W,r) &= \mathcal{W}_{2\kappa}(V + \Theta,r) \\
& = \frac{1}{r^{4\kappa}} [D_0(V+\Theta,r) - \kappa H(V+\Theta,r)] - \frac{1}{r^{4\kappa + 2}} \int_0^{r^2} t \int_{\mathbb{R}^N\times\{0\}} q (v+\vartheta)^2 \,d\mu_t(x,0)dt \\
& = \mathcal{W}_{2\kappa}(V,r) + \frac{1}{r^{4\kappa}} [D_0(\Theta,r) - \kappa H(\Theta,r)] + \frac{1}{r^{4\kappa+2}} \int_0^{r^2} \int_{\mathbb{R}_+^{N+1}} [t \nabla \Theta \cdot \nabla V - \kappa \Theta V] d\mu_t dt \\
& \quad - \frac{1}{r^{4\kappa + 2}} \int_0^{r^2} t \int_{\mathbb{R}^N\times\{0\}} q \vartheta \left(\vartheta + 2v \right) \,d\mu_t(x,0)dt \\
& = \mathcal{W}_{2\kappa}(V,r) - \frac{1}{r^{4\kappa + 2}} \int_0^{r^2} t \int_{\mathbb{R}^N\times\{0\}} q \vartheta \left(\vartheta + 2v \right) \,d\mu_t(x,0)dt,
\end{aligned}
\]
where we have used that $\Theta$ is parabolically $2\kappa$-homogeneous (see  Theorem \ref{Theorem:MonAlmgrenPoon}) and it is a re-scaled eigenfunction with eigenvalue $\kappa$ (see  Theorem \ref{SPECTRALTHEOREMEXTENDEDORNUHL1}). Finally, using the boundedness of $q$, \eqref{eq:EstthMonneau} and \eqref{eq:EstvthMonneau} , it follows
\[
\bigg| \int_0^{r^2} t \int_{\mathbb{R}^N\times\{0\}} q \vartheta \left(\vartheta + 2v \right) \,d\mu_t(x,0)dt \bigg| \leq Cr^{4\kappa + 3 - a},
\]
for some suitable constant $C > 0$ and all $r \in (0,r_0)$, and our claim follows.

\emph{Step 3: Conclusions.} Now, using \eqref{eq:WeissPropWV} and recalling that $r \to W_{2\kappa}(W,r) + C'r^{1-a}$ is nondecreasing in view of Theorem \ref{Lemma:WeissFormula}, \eqref{eq:EstMonneauDer1} becomes
\[
\begin{aligned}
\frac{d}{dr} \mathcal{M}_{2\kappa}(W,\Theta,r) &\geq \frac{4}{r} \left[ \mathcal{W}_{2\kappa}(W,r) + C'r^{1-a} \right] - (4C' + C) r^{-a} \\
& \geq \frac{4}{r} \mathcal{W}_{2\kappa}(W,0^+) - C'' r^{-a} \geq - C'' r^{-a},
\end{aligned}
\]
for all $r \in (0,r_0)$, where we have used that $\mathcal{W}_{2\kappa}(W,0^+) = 0$. This concludes the proof. $\Box$
\begin{rem}\label{Rem:WeakAssPotential} Notice that in the proof above we have used (for the first time) the fact that $q$ is defined and bounded in the whole space (together with its derivatives of suitable order), see  \eqref{eq:ASSUMPTIONSPOTENTIAL}. They are needed since the function $V := W - \Theta$ has not compact support, whilst $W$ does. In particular, it follows that the statements of Theorem \ref{Theorem:GeneralizedFrequency}, Theorem \ref{THEOREMBLOWUP1New}, Theorem \ref{CONVERGENCELINFINITYLOC} and Theorem \ref{THMDIMENSIONALESTIMATENODALSET}, can be proven under the weaker assumption
\[
\begin{cases}
\|q\|_{C^1(Q_1)} \leq K \quad &\text{ if } 1/2 \leq s < 1 \\
\|q\|_{C^2(Q_1)}, \; \|x\cdot\nabla_x q\|_{L^{\infty}(Q_1)} \leq K \quad &\text{ if } 0 < s < 1/2,
\end{cases}
\]
where $K > 0$ is a universal constant and $Q_1 := B_1\times(0,1)$.
\end{rem}	
Theorem \ref{Lemma:WeissFormula} and Theorem \ref{Lemma:MonneauFormula} allow us to prove the main result concerning the asymptotic behaviour and the differentiability of solutions near nodal points belonging to $\Gamma_{2\kappa}(W)$. Furthermore, we establish the continuity of the tangent map $p_0 \to \Theta_{p_0}$, seen as a function from $\Gamma_{2\kappa}(W)$ to $\mathfrak{B}_{2\kappa}(W)$.
\begin{thm}\label{CONTINUOUSDEPOFBLOWUPLIMITS} (Continuous dependence of the blow-up limits)
Let $W$ be a weak solution to \eqref{eq:TruncatedExtensionEquationLocalBackward}, with $q$ satisfying \eqref{eq:ASSUMPTIONSPOTENTIAL}. Let $p_0 \in \Gamma_{2\kappa}(W)$ and $\Theta_{p_0}  \in \mathfrak{B}_{2\kappa}(W)$ the tangent map of $W$ at $p_0$ (see  Definition \ref{DEFBLOWUPSET}), for some $\kappa \in \widetilde{\mathcal{K}}$ with $\kappa > 0$. Then the following two assertions hold:
	
	(i) We have as $\|(X,t)\|^2 := |X|^2 + |t| \to 0^+$
	\[
	W_{p_0}(X,t) = \Theta_{p_0}(X,t) + o(\|(X,t)\|^{2\kappa}).
	\]

	(ii) The map $p_0 \to \Theta_{p_0}$ from $\Gamma_{2\kappa}(W)$ to $\mathfrak{B}_{2\kappa}(W)$ is continuous.

	(iii) For any compact set $K \subset \Gamma_{2\kappa}(W)$ there exists a modulus of continuity $\sigma = \sigma_{K}$ with $\sigma(0^+) = 0$, such that as $\|(X,t)\|^2 := |X|^2 + |t| \to 0^+$,
	\[
	|W_{p_0}(X,t) - \Theta_{p_0}(X,t)| \leq \sigma(\|(X,t)\|) \, \|(X,t)\|^{2\kappa},
	\]
	for any $p_0 \in K$ (the main fact here is that $\sigma$ does not depend on $p_0$, but only on $K$).
\end{thm}
\emph{Proof.} Part (i) easily follows from Corollary \ref{Corollary:THMECONVERGENCELINFINITYLOCUNIQUENESS}. Indeed, assuming for simplicity $p_0 = O \in \Gamma_{2\kappa}(U)$, we have as $r \to 0^+$
\[
\begin{aligned}
o(1) &= \sup_{(X,t) \in \QQ_1^+} |W_r(X,t) - \Theta(X,t)| = \sup_{(X,t) \in \QQ_1^+} \left|\frac{W(rX,r^2t)}{r^{2\kappa}} - \frac{\Theta(rX,r^2t)}{r^{2\kappa}}\right| \\
&= r^{-2\kappa} \sup_{(X,t) \in \QQ_r^+} |W(X,t) - \Theta(X,t)|,
\end{aligned}
\]
where $\QQ_r^+ := \BB_r^+\times(0,r^2)$, which is equivalent to our first claim.

Let us prove part (ii). Fix $p_0 = O \in \Gamma_{2\kappa}(W)$. We apply Theorem \ref{THEOREMBLOWUP1New} (part (iii)) and use the homogeneity of $\Theta$ to infer that for any $\varepsilon > 0$, there is $r_{\varepsilon} > 0$ such that
\[
\mathcal{M}_{2\kappa}(W,\Theta,r_{\varepsilon}) := \frac{1}{r_{\varepsilon}^{4\kappa}} H(W-\Theta,r_{\varepsilon}) < \frac{\varepsilon}{2}.
\]
On the other hand, using the continuity of the map $\Gamma(U) \ni p \to U_p$, we deduce the existence of $\delta_\varepsilon > 0$ such that
\[
\mathcal{M}_{2\kappa}(W,W_p,r_{\varepsilon}) := \frac{1}{r_{\varepsilon}^{4\kappa}} H(W-W_p,r_{\varepsilon}) < \frac{\varepsilon}{2},
\]
for any $p \in \Gamma_{2\kappa}(W)$ satisfying $\|p\| < \delta_\varepsilon$, and so
\[
\mathcal{M}_{2\kappa}(W_p,\Theta,r_{\varepsilon}) \leq \mathcal{M}_{2\kappa}(W,W_p,r_{\varepsilon}) + \mathcal{M}_{2\kappa}(W,\Theta,r_{\varepsilon}) < \varepsilon.
\]
Thus, in view of Theorem \ref{Lemma:MonneauFormula}, we obtain
\[
\mathcal{M}_{2\kappa}(W_p,\Theta,r) + Cr^{1-a} \leq \mathcal{M}_{2\kappa}(W_p,\Theta,r_{\varepsilon}) + Cr_{\varepsilon}^{1-a} < \varepsilon + Cr_{\varepsilon}^{1-a},
\]
for some constant $C > 0$ depending $p$ and all $r \in (0,r_{\varepsilon})$ (notice that we are crucially using the arbitrariness of $\Theta \in \mathfrak{B}_{2\kappa}(W)$,  see  Theorem \ref{Lemma:MonneauFormula}). Consequently, taking the limit as $r \to 0$, we conclude
\[
\int_0^1 \int_{\RR_+^{N+1}} (\Theta_p - \Theta)^2 d\mu_t dt = \lim_{r \to 0^+} \mathcal{M}_{2\kappa}(W_p,\Theta,r) \leq \varepsilon + Cr_{\varepsilon}^{1-a}.
\]
Since $\mathfrak{B}_{2\kappa}(W)$ is a finite dimensional vector space and $r_{\varepsilon} > 0$ can be chosen arbitrarily small, the proof of part (ii) is completed.

Finally, we prove part (iii). Fix a compact set $K \subset \Gamma_{2\kappa}(W)$. We apply part (ii) (the last two formulas above) to deduce that for any fixed $p_0 \in K$
\[
\begin{aligned}
\int_0^{r^2} \int_{\RR_+^{N+1}} (W_p - \Theta_p)^2 d\mu_t dt &\leq \int_0^{r^2} \int_{\RR_+^{N+1}} (W_p - \Theta_{p_0})^2 d\mu_t dt + \int_0^{r^2} \int_{\RR_+^{N+1}} (\Theta_{p_0} - \Theta_p)^2 d\mu_t dt \\
& \leq (2\varepsilon + Cr_{\varepsilon}^{1-a})r^{4\kappa + 2},
\end{aligned}
\]
for all $p \in K$ satisfying $\|p_0 - p\| < \delta_\varepsilon$ and all $r \in (0,r_{\varepsilon})$. Equivalently, we can write
\[
\|W_{p,r} - \Theta_p \|_{L^2_{\mu_t}}^2 := \int_0^1 \int_{\RR_+^{N+1}} (W_{p,r} - \Theta_p)^2 d\mu_t dt \leq (2\varepsilon + Cr_{\varepsilon}^{1-a}),
\]
whenever $\|p_0-p\| < \delta_\varepsilon$ and all $r \in (0,r_{\varepsilon})$. Notice that a straightforward contradiction argument shows that the constant $C > 0$ can be made to depend only on $p_0$. Now, repeating the above considerations for any $p_0 \in K$, we obtain an open covering of $K$ made by cylinders of the type $\QQ_{\delta_{\varepsilon,p_0}}^+(p_0)$ such that
\[
\|W_{p,r} - \Theta_p \|_{L^2_{\mu_t}}^2 \leq C_{p_0}\varepsilon,
\]
for some suitable $C_{p_0} > 0$, all $p \in \QQ_{\delta_{\varepsilon,p_0}}^+(p_0) \cap K$ and all $r \in (0,r_{\varepsilon})$. Thus, extracting a finite covering $\{\QQ_{\delta_{\varepsilon,p_i}}^+(p_i)\}_{i \in \{ 1,\ldots,M\}}$ it follows
\[
\|W_{p,r} - \Theta_p \|_{L^2_{\mu_t}}^2 \leq C_K\varepsilon,
\]
for some suitable $C_K,r_{\varepsilon}^K > 0$, all $p \in K$ and all $r \in (0,r_{\varepsilon}^K)$. Now, from  \cite[Formula (5.7)]{BanGarofalo2017:art}, we have the $L^2-L^{\infty}$ bound
\begin{equation}\label{eq:L2LINFINITYBOUND}
\sup_{\QQ_{1/4}^+} |W_{p,r} - \Theta_{p}|^2 \leq C_0 \int_0^{1/4}\int_{\BB_{1/2}^+} y^a |W_{p,r} - \Theta_{p}|^2  \,dXdt,
\end{equation}
where $C_0 > 0$ is a universal constant. Using the fact that $W_{p,r} \to \Theta_{p}$ in $L^2(0,1;L^2_{\mu_t})$ as $r \to 0^+$, we deduce that the r.h.s. of \eqref{eq:L2LINFINITYBOUND} converges to zero as $r \to 0^+$ on sets of the type $\BB_{1/2}^+\cap\{y > A\} \times(\delta,1/4)$, for any fixed $\delta,A > 0$: this follows from the fact that the function $y^a \mathcal{G}_a$ is bounded above and below on $\BB_{1/2}^+\cap\{y > A\}\times(\delta,1)$, for any fixed $\delta,A > 0$. Consequently, using a standard Cantor diagonal argument, it is not difficult to see that
\[
\int_0^{1/4}\int_{\BB_{1/2}^+} y^a |W_{p,r} - \Theta_{p}|^2  \,dXdt \to 0
\]
as $r \to 0$, for any fixed $p \in K$. Further, using a covering argument as before, it follows
\[
\int_0^{1/4}\int_{\BB_{1/2}^+} y^a |W_{p,r} - \Theta_{p}|^2  \,dXdt \leq C_K\varepsilon,
\]
for some $C_K,r_{\varepsilon}^K > 0$ and all $r \in (0,r_{\varepsilon}^K)$. Consequently, in view of \eqref{eq:L2LINFINITYBOUND} and the homogeneity of $\Theta_p$, we obtain
\[
\sup_{\QQ_{r/4}^+} |W_p - \Theta_p|^2 \leq C_K\varepsilon r^{2\kappa},
\]
for all $r \in (0,r_{\varepsilon}^K)$, which is equivalent to our claim, thanks to the arbitrariness of $\varepsilon > 0$. $\Box$

\bigskip

Finally, to state the last result of this section, we must give two more definitions inspired by the work of Danielli et al. \cite[Definition 12.9]{DainelliGarofaloPetTo2017:book}.
\begin{defn}\label{DEFOFDIMENSIONSINGULARPOINTS}
	Let $W$ be a weak solution to \eqref{eq:TruncatedExtensionEquationLocalBackward}, with $q$ satisfying \eqref{eq:ASSUMPTIONSPOTENTIAL} and $p_0 \in \Gamma_{2\kappa}(W)$. We define the spatial dimension of $\Gamma_{2\kappa}(W)$ at $P_0$ as
	\[
	\begin{aligned}
	d_{p_0}^{2\kappa} := \dim \left\{\xi \in \RR^N: \xi \cdot \nabla_x \partial_x^{\alpha} \partial_t^j \Theta_{p_0} = 0, \text{ for any } \alpha \in \ZZ_{\geq 0}^N, \, j \in \NN \text{ with } |\alpha| + 2j = 2\kappa - 1 \right\},
	\end{aligned}
	\]
	where $\Theta_{p_0} \in \mathfrak{B}_{2\kappa}(U)$ is the blow-up limit of $W$ at $p_0$. Further, for any $d = 0,\ldots,N$, we define
	\[
	\Gamma_{2\kappa}^d(W) := \left\{p_0 \in \Gamma_{2\kappa}(W): d_{p_0}^{2\kappa} = d \right\}.
	\]
	Finally, if $u$ is a solution to \eqref{eq:NONLOCALHEPOTENTIAL}, we set $\Gamma_{2\kappa}^d(u) = \Gamma_{2\kappa}^d(W)$.
\end{defn}
\begin{lem}\label{LEMTIMELIKESINGPOINTS} (Time-like singular points)
	Let $W$ be a weak solution to \eqref{eq:TruncatedExtensionEquationLocalBackward}, with $q$ satisfying \eqref{eq:ASSUMPTIONSPOTENTIAL}. Then for any $\kappa \in \widetilde{\mathcal{K}}$ with $\kappa > 0$ and any $p_0 \in \Gamma_{2\kappa}^N(W)$, the tangent map $\Theta_{p_0} \in \mathfrak{B}_{2\kappa}(W)$ of $W$ at $p_0$ depends only on $y$ and $t$.
	
	Furthermore, if $\kappa$ is an integer then $\vartheta_{p_0}(x,t) = \Theta_{p_0}(x,0,t) = t^{\kappa}$, up to a multiplicative constant, while if $\kappa$ is a half-integer, then $\Gamma_{2\kappa}^N(W) = \emptyset$.
\end{lem}
\emph{Proof.} Assume $p_0 = O$ and write
\[
\Theta(x,y,t) = \sum_{\beta + |\alpha| + 2j = 2\kappa} c_{\alpha,\beta,j} x^\alpha y^\beta t^j = \sum_{|\alpha| + 2j = 2\kappa} c_{\alpha,0,j} x^\alpha t^j + \sum_{\substack{\beta + |\alpha| + 2j = 2\kappa \\ |\alpha| + 2j \leq 2\kappa - 1} } c_{\alpha,\beta,j} x^\alpha y^\beta t^j.
\]
Assuming $p_0 \in \Gamma_{2\kappa}^N(W)$ is equivalent to say that
\[
\partial_{x_i} \partial_x^{\alpha} \partial_t^j \Theta = 0, \text{ for any } \alpha \in \ZZ_{\geq 0}^N, \, j \in \NN \text{ with } |\alpha| + 2j = 2\kappa - 1,
\]
and all $i \in \{1,\ldots,N\}$. Consequently, from the above expression of $\Theta$ we obtain $c_{\alpha,0,j} = 0$ whenever $|\alpha| + 2j = 2\kappa$. In particular, $\Theta(x,0,t) = 0$ for all $(x,t) \in \RR^{N+1}$. Now, recalling the formulation of $\Theta$ in terms of Hermite and Laguerre polynomials (see  Theorem \ref{THEOREMBLOWUP1New} and Theorem \ref{SPECTRALTHEOREMEXTENDEDORNUHL1}), it follows
\[
0 = \Theta(x,0,t) = t^{\kappa} \sum_{\alpha,m} v_{\alpha,m} H_{\alpha}\left(\frac{x}{\sqrt{t}}\right) L_{(\frac{a-1}{2}),m}\left(\frac{y^2}{4t}\right) \Big|_{y=0},
\]
for all $(x,t) \in \RR^{N+1}$, which is impossible unless $\Theta$ is independent of $x$. Now, if $\kappa$ is an integer, using the explicit expression of $L_{(\frac{a-1}{2}),m}$ (see  \eqref{eq:LaguerreExpression}), we can write
\[
\Theta(y,t) = \sum_{\beta + 2j = 2\kappa} c_{\beta,j} y^\beta t^j = c_{0,\kappa} t^{\kappa}  +  \sum_{\substack{\beta + 2j = 2\kappa \\ j < \kappa } } c_{\beta,j} y^\beta t^j,
\]
with $c_{0,\kappa} \not= 0$, and so $\vartheta_{p_0}(x,t) = c_{0,\kappa} t^{\kappa}$. On the other hand, whenever $\kappa$ is a half integer $2\kappa = 2m +1$ for some integer $m \geq 1$. Consequently,
\[
\Theta(y,t) = \sum_{\substack{\beta + 2j = 2\kappa \\ \beta \text{ odd } }} c_{\beta,j} y^\beta t^j,
\]
which is impossible in view of the explicit expression of $L_{(\frac{a-1}{2}),m}$ (see  again with \eqref{eq:LaguerreExpression}), unless $c_{\beta,j} = 0$ for all $\beta$ and $j$. This in contradiction with the fact that $\kappa > 0$ and so $\Gamma_{2\kappa}^N(W) = \emptyset$. $\Box$
\begin{defn}\label{DEFSPACELIKETIMELIKE} (Definition 12.11 of \cite{DainelliGarofaloPetTo2017:book})
	\\
	A $(d+1)$-dimensional manifold $\mathcal{M} \subset \RR^N\times\RR$ (with $d = 0,\ldots,N-1$) is said to be space-like of class $\mathcal{C}^{1,0}$ if it can be locally represented as a graph of a $\mathcal{C}^{1,0}$ function $g:\RR^d\times\RR \to \RR^{N-d}$
	\[
	(x_{d+1},\ldots,x_N) = g(x_1,\ldots,x_d,t),
	\]
	up to rotation of coordinate axis in $\RR^N$.
	\\
	A $N$-dimensional manifold $\mathcal{M} \subset \RR^N\times\RR$ is said to be time-like of class $\mathcal{C}^1$ if it can be locally represented as a graph of a $\mathcal{C}^1$ function $g:\RR^N \to \RR$
	\[
	t = g(x_1,\ldots,x_N).
	\]
\end{defn}
We now state our last theorem. The proof is based on a suitable adaptation of some techniques developed by Danielli et al. \cite{DainelliGarofaloPetTo2017:book} for the local parabolic setting (see  also with \cite{GarofaloPetrosyan2008:art} for the elliptic one), where they combine the implicit function theorem and a parabolic version of the Whitney's extension theorem.
\begin{thm}\label{STRUCTUREOFSINGULARSET} (Structure of the singular set) Let $u$ a solution to \eqref{eq:NONLOCALHEPOTENTIAL} having finite vanishing order in $B_1\times(-1,0)$. Then for any $k \in \widetilde{\mathcal{K}}$ with $\kappa > 0$, the set $\Gamma_{2\kappa}^d(u)$ is contained in a countable union of $(d+1)$-dimensional space-like $\mathcal{C}^{1,0}$ manifolds for any $d = 0,\ldots,N-1$, while $\Gamma_{2\kappa}^N(u)$ is contained in a countable union of $N$-dimensional time-like $\mathcal{C}^1$ manifolds.
\end{thm}
\emph{Proof.} Let $W$ be a weak solution to \eqref{eq:TruncatedExtensionEquationLocalBackward}, with $q$ satisfying \eqref{eq:ASSUMPTIONSPOTENTIAL}, $p_0 \in \Gamma_{2\kappa}(W)$ and $\Theta_{p_0}$ the tangent map of $W$ at $p_0$. As always, we set $\vartheta_{p_0}(x,t) = \Theta_{p_0}(x,0,t)$.

\smallskip

\emph{Step1: Parabolic Whitney's extension.} Since, $\vartheta_{p_0}$ is a parabolically homogeneous polynomial of degree $2\kappa$, we can write it in the form
\begin{equation}\label{eq:POLYNOMIALEXPRESSIONBLOWUPS}
\vartheta_{p_0}(x,t) = \sum_{|\alpha| + 2j = 2\kappa} \frac{a_{\alpha,j}(x_0,t_0)}{\alpha!j!} x^{\alpha} t^j,
\end{equation}
where the coefficients $p_0 \to a_{\alpha,j}(p_0)$ are continuous on $\Gamma_{2\kappa}(u)$, thanks to part (ii) of Theorem \ref{CONTINUOUSDEPOFBLOWUPLIMITS}. Now, we define
\[
f_{\alpha,j}(x,t) :=
\begin{cases}
0 \quad & \text{ if } |\alpha| + 2j < 2\kappa \\
a_{\alpha,j}(x,t) \quad & \text{ if } |\alpha| + 2j = 2\kappa,
\end{cases}
\]
and we proceed by checking the compatibility conditions of Theorem \ref{PARABOLICWHITNEYEXTENSION}, i.e. proving the following claim:

\smallskip

\noindent CLAIM: Fix $h \in \NN$, define $E_h$ as in \eqref{eq:SETEJNOPENREGULARCLOSEDSINGULAR} and set $K = E_h$. Then for any $(x_0,t_0),(x,t) \in K$,
\begin{equation}\label{eq:FJALPHAWITHREMAINDER}
f_{\alpha,j}(x,t) = \sum_{|\beta| + 2i \leq 2\kappa - |\alpha| -2j} \frac{f_{\alpha+\beta,j+i}(x_0,t_0)}{\beta!i!} (x-x_0)^{\beta}(t-t_0)^j + R_{\alpha,j}(x,t;x_0,t_0),
\end{equation}
with
\begin{equation}\label{eq:REMAINDERDECAYING}
|R_{\alpha,j}(x,t;x_0,t_0)| \leq \sigma_{\alpha,j}(\|(x-x_0,t-t_0)\|) \|(x-x_0,t-t_0)\|^{2\kappa - |\alpha| -2j},
\end{equation}
where $\sigma_{\alpha,j}$ are suitable moduli of continuity depending on $K$.

\smallskip

To prove it we first consider the case $|\alpha| + 2j = 2\kappa$. In this setting, the claim follows from the continuity of the functions $p_0 \to a_{\alpha,j}(p_0)$ on $\Gamma_{2\kappa}(u)$, by taking
\[
R_{\alpha,j}(x,t;x_0,t_0) = a_{\alpha,j}(x,t) - a_{\alpha,j}(x_0,t_0).
\]
The case $0 \leq |\alpha| + 2j < 2\kappa$ is more involved. According to Taylor expansion, we define
\[
\begin{aligned}
R_{\alpha,j}(x,t;x_{0},t_{0}) & := - \sum_{\substack{(\beta,i) \geq (\alpha,j): \\ |\beta| + 2i = 2\kappa}} \frac{a_{\beta,i}(x_0,t_0)}{(\beta-\alpha)!(i-j)!} (x-x_0)^{\beta-\alpha}(t-t_0)^{i-j} \\
& = - \partial_x^{\alpha}\partial_t^j \vartheta_{p_0}(x-x_0,t-t_0), \qquad t \geq t_0.
\end{aligned}
\]
Now, assume by contradiction, there is no modulus of continuity $\sigma_{\alpha,j} = \sigma_{\alpha,j}(\cdot)$ such that \eqref{eq:REMAINDERDECAYING} is satisfied, i.e., there are sequences $p_l := (x_l,t_l), p_{0l}:=(x_{0l},t_{0l}) \in K$,
\[
\delta_l := \max\{|x_l-x_{0l}|,\sqrt{|t_l-t_{0l}|}\} \to 0 \quad \text{ as } l \to +\infty,
\]
such that
\begin{equation}\label{eq:REDUCTIONTOABSURBUMASSWHITNEY}
\begin{aligned}
&|R_{\alpha,j}(x_l,t_l;x_{0l},t_{0l})|  \\
&=  |\partial_x^{\alpha}\partial_t^j \vartheta_{p_0}(x_l-x_{l0},t_l-t_{l0})| = \left|\sum_{\substack{(\beta,i) \geq (\alpha,j): \\ |\beta| + 2i = 2\kappa}} \frac{a_{\beta,i}(x_0,t_0)}{(\beta-\alpha)!(i-j)!} (x_l-x_{0l})^{\beta-\alpha}(t_l-t_{0l})^{i-j}\right| \\
&\geq \underline{\sigma} \|(x_l-x_{0l},t_l-t_{0l})\|^{2\kappa - |\alpha| -2j},
\end{aligned}
\end{equation}
for some $\underline{\sigma} > 0$ and all $l \geq 0$. Now we define the families
\[
u_{p_{0l},\delta_l}(x,t) = \frac{u_{p_{0l}}(\delta_lx,\delta_l^2t)}{\delta_l^{2\kappa}}, \qquad (\xi_l,\theta_l) := \left(\frac{x_l - x_{0l}}{\delta_l},\frac{t_l - t_{0l}}{\delta_l^2}\right),
\]
and we assume (up to subsequences) $(x_{0l},t_{0l}) \to (x_0,t_0) \in K$ and $(\xi_l,\vartheta_l) \to ( \xi_0,\theta_0) \in \partial(B_1\times(-1,1))$, and so, from Corollary \ref{Corollary:THMECONVERGENCELINFINITYLOCUNIQUENESS}, it follows
\[
\|u_{p_{0l},\delta_l} - \vartheta_{p_0}\|_{L^{\infty}_{loc}} \leq \|u_{p_{0l},\delta_l} - u_{p_0,\delta_l}\|_{L^{\infty}_{loc}} + \|u_{p_0,\delta_l} - \vartheta_{p_0}\|_{L^{\infty}_{loc}} \to 0,
\]
as $l \to \infty$. The same holds true for the sequence
\[
u_{p_l,\delta_l}(x,t) = \frac{u_{p_l}(\delta_lx,\delta_l^2t)}{\delta_l^{2\kappa}},
\]
so that
\begin{equation}\label{eq:LINFINITYCONVERGECESEQUENCESBLOWUP}
\begin{aligned}
u_{p_{0l},\delta_l} \to \vartheta_{p_0} \quad &\text{ in } L^{\infty}_{loc} \\
u_{p_l,\delta_l} \to \vartheta_{p_0} \quad &\text{ in } L^{\infty}_{loc},
\end{aligned}
\end{equation}
as $l \to \infty$ and, consequently,
\begin{equation}\label{eq:LINFINITYCONVERGECESEQUENCESBLOWUP2}
\begin{aligned}
\|u_{p_{0l},\delta_l}(\cdot + \xi_l,\cdot + \theta_l) - u_{p_l,\delta_l}(\cdot,\cdot)\|_{L^{\infty}_{loc}} &\to 0 \\
\|u_{p_l,\delta_l}(\cdot,\cdot) - u_{p_{0l},\delta_l}(\cdot - \xi_l,\cdot - \theta_l)\|_{L^{\infty}_{loc}} &\to 0,
\end{aligned}
\end{equation}
as $l \to \infty$. Now, as in \cite{DainelliGarofaloPetTo2017:book}, we proceed by splitting the remaining part of the proof in two cases:

(i) There are infinitely many indexes $l \in \NN$ such $\theta_l \leq 0$.

(ii) There are infinitely many indexes $l \in \NN$ such $\theta_l \geq 0$.

\noindent Let us prove (i). After passing to a subsequence,  we can assume $\theta_l \leq 0$ for any $l \in \NN$. So, since $\theta_l \leq 0$, for any $(x,t) \in Q_1$ we have $(x-\xi_l,t-\theta_l) \in Q_2 := B_2\times(0,4)$ and, furthermore,
\[
\begin{aligned}
\|\vartheta_{p_0}(\cdot,\cdot) - \vartheta_{p_0}(\cdot - \xi_0,\cdot - \theta_0)\|_{L^{\infty}(Q_1)} & \leq \|\vartheta_{p_0} - \overline{u}_{p_l,\delta_l}\|_{L^{\infty}(Q_1)} \\
& + \|\overline{u}_{p_l,\delta_l}(\cdot,\cdot) - \overline{u}_{p_{0l},\delta_l}(\cdot - \xi_0,\cdot - \theta_0)\|_{L^{\infty}(Q_1)} \\
& + \|\overline{u}_{p_{0l},\delta_l}(\cdot - \xi_0,\cdot - \theta_0) - \vartheta_{p_0}(\cdot - \xi_0,\cdot - \theta_0)\|_{L^{\infty}(Q_1)} \to 0,
\end{aligned}
\]
as $l \to +\infty$, thanks to \eqref{eq:LINFINITYCONVERGECESEQUENCESBLOWUP} and \eqref{eq:LINFINITYCONVERGECESEQUENCESBLOWUP2}. Consequently, using the real analyticity of the polynomial $\vartheta_{p_0} = \vartheta_{p_0}(x,t)$, it follows
\[
\vartheta_{p_0}(x + \xi_0,t + \theta_0) = \vartheta_{p_0}(x,t) \quad \text{ for all } (x,t) \in \RR^N\times\RR,
\]
and so
\[
\partial_x^{\alpha}\partial_t^j\vartheta_{p_0}(\xi_0,\theta_0) = \partial_x^{\alpha}\partial_t^j\vartheta_{p_0}(0,0) = 0 \quad \text{ for all } |\alpha| + 2j < 2\kappa.
\]
On the other hand, multiplying \eqref{eq:REDUCTIONTOABSURBUMASSWHITNEY} by $\delta_l^{-2\kappa + |\alpha| + 2j}$, using the homogeneity of $\partial_x^{\alpha}\partial_t^j \vartheta_{p_0}$ and taking the limit as $l \to +\infty$, we obtain
\[
|\partial_x^{\alpha}\partial_t^j\vartheta_{p_0}(\xi_0,\theta_0)| = \left|\sum_{\substack{(\beta,i) \geq (\alpha,j): \\ |\gamma| + 2i = 2\kappa}} \frac{a_{\beta,i}(x_0,t_0)}{(\beta-\alpha)!(i-j)!} \xi_0^{\beta-\alpha}\theta_0^{i-j}\right| \geq \underline{\sigma} > 0,
\]
in contradiction with the computation above. The proof of case (ii) is almost identical to the previous one and we skip it (we just have to using the second convergence in \eqref{eq:LINFINITYCONVERGECESEQUENCESBLOWUP2} instead of the first one).

As a consequence, we deduce that the assumptions of Whitney's theorem are satisfied and so, there is a smooth function $F$ defined in the whole $\RR^N\times\RR$ such that
\[
\partial_x^{\alpha}\partial_t^j F = f_{\alpha,j} \quad \text{ in } K,
\]
for all $|\alpha| + 2j \leq 2 \kappa$.

\smallskip

\emph{Step2: Implicit function theorem.} Let $(x_0,t_0) \in \Gamma_{2\kappa}^d(u) \cap K$ and $d = 0,\ldots,N$. As in \cite[Theorem 12.12]{DainelliGarofaloPetTo2017:book} we consider two subcases.

Let us begin by assuming $d = 0,\ldots,N-1$. For these choices of the dimension, we have that there are multi-indexes $\alpha_i$ and nonnegative integers $j_i$ with $|\alpha_i| + j_i = 2\kappa-1$ such that
\[
v_i = \nabla_x \partial_x^{\alpha_i} \partial_t^{j_i} \vartheta_{p_0} = \nabla_x \partial_x^{\alpha_i} \partial_t^{j_i} F(x_0,t_0) \quad \text{ for } \; i = 1,\ldots,N-d,
\]
are linearly independent vectors. At the same time, in view of Theorem \ref{CONTINUOUSDEPOFBLOWUPLIMITS} (part (i)) and the Whitney's extension theorem, we have
\[
\Gamma_{2\kappa}^d(u) \cap K \subset \bigcap_{i = 1}^{N-d} \big\{(x,t) \in \RR^N\times\RR : \partial_x^{\alpha_i} \partial_t^{j_i} F(x,t) = 0 \big\},
\]
and so, using the linear independence of the $v_i$'s and the implicit function theorem, we immediately conclude that $\Gamma_{2\kappa}^d(u) \cap K$ is contained in a $(d+1)$-dimensional space-like $\mathcal{C}^{1,0}$ manifold (see  Definition \ref{DEFSPACELIKETIMELIKE}). Finally, recalling that we have chosen $K = E_h$ (for some arbitrary $h \in \NN$) and since $\Gamma_{2\kappa}(u) = \cup_{h\in\NN} E_h$ (see  Lemma \ref{LEMMAGAMMAKISFSIGMA}), the proof in the case $d = 0,\ldots,N-1$ is ended.

Assume now $d = N$ and $(x_0,t_0) \in \Gamma_{2\kappa}^N(u) \cap K$. In view of Lemma \ref{LEMTIMELIKESINGPOINTS}, we can assume that $\kappa$ is an integer. From the same lemma, we deduce
\[
\partial_t^{\kappa} F(x_0,t_0) = \partial_t^{\kappa} \vartheta_{p_0} \not = 0,
\]
while
\[
\Gamma_{2\kappa}^N(u) \cap K \subset \big\{(x,t) \in \RR^N\times\RR : \partial_t^{\kappa - 1} F(x,t) = 0 \big\},
\]
and so, for the Implicit function theorem we obtain that $\Gamma_{2\kappa}^N \cap K$ is contained in a $N$-dimensional time-like $\mathcal{C}^1$ manifold. The proof is then completed. $\Box$

\section{Appendix}\label{Appendix}
In this final section, we study \eqref{eq:DEFINITIONPARABOLICEQUATION1}:
\[
\partial_t U + |y|^{-a}\nabla\cdot( |y|^a \nabla U) = 0 \quad \text{in } \RR^{N+1}\times(0,1).
\]
As mentioned in Section \ref{Section:BlowUp2}, the notion of weak solutions to \eqref{eq:DEFINITIONPARABOLICEQUATION1} is very similar to the one given in Definition \ref{def:WeakSolutions} and, under the additional assumption on the vanishing order of $U$ at $O$, we can prove the monotonicity of an Almgren-Poon quotient $r \to N_0(U,r)$ (see  Theorem \ref{Theorem:MonotonicityNewPoon} and Remark \ref{Rem:GrowthCond} for the definition of $N_0$ and the main results).
\begin{rem}
	Notice that the nodal properties, blow-up analysis and unique continuation type properties of solutions to \eqref{eq:DEFINITIONPARABOLICEQUATION1} are well known when $y$ is far away from the set $\Sigma = \{(X,t): y = 0\}$. We will thus focus on the local behaviour of solutions near $\Sigma$ by proving a complete blow-up analysis for weak solutions vanishing of finite order. Further properties (dimensional bounds on the nodal set, stratification results and so on) can be also investigated but, since this is a complementary section, we do not enter into details.
\end{rem}
\begin{rem}
We mention that, as in the above sections, the next results are completely local even if we are considering solutions defined in the whole space $\RR^{N+1}$. In other words, we could have considered weak solutions to
\[
\partial_t U + |y|^{-a}\nabla\cdot( |y|^a \nabla U) = 0 \quad \text{in } \BB_1\times(0,1),
\]
and then studying the local properties of $W = \zeta U$, where $\zeta$ is the cut-off defined in \eqref{eq:PropCutOff}, which satisfies \eqref{eq:DEFINITIONPARABOLICEQUATION1} with the right hand side $F := 2 \nabla \zeta \cdot \nabla U + \zeta |y|^{-a} \nabla \cdot (|y|^a \nabla\zeta) U$ (see  \eqref{eq:DefRightHS}). Since the main ideas have been treated above, we have opted for a less general setting.
\end{rem}

We begin by proving the corresponding version of Theorem \ref{THEOREMBLOWUP1New} for solutions to \eqref{eq:DEFINITIONPARABOLICEQUATION1}.

\begin{thm}\label{THEOREMBLOWUP1}
Let $a \in (-1,1)$, $p_0 = (x_0,0,t_0) \in \Gamma(U)\cap\Sigma$, and let $U$ be a nontrivial weak solution to \eqref{eq:DEFINITIONPARABOLICEQUATION1}, having finite vanishing order at $O$. Then there exist $n_0,m_0 \in \NN$ such that the following assertions hold:

(i) The Almgren-Poon quotient $N_0(U_{p_0},r)$ (see  Theorem \ref{Theorem:MonotonicityNewPoon}) satisfies
\[
\lim_{r \to 0^+} N_0(U_{p_0},r) = \kappa_{n_0,m_0} = \kk,
\]
where the admissible values for $\kappa_{n,m}$ are
\[
\kappa_{n,m} = \widetilde{\kappa}_{n,m} := \frac{n}{2} + m \qquad \text{ and } \qquad \kappa_{n,m} = \widehat{\kappa}_{n,m} := \frac{n}{2} + m + \frac{1-a}{2},
\]
are the eigenvalues of problem \eqref{eq:ORNSTUHLENFOREXTENSION1}, for any $m,n \in \NN$.

(ii) The vanishing order of $U$ at $p_0$ is $2\kk$ and, furthermore,
\begin{equation}\label{eq:LIMITOFTHEQUOTIENTFORSMALLTSTATEMENT1}
\lim_{r \to 0^+} \frac{H(U_{p_0},r)}{r^{4\kappa}} = L_0,
\end{equation}
for some constant $L_0 > 0$ (depending on $p_0$ and $U$), where the function $r \to H(U_{p_0},r)$ is the denominator of the Almgren-Poon quotient.

(iii) There exists a parabolically homogeneous polynomial of degree $2\kappa_{n_0,m_0}$ given by the expression
\[
\Theta_{p_0}(X,t) := t^{\kappa_{n_0,m_0}} \sum_{(\alpha,m) \in J_0} v_{\alpha,m} \overline{V}_{\alpha,m}\left(\frac{X}{\sqrt{t}}\right),
\]
such that, for all $T_{\ast} > 0$, it holds
\[
W_{p_0,r} \to \Theta_{p_0} \quad \text{ in } L^2(0,T_{\ast};H^1(\RR_+^{N+1},d\mu_t))\cap\mathcal{C}^0(0,T_{\ast};L^2(\RR_+^{N+1},d\mu_t)),
\]
as $r \to 0^+$, where $v_{\alpha,m}$ are suitable constants, the sum is done over the set of indexes
\[
J_0 := \left\{\right (\alpha,m) \in \ZZ_{\geq 0}^N\times\NN : |\alpha| = n \in \NN \text{ and } \kappa_{n,m} = \kappa_{n_0,m_0} \},
\]
and the integration probability measure is defined in \eqref{eq:EXTENDEDGAUSSIANMEASURE1}. Moreover,
\[
\overline{V}_{\alpha,m}(x,y) := \frac{V_{\alpha,m}(x,y)}{\|V_{\alpha,m}\|_{L_{\mu}^2}}, \qquad \alpha \in \ZZ_{\geq 0}^N, \; m \in \NN,
\]
are the normalized versions of the eigenfunctions $V_{\alpha,m} = V_{\alpha,m}(x,y)$ to problem \eqref{eq:ORNSTUHLENFOREXTENSION1} corresponding to the eigenvalue $\kappa_{n,m}$ and defined in the statement of Theorem \ref{SPECTRALTHEOREMEXTENDEDORNUHL1}.
\end{thm}
\emph{Skecth of the proof.} We fix $T_\ast \geq 1$, $T = \sqrt{T_{\ast}}$ and $p_0 = O$. We consider a weak solution $U$ to \eqref{eq:DEFINITIONPARABOLICEQUATION1} and the blow-up family $U_r$ defined as in \eqref{eq:BLOWUPPLUSTRANSSEQUENCE}. Notice that in this case it is still a solution to \eqref{eq:DEFINITIONPARABOLICEQUATION1}.

Then we pass to its re-scaled version
\[
\widetilde{U}_r(X,t) := U_r(\sqrt{t}X,t),
\]
which satisfies
\[
\begin{aligned}
t \int_{\mathbb{R}^{N+1}} \partial_t \widetilde{U}_r \eta \,d\mu = \int_{\mathbb{R}^{N+1}}  \nabla \widetilde{W}_r \cdot \nabla \eta \,d\mu = 0,
\end{aligned}
\]
for all $\eta \in L^2(0,1/r^2;H_{\mu}^1)$ and for all $t \in (0,1/r^2)$, see  \eqref{eq:RescaledBlowUp}.

\smallskip

\emph{Step 1: Uniform estimates in Gaussian spaces.} Proceeding as above by using the monotonicity of the quotient $r \to N_0(U,r)$, we find that for any fixed $t_\ast \in (0,T_\ast)$, $\{\widetilde{U}_r\}_{r \in (0,r_0)}$ is uniformly bounded in $L^2(0,T_{\ast};H^1_{\mu})$, while $\{\partial_t\widetilde{U}_r\}_{r \in (0,r_0)}$ is uniformly bounded in $L^2(t_\ast,T_\ast;L^2_\mu)$ and so, applying Lemma \ref{SIMONLEMMARECALL}, it follows that
\[
\{\widetilde{W}_r\}_{r \in (0,r_0)} \subset \mathcal{C}^0(t_\ast,T_\ast; L^2_\mu)
\]
is relatively compact.

\smallskip

\emph{Step 2: Compactness and properties of the limit.} So, for any $r_n \to 0^+$ and any fixed $0 < t_{\ast} < T_{\ast}$, we can extract a sub-sequence $r_{n_j} \to 0^+$ (that we rename $r_j := r_{n_j}$ by convenience) such that
\[
\widetilde{W}_{r_j} \to \widetilde{\Theta} \quad \text{ in } \mathcal{C}^0\left(t_{\ast},T_{\ast};L_{\mu}^2\right), \quad \text{as } j \to +\infty,
\]
where $\widetilde{\Theta} \in \bigcap_{t_{\ast} \in (0,T_{\ast})} \mathcal{C}^0 \left(t_{\ast},T_{\ast}; L_{\mu}^2 \right)$ and it satisfies also
\[
\begin{aligned}
\int_{t_1}^{t_2} t \int_{\mathbb{R}^{N+1}} \partial_t \widetilde{\Theta} \eta \,d\mu dt = \int_{t_1}^{t_2} \int_{\mathbb{R}^{N+1}}  \nabla \widetilde{\Theta} \cdot \nabla \eta \,d\mu dt,
\end{aligned}
\]
for all $\eta \in L^2(t_1,t_2,;H_{\mu}^1)$ and for all $t_\ast \leq t_1 \leq t_2 \leq T_\ast$, that is $\widetilde{\Theta}$ is a global weak solution to
\[
t\partial_t \widetilde{\Theta} + \mathcal{O}_a\widetilde{\Theta} = 0 \quad \text{ in } \mathbb{R}^{N+1} \times (0,\infty).
\]
Then, exactly as before, it is possible to prove that $\widetilde{U}_{r_j} \to \widetilde{\Theta}$ strongly in $L^2\left(t_{\ast},T_{\ast};H_{\mu}^1\right)$, as $j \to +\infty.$

\smallskip

\emph{Step 3: Improved convergence.} Exactly as before we obtain that $\widetilde{U}_{r_j} \to \widetilde{\Theta}$ strongly in $L^2\left(0,T_{\ast};H_{\mu}^1\right)$, as $j \to +\infty$. In this setting it is also possible to prove that
\[
\widetilde{U}_{r_j} \to \widetilde{\Theta} \quad \text{ in } C^0\left(0,T_{\ast};L_{\mu}^2\right)
\]
as $j \to +\infty$. This easily follows from the fact that in this case, the function $t \to h(U,t)$ is monotone nondecreasing (see  \ref{Lemma:HeightsVariations}).

\smallskip

\emph{Step 4: The limit $\Theta$ is a re-scaled eigenfunction.} This step is almost identical and we skip it.

\smallskip

\emph{Step 5/Step 6: Uniqueness of the blow-up.} The proof of this step is easier since in this case the terms involving $F$ and $q$ are zero. $\Box$

\begin{ex}
Again, we report some specific examples. Since the blow-up profiles corresponding to $\widehat{\kappa}_{n,m} = \frac{n}{2} + m + \frac{1-a}{2} = \widetilde{\kappa}_{n,m} + \frac{1-a}{2}$ satisfy:
\[
\begin{aligned}
\widehat{\Theta}_{\alpha,m}(x,y,t) &= t^{\widehat{\kappa}_{n,m}} V_{\alpha,m}\left(\frac{x}{\sqrt{t}},\frac{y}{\sqrt{t}}\right) = t^{\widehat{\kappa}_{n,m}} H_n\left(\frac{x}{\sqrt{t}}\right) \,\frac{y}{\sqrt{t}}\left|\frac{y}{\sqrt{t}}\right|^{-a} L_{(\frac{1-a}{2}),m}\left(\frac{y^2}{4t}\right) \\
& = t^{\widetilde{\kappa}_{n,m}} y|y|^{-a} H_n\left(\frac{x}{\sqrt{t}}\right)L_{(\frac{1-a}{2}),m}\left(\frac{y^2}{4t}\right),
\end{aligned}
\]
we easily see that
\[
\begin{aligned}
&\widehat{\kappa}_{0,0} = \frac{1-a}{2}          \qquad\qquad          \widehat{\Theta}_{0,0}(x,y,t) = y|y|^{-a}\\
&\widehat{\kappa}_{1,0} = \frac{2-a}{2}    \qquad\qquad         \widehat{\Theta}_{1,0}(x,y,t) = xy|y|^{-a} \\
&\widehat{\kappa}_{2,0} = \frac{3-a}{2}         \qquad\qquad        \widehat{\Theta}_{2,0}(x,y,t) = y|y|^{-a}(x^2-2t) \\
&\widehat{\kappa}_{0,1} = \frac{3-a}{2}          \qquad\qquad           \widehat{\Theta}_{0,1}(x,y,t) = y|y|^{-a}\left[\left(\frac{3-a}{2}\right)t - \frac{y^2}{4}\right],
\end{aligned}
\]
and so on. Consequently, the blow-ups of solutions to \eqref{eq:DEFINITIONPARABOLICEQUATION1} are obtained by combining all the possible blow-ups above and the ones corresponding to solutions to \eqref{eq:DEFINITIONPARABOLICEQUATION}.
\end{ex}
As above we show that also Theorem \ref{CONVERGENCELINFINITYLOC} has a counter part for solutions to \eqref{eq:DEFINITIONPARABOLICEQUATION1}. Notice that it generalizes the regularity results proved in \cite[Appendix A]{BanGarDanPetr2019:art}.
\begin{thm}\label{CONVERGENCELINFINITYLOC1}
Let $a \in (-1,1)$ and let $U$ be a weak solution to \eqref{eq:DEFINITIONPARABOLICEQUATION1}. Assume that $U$ is bounded in $\mathbb{Q}_1 := \mathbb{B}_1\times(0,1)$. Then for any $\nu \in (0,\nu_\ast)$, there exists a constant $C = C(\nu) > 0$ such that
\[
\|U\|_{\mathcal{C}^{2\nu,\nu}(\overline{\mathbb{Q}_{1/2}})} \leq C \|U\|_{L^{\infty}(\mathbb{Q}_1)},
\]
where we recall that $\nu_{\ast} = \frac{1}{2}\min\{1,1-a\}$.
\end{thm}
\emph{Skecth of the proof.} As in the proof of Theorem \ref{CONVERGENCELINFINITYLOC}, we consider a locally bounded solution $U$ and we show that for any $\nu \in (0,\nu_{\ast})$,
\[
[\eta U]_{\mathcal{C}^{2\nu,\nu}(\QQ_1)} \leq C,
\]
where $\QQ_1 := \BB_1\times(0,1)$ and $\eta = \eta(X,t)$ is a smooth function satisfying
\[
\begin{cases}
\eta(X,t) = 1 \quad          &\text{ for } (X,t) \in \QQ_{1/2} \\
0 < \eta(X,t) \leq 1 \quad   &\text{ for } (X,t) \in \QQ_1 \setminus \QQ_{1/2} \\
\eta(X,t) = 0 \quad          &\text{ for } (X,t) \in \partial \QQ_1,
\end{cases}
\]
where $\partial \QQ_1 := [\partial \BB_1 \times (0,1)] \cup [\BB_1\times\{1\}]$.

\smallskip

%\color{blue}
\emph{Step 0: Approximation.} This part of the proof works as in the proof of Theorem \ref{CONVERGENCELINFINITYLOC}. Given a solution $U$ to \eqref{eq:DEFINITIONPARABOLICEQUATION1}, we define $\overline{v}_j := \varrho_j \star U$, where $\{\varrho_j\}_{j \in \mathbb{N}}$ is a family of mollifiers in the variables $(x,t) \in \mathbb{R}^{N+1}$. Then, $\overline{v}_j \in L^{\infty}(\mathbb{Q}_1^+)$ uniformly w.r.t. $j \in \mathbb{N}$, it is smooth w.r.t. to $x$ and $t$, satisfies \eqref{eq:DEFINITIONPARABOLICEQUATION1} in the weak sense, and $\overline{v}_j \to U$ in $L^\infty(\mathbb{Q}_1^+)$ as $j \to +\infty$.

Then, we extend the trace $v_j(x,t) := \overline{v}_j(x,0,t)$ to a function $\overline{u}_j$ which is smooth in $\mathbb{R}^{N+1}_+\times(0,1)$ using \eqref{eq:DEFOFEXTENDEDVERSIONOFU}, and satisfies \eqref{eq:DEFINITIONPARABOLICEQUATION1} in $\mathbb{R}^{N+1}_+\times(0,1)$. Moreover, from the explicit expression of the Poisson kernel $P_y^a$ (see \eqref{eq:POISSONKERNELWITHANOTATION}), it follows $\overline{u}_j \in \mathcal{C}_{loc}^{2\nu,\nu}(\overline{\mathbb{R}^{N+1}_+} \times (0,1))$ for all $\nu \in (0,\nu_\ast)$ and all $j \in \mathbb{N}$. Consequently, the function $\overline{w}_j := \overline{u}_j - \overline{v}_j$, satisfies \eqref{eq:DEFINITIONPARABOLICEQUATION1} up to an odd reflection w.r.t. $y$.

Repeating the procedure of Theorem \ref{CONVERGENCELINFINITYLOC}, using that $\partial_t \overline{w}_j \in L^{\infty}(\mathbb{B}_1\times(0,1))$ and the elliptic regularity estimates in \cite{SireTerVita2020:art}, it follows that $\overline{w}_j$ is $\nu$-H\"older continuous w.r.t. $y$ for all $\nu \in (0,\nu_\ast)$ and so, $\overline{v}_j \in \mathcal{C}^{2\nu,\nu}(\overline{\mathbb{Q}_{1/2}^+})$, for all $j \in \mathbb{N}$, $\nu \in (0,\nu_\ast)$.

Again, since $\overline{v}_j \to U$ in $L^\infty(\mathbb{Q}_1^+)$, it is enough to show that \eqref{eq:UNIFORMHOLDERBOUNDBLOWUPSEQUENCETEST} holds for each $\overline{v}_j$, for some constant $C > 0$ independent of $j \in \mathbb{N}$ and then passing to the limit as $j \to +\infty$.
\normalcolor

\smallskip

\emph{Step 1.} As above we assume by contradiction that there exists $\nu \in (0,\nu_{\ast})$, a sequence $\{\overline{v}_j\}_{j \in \mathbb{N}}$ and a sequence $(X_{1,j},t_{1,j}),(X_{2,j},t_{2,j}) \in \QQ_1$ such that
\[
[\eta \overline{v}_j]_{\mathcal{C}^{2\nu,\nu}(\overline{\mathbb{Q}_1})} := \frac{|\eta(X_{1,j},t_{1,j}) \overline{v}_j(X_{1,j},t_{1,j}) - \eta(X_{2,j},t_{2,j}) \overline{v}_j(X_{2,j},t_{2,j})|}{r_j^{2\nu}} := L_j \to +\infty,
\]
as $j \to +\infty$ where we have defined
\[
r_j := (|X_{1,j}-X_{2,j}|^2 + |t_{1,j}-t_{2,j}|)^{\frac{1}{2}},
\]
for all $j \in \NN$. Exactly as above, we obtain
\begin{equation}%\label{eq:SUMOFDISTANCESIMPORTANTSEQUENCE}
\frac{\text{dist}((X_{1,j},t_{1,j}),\partial \QQ_1)}{r_j} + \frac{\text{dist}((X_{2,j},t_{2,j}),\partial \QQ_1)}{r_j} \to + \infty,
\end{equation}
as $j \to + \infty$.

\emph{Step2.} We introduce the sequences:
\[
\begin{aligned}
U_j(X,t) &:= \eta(\widehat{X}_j,\widehat{t}_j)\frac{\overline{v}_j(\widehat{X}_j + r_jX,\widehat{t}_j + r_j^2t)}{L_jr_j^{2\nu}}, \\
\overline{U}_j(X,t) &:= \frac{(\eta \overline{v}_j)(\widehat{X}_j + r_jX,\widehat{t}_j + r_j^2t)}{L_jr_j^{2\nu}},
\end{aligned}
\]
where $\widehat{P}_j = (\widehat{X}_j,\widehat{t}_j) \in \QQ_1$ have to be chosen and
\[
(X,t) \in \widehat{\mathbb{Q}}_j := (\QQ_1 - \widehat{P}_j)/r_j = \BB_{1/r_j}(\widehat{X}_j) \times I_j,
\]
for all $j \in \NN$, where we have set $I_j := (-\widehat{t}_j/r_j^2,(1-\widehat{t}_j)/r_j^2)$. Again we can show:
\[
\left[ \overline{U}_j \right]_{\mathcal{C}^{2\nu,\nu}(\widehat{\QQ}_j)} = 1,
\]
for all $j \in \NN$ and
\[
\partial_t U_j + \mathcal{L}_a^jU_j = 0 \quad \text{ in } \RR^{N+1} \times I_j
\]
in the weak sense, where
\[
\mathcal{L}_a^jU_j = \left|\widehat{y}_jr_j^{-1} + y\right|^{-a} \nabla\cdot  \left(\left|\widehat{y}_jr_j^{-1} + y\right|^a \nabla U_j \right),
\]
and $W_j(O) = \overline{W}_j(O)$ by definition.

As in the proof of Theorem \ref{CONVERGENCELINFINITYLOC} the two sequences are asymptotically equivalent on compact sets of $\RR^{N+1}\times\RR$, i.e.
\[
\|U_j - \overline{U}_j\|_{L^{\infty}(K)} \to 0,
\]
as $j \to +\infty$. Furthermore, using the H\"older regularity of $\overline{U}_j$ and that  $\overline{U}_j(O) = U_j(O)$, we deduce also the existence of a constant $C > 0$ (depending on the compact set $K \subset \RR^{N+1}\times\RR$) such that
\[
\sup_{(X,t) \in K } |U_j(X,t) - U_j(O)| \leq C.
\]

\emph{Step3.} As above we show that there exists $C > 0$ such that
\[
\frac{\text{dist}((X_{1,j},t_{1,j}), \QQ_1 \cap \Sigma)}{r_j} + \frac{\text{dist}((X_{2,j},t_{2,j}), \QQ_1 \cap \Sigma)}{r_j} \leq C,
\]
for $j \in \NN$ large enough. We assume by contradiction that
\[
\frac{\text{dist}((X_{1,j},t_{1,j}), \QQ_1 \cap \Sigma)}{r_j} + \frac{\text{dist}((X_{2,j},t_{2,j}), \QQ_1 \cap \Sigma)}{r_j} \to +\infty,
\]
as $j \to +\infty$ and we take $(\widehat{X}_j,\widehat{t}_j) = (X_{1,j},t_{1,j})$ in the definition of $U_j$ and $\overline{U}_j$. Using the properties of $(X_{1,j},t_{1,j})$ we obtain that
\[
\bigcup_{j\in\mathbb{N}} \RR^{N+1} \times I_j
\]
equals either $\RR^{N+1}\times\RR$ or $\RR^{N+1}\times(0,+\infty)$. Again, we assume to be in the case $\RR^{N+1}\times\RR$ (the other can be treated in very similar way). As before, we consider the new sequences
\[
\begin{aligned}
V_j(X,t) &:= U_j(X,t) - U_j(O), \\
\overline{V}_j(X,t) &:= \overline{U}_j(X,t) - \overline{U}_j(O).
\end{aligned}
\]
Now, since the sequence $\{\overline{V}_j\}_{j\in\NN}$ is locally uniformly bounded with uniformly bounded $\nu$-H\"{o}lder semi-norm, we obtain the existence of a continuous function $V \in \mathcal{C}_{loc}$, uniform limit of $\overline{V}_j$ as $j \to +\infty$, thanks to the Ascoli-Arzel\`a theorem. Exactly as before, $V_j \to V$ uniformly on compact sets of $\RR^{N+1}\times\RR$, too and $V \in \mathcal{C}^{2\nu,\nu}(\RR^{N+1}\times\RR)$, $\nu \in (0,\nu_{\ast})$. In particular,
\[
|V(X,t)| \leq C[1 + (|X|^2 + |t|)^{\nu}] \leq C[1+ d^{2\nu}(X,t)]
\]
for some $C > 0$ and all $(X,t) \in \RR^{N+1}\times\RR$. Furthermore, testing the equation of $V_j$ with a $\eta \in \mathcal{C}_0^{\infty}(\RR^{N+1}\times\RR)$ and integrating by parts, it follows
\[
\int_{t_1}^{t_2} \int_{\RR^{N+1}} \left[ |1 + r_j\widehat{y}_j^{-1}y|^a\partial_t\eta - \nabla\cdot\left(|1 + r_j\widehat{y}_j^{-1}y|^a\nabla \eta \right)\right] V_j \,dXdt = 0,
\]
for all $j \in \NN$ large enough. Hence, exploiting the fact that $|1 + r_j\widehat{y}_j^{-1}y|^a \to 1$ on any compact set of $\RR^{N+1}\times\RR$ ($r_j\widehat{y}_j^{-1} \to 0$ as before) and using the uniform convergence of the sequence $V_j$, we pass to the limit as $j \to +\infty$ in the above relation to conclude
\[
\int_{\RR} \int_{\RR^{N+1}}\left(\partial_t\eta - \Delta_{x,y} \eta \right) V \,dXdt = 0,
\]
for any $\eta \in \mathcal{C}_0^{\infty}(\RR^{N+1}\times\RR)$. From the Liouville theorem \ref{1LIOUVILLETYPETHEOREMFORWEIGHTEDHE}, it follows that $V$ is constant, which is not possible (this can be seen exactly as above).

\smallskip

\emph{Step4.} Thanks to \emph{Step 3}, we can choose $\widehat{P}_j = (x_{1,j},0,t_{1,j})$ (see  \cite{TerVerZil2017:art,TerVerZil2016:art}). Consequently, $V_j$ satisfies
\[
\int_{\RR} \int_{\RR^{N+1}} \left[|y|^a\partial_t\eta - \nabla\cdot\left( |y|^a\nabla \eta \right)\right] V_j \,dXdt = 0
\]
for all test functions $\eta$ with compact support in $\RR^{N+1} \times \RR$ and all $j \in \NN$. Taking the limit as $j \to + \infty$, it follows
\[
\int_{\RR} \int_{\RR^{N+1}} \left[|y|^a\partial_t\eta - \nabla\cdot\left( |y|^a\nabla \eta \right)\right] V_j \,dXdt = 0
\]
for all test functions $\eta$ with compact support in $\RR^{N+1} \times \RR$, which means that $V$ is a global solution to \eqref{eq:DEFINITIONPARABOLICEQUATION1}. Since it satisfies the grow condition with $\nu \in (0,\nu_{\ast})$, it must be constant in view of Lemma \ref{1LIOUVILLETYPETHEOREMFORWEIGHTEDHE} (part (iii)). Proceeding exactly as in \eqref{eq:WISNONCONSTANT}, it follows that this is not possible, obtaining the desired contradiction. $\Box$
\begin{cor}\label{Corollary:THMECONVERGENCELINFINITYLOC1}
	Let $a \in (-1,1)$, $p_0 \in \RR^N\times\{0\}\times(0,1)$ and let $U$ be a weak solution to \eqref{eq:DEFINITIONPARABOLICEQUATION1}, having finite vanishing order at $O$. Let $\{U_{p_0,r}\}_{r \in (0,1)}$ the blow-up family defined in \eqref{eq:BLOWUPPLUSTRANSSEQUENCE}. Then for any compact set $K \subset \RR^{N+1}\times(0,\infty)$ and any $\nu \in (0,\nu_{\ast})$, there exist  constants $C,r_0 > 0$ depending on $K$ and $\nu$ such that
	\[
	\|U_{p_0,r}\|_{\mathcal{C}^{2\nu,\nu}(K)} \leq C.
	\]
	In particular, for any $\nu \in (0,\nu_{\ast})$, it holds
	\[
	U_{p_0,r} \to \Theta_{p_0} \quad \text{ in } \mathcal{C}^{2\nu,\nu}_{loc},
	\]
	as $r \to 0$, where $\Theta_{p_0}$ is the blow-up limit of $U$ at $p_0$ found in Theorem \ref{THEOREMBLOWUP1} part (iii).
\end{cor}
\emph{Proof.} The proof follows by repeating the arguments of Corollary \ref{Corollary:THMECONVERGENCELINFINITYLOC} and Corollary \ref{Corollary:THMECONVERGENCELINFINITYLOCUNIQUENESS}. $\Box$
%
%
%
%
%

%%%%%%%%%%%%%%%%%%%%%%%%%%%%%%%%%%%%%%%%%%%%%%%%%%%%%%%%%%%%%%%%%%%%%%%%%%%%%%%%%%%%%%%%%%%%%%%%%%%%%

%%%%%%%%%%%%%%%%%%%%%%%%%%%%%%%%%%%%%%%%%%%%%%%%%%%%%%%%%%%%%%%%%%%%%%%%%%%%%%%%%%%%%%%%%%%%%%%%%%%%%
%
%
%
%
%

\bigskip

%
%
%
%
%
%%%%%%%%%%%%%%%%%%%%%%%%%%%%%%%%%%%%%%%%%%%%%%%%%%%%%%%%%%%%%%%%%%%%%%%%%%%%%%%%%%%%%%%%%%%%%%%%%%%%%%%%%%%%
%
%


\begin{thebibliography}{99}

\bibitem{AbramowitzStegun1972:art}
{\sc M. Abramowitz, I. A. Stegun.} ``Handbook of mathematical functions, with formulas, graphs, and mathematical tables'', National Bureau of Standards Applied Mathematics Series, \textbf{55}, 1972.
%

\bibitem{AltCafFried1984:art}
{\sc H. W. Alt, L. Caffarelli, A. Friedman.} \emph{Variational problems with two phases and their free boundaries}, Trans. Amer. Math. Soc. \textbf{282} (1984), 431--461.
% Monotonicity Formula 2 species

\bibitem{AthCaffaMilakis2016:art}
{\sc I. Athanasopoulos, L. Caffarelli, E. Milakis.} \emph{On the regularity of the Non-dynamic Parabolic Fractional Obstacle Problem}, J. Differential Equations \textbf{265} (2018), 2614--2647.
%

\bibitem{AA2018:book}
{\sc A. Audrito}. ``Nonlinear and nonlocal equations. Qualitative theory and asymptotic behaviour'', Ph.D. dissertation (2018).
%

\bibitem{BaeumerEtAl2005:art}
{\sc B. Baeumer, M. M. Meerschaert, J. Mortensen.} \emph{Space-time fractional derivative operators}, Proc. Amer. Math. Soc. \textbf{133}, 2273--2282.
% CTRW

\bibitem{Balakrishnan1960:art}
{\sc A. V. Balakrishnan}, \emph{Fractional powers of closed operators and the semigroups generated by them}, Pacific J. Math. \textbf{10} (1960), 419--437.
%

\bibitem{BanGarofalo2017:art}
{\sc A. Banerjee, N. Garofalo.} \emph{Monotonicity of generalized frequencies and the strong unique continuation property for fractional parabolic equations}, Adv. Math., \textbf{336} (2018), 149--241.
%

\bibitem{BanGarDanPetr2018:art}
{\sc A. Banerjee, N. Garofalo, D. Danielli, A. Petrosyan.} \emph{The structure of the singular set in the thin obstacle problem for degenerate parabolic equations}, ArXiv preprint, arXiv:1902.07457 (2019).
%

\bibitem{BanGarDanPetr2019:art}
{\sc A. Banerjee, N. Garofalo, D. Danielli, A. Petrosyan.} \emph{The regular free boundary in the thin obstacle problem for degenerate parabolic equations}, ArXiv preprint, arXiv:1906.06885 (2019).
%

\bibitem{Baouendi1967:art}
{\sc M. Baouendi.} \emph{Sur une classe d'op\'erateurs elliptiques deg\'n\'er\'es}, Bull. Soc Math. France \textbf{95} (1967), 45--87.
%

\bibitem{BarriosFigalliRosOton2016:art}
{\sc B. Barrios, A. Figalli, X. Ros-Oton.} \emph{Free boundary regularity in the parabolic fractional obstacle problem}, Comm. Pure Appl. Math. \textbf{71} (2018), 2129--2159.
%

\bibitem{BatheRoberto2003:art}
{\sc F. Bathe, C. Roberto.} \emph{Sobolev inequalities for probability measures on the real line}, Studia Math.\textbf{159} (2003), 481--497.
%

\bibitem{Beckner1989:art}
{\sc W. Beckner.} \emph{A Generalized Poincar\'e Inequality for Gaussian Measures}, Proc. Amer. Math. Soc. \textbf{105} (1989), 397--400.
%


\bibitem{Bell2017:art}
{\sc J. Bell.} \emph{Gaussian measures, Hermite polynomials, and the Ornstein-Uhlenbeck semigroup}, online, 2015.
%

%\bibitem{BernardisMRStingaTorrea2013:art}
%{\sc A. Bernardis, F.-J. Mart\'in-Reyes, P. R. Stinga, J. L. Torrea}. \emph{Maximum principles, extension problem and inversion for nonlocal one-sided equations}, J. Differential Equations \textbf{260} (2016), 6333--6362.
%

\bibitem{BonSirVaz2016:art}
{\sc M. Bonforte, Y. Sire, J. L. V\'azquez.} \emph{Optimal Existence and Uniqueness Theory for the Fractional Heat Equation}, Nonlinear Anal. \textbf{157} (2017), 142--168.
% Different fractional laplacians

\bibitem{CafFigalli2013:art}
{\sc  L. A. Caffarelli, A. Figalli.} \emph{Regularity of solutions to the parabolic fractional obstacle problem}, J. Reine Angew. Math. \textbf{680} (2013), 191--233.
%

\bibitem{CafKen1998:art}
{\sc  L. A. Caffarelli, C. E. Kenig.} \emph{Gradient estimates for varible coefficient parabolic equations and singular perturbation problems}, Amer. J. Math. \textbf{120} (1998), 391--349.
% Parabolic monotonicity formulas, Local

\bibitem{CafLin2008:art}
{\sc L. A. Caffarelli, F.-H. Lin.} \emph{Singularly perturbed elliptic systems and multi-valued har-
monic functions with free boundaries}, J. Amer. Math. Soc. \textbf{21} (2008), 847--862.
% Regularity segregated profiles, Local

%\bibitem{CafRoqueSire2010:art}
%{\sc  L. A. Caffarelli, J.-M. Roquejoffre, Y. Sire.} \emph{Variational problems for free boundaries for
%	the fractional Laplacian}, J. Eur. Math. Soc. (JEMS) \textbf{12} (2010), 1151--1179.
% Regularity free boundaries, Variational, Fractional

\bibitem{CafSal2005:book}
{\sc L.A. Caffarelli, S. Salsa.} ``A Geometric Approach to Free Boundary Problems'', Grad. Stud. Math. \textbf{68}, AMS, 2005.
% Parabolic monotonicity formulas, local

\bibitem{CaffSil2007:art}
{\sc L. A. Caffarelli, L. Silvestre.}  \emph{An extension problem related to the fractional Laplacian},
Comm. Partial Differential Equations \textbf{32} (2007), 1245--1260.
%

\bibitem{Chen1998:art}
{\sc X.-Y. Chen}. \emph{A strong unique continuation theorem for parabolic equations}, Math. Ann. \textbf{311} (1998), 603--630.
%

\bibitem{ChiarenzaSerapioni1985:art}
{\sc F. M. Chiarenza, R. P. Serapioni.} \emph{A remark on a Harnack inequality for degenerate parabolic equations}, Rend. Sem. Mat. Univ. Padova \textbf{73} (1985), 179–190.
%


%\bibitem{ConTerVer2003:art}
%{\sc M. Conti, S. Terracini, G. Verzini.} \emph{An optimal partition problem related to nonlinear
%	eigenvalues}, J. Funct. Anal. \textbf{198} (2003), 160--196.
% Optimal partition prob, local


%\bibitem{ConTerVer2005_1:art}
%{\sc M. Conti, S. Terracini, G. Verzini.} \emph{Asymptotic estimates for the spatial segregation of
%	competitive systems}, Adv. Math. \textbf{195} (2005), 524--560.
% Spacial segregation, Asymptotic Estimates, local


\bibitem{DainelliGarofaloPetTo2017:book}
{\sc D. Danielli, N. Garofalo, A. Petrosyan, T. To.} ``Optimal regularity and the free boundary in the parabolic Signorini problem'', Memoirs AMS \textbf{249}, 2017.
%

%\bibitem{DanWanZha2011:art}
%{\sc E. N. Dancer, K. Wang, Z. Zhang.} \emph{Uniform H\"{o}lder estimate for singularly perturbed parabolic systems of Bose-Einstein condensates and competing species}, J. Differential Equations \textbf{251} (2011), 2737--2769.
% parabolic Holder estimate, local

%\bibitem{DanWanZha2012_1:art}
%{\sc E. N. Dancer, K. Wang, Z. Zhang.} \emph{The limit equation for the Gross-Pitaevskii equations
%	and S. Terracini’s conjecture}, J. Funct. Anal. \textbf{262} (2012), 1087--1131.
% Parabolic, limits local

\bibitem{DipierroSavinVald2017:art1}
{\sc S. Dipierro, O. Savin, E. Valdinoci}. \emph{All functions are locally s-harmonic up to
a small error}, J. Eur. Math. Soc. \textbf{19} (2017), 957--966.
%

\bibitem{DipierroSavinVald2016:art1}
{\sc S. Dipierro, O. Savin, E. Valdinoci}. \emph{Local approximation of arbitrary functions by solutions of nonlocal equations}, J. Geom. Anal. \textbf{29}, 1428--1455.
%

\bibitem{DipierroSavinVald2016:art2}
{\sc S. Dipierro, O. Savin, E. Valdinoci}. \emph{Definition of fractional Laplacian for functions with polynomial growth}, Arxiv preprint, arXiv:1610.04663 (2016).
%

\bibitem{DuvantLions1972:art}
{\sc G. Duvaut, J.-L.Lions}. \emph{Les in\'equations en m\'ecanique et en physique}, (French) Travaux et Recherches Math\'ematiques, No. 21. Dunod, Paris, 1972.
%

\bibitem{EscauriazaFernandez2003:art}
{\sc L. Escauriaza, J. Fernandez.} \emph{Unique continuation for parabolic operators}, Ark. Mat. \textbf{41} (2003), 35--60.
%

\bibitem{EscauriazaFernandezVess2006:art}
{\sc L. Escauriaza, J. Fernandez, S. Vessella.} \emph{Doubling properties of caloric functions}, Appl. Anal. \textbf{85} (2006), 205--223.
%

\bibitem{FabesKenigSerapioni1982:art}
{\sc E. B. Fabes, C. E. Kenig, R. P. Serapioni.} \emph{The local regularity of
solutions of degenerate elliptic equations}, Comm. Partial Differential Equations \textbf{7} (1982), 77--116.
%

\bibitem{FallFelli2014:art}
{\sc M. Fall, V. Felli} \emph{Unique continuation property and local asymptotics of solutions to fractional elliptic equations}, Comm. Partial Differential Equations, \textbf{39} (2014), 354--397.
%

\bibitem{FelliPrimo1949:art}
{\sc V. Felli, A. Primo.} \emph{Classification of local asymptotics for solutions to heat equations with inverse-square potentials}, Discrete Contin. Dyn. Syst. \textbf{31} (2011), 65--107.
%

%\bibitem{FernandezRosOton2017:art}
%{\sc X. Fern\'andez-Real, X. Ros-Oton.} \emph{Regularity theory for general stable operators: Parabolic equations}, J. Funct. Anal. \textbf{272} (2017), 4165--4221.
%

%\bibitem{Garofalo2018:art}
%{\sc N. Garofalo.} \emph{Fractional thoughts}, ArXive preprint, arXiv:1712.03347v4 (2018).
%

\bibitem{GarofaloLin1986:art}
{\sc N. Garofalo, F. Lin.} \emph{Monotonicity properties of variational integrals, $A_p$ weights and unique continuation}, Indiana Univ. Math. J. \textbf{35} (1986), 245--268.
%

\bibitem{GarofaloLin1987:art}
{\sc N. Garofalo, F. Lin.} \emph{Unique continuation for elliptic operators: a geometric-variational approach}, Comm. Pure Appl. Math. \textbf{40} (1987), 347--366.
%

\bibitem{GarofaloPetrosyan2008:art}
{\sc N. Garofalo, A. Petrosyan.} \emph{Some new monotonicity formulas and the singular set in the lower
dimensional obstacle problem}, Invent. Math. \textbf{177} (2009), 415--461.
%

\bibitem{GarofaloTralli2019:art}
{\sc N. Garofalo, G. Tralli.} \emph{A class of nonlocal hypoelliptic operators and their extensions}, ArXive preprint, arXiv:1811.02968.
%

\bibitem{GutierrezWheeden1991:art}
{\sc  C. E. Guti\'errez, R. L. Wheeden.} \emph{Harnack’s inequality for degenerate parabolic equations}, Comm. Partial Differential Equations \textbf{16} (1991), 745--770.
%

\bibitem{HanLin1994:art}
{\sc Q. Han, F.-H. Lin}. \emph{Nodal sets of solutions of parabolic equations II}, Comm. Pure Appl. Math. \textbf{47} (1994), 1219--1238.
%

\bibitem{Hilfer2000:book}
{\sc R. Hilfer.} ``Applications of Fractional Calculus in Physics'', World Scientific, 2000.
% Frac Deriv

\bibitem{Hirschman1952:note}
{\sc I. I. Hirschman}. \emph{A note on the Heat Equation}, Washington University, 1952.
% Liouville parabolic local theorem

\bibitem{Ishige1999:art}
{\sc K. Ishige.} \emph{On the behavior of the solutions of degenerate parabolic equations}, Nagoya Math. J. \textbf{155} (1999), 1--26.
%

\bibitem{Jones1977:art}
{\sc F. Jones.} \emph{A fundamental solution for the heat equation which is supported in a strip}, J. Math. Anal. Appl. \textbf{60} (1977), 314--324.
%

\bibitem{KenkreEtAl:1973}
{\sc V. M. Kenkre, E. W. Montroll, M. F. Shlesinger}. \emph{Generalized Master Equations for Continuous-Time Random Walks}, J. Stat. Phys. \textbf{9}, 45--50.
% CTRW

\bibitem{Lieberman1996:art}
{\sc G. Lieberman.} ``Second order parabolic equations'', World Scientific Publishing Co., Inc., River Edge, NJ, 1996.
%

\bibitem{Lin1991:art}
{\sc F.-H. Lin}. \emph{Nodal sets of solutions of elliptic and parabolic equations}, Comm. Pure Appl. Math. \textbf{45} (1991), 287--308.
%

\bibitem{Marchaud1927:art}
{\sc A. Marchaud}. \emph{Sur les d\'eriv\'ees et sur les diff\'erences des fonctions de variables r\'eelles}, J. math. pures et appl. \textbf{9} (1927), 337--425.
%

\bibitem{MetzlerKlafter2000:art}
{\sc R. Metzler, J. Klafter.} \emph{The random walk's guide to anomalous diffusion: a fractional dynamics approach}, Phys. Rep. \textbf{339} (2000), 1--77.
% CTRW

\bibitem{Monneau2003:art}
{\sc R. Monneau.} \emph{On the number of singularities for the obstacle problem in two dimensions}, J. Geom. Anal. \textbf{13} (2003), 359--389.
%

\bibitem{Monneau2009:art}
{\sc R. Monneau.} \emph{Pointwise estimates for Laplace equation. Applications to the free boundary
of the obstacle problem with Dini coefficients}, J. Fourier Anal. Appl. \textbf{15} (2009), 279--335.
%

\bibitem{MontrollWeiss1965:art}
{\sc E. W. Montroll, G. H. Weiss.} \emph{Random walks on lattices. II}, J. Mathematical Phys. \textbf{6} (1965), 167-181.
% CTRW

\bibitem{Muckenhoupt1972:art}
{\sc B. Muckenhoupt}, \emph{Hardy inequalities with weights}, Studia Math. \textbf{44} (1972), 31--38.
%

%\bibitem{NorTavTerVer2010:art}
%{\sc B. Noris, H. Tavares, S. Terracini, G. Verzini.} \emph{Uniform H\"{o}lder bounds for nonlinear
%Schrödinger systems with strong competition}, Comm. Pure Appl. Math., \textbf{63} (2010), 267--302.
% Regularity Gross-Pita, local

\bibitem{NystromSande2016:art}
{\sc K. Nystr\"{o}m, O. Sande.} \emph{Extension properties and boundary estimates for a fractional heat
operator}, Nonlinear Anal. \textbf{140} (2016), 29--37.
%

\bibitem{Poon1996:art}
{\sc C. C. Poon.} \emph{Unique continuation for parabolic equations}, Comm. Partial Differential Equations \textbf{21} (1996), 521--539.
%

\bibitem{Riesz1938:art}
{\sc M. Riesz}. \emph{Int\'egrales de Riemann-Liouville et potentiels}, Acta Sci. Math. Szeged, \textbf{9} (1938), 1--42.
%

\bibitem{Riesz1949:art}
{\sc M. Riesz}. \emph{L'int\'egrale de Riemann-Liouville et le probl\'eme de Cauchy}, (French) Acta Math. \textbf{81} (1949), 1--223.
%

\bibitem{Ruland2017:art}
{\sc A. R\"{u}land}. \emph{On quantitative unique continuation properties of fractional Schrödinger equations: Doubling, vanishing order and nodal domain estimates}, Trans. Amer. Math. Soc. \textbf{369} (2017), 2311--2362.
%

\bibitem{Samko2001:book}
{\sc S. Samko}. ``Hypersingular Integrals and their Applications'', CRC Press, 2001.
%

\bibitem{SamkoKilbasMarichev:1993}
{\sc S. Samko, A. A. Kilbas, O. Marichev}. ``Fractional Integrals and Derivatives. Theory and Applications'', Edited and with a foreword by S. M. Nikol'skii. Translated from the 1987 Russian original. Revised by the authors. Gordon and Breach Science Publishers, Yverdon, 1993. xxxvi+976 pp.
%

\bibitem{Simon1983:book}
{\sc L. Simon.} ``Lectures on Geometric Measure Theory'', Proceedings of the Centre for mathematical analysis, Australian National University \textbf{3} (1983).
% Appendix A

\bibitem{Simon1987:art}
{\sc J. Simon.} \emph{Compact sets in the space $L^p(0,T;B)$}, Ann. Mat. Pura Appl. \textbf{146} (1987), 65-96.
%

\bibitem{SireTerTor2018:art}
{\sc Y. Sire, S. Terracini, G. Tortone.} \emph{On the nodal set of solutions to degenerate or singular elliptic equations with an application to $s$-harmonic functions}, ArXiv preprint, arXiv:1808.01851.
%

\bibitem{SireTerVita2019:art}
{\sc Y. Sire, S. Terracini, S. Vita.} \emph{Liouville type theorems and regularity of solutions to degenerate or singular problems part I: even solutions}, ArXiv preprint, arXiv:1904.02143.
%

\bibitem{SireTerVita2020:art}
{\sc Y. Sire, S. Terracini, S. Vita.} \emph{Liouville type theorems and regularity of solutions to degenerate or singular problems part II: odd solutions}, preprint (2020).
%

\bibitem{StingaTorrea2009:art}
{\sc P. R. Stinga, J. L. Torrea.} \emph{Extension problem and Harnack's inequality for some fractional operators}, Comm. Partial Differential Equations \textbf{35} (2010), 2092--2122.
%

\bibitem{StingaTorrea2011:art}
{\sc P. R. Stinga, J. L. Torrea.} \emph{Regularity theory for the fractional harmonic oscillator}, J. Funct. Anal. \textbf{260} (2011), 3097--3131.
%

\bibitem{StingaTorrea2015:art}
{\sc P. R. Stinga, J. L. Torrea.} \emph{Regularity theory and extension problem for fractional nonlocal parabolic equations and the master equation}, SIAM J. Math. Anal. \textbf{49} (2017), 3893--3924.
%

\bibitem{Szego1939:art}
{\sc G. Szeg\"{o}.} ``Orthogonal Polynomials'', Colloquium publications - AMS \textbf{23} (1939).
%

%\bibitem{TavTer2012:art}
%{\sc H. Tavares, S. Terracini.} \emph{Regularity of the nodal set of segregated critical configurations
%under a weak reflection law}, Calc. Var. Partial Differential Equations \textbf{45} (2012), 273--317.
% Regul. Nodal Set, local

\bibitem{TerTorVita2018:art}
{\sc S. Terracini, G. Tortone, S. Vita.} \emph{On $s$-harmonic functions on cones}, Anal. PDE \textbf{11} (2018), 1653--1691.
%

\bibitem{TerVerZil2017:art}
{\sc S. Terracini, G. Verzini, A. Zilio.} \emph{Uniform H\"{o}lder bounds for strongly competing systems
involving the square root of the laplacian},, J. Eur. Math. Soc. \textbf{18} (2016), 2865--2924.
%

\bibitem{TerVerZil2016:art}
{\sc S. Terracini, G. Verzini, A. Zilio.} \emph{Uniform H\"{o}lder regularity with small exponent in competition-fractional diffusion systems}, Discrete Contin. Dyn. Syst. \textbf{34} (2014), 2669--2691.
%

\bibitem{Tortone2018:book}
{\sc G. Tortone.} ``On the nodal set of solutions of degenerate-singular and nonlocal equations'', Ph.D. Dissertation (2018).
%

\bibitem{Valdinoci2009:art}
{\sc E. Valdinoci.} \emph{From the long jump random walk to the fractional Laplacian}, Bol. Soc. Esp. Mat.
Apl. \textbf{49} (2009), 33--44.
%

%\bibitem{VerziniZilio2014:art}
%{\sc G. Verzini, A. Zilio}. \emph{Strong competition versus fractional diffusion: the case of Lotka-Volterra interaction}, Comm. Partial Differential Equations, \textbf{39} (2014), 2284--2313.
%

\bibitem{Vita2018:book}
{\sc S. Vita.} ``Strong competition systems ruled by anomalous diffusion'', Ph.D. Dissertation (2018).
%

\bibitem{Weiss1999a:art}
{\sc G. S. Weiss.} \emph{A homogeneity improvement approach to the obstacle problem}, Invent. Math. \textbf{138} (1999), 23--50.
%

\bibitem{Weiss1999b:art}
{\sc G. S. Weiss.} \emph{Self-similar blow-up and Hausdorff dimension estimates for a class of parabolic free boundary problems}, SIAM J. Math. Anal. \textbf{30} (1999), 623--644.
%

\bibitem{Weyl1917:art}
{\sc H. Weyl.} \emph{Bemerkungen zum Begriff des Differentialquotienten gebrochener Ordnung}, Zürich. Naturf. Ges. \textbf{62}, (1917) 296--302.
%

\bibitem{WuZhang2017:art}
{\sc J. Wu, L. Zhang.} \emph{Backward uniqueness of parabolic equations with variable coefficients in a half space}, Commun. Contemp. Math. \textbf{18} (2016), 1550011 (38 pages).
%

\bibitem{Yu2017:art}
{\sc J. Wu, L. Zhang.} \emph{Unique continuation for fractional orders operators of elliptic equations}, Ann. Partial Differ. Equ. \textbf{3} (2017), (21 pages).
%

\bibitem{Zaslavsky1994:art}
{\sc  G. Zaslavsky.} \emph{Fractional kinetic equation for Hamiltonian chaos. Chaotic advection, tracer
dynamics and turbulent dispersion}, Phys. D \textbf{76} (1994), 110--122.
%
\end{thebibliography}
\end{document}